\documentclass[10pt,leqno]{article}

\usepackage{amssymb} 
\usepackage{amsmath} 

\title{Infinite dimensional primitive linearly compact Lie superalgebras} 
\author{{\sc Nicoletta Cantarini}\thanks{Dipartimento di Matematica
Pura ed Applicata, Universit\`a di Padova, Padova, Italy -
Partially supported by Progetto Giovani Ricercatori CPDG031245}
\and\setcounter{footnote}{6}
{\sc Victor G.\ Kac}\thanks{Department of Mathematics, MIT, Cambridge,
Massachusetts 02139 - Partially supported by NSF grants DMS0201017 and DMS0501395}}

\newtheorem{theorem}{Theorem}[section] 
\newtheorem{lemma}[theorem]{Lemma} 
\newtheorem{corollary}[theorem]{Corollary} 
\newtheorem{proposition}[theorem]{Proposition} 
\newtheorem{definition}[theorem]{Definition} 
\newtheorem{remark}[theorem]{Remark}
\newtheorem{example}[theorem]{Example}
\def\Z{\mathbb{Z}} 
\def\ZZ{\mathbb{Z}}
\def\g{\mathfrak{g}}
\def\fg{\mathfrak{g}}
 
\def\N{\mathbb{N}}
\def\C{\mathbb{C}}
\def\CC{\mathbb{C}}
\def\ZZ{\mathbb{Z}}

\numberwithin{equation}{section} 

\newcommand{\st}[1]{\ensuremath{^{\scriptstyle \textrm{#1}}}}

\makeatletter

\def\enumerate{%
  \ifnum \@enumdepth >\thr@@\@toodeep\else
    \advance\@enumdepth\@ne
    \edef\@enumctr{enum\romannumeral\the\@enumdepth}%
      \list
        {\csname label\@enumctr\endcsname}%
        {\usecounter\@enumctr
          \addtolength{\leftmargin}{-\leftmargin}
          \settowidth{\labelwidth}{(99)}
          \itemindent = \labelwidth
          \addtolength{\itemindent}{\labelsep}
        \listparindent=1em      
          \def\makelabel##1{{##1}\hfill}
          }%
  \fi}

\makeatother

\newcommand{\arabicparenlist}{
  \renewcommand{\theenumi}{\arabic{enumi}}%
  \renewcommand{\labelenumi}{(\theenumi)}%
}

\date{} 
\begin{document} 

\maketitle 
\date{} 
\begin{abstract} We classify open maximal subalgebras of all
 infinite-dimensional linearly
compact simple Lie superalgebras. This is applied to the classification
of infinite-dimensional Lie superalgebras of vector fields, acting
transitively and primitively in a formal neighborhood of a point
of a finite-dimensional supermanifold.
\end{abstract} 
\section*{Introduction}
A well-known theorem of E.\ Cartan \cite{Ca} states that an
infinite-dimensional Lie algebra $L$ of vector fields in a
neighborhood of a point $p$ of an $m$-dimensional manifold $M$ acting
transitively and primitively in this neighborhood, is formally
isomorphic to a member of one of the six series of Lie algebras of
formal vector fields:
\begin{itemize}
\item[$1.$] $W_m=\{\sum_{i=1}^m P_i\frac{\partial}{\partial x_i} ~|~
  P_i\in\C[[x_1, \dots, x_m]]\}$,
\item[$2.$] $S_m=\{X\in W_m ~|~ div(X)=0\}$,
\item[$2^\prime.$] $CS_m=\{X\in W_m ~|~ div(X)=const\}$,
\item[$3.$] $H_m=\{X\in W_m ~|~ X\omega_s=0\}$ $(m=2k)$, where
  $\omega_s=\sum_{i=1}^k dx_i\wedge dx_{k+i}$ is a symplectic form,
\item[$3^\prime.$] $CH_m=\{X\in W_m ~|~ X\omega_s=const \,\omega_s\}$
  $(m=2k)$,
\item[$4.$] $K_m=\{X\in W_m ~|~ X\omega_c=f\omega_c\} (m=2k+1)$, where
  $\omega_c=dx_m+\sum_{i=1}^k x_idx_{k+i}$ is a contact form and $f$
  is a formal power series (depending on $X$).
\end{itemize}
Recall that the primitivity of an action means that there are no
non-trivial $L$-invariant fibrations in $M$. The Lie algebra $L$ has a
canonical filtration $L\supset L_0\supset L_1\supset \dots$, where
$L_j$ consists of vector fields vanishing up to order $j$ at $p$, and
the formal isomorphism means the isomorphism of the completed Lie
algebras with respect to this filtration.
The transitivity of the action implies that $L_0$ contains no non-zero
ideals of $L$, and primitivity implies that $L_0$ is a maximal
subalgebra.

It is easy to see (cf. \cite{G1}) that, in fact, Cartan's theorem gives a classification of
infinite-dimensional linearly compact Lie algebras $L$, which admit a
maximal open subalgebra $L_0$ containing no non-zero ideals of
$L$ (recall that the {\it linear compactness}\, of $L$ means 
that $L$ is a topological Lie algebra whose underlying topological space
is a topological product of finite-dimensional vector spaces
with discrete topology). Such a Lie algebra $L$ is called {\it primitive}, the subalgebra
$L_0$ is called a fundamental maximal subalgebra, and the pair $(L, L_0)$ is called
a {\it primitive pair}. It is easy to see that  all $L$ from the  six series
contain  a unique fundamental maximal subalgebra. Also, the Lie algebras $W_m$,
$S_m$, $H_m$ and $K_m$ are simple, and the remaining two series $CS_m$
and $CH_m$ are the Lie algebras of derivations of $S_m$ and $H_m$
respectively, obtained by adding the Euler vector field
$E=\sum_ix_i\frac{\partial}{\partial x_i}$.

In the present paper we solve the problem of classification of
primitive pairs in the Lie superalgebra case.
This problem is much more difficult than in the Lie algebra case for several 
reasons. First, in the Lie algebra case, 
a primitive $L$ is contained between $S$ and $Der S$, where $S$ is
simple (cf.\ Theorem \ref{Theorem2}),
which instantly reduces the classification of primitive Lie algebras
$L$ to simple ones, but the situation is more complicated in the super case. Second, there are many more
simple linearly compact Lie superalgebras than in the Lie algebra case
(see \cite{K}). Third, in a sharp contrast to the Lie algebra case, almost
all infinite-dimensional simple linearly compact Lie superalgebras
contain more than one maximal open subalgebra. Most of the space
of the present paper deals with the problem of their classification.

The infinite-dimensional linearly compact simple Lie superalgebras
have been classified in \cite{K}. The list consists of ten series
($m\geq 1$): $W(m,n)$, $S(m,n)$ ($(m,n)\neq (1,1)$), $H(m,n)$ ($m$ even), $K(m,n)$ ($m$ odd), $HO(m,m)$ ($m\geq 2$), $SHO(m,m)$
($m\geq 3$),
$KO(m,m+1)$, $SKO(m,m+1;\beta)$ ($m\geq 2$), $SHO^\sim(m,m)$ ($m$ even), $SKO^\sim(m,m+1)$
($m\geq 3$, $m$ odd), and
five exceptional Lie superalgebras: $E(1,6)$, $E(3,6)$, $E(3,8)$,
$E(4,4)$, $E(5,10)$.

The main idea of \cite{K} is to pick a maximal open subalgebra $S_0$ of a simple linearly compact Lie superalgebra $S$,
which is invariant with respect to all inner automorphisms of $S$. The existence of such an
{\it invariant} subalgebra $S_0$ is a non-trivial fact, the proof of which uses characteristic
varieties (cf. \cite{G2}). Remarkably, an invariant subalgebra is unique in most,
though not all, of the cases. 
After that the classification procedure is more or less routine. One
constructs an irreducible Weisfeiler filtration \cite{W} associated
to the pair $(S, S_0)$ and shows, using ideas from \cite{G2},  
that the associated graded Lie superalgebra $Gr S=
\oplus_j\mathfrak{g}_j$ 
has the property that $[\mathfrak{g}_0, v]=\mathfrak{g}_{-1}$ for any even element 
$v \in \mathfrak{g}_{-1}$ 
(which does not hold for a random choice of
a maximal open subalgebra $S_0$). After that one is able to describe
all possibilities for the $\mathfrak{g}_0$-module $\mathfrak{g}_{-1}$
and the subalgebra $\oplus_{j\leq 0} \mathfrak{g}_j$ of $Gr S$ 
\cite{K}, and, after some further work, all the possibilities for $Gr S$ \cite{CK}. Finally, 
one finds all simple filtered deformations of all these $Gr S$ \cite{ChK}.

Recall that one has the following isomorphisms
(cf.\ \cite[Remark 6.6]{K}):
$$W(1,1)\cong K(1,2)\cong KO(1,2), ~~S(2,1)\cong HO(2,2), ~~SHO^\sim(2,2)\cong H(2,1).$$
Besides, $S(2,1)\cong SKO(2,3;0)$. 
Hence, when discussing $S(m,n)$, $K(m,n)$, $KO(m,m+1)$, $HO(m,m)$ and $SHO^\sim(m,m)$, 
we shall assume that $(m,n)\neq (2,1)$,
 $(m,n)\neq (1,2)$, $m\geq 2$, $m\geq 3$ and $m>3$, respectively. 
Also we shall assume that $n\geq 1$ since the Lie algebra case is trivial.

In the first part of the present paper we give a description of
semisimple artinian linearly compact Lie superalgebras in terms of
simple ones (Theorem \ref{Theorem1}), which is similar to that in the
finite-dimensional case (\cite{K2}, \cite{Ch}).
Next, we show that if an infinite-dimensional linearly compact Lie
superalgebra $L$ is primitive, then $L$ is artinian
semisimple and, moreover, contains an open ideal
isomorphic to $S\otimes\Lambda(n)$, where $S$ is a simple linearly
compact Lie superalgebra and $\Lambda(n)$ is the Grassmann algebra in
$n$ indeterminates, and is contained in $(Der
S)\otimes\Lambda(n)+1\otimes Der\Lambda(n)$, so that the projection
of $L$ on the second summand acts transitively on $\Lambda(n)$
(Theorem \ref{Theorem2}).

Next, Theorem \ref{Theorem4} gives a description of fundamental maximal subalgebras in $L$
in terms of those in $S$. In fact, the situation is slightly more complicated,
namely in general $Der S=S\rtimes\mathfrak{a}$, where either
$\mathfrak{a}\cong gl_2$ 
or $\dim\mathfrak{a}\leq 3$, and we need to classify all 
 maximal among open $\mathfrak{a}_0$-invariant subalgebras of $S$, for each subalgebra 
$\mathfrak{a}_0$
of $\mathfrak{a}$. 
We, thus, arrive at a problem of
classification of  maximal among $\mathfrak{a}_0$-invariant open subalgebras of each 
infinite-dimensional simple linearly compact Lie superalgebra $S$.

If $S=\prod_{j\geq -d}\g_j$ is an {\it irreducible}
grading of $S$, i.e.\ the $\g_0$-module $\g_{-1}$ is irreducible and
$\g_{-j}=\g_{-1}^j$ for all $j\leq -2$, then $S_0=\prod_{j\geq 0}\g_j$
is called a {\it graded} subalgebra of $S$. 
All irreducible gradings (apart for a few omissions) were described in \cite{CK} and \cite{S}, and in
the present paper we give a detailed proof that these are all. It
turns out by inspection (using Proposition \ref{basic} $(b)$) that for every irreducible grading of a 
simple infinite-dimensional linearly compact Lie superalgebra $S$, the corresponding graded
subalgebra $S_0=\prod_{j\geq 0}\g_j$ is maximal. 

 A surprising discovery of the present paper is a large number of new open 
maximal
subalgebras (which are not graded). The main result of the present
paper is a classification of all maximal open 
$\mathfrak{a}_0$-invariant subalgebras of all
infinite-dimensional linearly compact simple Lie superalgebras $S$,
up to conjugation by the group $G$ of inner automorphisms of $Der S$.
(The group $G$ can be thought of as the unity component of the group
of all automorphisms of $S$.) Unless otherwise specified, by conjugation we always
mean the conjugation by $G$.

An important part of this classification is the classification of all regular maximal
open $\mathfrak{a}_0$-invariant subalgebras of $S$. A subalgebra of $S$ is called {\it regular}
if it is invariant with respect to a maximal torus of $Der S$. 
By Theorem \ref{Theorem3}, all maximal tori in $Der S$ are
conjugate, hence fixing one ``standard'' torus $T$, and classifying all
$T$-invariant maximal open subalgebras 
we obtain all regular maximal open subalgebras of $S$, up to conjugation (by 
$G$).

The numbers $a$ of graded and $b$ of non-graded regular maximal open
subalgebras of $S$, up to conjugation,
are given in Table 1 (the case $\mathfrak{a}_0$=0).
Thus we see that, with the exception of $K(m,1)$, any simple infinite-di\-men\-sional
primitive Lie superalgebra which is not a Lie algebra
 contains more than one maximal open subalgebra.
It turned out that in all cases
except for $H(m,n)$ with $n$ positive even, all maximal open subalgebras are regular, but
$H(m,n)$ with $n=2h$ even, contains, up to 
conjugacy, $h(h+1)/2$ non-regular maximal open subalgebras. 

\bigskip

\begin{table}[htbp]
\begin{center}
\begin{tabular}{c|c|c|c|c}
$S$ & a & b & c & e \\
\hline
$W(1,1)$ & $2$ & 0 & 3 & 1\\
$W(m,n)$, $(m,n)\neq (1,1)$ & $n+1$ & 0 & $n+1$ & 0\\
$S(1,2)$  & $2$ & 0& $4$ & 2 \\
$S(m,n)$, $(m,n)\neq (1,2)$ & $n+1$ & 0& $n+1$ & 0 \\
$K(1,2h)$\!  & $h+1$ & 0& $h+2$ & 1 \\
$K(m,2h)$, $m> 1$ & $h+2$& 0 & $h+2$& 0 \\
$K(m,2h+1)$ & $h+1$& 0 & $h+1$& 0 \\
$HO(n,n)$, $n>2$ & $n$ & 0 &$n+1$&1\\
$SHO(3,3)$ \! & $2$ & 0 &$5$&3\\ 
$SHO(n,n)$, $n>3$  & $n$ & 0 &
$n+1$ & $1$\\
$H(m,2h)$& $h+2$ & $\frac{h}{2}(1+h)$ &$h^2+2h+2$&0\\
$H(m,2h+1)$ & $h+1$ & $(h+1)^2$ &$h^2+3h+2$&0\\
$KO(2,3)$  & $2$ & $2$&$4$&0\\
$KO(n,n+1)$, $n>2$ & $n$ & $n$ & $2n+2$ & $2$\\
$SKO(2,3;0)$  & 2 & 0 &2&0\\
$SKO(2,3;1)$ & $2$ & 1 & 3 &0\\
$SKO(2,3;\beta)$, $\beta\neq 0,1$ & $3$ & 1 & 5 &1\\
$SKO(3,4;\beta)$\!  & $3$ & $3$&$8+8\delta_{3\beta, 1}$&
$2+8\delta_{3\beta, 1}$\\ 
$SKO(n,n+1;\beta)$, $n>3$ & $n$ & $n$ &$2n+2$&2\\
$SHO^\sim(n,n)$, $n>2$ & 1 & $n-1$ &$n+1$&1\\
$SKO^\sim(n,n+1)$ & 0 & $2n-1$ &$2n+2$&3\\
 $E(1,6)$ & 4 & 0 &5&1\\
 $E(3,6)$ & 3 & 0 &5&2\\
 $E(5,10)$ & 4 & 0 &6&2\\
$E(4,4)$ & 1 & 3 &5&1\\
 $E(3,8)$ & 3 & 6 &18&9\\
%
\end{tabular}
\begin{center}
\textbf{Table 1.}
\end{center}
\end{center}
\end{table}

The main tool in the classification of maximal open subalgebras in 
non-exceptional
simple linearly compact Lie superalgebras is a formal analogue of the
Frobenius theorem (Theorem \ref{1}$(a)$), which implies that a maximal
open subalgebra of a transitive subalgebra of $W(m,n)$
consists of vector fields, leaving invariant a conjugate, by a change of variables,
of a standard ideal of $\Lambda(m,n)$, that is, an ideal generated by a
subspace of the span of all odd indeterminates.
This instantly solves the problem in question for $W(m,n)$,
but for other non-exceptional simple Lie superalgebras it requires  more
subtle arguments to show that a conjugate of a standard ideal of $\Lambda(m,n)$
can be replaced by a standard ideal. 

In the case of exceptional simple linearly compact Lie superalgebras $S$
we use the notions of growth and size of $S$ (which remain unchanged
when passing from $S$ to $Gr S$) in order to list possible $Gr S$. This allows
us to find all maximal open subalgebras of $S$ by analyzing 
its deviation from a
maximal open invariant subalgebra (which is unique
in all exceptional superalgebras $S$).

A posteriori, it follows from the present paper that an open
subalgebra of minimal codimension in a
linearly compact infinite-dimensional simple Lie superalgebra $L$
 is always invariant
under all inner automorphisms of $L$. Moreover, in all cases, but 
$S=W(1,1)$, $S(1,2)$, 
$SHO(3,3)$, and $SKO(3,4;1/3)$, such a subalgebra
is unique (hence invariant under all automorphisms), and in 
$S=W(1,1)$, $S(1,2)$, and 
$SHO(3,3)$ such subalgebras
are conjugate by (outer)
automorphisms of $S$. We denote by $S_0$ the intersection of all 
open subalgebras 
of minimal 
codimension in $S$, and call it the \emph{canonical subalgebra} of $S$.
The canonical subalgebra is, of course, invariant with respect to the group 
$Aut S$ of all continuous automorphisms of $S$.
Let $S_{-1}$ be a minimal subspace of $S$, properly containing $S_0$
and invariant with respect to the group $Aut S$, and let
$S=S_{-d} \supsetneq
S_{-d+1} \supset \cdots \supset S_{-1} \supset S_0 \supset
\cdots$ be the associated Weisfeiler filtration of
$S$.  
All members of the Weisfeiler  filtration associated to $S_0$
 are invariant with respect
to the group $Aut S$, 
hence we have the induced filtration on the
superspace $V: =S/S_0 = V_{-d} \supset \cdots \supset V_{-1}$,
and the induced action of $Aut S$ on $V$ preserving this filtration.
Note that $Gr V$ carries a canonical $\ZZ$-graded Lie
superalgebra structure, isomorphic to $\oplus^{-1}_{j=-d} \fg_j$.
A subspace $U$ of $V$ is called \emph{abelian} if $Gr U$ is an
abelian subalgebra of $Gr V$.

Now, it is easy to see that if $L_0$ is a (proper) open subalgebra of $S$, 
its image under the canonical map $S\rightarrow V$
is a purely odd abelian subspace of $V$, denoted by $\pi
(L_0)$.  Thus, we obtain an $Aut S$-equivariant map $\pi$ from the set
of all open subalgebras of $S$ to the set of abelian subspaces of
$V_{\bar{1}}$ (the odd part of $V$).

The $G$-orbit of $\pi (L_0)$ in $V_{\bar{1}}$ is called the
\emph{canonical invariant} of the open subalgebra $L_0$ of $S$.
A posteriori, it turns out that the canonical invariant uniquely
determines an open maximal subalgebra of $S$, so we have an
injective map $\Pi$ from the set of conjugacy classes (by $G$) of
maximal open subalgebras of $S$ to the set of $G$-orbits of
abelian subspaces of $V_{\bar{1}}$.  The number $c$ of
elements of the latter set along with the number~$e$ of those of
them which are not canonical invariants of any open maximal
subalgebra are given in Table~1.  Looking at this table, we see
that in many cases
$e=0$, i.e.,~the map $\Pi$ is
bijective, and in the remaining cases it is very close to being bijective.

The contents of the paper are as follows. In Section 1 we prove 
a formal analogue of the Frobenius theorem (Theorem \ref{1}),
establish some general
classificational results on artinian semisimple and primitive infinite-dimensional
linearly compact Lie superalgebras (Theorems \ref{Theorem1} and \ref{Theorem2}),
and reduce the classification of primitive pairs $(L,L_0)$ to the case
of simple $L$ (Theorem \ref{Theorem4}). We also prove conjugacy of
maximal tori in  artinian semisimple linearly compact Lie superalgebras
(Theorem \ref{Theorem3}),
and discuss the notions of growth and size. 

In Sections 2 through 10 we give a classification of all  maximal 
open subalgebras (and all $\mathfrak{a}_0$-invariant maximal open  subalgebras as well)
of all infinite-dimensional simple linearly compact
Lie superalgebras.
As an immediate application of this long and tedious work, we obtain the list of all 
irreducible graded
infinite-dimensional linearly compact Lie superalgebras which admit a non-trivial
simple filtered deformation.

In Section 11 we classify all maximal open subalgebras which are invariant with respect to all 
inner automorphisms 
and we
discuss the canonical invariant. An a priori proof that the canonical invariant determines a
maximal open subalgebra uniquely would considerably shorten the paper, but
we were unable to find such a proof.

In a subsequent paper \cite{CantaK} we use the canonical subalgebras
to describe automorphisms and real forms of all simple infinite-dimensional
linearly compact Lie superalgebras.
  
Throughout the paper all vector spaces and algebras, as well as tensor products,
are considered over the field of complex numbers $\C$.

\section{General results on semisimple and primitive linearly compact
Lie superalgebras}\label{victor}
Recall that a {\it linearly compact space} is a topological
vector space which is isomorphic to a topological product of
finite-dimensional vector spaces endowed with discrete topology. The
basic examples of linearly compact spaces are finite-dimensional
vector spaces with the discrete topology, and the space of formal
power series $V[[t]]$ over a finite-dimensional vector space $V$, with
the formal topology defined by taking as a fundamental system of
neighborhoods of $0$ the set $\{t^jV[[t]]\}_{j\in\Z_+}$. 
We recall Chevalley's principle (\cite{G1}): if $F_1\supset
F_2\supset\dots$ is a sequence of closed subspaces of a linearly
compact space such that $\cap_jF_j=0$ and $U$ is an open subspace,
then $F_j\subset U$ for $j>>0$.

A {\it linearly
compact superalgebra} is a topological superalgebra whose underlying
topological space is linearly compact. The basic example of an
associative linearly compact superalgebra is
$\Lambda(m,n)=\Lambda(n)[[x_1,\dots, x_m]]$, where $\Lambda(n)$
denotes the Grassmann algebra on $n$ anticommuting indeterminates
$\xi_1, \dots, \xi_n$, and the superalgebra parity is defined by
$p(x_i)=\bar{0}$, $p(\xi_j)=\bar{1}$, with the formal topology defined
by the following fundamental system of neighborhoods of $0$:
$\{(x_1, \dots, x_m, \xi_1,\dots,\xi_n)^j\}_{j\in\Z_+}$. The basic example of a
linearly compact Lie superalgebra is $W(m,n)=Der \Lambda(m,n)$, the
Lie superalgebra of all continuous derivations of the superalgebra
$\Lambda(m,n)$. One has:
$$W(m,n):=\{X=\sum_{i=1}^mP_i(x,\xi)\frac{\partial}{\partial
x_i}+\sum_{j=1}^nQ_j(x,\xi)\frac{\partial}{\partial\xi_j} ~|~ P_i,
Q_j\in\Lambda(m,n)\}.$$
Letting $\deg x_i=\deg \xi_j=1$, $\deg\frac{\partial}{\partial x_i}=
\deg\frac{\partial}{\partial \xi_j}=-1$, we obtain the {\em principal}
$\Z$-grading $W(m,n)=\prod_{j\geq -1} W(m,n)_j$. 
A subalgebra $L$ of $W(m,n)$ is called {\it transitive} if the
projection of $L$ on $W(m,n)_{-1}$  is onto. 

Given a subspace $U$ of the subspace $\sum_{i=1}^m\C x_i+\sum_{j=1}^n\C\xi_j$
of $\Lambda(m,n)$, denote by
$I_U$ the ideal of $\Lambda(m,n)$ generated by $U$. 
Let $W_U=\{X\in W(m,n)~|~XI_U\subset I_U\}$ be the corresponding
subalgebra of $W(m,n)$.
More generally, for any subalgebra $L$ of
$W(m,n)$, let $L_U=\{X\in L~|~XI_U\subset I_U\}$. We
shall call $I_U$ a {\em standard ideal} of $\Lambda(m,n)$ and $L_U$ a
{\em standard subalgebra} of $L$. 

\begin{theorem}\label{1} $(a)$ Let $L$ be a closed subalgebra 
of $W(m,n)$, let $V$ be the projection of $L$ on $W(m,n)_{-1}=
\sum_i\C\frac{\partial}{\partial x_i}+\sum_j\C\frac{\partial}{\partial\xi_j}$,
and let $V^*\subset \sum_i\C x_i+\sum_j\C\xi_j$ be the dual subspace of $V$.
Then there exists a continuous automorphism of $\Lambda(m,n)$ such that
the induced automorphism of $W(m,n)$ maps $L$ to the subalgebra
$W_{V^*}$ of $W(m,n)$.

$(b)$ The algebra
$\Lambda(m,n)$ has no non-trivial closed $L$-invariant ideals if and
only if $L$ is a transitive subalgebra.
\end{theorem}
{\bf Proof.}
Making a linear change of variables, we may assume that $V$ is the span
of $\frac{\partial}{\partial x_1}, \dots, \frac{\partial}{\partial x_p},
\frac{\partial}{\partial \xi_1}, \dots, \frac{\partial}{\partial \xi_q}$.
Also, we may assume that $L$ is invariant with respect to multiplication
by elements of $\Lambda(m,n)$. Indeed, $\Lambda(m,n)W_U=W_U$,
an ideal of $\Lambda(m,n)$ is $L$-invariant if and only if it is
$\Lambda(m,n)L$-invariant, and $L$ is transitive if and only if
$\Lambda(m,n)L$ is transitive.

We turn now to the proof of $(a)$.
If $p\geq 1$, then
$L$ contains a vector field $X_1=\frac{\partial}{\partial x_1}+D_1$,
where $D_1$ is an even operator such that $D_1(0)=0$. Making change of variables (cf.\ \cite[p.\ 12]{K}),
we may assume that $X_1=\frac{\partial}{\partial x_1}$. Consider
$X_2=\frac{\partial}{\partial x_2}+D_2\in L$, $D_2(0)=0$. Subtracting
$f\frac{\partial}{\partial x_1}$ from $X_2$, we may assume that $X_2$,
hence $D_2$, do not involve $\frac{\partial}{\partial x_1}$. Next, we
show that we may assume that all coefficients of $D_2$ do not involve $x_1$. Here we use
that $L$ is a subalgebra. Let $D_2=\sum_{j\geq 0}
x_1^j\overline{D}_j$. Since $[X_1, X_2]=[X_1, D_2]\in L$, we see that $\sum_{j\geq
  0}jx_1^{j-1}\overline{D}_j\in L$, hence, $x_1\sum_{j\geq
  0}jx_1^{j-1}\overline{D}_j\in L$, then, $D_2-\sum_{j\geq
  0}jx_1^{j}\overline{D}_j\in L$, and 
we can assume that $\overline{D}_1=0$. Repeating this procedure, since
$L$ is closed, we get  
$\frac{\partial}{\partial x_2}+\overline{D}_0\in L$, where
$\overline{D}_0(0)=0$ and $\overline{D}_0$ does not  depend on $x_1$.
 Making change of variables, we may assume that
$\frac{\partial}{\partial x_1}, \frac{\partial}{\partial x_2}\in L$. Continuing one gets
$\frac{\partial}{\partial x_1}, \dots, \frac{\partial}{\partial
  x_p}\in L$. 
Indeed, if  $q\geq 1$, let $Y_1$  be an odd vector field in $L$ whose
projection on $W(m,n)_{-1}$ is $\frac{\partial}{\partial \xi_1}$. Up
to a change of variables
we may assume that
$Y_1=\frac{\partial}{\partial \xi_1}+\xi_1D$, 
where $D$ is an even operator. Since $[Y_1, Y_1]=2D$,
$D$ lies in $L$, hence $\xi_1D\in L$ and  $\frac{\partial}{\partial
  \xi_1}\in L$. Then, arguing as above, we can assume that
$\frac{\partial}{\partial \xi_1}, \dots, \frac{\partial}{\partial
  \xi_q}$ lie in $L$. 
Hence  $L$ is generated, as a $\Lambda(m,n)$-module, by
$\frac{\partial}{\partial x_1}, \dots, \frac{\partial}{\partial x_p},
\frac{\partial}{\partial \xi_1}, \dots, \frac{\partial}{\partial
  \xi_q}$ and by vector fields $X_k$ which do not involve $\frac{\partial}{\partial x_1}, \dots, \frac{\partial}{\partial x_p},
\frac{\partial}{\partial \xi_1}, \dots, \frac{\partial}{\partial
  \xi_q}$ and such that $X_k(0)=0$. 
As above, we may assume that all coefficients of all $X_k$ do not
depend on $x_1, \dots, x_p, \xi_1. \dots, \xi_q$.
 Therefore the ideal of $\Lambda(m,n)$ generated by
$x_{p+1},\dots, x_m, \xi_{q+1}, \dots,\xi_n$ is $L$-invariant,
which proves $(a)$. 

Now we prove $(b)$. The transitivity of $L$ is equivalent to saying that $L$
contains elements $a_i=\frac{\partial}{\partial x_i}+X$ and
$b_j=\frac{\partial}{\partial \xi_j}+Y$ for some vector fields $X$ and $Y$
such that $X(0)=0$ and $Y(0)=0$, for every $i$ and $j$.
Let $I$ be an $L$-stable non-zero ideal of $\Lambda(m,n)$. Then $I$ contains
a non-zero element $P(x,\xi)\in\Lambda(m,n)$. Since $I$
is stable under the action of the vector fields $a_i$ and $b_j$, we may
assume that $P(0,0)=1$ and, since $I$ is an ideal, 
multiplying $I$ by $P^{-1}$, we find that $I$ contains $1$, i.e., $I=\Lambda(m,n)$. Conversely, if $L$ is not transitive, then $V^*\neq 0$, and we arrive
at a contradiction with $(a)$.
 \hfill$\Box$

\begin{remark}\label{frobenius}\em
Theorem \ref{1}$(a)$ is an analogue of
  the Frobenius theorem for the superalgebra $\Lambda(m,n)$.
Namely, if the projection $V$ of $L$ on $W(m,n)_{-1}$ has dimension
  $(p\,|\,q)$, then there exists a continuous automorphism $\varphi$ of
  $\Lambda(m,n)$ such that the ideal $J_V=(\varphi(x_{p+1}), \dots,
  \varphi(x_m), \varphi(\xi_{q+1}), \dots, \varphi(\xi_n))$ is
  $L$-invariant. 
  In the purely odd case this was proved in \cite{FKR}. 

Note that
  $J_V$ is maximal among $L$-invariant ideals. Indeed, up to automorphisms, this is
  equivalent to saying that, if 
 $V=\langle \frac{\partial}{\partial x_1}, \dots, \frac{\partial}{\partial x_p},
\frac{\partial}{\partial \xi_1}, \dots, \frac{\partial}{\partial
  \xi_q}\rangle$ then 
$J_V=(x_{p+1}, \dots,
  x_m, \xi_{q+1}, \dots, \xi_n)$ is maximal among
  $L$-invariant ideals. Indeed, if we 
add a polynomial $P$ to the ideal $J_V$, we may assume that $P$ depends only
on the variables $x_1, \dots, x_p, \xi_1, \dots, \xi_q$. Then, since
$\frac{\partial}{\partial x_i}$ and $\frac{\partial}{\partial \xi_j}$
  lie in $L_0$ for every $i=1, \dots, p$ and $j=1, \dots q$,  
adding $P$ to the ideal adds $1$ to it. 
\end{remark}

\begin{remark}\label{rome}\em
Let L be an infinite-dimensional 
linearly compact Lie superalgebra 
embedded in $W(m,n)$, and
let $L_0$ be a fundamental maximal subalgebra of $L$ such that
the projection of $L_0$ to $W(m,n)_{-1}$ does not contain the even derivations 
$\frac{\partial}{\partial x_i}$ for any $i=1,\dots,m$. 
Then, by \ Theorem 1.1(a), $L_0$  stabilizes an ideal $J$ of $\Lambda(m,n)$ which is, up to
changes of variables, a standard ideal
containing all even indeterminates $x_1,\dots , x_m$.
Besides, 
$J$ is maximal among the $\Lambda(m,n)L_0$-invariant ideals of $\Lambda(m,n)$
by \ Remark \ref{frobenius}.
Notice that an ideal $I$ of $\Lambda(m,n)$ is $L_0$-invariant if and
only if
it is  $\Lambda(m,n)L_0$-invariant. Therefore $J$ is also
maximal among the $L_0$-invariant ideals of $\Lambda(m,n)$.
It follows that
$L_0$ stabilizes an ideal $J=(x_1+f_1,\dots, x_m+f_m, \eta_1+g_1,..., \eta_r+g_r)$
for some
linear functions $\eta_j$ in odd indeterminates and even functions
$f_i$ and odd functions $g_j$ without constant and linear terms,
and that 
$J$ is maximal among the
$L_0$-invariant ideals of $\Lambda(m,n)$.
\end{remark}

Recall that a linearly compact Lie superalgebra $L$ is called {\it simple} if it is
not abelian and contains no closed ideals different from $0$ and $L$;
$L$ is called {\it semisimple} if it contains no non-zero abelian
ideals; $L$ is called {\it artinian} if any descending sequence of
closed ideals in $L$ stabilizes.

A subalgebra $L_0$ of $L$ is called {\it fundamental} if it is proper, open and contains no
non-zero ideals of $L$.
Due to Guillemin's theorem \cite{G1} a linearly compact Lie superalgebra is
artinian if and only if it contains a fundamental subalgebra (the
proof in \cite{G1} is given in the Lie algebra case, but it extends
verbatim to the super case).

Let $L_0$ be a fundamental subalgebra of a Lie superalgebra $L$ and let $L_{-1}$ be
an $ad\,L_0$-stable subspace of $L$ generating $L$ as a Lie superalgebra. The Weisfeiler
filtration (\cite{W}) associated to the triple $L\supset L_{-1}\supset L_0$ is the
filtration of $L$ inductively defined as follows:
for $s\geq 1$, 
$$L_{-(s+1)}=\left[ L_{-1}, L_{-s}\right]+L_{-s}, ~~L_s=\{a\in L_{s-1}~|~
\left[ a, L_{-1} \right]\subset L_{s-1}\}.$$
If $L_{-1}$ is a minimal $ad\, L_0$-stable subspace properly
containing $L_0$, then the Weisfeiler filtration is called
\emph{irreducible}.  
If $L_{-1}=L$,  the Weisfeiler filtration is called the 
\emph{canonical filtration}.

A linearly compact Lie superalgebra $L$ is called {\it primitive} if
it contains a fundamental subalgebra $L_0$ which is a maximal
subalgebra.   In this case, $(L,L_0)$ is called a \emph{primitive
pair}.  Note that for a primitive pair $(L,L_0)$ there exists an
irreducible Weisfeiler filtration whose $0$\st{th} term is $L_0$.

Given a filtered Lie superalgebra $L=L_{-d}\supset\dots\supset L_{-1}\supset L_0\supset L_1\supset\dots$, we shall denote by $Gr L$ the associated $\Z$-graded
Lie superalgebra:
$$Gr L=\oplus_{j\geq -d}~Gr_jL, ~~~Gr_jL=L_j/L_{j+1}.$$

If ${\mathfrak{g}}=\oplus_{j\geq -d}{\mathfrak{g}}_j$ is a graded Lie
superalgebra, we denote by $\overline{\mathfrak{g}}=\prod_{j\geq -d}\mathfrak{g}_j$ its 
completion. Then $\overline{\mathfrak{g}}$ has a natural filtration 
 given by the subspaces
$$\overline{\mathfrak{g}}_i=\prod_{j\geq i}\mathfrak{g}_j$$
for $i\geq -d$.
We shall call such a filtration a {\em graded} filtration (or, equivalently,
 a trivial filtered deformation of $\overline{\mathfrak{g}}$, cf.\ Section \ref{last}).

Let $L_0$ be a fundamental subalgebra of $L$, let
$L=L_{-d}\supset \cdots \supset L_{-1}\supset L_0 \supset L_1
\supset \cdots$ be a Weisfeiler filtration, and let $GrL=
\oplus_{j \geq -d} \fg_j$, where $\fg_j=Gr_jL$, be the associated
graded superalgebra.  Then (\cite{W}):
\begin{eqnarray}
  \label{eq:1.1}
  &&\fg_{-j}=\fg^j_{-1} \hbox{ for }j \geq 1\, , \\
 \label{eq:1.2}
&&\hbox{if } x \in \fg_j \, , \, j \geq 0 \hbox{ and }
[x,\fg_{-1}]=0 \, ,\, \hbox{ then } x=0\, .
\end{eqnarray}
If, in addition, the Weisfeiler filtration is irreducible, then
\begin{equation}
  \label{eq:1.3}
  \fg_{-1} \hbox{ is an irreducible } \fg_0
\hbox{-module}\, .
\end{equation}
A $\ZZ$-graded Lie superalgebra $\fg = \oplus_{j \geq -d}$
$\fg_i$ is called \emph{transitive} if properties (\ref{eq:1.1})
and (\ref{eq:1.2}) hold, and it is called \emph{irreducible} if,
in addition (\ref{eq:1.3}) holds.

\medskip

The following theorems describe the artinian semisimple and the
infinite-dimensional primitive linearly compact Lie superalgebras.
We denote by $Der S$ (resp.\ $Inder S$)
the Lie superalgebra of all (resp.\ all inner) continuos derivations
of a linearly compact Lie superalgebra $S$.

\begin{theorem}\label{Theorem1} Let $S_1, \dots, S_r$ ($r\in\N$) be simple linearly
compact Lie superalgebras, let $m_1, \dots, m_r, n_1, \dots, n_r$ be
non-negative integers and let
$S=\oplus_{i=1}^r$ $(S_i\hat{\otimes}\Lambda(m_i, n_i))$. Then
\begin{equation} Der S= \oplus_{i=1}^r((Der
S_i)\hat{\otimes}\Lambda(m_i, n_i)+1\otimes W(m_i, n_i))
\label{(1.1)}
\end{equation}
is a linearly compact Lie superalgebra and $S=Inder S$ canonically
embeds in $Der S$.
Let $L$ be an open subalgebra of $Der S$ containing $S$, and denote by
$F_i$ the projection of $L$ on $1\otimes W(m_i, n_i)$. Then
\begin{itemize}
\item[$(a)$] $L$ is semisimple if and only if $F_i$ is a transitive
subalgebra of $W(m_i, n_i)$ for all $i=1, \dots, r$.
\item[$(b)$] All artinian semisimple linearly compact Lie superalgebras
can be obtained as in (a).
\item[(c)] If $L$ is semisimple, then $Der L$ is the normalizer of $L$
in $Der S$ (and is semisimple).
\end{itemize}
\end{theorem}
{\bf Proof.} It follows traditional lines (cf.\ \cite{B},
\cite{Ch}). Let $L$ be an artinian semisimple
linearly compact Lie superalgebra, and let $I$ denote the sum of
all its minimal closed ideals. Since $L$ is
semisimple, for any (non-zero) minimal closed ideal $J$ one has $[J,J]=J$.
Using this, it is standard to show that $I$ is a direct sum of all
minimal closed ideals of $L$. Since $L$ is artinian,
it follows that it contains a finite number of (non-zero)
minimal closed ideals; denote them by $I_1,\dots, I_r$. We have
a homomorphism 
$\varphi :
L\rightarrow \oplus_j Der I_j$ defined by $\varphi(a)=\sum_j(ad
~a)|_{I_j}$. 
 The homomorphism $\varphi$ is injective since
$(\ker\varphi) \cap I =0$,
and therefore, by the artinian property, if $\ker\varphi$ is non-zero, it would
contain a (non-zero) minimal closed ideal different from all $I_j$'s.
Thus, we have the following inclusions:
\begin{equation}
\oplus_{j=1}^r I_j \subset L\subset \oplus_{j=1}^r Der I_j. 
\label{(1.2)}
\end{equation}
Next we use the super analogue of the Cartan-Guillemin theorem
(\cite{G1}, \cite{Bl}), established in \cite{FK}, according to which
$I_j\cong S_j\hat{\otimes}\Lambda(m_j, n_j)$, where $S_j$ is a simple
linearly compact Lie superalgebra and $m_j, n_j\in\Z_+$.

Next, the same argument as in \cite{B} or \cite{Ch} shows that 
$$Der I=(Der S_j)\hat{\otimes}\Lambda(m_j, n_j)+1\otimes W(m_j, n_j),$$
and that $L$ in (\ref{(1.2)}) is semisimple if and only if $\Lambda(m_j,
n_j)$ contains no non-trivial $F_j$-invariant ideals. Now $(a)$ and
$(b)$ follow from Theorem \ref{1}. The proof of $(c)$ is the same as in
\cite{B} or \cite{Ch}. \hfill$\Box$
\begin{theorem}\label{Theorem2} If $L$ is an infinite-dimensional primitive Lie
superalgebra, then $L$ is artinian semisimple, and, moreover,
$$S\otimes\Lambda(n)\subset L\subset (Der S)\otimes\Lambda(n)+1\otimes
W(0,n)$$
for some infinite-dimensional simple linearly compact Lie superalgebra
$S$ and $n\in\Z_+$, where the projection of $L$ on $W(0,n)$ is a
transitive subalgebra.
\end{theorem}
{\bf Proof.} By the above mentioned Guillemin's theorem, $L$ is
artinian. By another result of Guillemin (\cite[Proposition 4.1]{G1}),
whose proof works verbatim in the super case, any non-zero closed
ideal of $L$ has finite codimension.

In order to show that $L$ is semisimple, choose an irreducible Weisfeiler
filtration of $L$ associated with the  fundamental maximal subalgebra
$L_0$ of $L$, and let $\g=\oplus_{j\geq -d}\g_j$ be the associated
graded Lie superalgebra. Suppose that $L$ contains a non-zero closed
abelian ideal. Then the corresponding ideal $I=\oplus_{j\geq -d}I_j$
in $\g$ has finite codimension, and since $\dim\g=\infty$, we conclude
that $I_j\neq 0$ for some $j\geq 0$. By the transitivity of $\g$,
$I_0\neq 0$ and $I_{-1}\neq 0$, and by the irreducibility of the
$\g_0$-module $\g_{-1}$, $I_{-1}=\g_{-1}$. Hence $\left[I_0,
\g_{-1}\right]=0$ (since $I$ is an abelian ideal), which contradicts
the transitivity of $\g$.

Thus, by Theorem \ref{Theorem1}, $L$ contains the ideals
$S_i\hat{\otimes}\Lambda(m_i, n_i)$, $i=1,\dots, r$. Since $\dim
L=\infty$ and all non-zero ideals of $L$ have finite codimension, we
conclude that $r=1$ and 
$$S\hat{\otimes}\Lambda(m, n)\subset L\subset (Der
S)\hat{\otimes}\Lambda(m, n)+1\otimes W(m,n) \, ,$$
where $S$ is a simple linearly compact Lie superalgebra and the
projection $F$ of $L$ on $W(m,n)$ is a transitive subalgebra. It
remains to show that $m=0$.

Since $L_0$ is a maximal subalgebra of $L$, and
$S\hat{\otimes}\Lambda(m, n)$ is an ideal, we conclude that
\begin{equation}
L=L_0+(S\hat{\otimes}\Lambda(m, n)).
\label{(1.3)}
\end{equation}
Suppose that $m\geq 1$. Since $L_0$ is an open subalgebra, by
Chevalley's principle,
\begin{equation}
S\hat{\otimes}(x_1, \dots, x_m)^j\Lambda(m, n)\subset L_0
~~~\mbox{for}~ j>>0.
\label{(1.4)}
\end{equation}
By transitivity of $F$ and (\ref{(1.3)}), the projection of $L_0$ on
$1\otimes W(m,n)_{-1}$ is surjective. Hence it follows from
(\ref{(1.4)}) that $S\hat{\otimes}\Lambda(m, n)\subset L_0$, a
contradiction since $L_0$ contains no non-zero ideals of $L$. \hfill$\Box$

\medskip

An $ad$-diagonalizable subalgebra $T$ of a linearly compact Lie
superalgebra $L$ is called a \emph{torus} of $L$.
The following proposition allows one to construct maximal tori.

\begin{proposition}
  \label{prop:1.5}
Let $L$ be a linearly compact Lie superalgebra with trivial
center and let $L=L_{-d} \supset \cdots \supset L_0 \supset L_1
\supset \cdots$ be a filtration of $L$ such that $L_0$ contains
all $ad$-exponentiable elements of $L$.  Then any torus $T$ of
$L$ lies in $L_0$ and $T$ is a maximal torus in $L$ if and only if its image
in $L_0 /L_1$ is a maximal torus.  Any maximal torus of $L_0
/L_1$ can be lifted to that of $L$.
\end{proposition}

{\bf Proof.}
Since all elements of $T$ are exponentiable, $T\subset L_0$.
Since, obviously, $T \cap L_1 =0$, $T$ is a maximal torus of $L_0$ (and
hence of $L$) if and only if its image is a maximal torus of $L_0 /L_1$.
 \hfill$\Box$

\bigskip

We do not know examples for which the maximal tori are not
conjugate, but we can prove their conjugacy only for the artinian
semisimple $L$ (which we shall apply to primitive~$L$).

\begin{theorem}\label{Theorem3} If $L$ is an artinian semisimple
linearly compact Lie superalgebra, then all maximal tori of $L$ are
conjugate by inner automorphisms of $L$.
\end{theorem}

{\bf Proof.} We may assume that $\dim L=\infty$. The socle
$\oplus_{i=1}^r S_i\otimes\Lambda(n_i)$ of $L$ (see Theorem
\ref{Theorem1}) is invariant with respect to all automorphisms of
$L$. But due to \cite{K}, each $Der S_i$ contains a fundamental
subalgebra $S_i^0$, which is proper  if $\dim
S_i=\infty$, and which contains all exponentiable elements of $Der
S_i$. 

Consider the Lie superalgebra
$$\tilde{L}=\oplus_{i=1}^r((Der S_i)\hat{\otimes}\Lambda(m_i,
n_i))\oplus (1\otimes W(m_i, n_i))$$
containing $L$. Take the canonical filtration of $Der S_i$ defined by
$S_i^0$ and tensor it with the filtration of $\Lambda(m_i, n_i)$ whose
$j$-th member is $(x_1, \dots, x_m, \xi_1, \dots, \xi_n)^j$; this
defines a filtration of $(Der S_i)\hat{\otimes}\Lambda(m_i, n_i)$ all
of whose exponentiable elements lie in the $0$-th member of the
filtration. These and the principal filtration of $W(m_i, n_i)$ for
each $i$ add up to produce a filtration of $\tilde{L}$. Intersecting
the members of this filtration with $L$, we get a filtration of $L$ by
open subspaces $L\supset L_0\supset L_1\supset\dots$, such that $L_0$
contains all exponentiable elements of $L$.
In particular, $L_0$ contains any two maximal tori $T$ and
$T^{\prime}$ of $L$. But  $T$ and
$T^{\prime}$ are conjugate in $L_0$ mod $L_N$ for each $N\geq 1$ by
the conjugacy of maximal tori in any finite-dimensional Lie
superalgebra. Taking the limit as $N\rightarrow\infty$, we obtain that $T$ and
$T^{\prime}$ are conjugate in $L_0$. \hfill$\Box$

\medskip

We shall use the following (corrected) explicit description of the Lie
superalgebras $Der S$ for all simple linearly compact Lie superalgebras $S$,
given in \cite{K}. 
\begin{proposition}\cite[Proposition 6.1]{K}\label{derivations}
Let $S$ be a simple infinite-dimensional linearly compact Lie superalgebra.
Then $Der S = S \rtimes \mathfrak{a}$, where $\mathfrak{a}$ is a
finite-dimensional subalgebra, described below:
\begin{itemize}
\item[$(a)$] If $S$ is one of the Lie superalgebras $W(m,n)$, $SHO^\sim(m,m)$,
$K(m,n)$, $KO(m,m+1)$,  $SKO^\sim(m,m+1)$,
$E(4,4)$, $E(1,6)$, $E(3,6)$, $E(3,8)$, then $\mathfrak{a}=0$.
\item[$(b)$] If $S$ is one of the Lie superalgebras $S(m,n)$ with $m\geq 2$,
$(m,n)\neq (2,1)$,
$H(m,n)$, $HO(m,m)$ with $m\geq 3$, $SKO(m,m+1;\beta)$ with
$m\geq 2$ and $\beta\neq 1, (m-2)/m$,
$E(5,10)$, then $\mathfrak{a}$ is a one-dimensional torus of $Der S$.
\item[$(c)$] If $S$ is one of the Lie superalgebras $S(1,n)$ with $n\geq 3$,
 $SKO(m,m+1;(m-2)/m)$ with $m\geq 2$,
$SKO(m,m+1;1)$ with $m>2$,
then $\mathfrak{a}=\mathfrak{n} \rtimes \mathfrak{t}_1$, where 
$\mathfrak{t}_1$ is a
one-dimensional torus of $Der S$ and $\mathfrak{n}$ is a one-dimensional subalgebra such that
$[\mathfrak{t}_1,\mathfrak{n}]=\mathfrak{n}$.
\item[$(d)$] If $S=SHO(m,m)$ with $m\geq 4$, then $\mathfrak{a}=
\mathfrak{n}\rtimes\mathfrak{t}_2$, where $\mathfrak{t}_2$
is a two-dimensional torus of $Der S$ and $\mathfrak{n}$ is a one-dimensional
subalgebra such that $[\mathfrak{t}_2, \mathfrak{n}]=\mathfrak{n}$.
\item[$(e)$] If $S=S(1,2)$ or $S=SKO(2,3;1)$, then
$\mathfrak{a}\cong sl_2$. 
\item[$(f)$] If $S=SHO(3,3)$, then $\mathfrak{a}\cong gl_2$.
\end{itemize}
The subalgebra $\mathfrak{a}$ of $Der S$ is called the {\em subalgebra
of outer derivations} of $S$. 
\end{proposition}

The following theorem describes all
primitive pairs in terms of simple ones.
\begin{theorem}\label{Theorem4} $(a)$ Let $L=S\otimes\Lambda(n)\rtimes
  F$, where $S$ is a linearly compact Lie superalgebra and $F$ is a
  transitive subalgebra of $W(0,n)$. Then any fundamental maximal
  subalgebra $L_0$ of $L$ is of the form $(S_0\otimes
  \Lambda(n))\rtimes F$, where $S_0$ is a fundamental maximal
  subalgebra of $S$.

$(b)$ Let $S$ be a simple
infinite-dimensional linearly compact Lie superalgebra. 
Let $\mathfrak{a}_0$ be a subalgebra of the subalgebra $\mathfrak{a}$
of outer derivations of $S$ and let $L_0$
be a fundamental maximal subalgebra of $S\rtimes\mathfrak{a}_0$. Then $L_0=S_0\rtimes
\mathfrak{a}_0$, where $S_0$ is a maximal among open $\mathfrak{a}_0$-invariant
subalgebras of $S$. Thus, all fundamental maximal subalgebras of
$S\rtimes\mathfrak{a}_0$ are $S_0\rtimes\mathfrak{a}_0$, where
$S_0$ is a maximal among open $\mathfrak{a}_0$-invariant subalgebras of $S$.

$(c)$ Let $S$ be a simple
infinite-dimensional linearly compact Lie superalgebra.
  Let
$F$ be a subalgebra of
$(\mathfrak{a}\otimes\Lambda(n))\rtimes W(0,n)$ containing
elements $f_i$, for $i=1,\dots,n$, such that 
$f_i(0)=\frac{\partial}{\partial\xi_i}$.
Let $L=(S\otimes\Lambda(n))\rtimes
F$. Then these $L$ exhaust, up to automorphisms,
 all that occur in a primitive pair. All
possible fundamental maximal subalgebras $L_0$ in $L$ can be obtained as follows. Let
$\mathfrak{a}_0=\{a(0) ~|~ a(\xi)\in {\mbox{projection of}}~ F
~{\mbox{on}}~ \mathfrak{a}\otimes\Lambda(n)\}\subset
\mathfrak{a}$. Let $S_0$ be a  maximal among
$\mathfrak{a}_0$-invariant  subalgebras of $S$. Then
$L_0=(S_0\otimes\Lambda(n))\rtimes F$. 
\end{theorem} 
{\bf Proof.} $(a)$ First, we show that $F\subset L_0$. In
 the contrary case, consider an irreducible Weisfeiler filtration of
 $L$ associated to $L_0$. Then we have:
$Gr L=Gr(S\otimes\Lambda(n))\rtimes Gr F.$ Since $Gr_{-1}L$ is
 irreducible with respect to $Gr_0L$ and $Gr_{-1}(S\otimes\Lambda(n))$ is a
 submodule of $Gr_{-1}L$, we conclude that
 $Gr_{-1}(S\otimes\Lambda(n))=0$, i.e.,
 $Gr_{-1}L=Gr_{-1}F$, hence $Gr_{<0}L=Gr_{<0}F$. It follows that
 $S\otimes\Lambda(n)\subset L_0$, which is impossible since $S\otimes
 \Lambda(n)$ is an ideal of $L$.

We write elements of $S\otimes\Lambda(n)$ in the form
 $s(\xi)=\sum_Is_I\xi^I$, where $I=\{i_1, \dots, i_r\}\subset \{1,
 \dots, n\}$, $s_I\in S$, $\xi^I=\xi_{i_1}\dots\xi_{i_r}$. Let
 $S_I=\{s_I ~|~ s(\xi)\in L_0\}$; then $S_0:=S_{\emptyset}$ is a
 subalgebra of $S$. Due to the transitivity of $F$, we
 conclude that  $S_I\subset S_0$ for all $I$, hence
 $S_0\otimes\Lambda(n)\supset L_0\cap(S\otimes \Lambda(n))$. Since
 $F\subset L_0$, we deduce that
 $(S_0\otimes\Lambda(n))+F\supset L_0$. Hence these two subalgebras
 coincide due to the maximality of $L_0$. Since $L_0$ is a fundamental
 maximal subalgebra of $L$, $S_0$ is a fundamental maximal subalgebra
 of $S$.

 $(b)$ The same argument as in $(a)$ shows that $\mathfrak{a}_0\subset
  L_0$.
Therefore $L_0=(L_0\cap S)\rtimes\mathfrak{a}_0$, and $L_0\cap S$ is an
  $\mathfrak{a}_0$-invariant subalgebra of $S$. By the maximality of
  $L_0$ it follows that $L_0\cap S$ is maximal among the
  $\mathfrak{a}_0$-invariant subalgebras of $S$.

$(c)$ Let $(L,L_0)$ be a primitive pair. Then, by Theorem \ref{Theorem2},
$L=(S\otimes\Lambda(n))\rtimes F$, where $S$ is a  simple Lie superalgebra,
$\mathfrak{a}$ is the subalgebra of outer derivations of $S$,
 and $F$  is a subalgebra of 
$(\mathfrak{a}\otimes\Lambda(n))\rtimes W(0,n)$ with transitive projection
on $W(0,n)$. Since
the projection of $F$ on $W(0,n)$ is transitive, we may assume, up
to automorphisms, that $F$ contains some elements $f_i$, for every
$i=1,\dots,n$, such that $f_i(0)=\frac{\partial}{\partial \xi_i}$.
Indeed, if $g_i$ lies in $F$,
 such that $g_i(0)=\frac{\partial}{\partial \xi_i} +a_i$ for some 
$a_i\in\mathfrak{a}$, the automorphism $1+ad(a_i\xi_i)$
brings $g_i$ to an element $f_i$ such that $f_i(0)=\frac{\partial}{\partial \xi_i}$.

The same argument as in $(a)$ shows that
$F\subset L_0$, hence $L_0=L_0\cap (S\otimes\Lambda(n))\rtimes F$.
Let us write the elements of $S\otimes\Lambda(n)$ in the form
$s(\xi)$ as in $(a)$, let $S_I$ be defined as in $(a)$, and let
$S_0=\{s(0) ~|~ s(\xi)\in L_0\}$. Then $S_0$ is a subalgebra of $S$
and, since $f_i\in L_0$ for every $i=1,\dots,n$, $S_I\subset S_0$ for all
$I$. It follows that $L_0\subset (S_0\otimes\Lambda(n))\rtimes F$, hence,
by the maximality of $L_0$, equality holds.

Likewise, let us write the elements of $\mathfrak{a}\otimes\Lambda(n)$ in
the form $a(\xi)=\sum a_I\xi^I$, let $\mathfrak{a}_0$ be as in
the statement and let $\mathfrak{a}_I=\{a_I ~|~ a(\xi)\in {\mbox{projection of}}~ F
~{\mbox{on}}~ \mathfrak{a}\otimes\Lambda(n)\}$. 
Then $\mathfrak{a}_0$ is a subalgebra of $\mathfrak{a}$ and, since $L_0$
contains the elements $f_i$, $\mathfrak{a}_I\subset \mathfrak{a}_0$ for
all $I$. It follows that $L_0\subset S_0\otimes\Lambda(n)+ \mathfrak{a}_0\otimes\Lambda(n)+F'$, where $F'$ is the projection of $F$ on $W(0,n)$.
Since $S_0$ is $\mathfrak{a}_0$-invariant, the maximality of $S_0$
among the $\mathfrak{a}_0$-invariant subalgebras of $S$ follows from
the maximality of $L_0$.  \hfill$\Box$

\bigskip

Recall that the growth of an artinian linearly
compact Lie superalgebra $L$ is defined as follows. Choose a
fundamental subalgebra $L_0$ of $L$ and construct a Weisfeiler
filtration $L=L_{-d}\supset\dots \supset L_0\supset L_1\supset\dots$,
containing $L_0$ as its $0$-th member, for some choice of $L_{-1}$
containing $L_0$ and generating $L$. Consider the function $F(j)=\dim
L/L_j$. It depends on the choice of $L_0$ and on the Weisfeiler
filtration, but it is easy to show (see \cite{BDK}, \cite{FK}), that the leading term
of $F(j)$ is independent of these choices. Namely, there exist unique
positive real numbers $a$ and $g$ such that
$\overline\lim_{j\rightarrow\infty}\frac{F(j)}{j^g}=a$. The number $g$ is
    called the {\it growth} of $L$, and is denoted by $g(L)$.

It is easy to see from the classification, that, if $L$ is simple, then $g(L)$ is a positive
integer and, moreover, $s(L):=a g(L)!$ is a positive integer
. The number $s(L)$ is called the {\em size} of $L$. One can think
of the growth (resp.\ size) of $L$ as the minimal number of even
variables (resp.\ minimal number of functions in these variables) involved
in vector fields from $L$.
It is also easy
to see from the classification that if $L$ is simple and is not a Lie algebra, then
$s(L)$ is an even integer, and, moreover, the sizes of the even and
the odd parts of $L$ are $\frac{1}{2}s(L)$ (of course, their growths
are both equal to $g(L)$). Due to Theorem \ref{Theorem2}, any
primitive $L$ contains $S\otimes\Lambda(n)$ as an open ideal, hence
$g(L)=g(S)$ and $s(L)=2^ns(S)$.

If a simple $L$ is of type $X(m,n)$, then $g(L)=m$. The sizes are
given in the following table:

\medskip

{\center \begin{tabular}{c|c||c|c||c|c}
$L$ & $s$ & $L$ & $s$ &  $L$ & $s$\\
\hline
$W(m,n)$ & $(m+n)2^{n}$ & $SHO(n,n)$ & $2^{n-1}$  &E(1,6) &  32\\
$S(m,n)$ & $(m+n-1)2^{n}$  & $KO(n,n+1)$ & $2^{n+1}$  &$E(3,6)$ & 12\\
$H(m,n)$  & $2^{n}$ &$SKO(n,n+1;\beta)$ & $2^{n}$  & $E(3,8)$ & 16\\
$K(m,n)$ & $2^{n}$ & $SHO^\sim(n,n)$ & $2^{n-1}$ & $E(4,4)$ & 8\\
$HO(n,n)$ & $2^{n}$ & $SKO^\sim(n,n+1)$  & $2^{n}$ & $E(5,10)$ & 8 
\end{tabular}}

\begin{center} Table 2.
\end{center}

\begin{remark}\label{rem:1.9}\em
  If $(L,L_0)$ is a primitive pair, and $Gr L$ is  its associated
graded superalgebra for a Weisfeiler filtration, 
then $g(L)=g(\overline{Gr L})$ and $s(L)=s(\overline{Gr L})$. This puts
stringent restrictions on the possibilities for $Gr L$ for the
given primitive pair $(L,L_0)$.
\end{remark}

The following proposition allows one to construct graded maximal
subalgebras.
\begin{proposition}\label{basic} 
Let ${\mathfrak{g}}=\oplus_{j\geq -d}{\mathfrak{g}}_j$ be a
$\Z$-graded Lie superalgebra and let 
${\mathfrak{g}}_{\geq 0}=\oplus_{j\geq 0}{\mathfrak{g}}_j$, 
${\mathfrak{g}}_{\pm}=\oplus_{j>0}{\mathfrak{g}}_{\pm j}$. 

\medskip

\noindent 
$(a)$ If ${\mathfrak{g}}_{\geq 0}$ is a maximal subalgebra of $\g$, then: 
\begin{itemize} 
\item[$(i)$] ${\mathfrak{g}}_{-1}$ is an irreducible 
${\mathfrak{g}}_0$-module; 
\item[$(ii)$] ${\mathfrak{g}}_-$ is generated by ${\mathfrak{g}}_{-1}$; 
\item[$(iii)$]${\mathfrak{g}}_-$ contains no ideals of $\g$ different from 
${\mathfrak{g}}_-$ or zero.
\end{itemize} 

\noindent 
$(b)$ If $(i)$ and $(ii)$ hold and, in addition, 
\begin{itemize} 
\item[$(iii)^{\prime}$]$[a,\mathfrak{g}_1]\neq 0$ for any non-zero $a\in \mathfrak{g}_j$, $j<-1$,
\end{itemize} 
then $\g_{\geq 0}$ is a maximal 
subalgebra of $\g$. 
\end{proposition} 
{\bf Proof.} $(a)$ $(i)$ If $V$ is a ${\mathfrak{g}}_0$-submodule of 
${\mathfrak{g}}_{-1}$ and $V_-$ is the subalgebra of $\g$
generated by $V$ then $V_-+{\mathfrak{g}}_{\geq 0}$ is a subalgebra 
of $\mathfrak{g}$. 

\noindent 
$(ii)$ If $\g_-^{\prime}$ is the subalgebra of $\g$ generated by 
$\g_{-1}$, then $\g_-^{\prime}+\g_{\geq 0}$ is a subalgebra of $\g$. 

\noindent
$(iii)$ If $I$ is such an ideal, then $I + \g_{\geq 0}$ is a subalgebra of $\g$.

\noindent 
$(b)$ Suppose that $\g_{\geq 0}$ is properly contained in a 
subalgebra $\g^{\prime}$ of $\g$. It follows that there
exists a non-zero element 
$a\in\g_-\cap \g^{\prime}$. Now $(iii)^{\prime}$ implies that 
$[a,\g_1]\neq 0$. It follows that $\g_{-1}\cap\g^{\prime}\neq \{0\}$, 
therefore, due to $(i)$, $\g_{-1}\subset\g^{\prime}$ and, due to 
(ii), $\g^{\prime}=\g$. 
~$\hfill\Box$ 


\begin{corollary}\label{cor}
If $L$ is a filtered Lie superalgebra such that $Gr L$ has
properties $(i)$, $(ii)$, $(iii)^{\prime}$ of Proposition \ref{basic}, then $L_0$ is a maximal
subalgebra of $L$. 
\end{corollary}

%
%
\begin{remark}\label{gd}\em If $\g=\oplus_{i\geq -d}\g_i$ is simple then $\g_{-d}$ 
is irreducible. Indeed if $V$ is a $\g_0$-stable subspace of $\g_{-d}$ then $V+(\oplus_{i>-d}\g_i)$
is an ideal of $\g$. In particular any $\Z$-grading of depth 1 of a simple Lie superalgebra is irreducible.
\end{remark}
\begin{definition} 
Let $T$ be a maximal torus in ${\mbox Der} L$. We call an open  subalgebra of $L$ 
regular if it is $T$-invariant. 
\end{definition}
\begin{remark}\label{regularideal}\em Let $L$ be a 
subalgebra of $W(m,n)$ and let $I_U$ 
be a standard ideal of $\Lambda(m,n)$. 
If $I_U$ is 
stabilized by a maximal torus $T$ of $Der L$ then the standard subalgebra
$L_U$ is regular.
\end{remark}
\section{Maximal open subalgebras of $\boldsymbol{W(m,n)}$, \break $\boldsymbol{S(m,n)}$, $\boldsymbol{K(m,n)}$, 
$\boldsymbol{HO(n,n)}$ and
$\boldsymbol{SHO(n,n)}$}\label{classical}
{\bf {\em The Lie superalgebras $\boldsymbol{W(m,n)}$ and $\boldsymbol{S(m,n)}$, $\boldsymbol{m\geq 1}$}.}
In Section \ref{victor} we introduced the Lie superalgebra $W(m,n)$ of continuous derivations of the Lie
superalgebra $\Lambda(m,n)$.
We shall  assume  $m\geq 1$ (note that $\dim W
(0,n)<\infty$). For every vector field $X$ in $W(m,n)$, 
$X=\sum_{i=1}^mP_i\frac{\partial}{\partial x_i}+\sum_{j=1}^nQ_j\frac{\partial}{\partial \xi_j}$, we shall set $X(0)=\sum_{i=1}^mP_i(0)\frac{\partial}{\partial x_i}
+\sum_{j=1}^nQ_j(0)\frac{\partial}{\partial \xi_j}$.  Let us fix the standard maximal torus $T=\langle x_i\frac{\partial}{\partial x_i}, 
\xi_j\frac{\partial}{\partial \xi_j} ~|~ i=1,\dots, m, j=1,\dots,n\rangle$ of
$W(m,n)$.

The simple Lie superalgebras $L$ considered in this section and in the
following three,  are
subalgebras of $W(m,n)$ such that $Der L\subset W(m,n)$ and $T\cap Der
L$ is a maximal torus of $Der L$. 
Such a maximal torus
of $Der L$ will be called {\em standard}.
\begin{remark}\label{standard}\em 
By Theorem \ref{Theorem3} each regular subalgebra of $L$
is conjugate by $G$
to a subalgebra which is invariant with
respect to the standard torus of $L$. Thus, in order to classify
regular subalgebras up to conjugation by $G$, it suffices to consider
the ones that contain $T$. In what follows, 
conjugation will always mean conjugation by $G$, unless otherwise
specified. We will often use automorphisms of $L$ defined by changes of
variables; each time it will not be difficult to check that they are
inner, hence lie in $G$. Note that when the linear part of a change of 
variables is
the identity then this is always an inner automorphism (cf. \cite{R}).
\end{remark}
\begin{remark}\label{gradingsofW}\em A $\Z$-grading, called the grading of type 
$(a_1,\dots, a_m|b_1,\dots, b_n)$, can be defined on $W(m,n)$ 
by setting $a_i=\deg x_i=-\deg\frac{\partial}{\partial x_i}\in\N$ and
$b_i=\deg\xi_i=-\deg\frac{\partial}{\partial \xi_i}\in\Z$ 
(cf.\ \cite[Example 4.1]{K}).
The $\Z$-grading of type $(1,\dots,1|1,\dots, 1)$ is the
principal grading of $W(m,n)$.
In this grading $W(m,n)$ has depth 1 with $0$-th graded component isomorphic to the Lie
superalgebra $gl(m,n)$ and $-1$-st graded component isomorphic to the
standard $gl(m,n)$-module $\C^{m|n}$. 
Since $W(m,n)$ is simple for every $(m,n)\neq (0,1)$, under our hypotheses the principal
grading of $W(m,n)$ is irreducible by Remark \ref{gd}. 
More  generally, 
the gradings of type $(1,\dots,1|1,\dots, 1, 0,\dots,0)$ with $k$ zeros, are
irreducible for every $k=0,\dots, n$ and satisfy the hypotheses of
Proposition \ref{basic}$(b)$. It follows, by Proposition
\ref{basic}$(b)$, that the corresponding subalgebras $\prod_{j\geq 0}W(m,n)_j$ of
$W(m,n)$ are maximal.
The $\Z$-grading of $W(m,n)$ of type $(1,\dots,1|0,\dots,0)$
is called {\em subprincipal}.
\end{remark}

\begin{theorem}\label{W(m,n)}
Let $W=W(m,n)$ with $m\geq 1$. Then all
maximal open subalgebras of\, $W$ are, up to conjugation, the
graded subalgebras of type
$(1,\dots,1|1,\dots,1,$ $0,\dots,0)$ with $k$ zeros, for  
$k=0,\dots, n$.
\end{theorem}
{\bf Proof.} Let $L_0$ be a maximal open subalgebra of $W$.
Since the vector fields $\frac{\partial}{\partial x_i}$ are not
exponentiable, $L_0$
does not contain any vector field
of the form $\sum\alpha_i\frac{\partial}{\partial x_i}+X+Y$ for
any non-zero linear combination $\sum\alpha_i\frac{\partial}{\partial x_i}$,
any $X\in W$ such that $X(0)=0$ and any $Y\in W(0,n)$.
By Theorem \ref{1}$(a)$,  $L_0$ is conjugate to the
subalgebra $W_U$ for some subspace 
$U=\langle x_1, \dots, x_m, \xi_1,\dots, \xi_k\rangle$
of $\Lambda(m,n)$, with $0\leq k\leq
n$. The subalgebra  $W_U$ is in fact
the graded subalgebra of type $(1,\dots,1|1,\dots,1,0,\dots,0)$ with
$n-k$ zeros.
\hfill $\Box$

\begin{definition}\label{divergence} Let $L$ be a subalgebra of $W(m,n)$.
A linear map $Div : L\rightarrow \Lambda(m,n)$
is called a {\em divergence} if the
action of $L$ on the space $\Lambda(m,n)v$
given by:
\begin{equation}
X(fv)=(Xf)v + (-1)^{p(X)p(f)} fDiv(X)v
\label{Div}
\end{equation}
is a representation of $L$.
The symbol $v$ is called the {\em volume form}
attached to the divergence $Div$.
Then $S'L:=\{X \in L~|~ Div(X)=0\}$ is a subalgebra of $L$.
Moreover, $Div$ is a homomorphism of $S'L$-modules.
\end{definition}
\begin{definition}\label{twisted}
If we have a representation of $L\subset W(m,n)$ on a vector
space $V$, 
which is also a left module over
$\Lambda(m,n)$, compatible with the action of $L$,
and $v$ is a volume form for $L$, then, for any complex number
$\lambda$, $L$ acts   on the space  $V^{\lambda}:=v^{\lambda} V$,
by the {\em twisted} action defined as follows: 
$$X(v^{\lambda} u)=\lambda v^{\lambda} Div(X)u +  v^{\lambda} Xu.$$
\end{definition}
\begin{remark}\label{annihilating}\em The subalgebra $S'L$ 
consists of vector fields $X$ in $L$ such that $Xv=0$.
Likewise,
$CS'L:=\{X\in L~|~ Div(X)\in \C\}$ is the subalgebra of $L$ 
consisting of vector fields $X$ in $L$ such that $Xv=cv$
with $c\in\C$. 
\end{remark}
\begin{remark}\label{tilde}\em
If $Div$ is a divergence and $F$ is an even invertible function in
$\Lambda(m,n)$, then the map $Div_F: L \rightarrow \Lambda(m,n)$
defined by: $$Div_F(X)=X(F)F^{-1}+Div(X)$$ 
is also a divergence. If $v$ is the volume form attached to $Div$,
then  $Fv$ is the volume form attached to $Div_F$.
\end{remark}
\begin{example}\label{S'W}\em  The function $div: W(m,n)\rightarrow \Lambda(m,n)$
defined by $$div(\sum_{i=1}^m P_i\frac{\partial}{\partial x_i}+
\sum_{j=1}^n Q_j\frac{\partial}{\partial\xi_j})=\sum_{i=1}^m
\frac{\partial P_i}{\partial x_i}+\sum_{j=1}^n(-1)^{p(Q_j)}
\frac{\partial Q_j}{\partial \xi_j}$$ is a divergence.
We will refer to it as the {\it usual divergence}. It follows,
according to Definition \ref{divergence}, that the set
$S^{\prime}(m,n):=S'W(m,n)=\{X\in W(m,n) ~|~ div(X)=0\}$
is a subalgebra of $W(m,n)$ (cf.\ \cite[Example
4.2]{K}). Moreover,
$CS'(m,n)=S'(m,n)+\mathbb{C}\sum_{i=1}^mx_i\frac{\partial}{\partial x_i}$.
\end{example}
\begin{remark}\label{usual}\em Let $div$ be the usual divergence 
(see  Example
\ref{S'W}). Then, for every $X\in W(m,n)$ and any even invertible function
$F\in \Lambda(m,n)$, $div(FX)=X(F)+Fdiv(X)$. Therefore
$div_F(X)=0$ if and only if $div(FX)=0$.
\end{remark}

Let $S(m,n)=[S^{\prime}(m,n), S^{\prime}(m,n)]$. We recall that
 if $m>1$ then
$S(m,n)=S^{\prime}(m,n)$ is simple. Besides,
$S^{\prime}(1,n)=S(1,n)+\C\xi_1\dots \xi_n\frac{\partial}{\partial
x_1}$ and
$S(1,n)$ is simple if and only if $n\geq 2$ (cf.\ \cite[Example
4.2]{K}). Since $S(2,1)\cong SKO(2,3;0)$, when talking about $S(m,n)$
we shall always assume $(m,n)\neq (2,1)$.

\begin{remark}\em Every $\Z$-grading of $W(m,n)$ induces a grading on $S(m,n)$. In
particular the
$\Z$-gradings of
type $(1,\dots,1|1,\dots, 1,0,\dots,0)$, with $k$ zeros, induce on
$S(m,n)$, by
Remark \ref{gd}, irreducible
gradings for $m>1$ or $m=1$ and $n\geq 2$. As in Remark \ref{gradingsofW}, the corresponding subalgebras $\prod_{j\geq 0}S(m,n)_j$ of
$S(m,n)$ are maximal. 
The $\Z$-grading of $S(m,n)$ of type $(1,\dots,1|0,\dots,0)$
is called {\em subprincipal}.
\end{remark}

\begin{theorem}\label{newS} Let $S=S(m,n)$ or $S=S'(m,n)$ or $S=CS'(m,n)$
 with $m>1$ or $m=1$ and $n\geq 2$.
Then every 
maximal open subalgebra of $S$ is 
regular.
\end{theorem}
{\bf Proof.} Let $L_0$ be a maximal open subalgebra of $S=S'(m,n)$ and let
$U=\langle x_1,\dots, x_m, \xi_{1},\dots, \xi_k\rangle$ with $0\leq
k\leq n$. Then, by Theorem \ref{1}, there exists a continuous
automorphism $\varphi$ of $\Lambda(m,n)$ such that $L_0=S\cap\varphi
W_{U}\varphi^{-1}$.

Let $\omega$ be the volume form attached to the divergence $div$. Then:
$$\varphi^{-1}S\varphi=\{\varphi^{-1}X\varphi ~|~ X\omega=0\}=\{Y~|~
\varphi Y \varphi^{-1}(\omega)=0\}=
\{Y~|~ Y(f\omega)=0,
~\mbox{for}$$

\noindent
~some invertible
$f\in\Lambda(m,n)\}=\{Y~|~ f^{-1}Yf\omega=0\}=fSf^{-1}.$

\medskip

It follows that:
$$\varphi^{-1}S\varphi\cap W_{U}=fSf^{-1}\cap W_{U}=\{fXf^{-1}~|~ X\in S,
~fXf^{-1}(I_{U})\subset I_{U}\}=$$
$$=\{fXf^{-1} ~|~ X\in S, ~X(I_{U})\subset I_{U}\}=f(S\cap
W_{U})f^{-1}.$$
Therefore $L_0=S\cap\varphi
W_{U}\varphi^{-1}=\varphi f(S\cap W_{U})f^{-1}\varphi^{-1}$.
Since  $S\cap W_{U}$ is a regular subalgebra of $W(m,n)$, 
its image under an automorphism of $W(m,n)$ is again a regular subalgebra of $W(m,n)$.

The same argument holds if we replace $S'(m,n)$ by $S(m,n)$ or by $CS'(m,n)$.
 \hfill$\Box$

\begin{remark}\label{outerS(1,2)}\em
We recall that $Der~S(1,2)=S(1,2)+\mathfrak{a}$ with
$\mathfrak{a}\cong sl_2$ 
(cf.\ Proposition \ref{derivations}). Let us denote by $e, f, h$
the standard basis of $\mathfrak{a}\cong sl_2$ defined in \cite[Lemma 5.9]{FK}.
Let $S=\prod_{j\geq -2}S_j$ denote
 the Lie superalgebra $S(1,2)$ with respect to the grading of type
 $(2|1,1)$. Then $S_0\cong gl_2$ and $S_{-1}$ is isomorphic, as an
 $S_0$-module, to the direct sum of two copies of the standard
 $gl_2$-module. It follows that, for every irreducible $gl_2$-submodule
 $U$ of $S_{-1}$, $S_U:=U+\prod_{j\geq 0}S_j$ is a maximal open
 subalgebra of $S$. In particular, if $U=\langle
 \xi_i\frac{\partial}{\partial x} ~|~ i=1,2\rangle$ or $U=\langle
 \frac{\partial}{\partial \xi_i}~|~ i=1,2\rangle$, then $S_U$ is the maximal
 graded subalgebra of type $(1|1,1)$ or $(1|0,0)$, respectively.
The subalgebras $S_U$ are not conjugate by inner automorphisms of $S$,
but they are conjugate by inner automorphisms of $Der S$,
since the subalgebra $\mathfrak{a}$ of outer derivations of $S$ permutes the subspaces $U$. In particular 
the graded subalgebras of principal and subprincipal type 
are conjugate by the (outer) automorphism $\exp(e)\exp(-f)\exp(e)\in G$.
\end{remark}

\begin{theorem}\label{S(m,n)} $(a)$ Let $S=S(m,n)$  with $m>1$ or $m=1$ and $n\geq 3$.
Then all
maximal open subalgebras of $S$ are, up to conjugation, the
graded subalgebras of type
$(1,\dots,1|1,\dots,1,0,\dots,0)$ with $k$ zeros, for  
$k=0,\dots, n$.

$(b)$ All maximal open subalgebras of $S(1,2)$ are, up to conjugation, the graded
subalgebras of type $(1|1,1)$ and $(1|1,0)$.
\end{theorem}
{\bf Proof.} Let $L_0$ be a maximal open  subalgebra of $S$. Then, by
Theorem \ref{newS}, $L_0$ is regular and 
 we can assume, by Remark \ref{standard}, that it is invariant with
respect to the standard torus $T$ of
$W(m,n)$. In particular $L_0$ decomposes into the direct product of weight spaces
with respect to $T$. Note that $\mathbb{C}\frac{\partial}{\partial
  x_i}$, $\mathbb{C}\frac{\partial}{\partial
  \xi_i}$, $\mathbb{C}\xi_{j_1}\dots\xi_{j_h}\frac{\partial}{\partial
  x_i}$, $\mathbb{C}\xi_{j_1}\dots\xi_{j_h}\frac{\partial}{\partial
  \xi_k}$ with $k\neq j_1, \dots, j_h$, are one-dimensional weight spaces. Besides,
the vector fields $\frac{\partial}{\partial
  x_i}$
cannot lie in
$L_0$ since they are not exponentiable (cf.\ \cite[Lemma
1.2]{K}). 
We may thus assume that one of the following situations occurs:
\begin{enumerate}
\item[1)] no element $\frac{\partial}{\partial \xi_i}$ lies in $L_0$. It
follows that the $T$-invariant complement of $L_0$ contains the 
$T$-invariant complement of the
subalgebra of type $(1,\dots, 1|1,\dots,1)$, thus $L_0$ coincides
with the
subalgebra of type $(1,\dots, 1|1,\dots,1)$, since it is maximal;
\item[2)] the elements $\frac{\partial}{\partial \xi_{k+1}}, \dots, \frac{\partial}{\partial
\xi_n}$ lie in $L_0$ for some $k=0,\dots,n-1$, and
$\frac{\partial}{\partial \xi_1}, \dots, \frac{\partial}{\partial
\xi_k}$ do not. Then the elements $\xi_i\frac{\partial}{\partial
x_j}$ and $\xi_i\frac{\partial}{\partial\xi_h}$ cannot lie in $L_0$
for any $j=1,\dots, m$, any $i=k+1,\dots, n$ and any $h=1,\dots, k$,
since $\left[\frac{\partial}{\partial \xi_i}, \xi_i\frac{\partial}{\partial x_j}
\right]=\frac{\partial}{\partial x_j}$ and $\left[\frac{\partial}{\partial \xi_i}, \xi_i\frac{\partial}{\partial \xi_h}
\right]=\frac{\partial}{\partial \xi_h}$. Similarly, the elements
$P\frac{\partial}{\partial x_j}$ and
$P\frac{\partial}{\partial \xi_h}$, with $P\in\Lambda(\xi_{k+1},\dots,
\xi_n)$, cannot lie
in $L_0$ for any $j=1,\dots, m~$ and any $h=1,\dots, k$.
It follows that $L_0$ is contained
in the graded subalgebra of $S$ of type
$(1,\dots,1|1,\dots,1,0,\dots,0)$ with $n-k$ zeros and thus coincides
with it since $L_0$ is maximal.
\end{enumerate}
By Remark \ref{outerS(1,2)}, when $m=1$ and $n=2$, the subalgebras
of principal and subprincipal type
 are conjugate by an element of $G$.
\hfill$\Box$

\begin{corollary} $(a)$ All irreducible $\Z$-gradings of $W(m,n)$ with $m\geq 1$, and of $S(m,n)$
with $m>1$ or $m=1$ and $n\geq 3$, are, up to conjugation, the gradings of type
$(1,\dots, 1|1,\dots,1,0,\dots,0)$ with $k$ zeros, for  
$k=0,\dots, n$. 

$(b)$ All irreducible $\Z$-gradings of $S(1,2)$ are, up to conjugation, the
gradings of type $(1|1,1)$ and $(1|1,0)$.
\end{corollary}

%


\begin{theorem} Let $S=S(m,n)$ with $m>1$, so that $S(m,n)=S'(m,n)$
and  $Der S=CS'(m,n)=S(m,n)+\C\sum_{i=1}^mx_i\frac{\partial}{\partial x_i}$.
 Then
all maximal among open 
$\sum_{i=1}^mx_i\frac{\partial}{\partial x_i}$-invariant subalgebras
of $S$ are, up to conjugation, the subalgebras of $S$ listed in
Theorem \ref{S(m,n)} (a).
\end{theorem}
{\bf Proof.} Let $L_0$ be a maximal among open
$\sum_{i=1}^m x_i\frac{\partial}{\partial x_i}$-invariant subalgebras
of $S$. Then $L_0+\C\sum_{i=1}^m x_i\frac{\partial}{\partial x_i}$
is a maximal open subalgebra of $CS'(m,n)$, hence it is regular
by Theorem \ref{newS}. Then one uses the same arguments as in the proof
of Theorem \ref{S(m,n)}. \hfill$\Box$

\bigskip

We recall that if  $L=S(1,n)$, with  $n\geq 3$, then
 $Der L=CS'(1,n)=\C E+S^\prime(1,n)$ where 
$E=x\frac{\partial}{\partial x}+\sum_{i=1}^n\xi_i\frac{\partial}
{\partial\xi_i}$ is the Euler operator and
$S^\prime(1,n)=S(1,n)+\C\xi_1\dots\xi_n\frac{\partial}{\partial x}$
(cf.\ Proposition \ref{derivations}). 
We are now interested in the subalgebras of $S(1,n)$ which are
maximal among the $\mathfrak{a}_0$-invariant subalgebras of
$S(1,n)$, for every subalgebra $\mathfrak{a}_0$ of the subalgebra
$\mathfrak{a}$ of outer derivations of $S(1,n)$
(cf.\ Theorem \ref{Theorem4}$(b)$). 

\begin{remark}\label{mfS'(1,n)}\em 
By Theorem \ref{newS} every maximal open subalgebra
of $S'(1,n)$ or $CS'(1,n)$, for every $n\geq 2$, is regular. 
Therefore the same
argument as in the proof of Theorem \ref{S(m,n)}
shows that all fundamental maximal subalgebras of $S'(1,n)$
or $CS'(1,n)$
(i.e., fundamental among maximal subalgebras) 
are, up to conjugation, the graded subalgebras 
of type $(1|1,\dots,1,0,\dots,0)$
with $k$ zeros, for $k=0,\dots,n-1$.
Indeed, the graded subalgebra of $S'(1,n)$ (resp.\ $CS'(1,n)$) of 
type $(1|0,\dots,0)$ is not  maximal, since it is contained 
in $S(1,n)$ (resp. $S(1,n)+\C E$). Notice that the graded subalgebras
of principal and subprincipal type of $S'(1,2)$ (resp.\ $CS'(1,2)$)
are not conjugate.
By the same arguments,
all maximal fundamental subalgebras of $S'(1,n)$ and $CS'(1,n)$
(i.e.,  maximal among fundamental subalgebras) are,
up to conjugation, the graded subalgebras of
type $(1|1,\dots,1,0,\dots,0)$ with $k$ zeros, for $k=0,\dots,n$.

In order to distinguish, when needed, a subalgebra 
of type 
$(a|b_1,\dots ,b_n)$ of
$S(1,n)$ from the graded subalgebra of $S'(1,n)$ or $CS'(1,n)$ of
the same type, we shall use subscripts. For example
$(1|1,\dots,1)_{S'(1,n)}$ will
denote the graded subalgebra of $S'(1,n)$ of principal type,
so that $(1|1,\dots,1)_{S'(1,n)}=(1|1,\dots,1)_{S(1,n)}+\C\xi_1\dots\xi_n
\frac{\partial}{\partial x}$.
\end{remark}

\begin{theorem}\label{derS(1,n)} Let $L=S(1,n)$ with $n\geq 3$.

\noindent
(i) All maximal among open $E$-invariant  subalgebras of $L$
are, up to conjugation,
the graded subalgebras
of type $(1|1,\dots, 1, 0,\dots, 0)$ with $k$ zeros, for some $k=0,\dots, n$.

\noindent
(ii) If $\mathfrak{a}_0=\C\xi_1\dots\xi_n\frac{\partial}{\partial x}$
or $\mathfrak{a}_0=\mathfrak{a}$,
then all maximal among open $\mathfrak{a}_0$-invariant subalgebras of $L$ 
are, up to conjugation,
the graded subalgebras
of type $(1|1,\dots, 1,$ $0,\dots, 0)$ with $k$ zeros, for some $k=0,\dots, n-1$.
%
\end{theorem} 
{\bf Proof.} Let $L_0$ be a maximal among open $E$-invariant subalgebras of
$S(1,n)$. Then $L_0+\C E$ is a fundamental subalgebra of $CS'(1,n)$, hence
it is contained in a maximal fundamental subalgebra of $CS'(1,n)$,
i.e., by Remark \ref{mfS'(1,n)}, 
in a conjugate of the graded subalgebra of
$CS'(1,n)$ of type $(1|1,\dots,1,0,\dots,0)$ with $k$ zeros, for some $k=0,\dots,n$.
Suppose $L_0+\C E\subset \varphi((1|1,\dots,1)_{CS'(1,n)})=
\varphi((1|1,\dots,1)_{S(1,n)})+\C\varphi(E)+\C\varphi(\xi_1\dots\xi_n
\frac{\partial}{\partial x})$ for some inner automorphism
$\varphi$ of $CS'(1,n)$.
Since $E$ is contained in $\varphi((1|1,..,1)_{CS'(1,n)})$, 
$\varphi((1|1,..,1)_{S(1,n)})$ is an $E$-invariant subalgebra of $S(1,n)$,
hence $L_0=\varphi((1|1,..,1)_{S(1,n)})$ by maximality.
If $L_0+\C E$ is contained in a conjugate of the subalgebra
of type $(1|1,\dots,1,0,\dots,0)$ with $k$ zeros, for some $k=1,\dots,n$,
the argument is similar. 

Now let $S_0$ be a maximal among open $\xi_1\dots\xi_n\frac{\partial}{\partial x}$-invariant subalgebras of $S(1,n)$. 
Then
$S_0+\mathbb{C}\xi_1\dots\xi_n\frac{\partial}{\partial x}$ is
a fundamental maximal subalgebra of $S'(1,n)$ containing  $\xi_1\dots\xi_n\frac{\partial}{\partial x}$.
Likewise, if $S_0$ is a maximal among open $\mathfrak{a}$-invariant 
subalgebras 
 of $S(1,n)$, then
 $S_0+\mathbb{C}E+\mathbb{C}\xi_1\dots\xi_n\frac{\partial}{\partial
 x}$  is a fundamental maximal  
subalgebra of $CS'(1,n)$ containing $E$ and $\xi_1\dots\xi_n\frac{\partial}{\partial
 x}$. Then statements $(ii)$ and $(iii)$ follow from Remark \ref{mfS'(1,n)}.
\hfill$\Box$

\begin{theorem}\label{a-invS(1,2)}
Let $S=S(1,2)$ and let $\mathfrak{b}$ be the
$2$-dimensional subalgebra of $\mathfrak{a}$ spanned by $e$ and $h$.

\noindent
(i) If $\mathfrak{a}_0$ is a one-dimensional subalgebra of 
$\mathfrak{a}$, then
all maximal among open $\mathfrak{a}_0$-invariant subalgebras of $S(1,2)$
are, up to conjugation, the subalgebras of type $(1|1,1)$
and  $(1|1,0)$.

\noindent
(ii) The graded
subalgebra of type $(1|1,1)$  is, up to conjugation, the only 
maximal among open $\mathfrak{b}$-invariant 
subalgebras of $S(1,2)$, which is not
invariant with respect to $\mathfrak{a}$.

\noindent
(iii) All  maximal open among $\mathfrak{a}$-invariant subalgebras
of $S(1,2)$ are, up to conjugation, the subalgebras of type $(2|1,1)$
and  $(1|1,0)$.
\end{theorem}
{\bf Proof.} By Remark \ref{mfS'(1,n)}, the proof of $(i)$ is the same as the 
proof of $(i)$ and
$(ii)$ in
Theorem \ref{derS(1,n)}. Recall that the graded
subalgebras of principal and subprincipal type of $S(1,2)$ are
conjugate.
  
Now, using \cite[Lemma 5.9]{FK} one can check that
 the graded subalgebras  of $S(1,2)$ of type $(1|1,0)$
and $(2|1,1)$ are invariant with respect to $\mathfrak{a}$.
On the other hand,
the graded subalgebra $L_0$ of type $(1|1,1)$  is 
invariant with respect to $\mathfrak{b}$
but it is not $\mathfrak{a}$-invariant.
Indeed, one has: $\xi_i\frac{\partial}{\partial x}\in L_0$, $\frac{\partial}{\partial \xi_j}\notin L_0$ and
 $f(\xi_i\frac{\partial}{\partial x})=\pm \frac{\partial}{\partial
 \xi_j}$ with $j\neq i$.  
Let $S_0$ be a maximal among
$\mathfrak{b}$-invariant subalgebras of $S(1,2)$. Then
$S_0+\mathbb{C}\sum_{i=1}^2\xi_i\frac{\partial}{\partial\xi_i}
+\mathbb{C}\xi_1\xi_2\frac{\partial}{\partial x}$ is a fundamental
 maximal subalgebra of $CS'(1,2)$ containing
$\sum_{i=1}^2\xi_i\frac{\partial}{\partial\xi_i}$ and $\xi_1\xi_2\frac{\partial}{\partial x}$,
hence, by Remark \ref{mfS'(1,n)}, $S_0$ is
conjugate either to the graded
subalgebra of $S(1,2)$ of type $(1|1,1)$ or to the graded subalgebra of type
$(1|1,0)$.


Now suppose that  $\tilde{S}$ is a maximal among open $\mathfrak{a}$-invariant subalgebras
of $S(1,2)$. Then $\tilde{S}$ is invariant with respect to  
$\mathfrak{b}$, hence $\tilde{S}+\mathfrak{b}$ is contained 
in a maximal fundamental subalgebra of $CS'(1,2)$
containing $\mathfrak{b}$.
It follows that either $\tilde{S}$ is contained in a conjugate of the subalgebra of $S(1,2)$
of type $(1|1,0)$, thus coincides with it by maximality, or
it is contained in a conjugate $S_U$ of the subalgebra of principal type.
As we noticed in Remark \ref{outerS(1,2)}, $S_U$ is conjugate to the subalgebra
of principal type by an automorphism $\varphi=\exp(ad~ a)$ 
for some $a\in\mathfrak{a}$. Since $\tilde{S}$ is $\mathfrak{a}$-invariant,
$\varphi(\tilde{S})=\tilde{S}$, therefore
$\tilde{S}$ is contained in the intersection of $S_U$ with the subalgebra
of principal type, i.e., in the graded subalgebra of type $(2|1,1)$. Since
the subalgebra of type $(2|1,1)$ is $\mathfrak{a}$-invariant,
$\tilde{S}$ coincides with it by maximality. 
 \hfill$\Box$

\bigskip

\noindent 
{\bf {\em The Lie superalgebra $\boldsymbol{K(2k+1,n)}$}.} Let $k\in\Z_+$  
and let $t$, $p_1,\dots, p_k$, 
$q_1,\dots, q_k$, be $2k+1$ even indeterminates and $\xi_1,\dots,\xi_n$  be 
$n$ odd indeterminates. 
Consider the differential form  
$\tau=dt+\sum_{i=1}^k(p_idq_i-q_idp_i)+\sum_{j=1}^n\xi_jd\xi_{n-j+1}$.
The contact Lie superalgebra is defined as follows (\cite[Example 4.4]{K}):
$$K(2k+1,n):=\{X\in W(2k+1,n)~|~X\tau=f\tau ~{\mbox{for some}}~
f\in\Lambda(2k+1,n)\}.$$ 
The algebra $K(2k+1,n)$ is simple for every $k,n$. Recall that we  assumed $(k,n)\neq (0,2)$.  

Consider the Lie superalgebra $\Lambda(2k+1,n)$ with the
following bracket: 
\begin{equation}
\left[f,g\right]=(2-E)f\frac{\partial g}{\partial t}-\frac{\partial
f}{\partial t}(2-E)g+\!\sum_{i=1}^k(\frac{\partial f}{\partial
p_i}\frac{\partial g}{\partial q_i}-\frac{\partial f}{\partial
q_i}\frac{\partial g}{\partial
p_i})+(-1)^{p(f)}\!\!\sum_{i=1}^n\frac{\partial f}{\partial
\xi_i}\frac{\partial g}{\partial \xi_{n-i+1}},
\label{Leibniz!!}
\end{equation}
where $E=\sum_{i=1}^k(p_i\frac{\partial}{\partial
p_i}+q_i\frac{\partial}{\partial q_i})+\sum_{i=1}^n\xi_i
\frac{\partial}{\partial \xi_i}$ is the Euler operator.
Then the map $\varphi: \Lambda(2k+1,n)\longrightarrow K(2k+1,n)$ given by:
$$f\longmapsto X_f=(2-E)f\frac{\partial }{\partial t}+\frac{\partial
f}{\partial t}E+\sum_{i=1}^k(\frac{\partial f}{\partial
p_i}\frac{\partial}{\partial q_i}-\frac{\partial f}{\partial
q_i}\frac{\partial}{\partial
p_i})+(-1)^{p(f)}\sum_{i=1}^n\frac{\partial f}{\partial
\xi_i}\frac{\partial}{\partial \xi_{n-i+1}}$$
is an isomorphism of Lie superalgebras (cf.\ \cite[\S 1.2]{CK}, \cite{CCK}).
We will therefore identify
$K(2k+1,n)$ with $\Lambda(2k+1,n)$. Fix the maximal torus
$T=\langle t, p_iq_i, \xi_j\xi_{n-j+1} ~|~ i=1,\dots, k, ~j=1,\dots,
[n/2]\rangle$.

\begin{remark}\label{Leibniz1}\em Bracket (\ref{Leibniz!!}) satisfies the
following rule:
$$[f,gh]=[f,g]h+(-1)^{p(f)p(g)}g[f,h]+2\frac{\partial f}{\partial t}gh.$$
Besides
we have: $X_f(g)=[f,g]+2\frac{\partial f}{\partial t}g$.
It follows, in particular, that an ideal $I=(f_1, \dots, f_r)$ of
$\Lambda(2k+1,n)$ is stabilized by a function $f$ in $K(2k+1,n)$ if
and only if $[f,f_i]$ lies in $I$ for every $i=1,\dots, r$.

Notice that, if $f$ is an even function independent of $t$, then $\varphi=\exp ad(f)$ is an
automorphism of $\Lambda(2k+1,n)$ with respect to both the Lie bracket
and the usual product of polynomials.
It follows that a subalgebra $L_0$ of $K(2k+1,n)$ stabilizes an ideal
$I=(f_1,\dots, f_r)$ of $\Lambda(2k+1,n)$ if and only if the subalgebra $\varphi(L_0)$ stabilizes
the ideal $J=(\varphi(f_1), \dots, \varphi(f_r))$.
\end{remark}

\begin{remark}\label{gradingsofK}\em
A $\Z$-grading of $W(2k+1,n)$ induces a $\Z$-grading on $K(2k+1,n)$
if and only if the differential form $\tau$ is homogeneous.
It follows that,
for every $s=0,\dots, [n/2]$, the $\Z$-grading of
$W(2k+1,n)$ of type $(2,1,\dots,1|2,\dots,
2,1,\dots,1,$ $0,\dots, 0)$, with $s+1$ 2's and $s$ zeros,  induces on
$K(2k+1,n)$ a
$\Z$-grading of depth 2, 
where $\g_0\cong
cspo(2k,n-2s)\otimes\Lambda(s)+W(0,s)$, $\g_{-1}\cong
\C^{2k|n-2s}\otimes\Lambda(s)$ and $\g_{-2}=[\g_{-1}, \g_{-1}]\cong\C\otimes\Lambda(s)$. 
This grading is thus irreducible for every $s=0,\dots, [n/2]$ when
$k=0$ and $n$ is odd, or $k>0$, and it is irreducible for every
$0\leq s <\frac{n-2}{2}$ when 
$k=0$ and $n$ is even. Besides, when $k=0$ and $n$ is
even the grading of type $(1|1,\dots,1,0,\dots,0)$ with
$n/2$ zeros,
is also irreducible by Remark \ref{gd}. One can verify that these
irreducible gradings satisfy the hypotheses of Proposition
\ref{basic}$(b)$, therefore the corresponding graded subalgebras of
$K(2k+1,n)$ are maximal. 
\end{remark}

The grading of $K(2k+1,n)$ of type 
$(2,1\dots,1|1,\dots,1)$
is called {\em principal}. The grading of $K(2k+1,2h)$ of type
$(2,1\dots,1|2,\dots,2,0,\dots,0)$, with $h$ zeros, is called {\em subprincipal}.

\begin{remark}\label{K(1,n)}\em Notice that when $k=0$, $n$ is even and
  $s=\frac{n-2}{2}$, then the grading of
$W(1,n)$ of type $(2|2,\dots,
2,1,\dots,1,0,\dots, 0)$, with $s+1$ 2's and $s$ zeros,   induces on
$K(1,n)$ a grading which is not irreducible. In particular, 
the subalgebra $\prod_{j\geq 0}K(1,n)_j$ of $K(1,n)$
  corresponding to this grading is contained in the maximal subalgebra
  of type $(1|1,\dots,1,0,\dots,0)$ with $n/2$ zeros.
\end{remark}
\begin{remark}\label{notconjofK}\em The group of inner automorphisms 
that preserve the principal grading of
$L=K(2k+1,n)$ is isomorphic to $\mathbb{C}^{\times}(Sp(2k)\times SO(n))$. It follows that
when $k=0$ and $n$ is even the graded subalgebras of $L$ of type
$(1|1,\dots,1,0,\dots,0)$ and $(1|1,\dots,1,0,1,0\dots,0)$ with
$n/2$ zeros, are not conjugate by an inner automorphism of $L$.
Likewise, when $k>0$ and $n$ is even the subalgebra of $L$ of subprincipal
type is not conjugate by an inner automorphism
 to the subalgebra of type $(2,1,\dots,1|2,\dots,2,
0,$ $2,0,\dots,0)$ with $n/2$ zeros. 
\end{remark}

\begin{remark}\label{valuation1}\em One can define a valuation
$\nu$ on $\Lambda(m,n)$ (and the induced valuation on $\Lambda(m,n)/\C 1$) with values in 
$\Z$, by assigning the values of $\nu$ on the generators $\{x_i, \xi_j
~|~ i=1, \dots, m; j=1,\dots, n\}$ of $\Lambda(m,n)$ as an associative
algebra, and by extending $\nu$ to $\Lambda(m,n)$ through the usual two rules:
\begin{enumerate}
\item[$a)$] $\nu(f\cdot g)=\nu(f)+\nu(g)$;
\item[$b)$] $\nu(f+g)=\min(\nu(f), \nu(g))$.
\end{enumerate}
\end{remark}
\begin{example}\label{step2}\em
Consider the
symmetric bilinear form $(\xi_i,\xi_j)=\delta_{i, n-j+1}$
on the vector space $V=\langle \xi_1,\dots, \xi_n\rangle$.
Given a subspace $U$ of $V$ we shall denote by
$U^0$ 
the kernel of the restriction of the bilinear form $(\cdot, \cdot)$
to $U$. Then $U=U^0\oplus U^1$ where $U^1$ is a maximal
subspace of $U$ with non-degenerate metric.
Let $(U^1)^\perp$ be the orthogonal complement of $U^1$ in $V$.
Then $(U^1)^\perp$ contains $U^0$
and a subspace $(U^0)'$ non-degenerately paired with $U^0$. Let
us denote by $(U^1)'$ the orthogonal complement of $U^0+(U^0)'$ in
$(U^1)^\perp$.

Now suppose that $U$ is a coisotropic subspace of $V$ and
consider the ideal 
$I_{\cal U}=(t, p_1,\dots, p_k, q_1, \dots, q_k, U)$ of $\Lambda(2k+1,n)$.
 We define 
a valuation $\nu$ on $\Lambda(2k+1,n)$ as follows:
$$\nu(t)=2; ~~\nu(p_i)=\nu(q_i)=1;$$
$$\nu(x)=1 ~\mbox{for}~ x\in U^1;$$
$$\nu(x)=2 ~\mbox{for}~ x\in U^0; ~~\nu(x)=0 ~\mbox{for}~ x\in (U^0)^\prime.$$

Then the subspaces 
$$L_j(U)=\{x\in\Lambda(2k,n) ~|~ \nu(x)\geq j+2\}$$
define a filtration of $K(2k+1,n)$ where
$L_{-1}=I_{\cal U}$ and $L_0=Stab(I_{\cal U})$. 
If $n$ is not even or $n$ is even and $\dim(U^0)<n/2$, 
this is in fact the graded filtration of
 $K(2k+1,n)$ associated, up to conjugation, to
the grading of type $(2,1,\dots,1|2,\dots,2,1,\dots,1,
0,\dots,0)$ with $s+1$ 2's and $s$ 0's, $s$ being the dimension
of $U^0$. If $n$ is even, $k=0$ (resp.\ $k>0$), and $U=U^0$ is a maximal isotropic subspace of
$V$, then $L_0$ is conjugate either to the graded subalgebra of $L$
of type $(1|1,\dots,1,0,\dots,0)$ (resp.\ $(2,1,\dots,1|2,\dots,2,0,\dots,0)$)
or to the graded subalgebra of type $(1|1,\dots,1,0, 1,0,\dots,0)$
(resp.\ $(2,1,\dots,1|2,\dots,2,0,2,0\dots,0)$) with $n/2$ zeros. 
\end{example}

\begin{remark}\label{2-1}\em  If $k=0$ and $n=2h$ then the maximal graded 
subalgebra 
  of $K(1,n)$ of type $(2|2,\dots,2,0,\dots,0)$ is not irreducible
since its component of degree $-1$ does not generate its negative part. 
Notice that the non-negative part of the irreducible grading of type
 $(1|1,\dots,1,0,\dots,0)$ with $h$ zeros coincides
with the non-negative part of the grading of type $(2|2,\dots,2,0,\dots,0)$.
In particular, it
  stabilizes the ideal $I_{\cal{U}}=(t, p_1, \dots,$ $p_k, q_1, \dots, q_k, U^0)$
  where $U^0=\langle \xi_1, \dots, \xi_h\rangle$. 
\end{remark}

\begin{lemma}\label{K1} In the associative superalgebra $\Lambda(2k+1,n)$, 
let us consider an ideal $J=(t, p_1,\dots, p_k, q_1,\dots, q_k, h_1,...,h_r)$
where $h_1, \dots, h_r\in \Lambda(0,n)$.
Suppose that 
$h_1=\eta_1+F$ and
$h_2=\eta^{\prime}_1+G$
where $\eta_1, \eta^{\prime}_1$  are non-degenerately paired, distinct  elements
of $V=\langle\xi_1, \dots,\xi_n\rangle$ and $F, G$ contain no constant and linear terms. Then $J$ is
conjugate
 to the ideal $K=(t+T, p_1,\dots, p_k, q_1,\dots, q_k, \eta_1, \eta^{\prime}_1,$ $f_1,...,f_{r-2})$
for some functions $T, f_i\in\Lambda(U)$ where $U$ is the orthogonal 
complement of $\langle \eta_1, \eta^{\prime}_1\rangle$ in $V$.
\end{lemma}
{\bf Proof.} Multiplying $h_1$  by some invertible function
 we can assume that $F$ does not depend on $\eta_1$, i.e., 
$\eta_1+F=\eta_1+f_1\eta^{\prime}_1+f_2$ where $f_1, f_2$ lie in $\Lambda(U)$.
Also, we can assume that $G$ lies in $\Lambda(U_1)$ where
$U_1=\langle U, \eta_1\rangle$.
Notice that $f_1\eta^{\prime}_1+f_1G$ lies in $J$, therefore
$J=(t, p_1,\dots, p_k, q_1,\dots, q_k, \eta_1+f_2-f_1G, \eta^{\prime}_1+G, h_3,\dots, h_r)$ where
$f_2-f_1G\in\Lambda(U_1)$. Therefore, multiplying
$\eta_1+f_2-f_1G$ by an invertible function, we can
write $J=(t, p_1,\dots, p_k, q_1,\dots, q_k, \eta_1+F', \eta^{\prime}_1+G, h_3,\dots, h_r)$ where
$F'\in\Lambda(U)$.

Now (see Remark \ref{Leibniz1}) the automorphism $\exp ad(\eta^{\prime}_1F')$
maps $J$ to
the ideal $J'=(t+T_1,$ $p_1,\dots, p_k, q_1,\dots, q_k,\eta_1, \eta^{\prime}_1+H, h'_3,...,h'_r)$ where
$T_1$ and  the functions
$h'_i$'s lie in $\Lambda(0,n)$, and 
$H\in\Lambda(U)$. Then, similarly as above, 
the automorphism $\exp ad(\eta_1H)$, maps $J'$ to the ideal
$K=(t+T_2, p_1,\dots, p_k, q_1,\dots, q_k, \eta_1, \eta^{\prime}_1,$ 
$f_1,...,
f_{r-2})$,
for some $T_2, f_1, \dots f_{r-2}\in\Lambda(0,n)$. Since $\eta_1, \eta^{\prime}_1$ lie in $K$,
it is immediate to see that we can assume $T_2, f_1,...,f_{r-2}\in\Lambda(U)$.
\hfill$\Box$

\begin{lemma}\label{K1+}
In the associative superalgebra $\Lambda(2k+1,n)$, 
let us consider an ideal $J=(t, p_1,\dots, p_k, q_1,\dots, q_k,  h_1,...,h_r)$
where $h_1, \dots, h_r\in \Lambda(0,n)$.
Suppose that 
$h_1=\eta_1+F$ 
where $\eta_1$  is an element of $V$ non-degenerately paired with itself,
 and $F$ contains no constant and linear terms. Then $J$ is
conjugate
 to the ideal $K=(t+T, p_1,\dots, p_k, q_1,\dots, q_k, \eta_1, f_1,...,f_{r-1})$
for some functions $T, f_i\in\Lambda(U)$ where $U$ is the orthogonal 
complement of $\langle \eta_1\rangle$ in $V$.
\end{lemma}
{\bf Proof.} One uses the same argument as in the first part of the proof of
Lemma \ref{K1}. \hfill$\Box$

\begin{lemma}\label{max} Let $L_0$ be a maximal open subalgebra of $L=K(m,n)\cong\Lambda(m,n)$
and let $I$ be an ideal of $\Lambda(m,n)$ stabilized by $L_0$.
Suppose that $I$ is maximal among the $L_0$-invariant ideals.
Then $L_0\subset I$.
\end{lemma}
{\bf Proof.} Let $(L_0)$ be the ideal generated by $L_0$. 
Every invertible element of $\Lambda(m,n)$ is not exponentiable
therefore $(L_0)$ contains no invertible element of $\Lambda(m,n)$.
It follows that $(L_0)+I$ is a proper ideal of $\Lambda(m,n)$ containing
$I$, and it is $L_0$-invariant. By the maximality of $I$ it follows $L_0\subset I$. ~$\Box$

\begin{lemma}\label{victorshelp} Let $J=(t+f_0, p_1+f_1, q_1+h_1, \dots, 
p_k+f_k, q_k+h_k)$
 be an ideal of $\Lambda(2k+1,n)$,
for some even functions $f_i, h_j$ containing
no linear and constant terms.
Then $J=(t+f'_0, p_1+f'_1, q_1+h'_1, \dots, p_k+f'_k, q_k+h'_k)$
with $f'_i, h'_j$ in $\Lambda(0,n)$. 
\end{lemma}
{\bf Proof.} Suppose 
$f_0=t+t\phi_1+\phi_2$ with $\phi_2$ independent of $t$ and 
$n_2=\deg(\phi_2)>1.$ Then 
$f_0-f_0\phi_1=
t-t\phi_1^2-\phi_2\phi_1+\phi_2$. Then the coefficients
 of $t$ in the second and in the third term have degree $2n_1$ and 
$n_1+n_2-1>n_1$, respectively, where $n_1=\deg\phi_1$.
Hence in the limit we get $t+\psi$ for some function $\psi$ independent of $t$.
Similarly we can make $f_0$, and $f_j, h_j$ 
independent of all even variables for every $j=1,\dots,k$. \hfill$\Box$

\begin{theorem}\label{standardforK} Let $L_0$ be a maximal open subalgebra of $L=K(2k+1,n)$.
Then $L_0$ is conjugate to the standard subalgebra $L_{\cal U}$ of
$L$ stabilizing the ideal $I_{\cal U}=(t, p_1, q_1, \dots, p_k, q_k,
U)$ of $\Lambda(2k+1,n)$, for some coisotropic subspace $U$ of
$V=\langle \xi_1,\dots, \xi_n\rangle$.
\end{theorem}
{\bf Proof.} By Remark \ref{rome} $L_0$ stabilizes an ideal of the form
$$J=(t+f_0, p_1+f_1, q_1+h_1, \dots, p_k+f_k, q_k+h_k,
 \nu_1+g_1,\nu_2+g_2, \dots, \nu_s+g_s)$$
for some
linear functions $\nu_j$ in odd indeterminates, and even functions
$f_i, h_i$ and odd functions $g_j$ without constant and linear terms,
 and $J$ is maximal among the $L_0$-invariant
ideals of $\Lambda(2k+1,n)$.

By Lemma \ref{victorshelp}
we can assume  $f_0$ and, similarly, $f_i$, $h_i$
in $\Lambda(0,n)$ for every $i$. 
Therefore the automorphism $\exp(ad(f_1q_1))$  maps $J$ to  
$$J_1=(t+f'_0, p_1, q_1+h'_1, p_2+f'_2, q_2+h'_2, \dots, p_k+f'_k, q_k+h'_k, \nu_1+g'_1, \nu_2+g'_2, \dots, \nu_s+g'_s).$$
As above we can make $h'_1$ independent of even variables.
It follows that the automorphism $\exp ad(-h'_1p_1)$ maps $J_1$
to $J_2=(t+f^{\prime\prime}_0, p_1, q_1, p_2+f^{\prime\prime}_2, q_2+h^{\prime\prime}_2, \dots, p_k+f^{\prime\prime}_k, q_k+h^{\prime\prime}_k, \nu_1+g^{\prime\prime}_1, \nu_2+g^{\prime\prime}_2, \dots, \nu_s+g^{\prime\prime}_s)$.
The same procedure applied to all  generators $p_i+f^{\prime\prime}_i$ and
$q_j+h^{\prime\prime}_j$ shows that $J$ is in fact conjugate to the ideal
$$I=(t+T_0, p_1, \dots, p_k, q_1, \dots, q_k, \nu_1+\ell_1, \nu_2+\ell_2, \dots, \nu_s+\ell_s)$$
where 
$\nu_1, \dots, \nu_s$ are linearly independent vectors in $V$ and 
$T_0$, $\ell_1, \dots, \ell_s$ are functions in $\Lambda(0,n)$
without constant and  linear terms. 

Let $U=\langle \nu_1,\dots, \nu_s\rangle$. 
Then, using the notation introduced in Example \ref{step2}, by Lemmas \ref{K1}
and \ref{K1+}, 
$$I=(t+T_1, p_1,\dots, p_k, q_1, \dots, q_k, U^1, \eta_1+\ell_1,\dots, \eta_r+\ell_r)$$
where $U^0=\langle \eta_1, \dots, \eta_r\rangle$ and 
$T_1, \ell_1, \dots, \ell_r\in\Lambda((U^1)^\perp)$. 
Let
$(U^0)'=\langle \eta'_1, \dots, \eta'_r\rangle$ with $(\eta_i, \eta'_j)=\delta_{i,j}$.

Denote by
$I'$ the ideal $I'=(t+T_1,\eta_1+\ell_1,
\dots, \eta_r+\ell_r)\subset I$.
Then, each function $f$ in $L_0$ (thus stabilizing $I$)
stabilizes the ideal
$K=(I, [I',I'])$. Indeed, for every $g, h\in I'$
we have:
$$[f,[g,h]]=[[f,g],h]\pm[g,[f,h]]\in[I, I']$$
and $[I, I']\subset K$ since $T_1$ and all odd generators of $I'$ are orthogonal
to $U^1$.
Notice that $K$ is generated by the generators of $I$ and by
the commutators between every pair of generators of $I'$.
Therefore $K$ is a proper ideal of
$\Lambda(2k+1,n)$ since among its generators there is
no invertible element. By the maximality of
$I$ among the ideals stabilized by $L_0$ we have 
$I=K$.

Let us rewrite the ideal $I$ as follows:
$$I=(t+T_1, p_1, q_1, \dots, p_k, q_k, U^1,\eta_1,\dots, \eta_{h-1},
\eta_h+\ell_h,\dots, \eta_r+\ell_r)$$ where $h=\min\{i=1,\dots, r ~|~ \ell_i\neq 0\}$.

We first show that the functions $\ell_h$ can be made 
independent of  $\eta'_1,\dots, \eta'_{h-1}$.
Indeed, let
$\eta_h+\ell_h=\eta_h+\eta'_1\phi_1+\phi_2$ where 
$\phi_1, \phi_2$ do not depend on $\eta'_1$.
Then $\phi_1=[\eta_1, \eta_h+\ell_h]\in [I^{\prime}, I^{\prime}]\subset
 K=I$, thus 
$I=(t+T_1, p_1, q_1, \dots, p_k, q_k, U^1,\eta_1,\dots, \eta_{h-1},$ $\eta_h+\phi_2, \eta_{h+1}+\ell^{\prime}_{h+1}, \dots, \eta_r+\ell^{\prime}_r)$,
where $\phi_2\in \Lambda((U^1)^\perp)$ does not
depend on 
$\eta'_1$. Arguing in the same way for the
variables $\eta'_2, \dots, \eta'_{h-1}$, we get
$$I=(t+T_1, p_1, q_1, \dots, p_k, q_k, U^1,\eta_1,\dots, \eta_{h-1}, \eta_h+\phi, \eta_{h+1}+\ell_{h+1}, \dots, \eta_r+\ell_r)$$
where $\phi$ does not depend on $\eta'_1,\dots, \eta'_{h-1}$.
Besides, multiplying $\eta_h+\phi$ by an invertible function, we can
assume that $\phi$ does not depend on $\eta_h$.
Now we can write $\phi=\eta'_h\psi_1+\psi_2$ with $\psi_1, \psi_2$
independent  of $\eta'_1,\dots, \eta'_{h}$.
Therefore, applying  
the automorphism $\exp(ad(\psi_2\eta'_h))$ we can write
$$I=(t+T_2, p_1, q_1, \dots, p_k, q_k, U^1,\eta_1,\dots, \eta_{h-1}, \eta_h+\eta'_h\psi_1, \eta_{h+1}+\ell^{\prime}_{h+1}, \dots, \eta_r+\ell^{\prime}_r)$$
for some $T_2, \ell^{\prime}_j\in \Lambda(0,n)$.
  Then $\psi_1=1/2[\eta_h+\eta'_h\psi_1, \eta_h+\eta'_h\psi_1]\in [I^{\prime},I^{\prime}]\subset K=I$.
Therefore, up to conjugation,
$$I=(t+T_2, p_1, q_1, \dots, p_k, q_k, U^1, \eta_1,\dots, \eta_{h-1}, \eta_h,
\eta_{h+1}+\ell^{\prime}_{h+1},\dots, \eta_r+\ell^{\prime}_r).$$
Arguing as above for $\ell^{\prime}_{h+1}, \dots, \ell^{\prime}_s$, we 
can assume, up to conjugation, that $I$ has the following form:
$$I=(t+f, p_1, q_1, \dots, p_k, q_k, U^1, \eta_1,\dots, \eta_r)=
(t+T, p_1, q_1, \dots, p_k, q_k, U)$$
for some function $f$ in $\Lambda(0,n)$. Notice that, in fact, we can assume
$f\in \Lambda((U^0)'\oplus (U^1)')$, since $U=U^0\oplus U^1\subset I$.
Now suppose $f=\eta'_1\varphi_1+\varphi_2$ with $\varphi_1$, $\varphi_2$
independent of $\eta'_1$. Then $[t+f,\eta_1]=-\eta_1+\varphi_1\in
[I',I']\subset K=I$, thus we can replace $t+f$ with $t+\varphi_2$
(here $I'=(t+f, U^0)$).
Similarly, we can make $f$ independent of $\eta'_2, \dots, \eta'_r$,
i.e., $f\in \Lambda((U^1)')$ with no linear and
constant terms.
In particular, if $U$ is coisotropic then $f=0$ and $I$ is a standard
ideal.

Suppose that $U$ is not coisotropic and
consider the ideal $Y=(t, p_1, \dots, p_k,$ $q_1, \dots, q_k,
U+(U^1)')$. Note that if $U$ is not coisotropic then the subspace
$U+(U^1)'$ is coisotropic.
Let
$L'_{-2}\supset L'_{-1}\supset\dots$ be the filtration of $K(2k+1,n)$
associated to the ideal $Y$
as in Example \ref{step2}, where $L'_0=Stab(Y)$.
Then the completion of the  graded superalgebra
$Gr L$ associated to this filtration
is isomorphic to $K(2k+1,n)$ with respect to the grading
of type $(2|2,\dots,2,1,\dots,1,0,\dots,0)$ with $s+1$ 2's and $s$ 0's,
where $s=\dim(U^0)$.
In particular we have:
$$Gr_{-2}(L)=\Lambda((U^0)'), 
~~Gr_{-1}(L)=(\langle p_i, q_i\rangle \oplus U^1\oplus (U^1)')\otimes\Lambda((U^0)').$$

\medskip

We want to show that $L_0$ is contained in $L'_0$.
Suppose that $X\in K(2k+1,n)$ stabilizes $I$. Then we can write
$$X=X_{-2}+X_{-1}+X_0$$ with $X_{-2}\in Gr_{-2}L$,
$X_{-1}\in Gr_{-1}L$ and $X_0\in \prod_{j\geq 0}Gr_j L$. In fact, since $L_0$
is open, $X_{-2}\in\Lambda((U^0)')/\mathbb{C}$.

Note that $t\in Gr_0L$, $f\in \prod_{j\geq 0}Gr_j L$, $U^0\subset Gr_0L$,
$U^1\subset Gr_{-1}L$ and $p_i, q_i\in Gr_{-1}L$. It follows that $I\subset Gr_{\geq -1}=L'_{-1}$.

Now, since $X\in Stab(I)$, we have:
$$[X,U^0]\subset I\subset L'_{-1} \Rightarrow [X_{-2}, U^0]=0$$
and, since $X_{-2}\in \Lambda((U^0)')/\mathbb{C}$, it follows $X_{-2}=0$.

Similarly,
$$[X, U^1]\subset I \Rightarrow [X_{-1}, U^1]=0$$
and 
$$[X, p_i]\subset I \Rightarrow [X_{-1}, p_i]=0, ~~[X, q_i]\subset I \Rightarrow [X_{-1}, q_i]=0,~~~\forall i=1,\dots,k,$$
hence
$X_{-1}\in (U^1)'\otimes\Lambda((U^0)')$.
Therefore  $X=X_{-1}+X_0\in L_0$ with $X_{-1}\in (U^1)'\otimes\Lambda((U^0)')\subset Gr_{-1}L$ and $X_0\in \prod_{j\geq 0}Gr_j L$.
By Lemma \ref{max}, $L_0\subset I$, hence $X\in I$. It follows $X_{-1}=0$,
i.e., $L_0\subset \prod_{j\geq 0}Gr_j L$. Indeed if $X_{-1}\neq 0$ then $\nu(X)=1$
but $I\cap\{x~|~\nu(x)=1\}=(\langle p_i, q_i\rangle +U^1)\otimes
\Lambda((U^0)')+\prod_{j\geq 0}Gr_j L$.
By the maximality of $L_0$ the statement follows.
 ~$\Box$

\begin{theorem}\label{K(2n+1)}

$(i)$ All maximal open subalgebras of $K(1,2h)$ ($h>1$) are, up to conjugation,
the graded subalgebras of type $(1|1,\dots,1,0,\dots,0)$
and $(1|1,\dots,1,0,$ $1,0,\dots,0)$
with $h$ zeros, 
and the graded subalgebras of type $(2|2,\dots,
2,1,\dots,1,$ $0,\dots, 0)$ with $s+1$ 2's and $s$ zeros, for  
$s=0,\dots, h-2$.

\medskip
\noindent
$(ii)$ If $k>0$ and $n$ is even, all  maximal open subalgebras
 of $K(2k+1,n)$ are, up to conjugation,
 the graded subalgebras of type $(2,1,\dots,1|2,\dots,
2,1,\dots,1,$ $0,\dots,0)$ with $s+1$ 2's and $s$ zeros, for
$s=0,\dots,n/2$ and the graded subalgebra of type
$(2,1,\dots,1|2,\dots,2,0,2,0,\dots,0)$ with $n/2$ zeros.

\medskip
\noindent
$(iii)$ If $n$ is odd, all  maximal open subalgebras
 of $K(2k+1,n)$ are, up to conjugation,
 the graded subalgebras of type $(2,1,\dots,1|2,\dots,
2,1,\dots,1,0,\dots,0)$ with $s+1$ 2's and $s$ zeros, for
$s=0,\dots,\left[n/2\right]$.
\end{theorem}
{\bf Proof.} By Theorem \ref{standardforK} every maximal open
subalgebra of $K(2k+1,n)$ is conjugate to the standard subalgebra 
associated to the ideal $I_{\cal
  U}=(t, p_1, \dots, p_k,$ $q_1,\dots, q_k, U)$ of $\Lambda(2k+1,n)$, for some
coisotropic subspace $U$ of $V=\langle \xi_1, \dots,
\xi_n\rangle$. Now the statement follows from Example \ref{step2} and
 Remarks \ref{gradingsofK}, 
\ref{K(1,n)}, \ref{notconjofK} and
\ref{2-1}. \hfill$\Box$

\begin{corollary} 

$(i)$ All irreducible $\Z$-gradings of $K(1,2h)$ are, up
  to conjugation,
the grading of type $(2|2,\dots,
2,1,\dots,1,0,\dots, 0)$, with
$s+1$ 2's and $s$ zeros, for $s=0,\dots, h-2$ and 
the gradings of type $(1|1,\dots,1,0,\dots,0)$ and
$(1|1,\dots,1,0,1,$ $0,\dots,0)$ 
with $h$ zeros.

\medskip
\noindent
$(ii)$ All irreducible $\Z$-gradings of $K(2k+1,n)$ where $k>0$ and $n$ is even
are, up to conjugation,
the gradings of type $(2,1,\dots,1|2,\dots,
2,1,\dots,1,$ $0,\dots, 0)$ with
$s+1$ 2's and $s$ zeros, for 
$s=0,\dots,n/2$ and the grading of type
$(2,1,\dots,1|2,\dots,2,$ $0,2,0,\dots,0)$ with $n/2$ zeros. 

\medskip
\noindent
$(iii)$ All irreducible $\Z$-gradings of $K(2k+1,n)$ where $n$ is odd 
are, up to conjugation,
the gradings of type $(2,1,\dots,1|2,\dots,
2,1,\dots,1,0,\dots, 0)$ with
$s+1$ 2's and $s$ zeros, for 
$s=0,\dots,\left[n/2\right]$. 
\end{corollary}

We take the opportunity here to describe the embedding of the Lie superalgebra
$S(1,2)$ in $K(1,4)$ and to correct Proposition 4.1.2
in \cite{CK}.

\begin{remark}\em Consider the Lie superalgebra
 $K(1,4)$ with its principal grading. Then $\mathfrak{g}_0=cso_4$ and we want to study
${\mathfrak g}_1$ as a ${\mathfrak g}_0$-module. 
${\mathfrak g}_1$ is spanned by the elements $t\xi_i$, for $i=1,\dots, 4$,
and  $\xi_i\xi_j\xi_k$ for $i,j,k=1, \dots, 4$, $i\neq j\neq k$,
thus it is the direct sum of two isomorphic irreducible
 representations  
of $so_4$: $V=\langle t\xi_i\rangle$ and $W=\langle\xi_i\xi_j\xi_k\rangle$, 
each of which is isomorphic to
the standard $so_4$-module.
Note that 
$[W,W]=\langle\xi_1\xi_2\xi_3\xi_4\rangle$ and $[{\mathfrak g}_{-2}, W]=0$. 
It follows that  ${\mathfrak g}_{-2}+{\mathfrak g}_{-1}+{\mathfrak g}_0+W+[W,W]$  is isomorphic to
$\hat{H}(0,4)+Ct$ where $t$ is the grading operator (see \cite[\S 1.2]{CK}).
As it was noticed in \cite[Proposition 4.1.2]{CK},
for every $\lambda\in\mathbb{C}$, the subspace
$V_{\lambda}=\langle\xi_1t+\lambda\xi_1\xi_2\xi_3, 
\xi_2t+\lambda\xi_2\xi_1\xi_4, \xi_3t+\lambda\xi_4\xi_1\xi_3,
\xi_4t-\lambda\xi_4\xi_2\xi_3\rangle$ is an irreducible 
${\mathfrak g}_0$-submodule of
${\mathfrak g}_1$, but, 
differently from what is claimed in \cite[Proposition 4.1.2]{CK},
$\dim ([V_{\lambda}, V_{\lambda}])= 1$ for every $\lambda\in\mathbb{C}$.
Besides, for every $\lambda\neq \pm 1$,
$[{\mathfrak g}_{-1}, V_{\lambda}]=cso_4={\mathfrak g}_0$, while, if $\lambda=1$ 
or
$\lambda=-1$, then $[{\mathfrak g}_{-1}, V_{\lambda}]=gl_2$.
Therefore, for every $\lambda\neq \pm 1$,
${\mathfrak g}_{-2}+{\mathfrak g}_{-1}+{\mathfrak g}_0+V_{\lambda}+[V_{\lambda}, V_{\lambda}]$
is a simple, 17-dimensional Lie superalgebra,
isomorphic to the Lie superalgebra $D(2,1;\alpha)$ for some $\alpha$
(cf.\ \cite[Remark 2.5.7]{K2}). 
If $\lambda=1$ or $\lambda=-1$, the Lie superalgebra
$L:={\mathfrak g}_{-2}+{\mathfrak g}_{-1}+gl_2+V_{\lambda}+V_{\lambda}^2$ has dimension 14 and it is
isomorphic to $sl(2,2)/\CC1$, and the copy of $sl_2$, lying in ${\mathfrak g}_0$
and outside of $L$, acts on $L$ by outer derivations.

Now consider the Lie superalgebra
 $S(1,2)=\sum_{j\geq -2} {\mathfrak h}_j$ with respect to the grading of type $(2|1,1)$. 
Then
 the positive part of this grading is not generated by
${\mathfrak h}_1$. On the contrary, $({\mathfrak h}_1)^2$ has dimension 1 and 
${\mathfrak h}_{-2}+{\mathfrak h}_{-1}+{\mathfrak h}_0+{\mathfrak h}_1+({\mathfrak h}_1)^2\cong sl(2,2)/\mathbb{C}1$.
We have the following
embedding of $S(1,2)$ in $K(1,4)$:
%
$$S(1,2)\cong\mathbb{C}[t+\xi_1\xi_4+\xi_2\xi_3]\Lambda(\xi_1,\xi_2)+
\mathbb{C}[t-\xi_1\xi_4-\xi_2\xi_3]\Lambda(\xi_3,\xi_4).$$
\end{remark}

\medskip

\noindent
{\bf {\em The Lie superalgebras $\boldsymbol{HO(n,n)}$ and $\boldsymbol{SHO(n,n)}$}.} Let $x_1,\dots, x_n$ be $n$ even indeterminates and $\xi_1,
\dots, \xi_n$ be $n$ odd indeterminates, and
let us consider the differential form 
$\sigma=\sum_{i=1}^ndx_id\xi_i$. The odd Hamiltonian
superalgebra is defined as follows (cf.\ \cite{ALS}):
$$HO(n,n):=\{X\in W(n,n) ~|~ X\sigma=0\}.$$
It is a simple Lie superalgebra if and only if $n\geq 2$. 
The Lie superalgebra $HO(n,n)$ contains the subalgebra
$$SHO^{\prime}(n,n):=S'HO(n,n)=\{X\in HO(n,n) ~|~ div(X)=0\}$$ 
 (see Definition \ref{divergence}
and Example \ref{S'W}).

Its derived algebra $SHO(n,n)=[SHO^{\prime}(n,n), SHO^{\prime}(n,n)]$
is simple if and only if $n\geq 3$.

The Lie superalgebra $HO(n,n)$ can be realized as follows
(cf. \cite[\S 1.3]{CK}): in $\Lambda(n,n)$ one can consider the Lie
superalgebra structure defined by  the Buttin bracket:
$$[f,g]:=\sum_{i=1}^n(\frac{\partial f}{\partial x_i}\frac{\partial g}{\partial \xi_i}+(-1)^{p(f)}\frac{\partial f}{\partial \xi_i}\frac{\partial g}{\partial x_i}).$$
Then the map $\Lambda(n,n)\rightarrow HO(n,n)$ given by:
$$f\longmapsto\sum_{i=1}^n(\frac{\partial f}{\partial x_i}\frac{\partial }{\partial \xi_i}+(-1)^{p(f)}\frac{\partial f}{\partial \xi_i}\frac{\partial }{\partial x_i})$$
defines a surjective homomorphism of Lie superalgebras whose kernel consists
of constant functions. Hence we will identify
$HO(n,n)$ with
$\Lambda(n,n)/\C 1$ with the Buttin bracket, with reversed
parity. Under this identification
$$SHO^{\prime}(n,n)=\{f\in\Lambda(n,n)~|~ \Delta(f)=0\}/\C 1=:\Lambda^{\Delta}(n,n)/\C 1,$$ where
$\Delta=\sum_{i=1}^n \frac{\partial^2}{\partial x_i\partial\xi_i}$ is
the odd Laplace operator, and $SHO(n,n)$ is the span of all monomials
in $SHO'(n,n)$ except for $\xi_1\dots\xi_n$.

Since $HO(2,2)\cong S(2,1)$ and since $SHO(n,n)$ is simple if and only if
$n\geq 3$, when talking about $HO(n,n)$ and $SHO(n,n)$ we shall assume $n\geq 3$. 
Fix the standard torus $T=\langle x_i\xi_i ~|~
i=1,\dots,n\rangle$ of $HO(n,n)$. 
Recall that $Der HO(n,n)=HO(n,n)+\C E$ 
 where $E=\sum_{i=1}^n(x_i\frac{\partial}{\partial x_i}+
\xi_i\frac{\partial}{\partial\xi_i})$ is the Euler operator. Besides,
if $n\geq 4$ 
then $Der SHO(n,n)=SHO^\prime(n,n)+\C E+\C \Phi$ where 
$\Phi=\sum_{i=1}^n x_i\xi_i$ (with $\sum_{i=1}^n(-x_i\frac{\partial}
{\partial x_i}+\xi_i\frac{\partial}{\partial\xi_i})$ the 
corresponding vector field)
(cf.\ Proposition \ref{derivations}).
Finally, $Der SHO(3,3)=SHO(3,3)+\mathfrak{a}$ where $\mathfrak{a}\cong gl_2$
and a maximal torus of $\mathfrak{a}$ is spanned by $E$ and $\Phi$
(cf.\ \cite[Remark 4.4.1]{CK}).


\begin{remark}\label{gradingsofHO}\em The
$\Z$-grading of type 
$(1,\dots,1|0,\dots,0)$ of $W(n,n)$ induces on $HO(n,n)$ (resp.\ $SHO(n,n)$) a grading
of depth 1 (called the {\em subprincipal} grading) which is irreducible
by Remark \ref{gd}. 

Consider the gradings induced on
$HO(n,n)$ (resp.\ $SHO(n,n)$) by the $\Z$-gradings of type $(1,\dots,
1, 2,\dots, 2|1,\dots,
1, 0, \dots, 0)$ of  $W(n,n)$, with
$k$ 2's and $k$ zeros. For any fixed $k$, $0\leq k\leq n-2$,
the $0$-th graded component of $HO(n,n)$ (resp.\ $SHO(n,n)$)
with respect to this grading is isomorphic to the Lie superalgebra
$\tilde{P}(n-k)\otimes\Lambda(k)+W(0,k)$ (resp.\ $P(n-k)\otimes\Lambda(k)+W(0,k)$)  and the $-1$-st graded component is
isomorphic to  $\C^{n-k|n-k}\otimes\Lambda(k)$  
where $\C^{n-k|n-k}$ is the standard $P(n-k)$-module (cf.\ \cite{K2}).
Therefore for every $k=0, \dots, n-2$ these are irreducible
gradings of $HO(n,n)$ (resp.\ $SHO(n,n)$). 
If $k>0$ then the grading of type $(1,\dots,
1, 2,\dots, 2|1,\dots,
1, 0, \dots, 0)$ with $k$ 2's and $k$ zeros, has depth 2, its $-1$-st graded component 
generates its negative part and property $(iii)^\prime$ of Proposition 
\ref{basic}$(b)$ is satisfied. It follows that the subalgebras of
$HO(n,n)$ (resp.\ $SHO(n,n)$) of type $(1,\dots,1|$ $0,\dots,0)$ and
$(1,\dots,1, 2,\dots, 2$ $|1,\dots, 1, 0, \dots, 0)$ with
$k$ 2's and $k$ zeros, for $k=0,\dots,n-2$, are maximal. 
(All claims  hold also for the Lie superalgebra $HO(2,2)$).

The $\Z$-grading induced on $HO(n,n)$ (resp. $SHO(n,n)$) by the principal
grading of $W(n,n)$ is also called {\em principal}.
\end{remark}


\begin{remark}\label{irrenotHO(n,n)}\em
The $\Z$-grading of type
$(1,2,\dots,2|1, 0,\dots, 0)$ of 
$HO(n,n)$  is not irreducible. Indeed the $0$-th graded component
of $HO(n,n)$ with respect to this grading is the subspace $\langle x_1^2, x_1\xi_1, 
x_i ~|~ i=2,\dots, n\rangle\otimes\Lambda(\xi_2, \dots, \xi_{n})$ and
its $-1$-st graded component is $\langle x_1, \xi_1\rangle 
\otimes\Lambda(\xi_2, \dots, \xi_{n})$. It follows that
the subspace $\langle x_1\rangle 
\otimes\Lambda(\xi_2, \dots, \xi_{n})$ of $HO(n,n)_{-1}$ is 
$HO(n,n)_{0}$-stable.

Likewise, the $\Z$-grading of type
$(1,2,\dots,2|1,0,\dots, 0)$ of 
$SHO(n,n)$ fails to be irreducible. Indeed, with respect to this grading,
$HO(n,n)_{-1}=SHO(n,n)_{-1}$.

The subalgebra of type $(1,2,\dots,2|1,0,\dots, 0)$ of 
$HO(n,n)$ (resp.\ $SHO(n,n)$) is contained in the maximal subalgebra
of type $(1,\dots, 1|0,\dots,0)$.
\end{remark}

%

\begin{remark}\label{stabHO}\em Let $L=HO(n,n)$ or $L=SHO(n,n)$. Then
the graded subalgebra $L_k$ of $L$ of type $(1,\dots,1, 2,\dots, 2$ 
$|1,\dots, 1, 0, \dots, 0)$ with $n-k$ 2's and $n-k$ zeros, is,
for every $k=1,\dots, n$, the standard
subalgebra $L_U$ of $L$ stabilizing the ideal 
$I_U=(x_1,\dots, x_n, \xi_1, \dots,
\xi_k)$. Indeed,  for every $k$, $L_k\subset L_U$
since $L_k$
is contained in the graded subalgebra of $W(n,n)$ of type $(1,\dots,1|1,\dots,1,0,\dots,0)$ with $n-k$ zeros, which stabilizes $I_U$ (cf.\ the proof of Theorem
\ref{W(m,n)}). If $k\neq 1$,  then, by Remark \ref{gradingsofHO}, $L_k$ is a maximal  
subalgebra of $L$, thus 
$L_k=L_U$ for every $k\neq 1$.

Now suppose $k=1$. 
By Remark \ref{regularideal}, $L_U$
is regular and, up to conjugation, we can assume that it is invariant
with respect to the standard torus $T+\C E$ of $Der HO(n,n)$.
Therefore $L_U$ decomposes into the direct
product of weight spaces with respect to $T+\C E$. 
Consider the $\Z$-grading of $L$ of type
$(1,2,\dots,2|1, 0,\dots,0)$. Then the negative part of this
grading is $\mathfrak{g}_{-}=$
$(<1, x_1, \xi_1>\otimes\Lambda(\xi_2, \dots, \xi_n))/
\mathbb{C}1$. Notice that $\C\xi_{i_1}\dots\xi_{i_h}$, with
$i_1\neq\dots\neq i_h$, and
$\C x_1\xi_{j_1}\dots \xi_{j_h}$, with $1\neq j_1\neq\dots\neq j_h$,
are one-dimensional weight spaces with respect to $T+\C E$.
Therefore,  in order to prove that $L_U$ is contained in $L_1$
(hence $L_U=L_1$)
it is sufficient to show that, for every $f\in\mathfrak{g}_{-}$,
$f$ does not lie in $L_U$.
Notice that $L_U$ contains the elements $x_2, \dots, x_n$
but it does not contain neither the elements $\xi_i$ 
for any $i=1, \dots, n$,
nor the element $x_1$, since these elements do not stabilize the ideal
$I_U$.
It follows that the elements $\xi_i\xi_j$ cannot lie in
$L_U$ for any $i\neq j$. Indeed, $[x_j, \xi_i\xi_j]=-\xi_i$.
Likewise, by induction on $k=1, \dots, n$, 
the elements $\xi_{i_1}\dots \xi_{i_k}$ cannot lie in
$L_U$ for any $k=1, \dots, n$.
Now, suppose that $x_1\xi_j$ lies in $L_U$ for some $j\neq 1$.
Then $L_U$ contains the element $[x_j, x_1\xi_j]=x_1$ 
and this contradicts our assumptions. It
follows that $L_U$ cannot contain the elements $x_1\xi_j$
and, similarly, the elements $x_1\xi_{j_1}\dots\xi_{j_k}$
for any $j_1\neq\dots\neq j_k\neq 1$.
$L_U$ is therefore contained in $L_1$, 
hence $L_U=L_1$.

Finally, the graded subalgebra of $L$ of type $(1, \dots,1|0,\dots,0)$ 
is the standard subalgebra of $L$ stabilizing the ideal $(x_1,\dots, x_n)$.
\end{remark}

\begin{remark}\label{SHO(3,3)}{\em We recall 
that
$Der SHO(3,3)=SHO(3,3)+\mathfrak{a}$ with $\mathfrak{a}\cong gl_2$
(cf.\ Proposition \ref{derivations}, \cite[Remark 4.4.1]{CK}).
The subalgebra $\mathfrak{a}$ of outer derivations
 is generated by the Euler operator $E$
and by  a copy of $sl_2$ 
with Chevalley basis $\{e, h, f\}$ where
$e=ad(\xi_1\xi_3\frac{\partial}{\partial x_2}-
\xi_2\xi_3\frac{\partial}{\partial x_1}-
\xi_1\xi_2\frac{\partial}{\partial x_3})$ and 
$h=ad(\frac{2}{3}\sum_{i=1}^3(\xi_i\frac{\partial}{\partial\xi_i}-
x_i\frac{\partial}{\partial x_i}))$.
In order to describe the action of the derivation $f$ one can proceed as
in \cite[Lemma 5.9]{FK}. Here it is convenient, as before,
 to identify $SHO'(3,3)$ with the set of elements $g$ in
$\Lambda(3,3)/\mathbb{C}1$ such that $\Delta(g)=0$, and
$SHO(3,3)$ with the subspace consisting of elements not containing
the monomial $\xi_1\xi_2\xi_3$.
Under this identification $e=ad(\xi_1\xi_2\xi_3)$ and $h=ad(\frac{2}{3}\sum_{i=1}^3 x_i\xi_i)$. 
Let us consider $SHO(3,3)$ with its  principal grading.
With respect to this grading,  $SHO(3,3)_j=(SHO(3,3)_1)^j$,
for $j>1$, therefore it is sufficient to
define the derivation $f$ on the local part $SHO(3,3)_{-1}\oplus SHO(3,3)_0\oplus SHO(3,3)_1$ of $SHO(3,3)$.
One has:
$f(\xi_1\xi_2)=-\frac{4}{3}x_3$, $f(\xi_1\xi_3)=\frac{4}{3}x_2$, 
$f(\xi_2\xi_3)=-\frac{4}{3}x_1$;
$f(x_1\xi_2\xi_3)=-\frac{1}{3}x_1^2$,
$f(x_2\xi_1\xi_3)=\frac{1}{3}x_2^2$,
$f(x_3\xi_1\xi_2)=-\frac{1}{3}x_3^2$,
$f(x_1\xi_1\xi_2-x_3\xi_3\xi_2)=-\frac{2}{3}x_1x_3$,
$f(x_2\xi_2\xi_1-x_3\xi_3\xi_1)=\frac{2}{3}x_2x_3$,
$f(x_1\xi_1\xi_3-x_2\xi_2\xi_3)=\frac{2}{3}x_1x_2$, 
and $f=0$ elsewhere on $SHO(3,3)_{-1}\oplus SHO(3,3)_0\oplus SHO(3,3)_1$.}
\end{remark}

\begin{remark}\label{inf.manySHO}\em Let $S=\prod_{j\geq -2}S_j$ denote
 the Lie superalgebra $SHO(3,3)$ with respect to the grading of type
 $(2,2,2|1,1,1)$. Then $S_0\cong sl_3$ and $S_{-1}$ is isomorphic, as an
 $S_0$-module, to the direct sum of two copies of the standard
$sl_3$-module. It follows that, for every irreducible $sl_3$-submodule
 $U$ of $S_{-1}$, $S_U:=U+\prod_{j\geq 0}S_j$ is a maximal open
 subalgebra of $S$. In particular, if $U=\langle
 \xi_i\xi_j ~|~ i,j=1,2,3\rangle$ or $U=\langle
 x_i ~|~ i=1,2,3\rangle$, then $S_U$ is the maximal
 graded subalgebra of type $(1,1,1|1,1,1)$ or $(1,1,1|0,0,0)$, respectively.
The subalgebras $S_U$ are not conjugate by inner automorphisms of $S$,
but they are conjugate by inner automorphisms of $Der S$,
since the copy of $sl_2$ of outer derivations of $S$
described in Remark \ref{SHO(3,3)}, permutes the subspaces $U$. In particular the graded subalgebras of principal and subprincipal
 type 
are conjugate by the automorphism $\exp(e)\exp(-3/4f)\exp(e)\in G$.
\end{remark}

\begin{remark}\label{change1}\em Let $1\leq i<j\leq n$. Then the change
of indeterminates that exchanges $x_i$ with $x_j$ and $\xi_i$ with $\xi_j$
preserves the form $\sigma$. 
\end{remark}
\begin{remark}\label{change2}\em
Let $\eta=\alpha_{i_1}\xi_{i_1}+\dots+\alpha_{i_k}\xi_{i_k}$
for some $k\leq n$, with $\alpha_{i_j}\in\mathbb{C}$, $\alpha_{i_j}\neq 0$.
According to Remark \ref{change1} we can assume 
$\eta=\alpha_{1}\xi_{1}+\dots+\alpha_{k}\xi_{k}$ with
$\alpha_i\neq 0$ for $i=1,\dots,k$.
Then the following change of indeterminates preserves the form
$\sigma$:
$$\begin{array}{ll}
x'_1=\frac{1}{\alpha_1}x_1 & \xi'_1=\eta\\
x'_2=x_2-\frac{\alpha_2}{\alpha_1}x_1 & \xi'_2=\xi_2\\
\vdots & \vdots\\
x'_k=x_k-\frac{\alpha_k}{\alpha_1}x_1 & \xi'_k=\xi_k\\
x'_i=x_i & \xi'_i=\xi_i ~~\forall i>k.
\end{array}$$
\end{remark}

\begin{theorem}\label{standardforHO} Let $L=HO(n,n)$ 
and let $L_0$ be a maximal open  subalgebra of $L$.
Then $L_0$ is conjugate to a standard subalgebra of $L$.
\end{theorem}
{\bf Proof.} Let $L=HO(n,n)$. By Remark \ref{rome} $L_0$ stabilizes an ideal of the form
$$J=(x_1+f_1, \dots, x_n+f_n,
 \eta_1+g_1, \dots, \eta_s+g_s)$$
for some
linear functions $\eta_j$ in odd indeterminates, and even functions
$f_i$ and odd functions $g_j$ without constant and linear terms,
 and $J$ is maximal among the $L_0$-invariant
ideals of $\Lambda(n,n)$.
By Remark \ref{change2}, up to changes of indeterminates, we can write
$$J=(x_1+F_1, \dots, x_n+F_n,
 \xi_1+G_1,\dots, \xi_s+G_s)$$
for some even functions
$F_i$ and odd functions $G_j$ without constant and linear terms,
where we can assume $G_j$ independent of $\xi_1, \dots, \xi_s$ for
every $j=1,\dots, s$.

Suppose that $x_1+F_1=x_1+\xi_1F'_1+F^{\prime\prime}_1$ with
$F'_1$ and $F^{\prime\prime}_1$ independent of $\xi_1$.
Then we can replace $x_1+F_1$ by $x_1+F_1-(\xi_1+G_1)F'_1=x_1+H_1$
with $H_1$ independent of $\xi_1$. Similarly we can make
every function 
$F_i$ independent of $\xi_j$ for every $j=1,\dots,s$.
Besides, as in Lemma \ref{victorshelp}, since the ideal $J$
is closed, we can make the functions $F_i$ and $G_i$ independent of
all even variables, i.e., $F_i, G_i\in\Lambda(0,n)$.
It follows that the automorphism $\exp(ad(\xi_1F_1))$ maps $J$
to the ideal 
$$I=(x_1, x_2+F'_2, \dots, x_n+F'_n,
 \xi_1+G_1, \xi_2+G_2,\dots, \xi_s+G_s).$$
Arguing in the same way for every function $F'_j$ with
$1\leq j\leq s$, we have, up
to automorphisms,
$$I=(x_1, \dots, x_s, x_{s+1}+h_{s+1}, \dots, x_n+h_n,
 \xi_1+G_1, \dots, \xi_s+G_s)$$
for some functions $h_i\in\Lambda(\xi_{s+1}, \dots, \xi_{n})$ with
no constant and linear terms.
Now
the automorphism
$\exp(ad(-x_1G_1))$ sends $I$ to the ideal
$$I_1=(x_1, \dots, x_s, x_{s+1}+h_{s+1}, \dots, x_n+h_n,
 \xi_1, \xi_2+G_2, \dots, \xi_s+G_s).$$
Analogous automorphisms for $G_i$, $i=1,\dots,s$, yield to the ideal
$$Y=(x_1, \dots, x_s, x_{s+1}+h_{s+1}, \dots,  x_n+h_n,
 \xi_1, \xi_2,\dots, \xi_s).$$

Consider the ideal $Y'=(x_{s+1}+h_{s+1}, \dots,  x_n+h_n)\subset Y$.
Then, each function $f$ in $L_0$ (thus stabilizing $Y$)
stabilizes the ideal
$K=(Y, [Y',Y'])$, i.e., the ideal
generated by the generators of $Y$ and 
by
the commutators between every pair of generators of $Y'$. Indeed, for every $g, h\in Y'$
we have:
$$[f,[g,h]]=[[f,g],h]\pm[g,[f,h]]\in[Y, Y']$$
and $[Y, Y']\subset K$ since all generators of $Y$ outside $Y'$ commute
with the generators of $Y'$.
Notice that $K$ is a proper ideal of
$\Lambda(2k+1,n)$ since among its generators there is
no invertible element. By the maximality of
$J$ among the ideals stabilized by $L_0$ we have 
$Y=K$.

Suppose that $h_{s+1}=\xi_{s+1}\psi_1+\psi_2$ with 
$\psi_1$ and $\psi_2$ independent of $\xi_{s+1}$. Then,
applying the automorphism $\exp(ad(\xi_{s+1}\psi_2))$, we can assume
$$Y=(x_1, \dots, x_s, x_{s+1}+\xi_{s+1}\psi_1, \dots,  x_n+h'_n,
 \xi_1, \xi_2,\dots, \xi_s).$$
Now $\psi_1=1/2[x_{s+1}+\xi_{s+1}\psi_1, x_{s+1}+\xi_{s+1}\psi_1]\in
[Y',Y']\subset K=Y$, therefore
$$Y=(x_1, \dots, x_s, x_{s+1}, x_{s+2}+h'_{s+2}, \dots,  x_n+h'_n,
 \xi_1, \xi_2,\dots, \xi_s).$$
Arguing in the same way for every function $h'_j$ we end up with a standard
ideal.
%
\hfill $\Box$

\begin{theorem}\label{HO(n,n)} $(a)$ Let $L=HO(n,n)$, or $SHO(n,n)$ with $n>3$. Then
all maximal open subalgebras of $L$ are, up to conjugation,
the graded subalgebras of type
$(1,\dots,
1,2,\dots, 2|1,\dots,
1, 0, \dots, 0)$ with
$k$ 2's and $k$ zeros, for $k=0,\dots,n-2$ and the graded
subalgebra of type $(1,\dots,1|0,\dots,0)$.

$(b)$ All maximal open subalgebras of $SHO(3,3)$ are, up to
conjugation, the graded subalgebras of type $(1,1,1|1,1,1)$ and
$(1,1,2|1,1,0)$.
\end{theorem}
{\bf Proof.} Let $L=HO(n,n)$ and let $L_0$ be a maximal open
subalgebra of $L$. By Theorem \ref{standardforHO},  $L_0$ is,
up to conjugation, the standard subalgebra of $L$
stabilizing the ideal $I_{\cal{U}}=(x_1, \dots, x_n, \xi_1, \dots, \xi_s)$
for some $s=0, \dots, n$. The statement then follows using Remarks 
\ref{gradingsofHO}, \ref{irrenotHO(n,n)} and \ref{stabHO}. 

Let now $L=SHO(n,n)$ and let $L_0$ be a maximal open subalgebra of $L$.
The same argument as in Theorem \ref{newS} shows that $L_0$ is regular
and we can assume, by Remark \ref{standard}, that it is invariant 
with respect to the standard torus $T+\C E$ of $Der HO(n,n)$. It follows that
$L_0$ decomposes into the direct product of weight spaces
with respect to $T+\C E$. As we noticed in Remark \ref{stabHO},
$\C\xi_{i_1}\dots\xi_{i_h}$, with
$i_1\neq\dots\neq i_h$, and
$\C x_1\xi_{j_1}\dots \xi_{j_h}$, with $1\neq j_1\neq\dots\neq j_h$,
are one-dimensional weight spaces with respect to $T+\C E$.
Besides,
note that the elements
$\xi_i$ cannot lie in $L_0$ since they are not exponentiable.

We may assume that one of the following 
two cases holds:

\noindent
1) $x_1, \dots, x_n$ lie in $L_0$.
Since $[x_i, \xi_i\xi_h]=\xi_h$, 
it follows that the $(T+\C E)$-invariant complement of $L_0$
contains the elements $\xi_i\xi_h$ for every $i, h=1,\dots,n$. 
Arguing inductively,
since $[x_{i_1}, \xi_{i_1}\dots \xi_{i_h}]=\xi_{i_2}\dots\xi_{i_h}$, 
one shows that $L_0$ cannot contain any element 
lying in the negative part of the grading of type
$(1,\dots,1|0,\dots,0)$, therefore
$L_0$ is contained in
the maximal graded subalgebra of $L$ of type
$(1,\dots,1|0,\dots,0)$, thus $L_0$ coincides with this graded subalgebra, 
by maximality;

\noindent
2) $x_1,\dots, x_k$ do not lie in $L_0$ for some $k=2, \dots, n$,
 and 
$x_{k+1},\dots, x_n$ lie in $L_0$.
Then the $(T+\C E)$-invariant complement of $L_0$ contains
the elements $\xi_hP$ for $h=1, \dots, n$ and $P\in\Lambda(\xi_{k+1},
\dots, \xi_n)$. Likewise,
since $[x_{i_1}, x_j\xi_{i_1}\dots \xi_{i_h}]=
x_j\xi_{i_2}\dots\xi_{i_h}$, the
$(T+\C E)$-invariant complement of $L_0$ contains the elements
$x_jP$, for $j=1, \dots, k$ and 
$P\in\Lambda(\xi_{k+1},
\dots, \xi_n)$.
  Therefore the $(T+\C E)$-invariant complement of $L_0$ contains
the $(T+\C E)$-invariant complement of the graded subalgebra of type
$(1, \dots, 1,2,\dots,2|1,\dots,1,0,\dots,0)$ with $n-k$ 2's and $n-k$ zeros.
Hence $L_0$ coincides with this graded subalgebra of $L$.

Note that any open regular subalgebra of $L$ containing $x_2,
\dots, x_{n}$ 
 and not containing $x_1$, is not a maximal  subalgebra
of $L$.  Indeed any such a subalgebra is contained in the graded subalgebra
of type $(1,2,\dots,2|1,0,\dots,0)$ which is not maximal by Remark
\ref{irrenotHO(n,n)}.  

By Remark \ref{inf.manySHO},  the subalgebras of principal and subprincipal
type of \break $SHO(3,3)$ are
conjugate by an element of $G$.
\hfill$\Box$

\begin{corollary}\label{HO(n,n)-gradings} $(a)$ All irreducible $\Z$-gradings
of $HO(n,n)$ and of
$SHO(n,n)$ with $n>3$, are, up to conjugation,
the gradings
of type
$(1,\dots,
1,2,\dots, 2|$ $1,\dots,
1, 0, \dots, 0)$ with
$k$ 2's and $k$ zeros, for $k=0,\dots,n-2$ and the grading
of type $(1,\dots,1|0,\dots,0)$.

$(b)$ All irreducible $\Z$-gradings of $SHO(3,3)$ are, up
to conjugation, the gradings of type $(1,1,1|1,1,1)$ and
$(1,1,2|1,1,0)$. 
\end{corollary}

\begin{remark}\label{DerHO}\em 
By Remark \ref{rome}, the proof
of Theorem \ref{standardforHO} works verbatim if we replace $L=HO(n,n)$ with
$Der L$ and $L_0$ with a
 fundamental maximal subalgebra of $Der L$.
Therefore every fundamental maximal subalgebra of $Der L$ is conjugate
to the standard subalgebra of $Der L$ stabilizing the ideal 
$I_{\cal U}=(x_1, \dots, x_n, \xi_1, \dots, \xi_s)$ of $\Lambda(n,n)$, 
for some $s=0,\dots,n$.
\end{remark}

\begin{theorem} Let $L=HO(n,n)$. Then all maximal 
among $E$-invariant subalgebras of $L$
are, up to conjugation, the 
subalgebras of $L$ listed in Theorem \ref{HO(n,n)}$(a)$.
\end{theorem}
{\bf Proof.} By Remark \ref{DerHO} every
 fundamental maximal subalgebra
of $Der L$ is conjugate to the standard subalgebra of $Der L$ stabilizing the
ideal $I_{\cal U}=(x_1, \dots, x_n,$ $\xi_1, \dots, \xi_s)$ of
$\Lambda(n,n)$, for some $s=0, \dots, n$.
Therefore, by Remarks
\ref{gradingsofHO}, \ref{irrenotHO(n,n)} and \ref{stabHO},
all fundamental maximal subalgebras
of $Der L$ are, up to conjugation, the subalgebras $L_0+\mathbb{C}E$ where $L_0$ is one of the
maximal open subalgebras of $L$ listed in Theorem \ref{HO(n,n)}$(a)$.
If $S_0$ is a maximal among open $E$-invariant subalgebras of $L$,
then $S_0+\C E$ is a fundamental maximal subalgebra of $Der L$ and the 
thesis follows. 
 \hfill$\Box$



\medskip

As in the case of the Lie superalgebra $S(1,n)$, we are now interested in the 
subalgebras of $SHO(n,n)$ which are maximal among its
$\mathfrak{a}_0$-invariant subalgebras, for any subalgebra $\mathfrak{a}_0$
 of the subalgebra $\mathfrak{a}$ of outer 
derivations of $SHO(n,n)$.  

\begin{remark}\label{SHO'}\em The same arguments as in Theorem
\ref{newS} show that 
every maximal open subalgebra of $SHO'(n,n)$ and $CSHO'(n,n)$
is regular. Therefore 
the same arguments as for $SHO(n,n)$
 in Theorem \ref{HO(n,n)}, show that all fundamental among maximal
subalgebras of $SHO'(n,n)$ or $CSHO'(n,n)$ with $n\geq 3$, are,
up to conjugation,
the graded subalgebras of type 
$(1,\dots,1,2,\dots,2|1,\dots,1,0,\dots,0)$ with $k$ 0's and $k$ 2's, for some
$k=0,\dots, n-2$. 
Indeed the graded subalgebra of $SHO'(n,n)$ (resp.\ $CSHO'(n,n)$)
of type $(1, \dots,1|0,\dots,0)$ 
is not maximal, since it is contained in $SHO(n,n)$
(resp.\ $SHO(n,n)+\C\Phi+\C E$). Notice that the graded
subalgebras of  principal and subprincipal type of $SHO'(3,3)$ (resp.\ 
$CSHO'(3,3)$) are not conjugate.
By the same arguments, all maximal among fundamental subalgebras of
$SHO'(n,n)$ and $CSHO'(n,n)$ are, up to conjugation,
the graded subalgebra of type $(1,\dots,1|0,\dots,0)$ and
the graded subalgebras of type 
$(1,\dots,1,2,\dots,2|1,\dots,1,0,\dots,0)$ with $k$ 0's and $k$ 2's, for some
$k=0,\dots, n-2$.
\end{remark}

\begin{theorem}\label{derSHO} Let $L=SHO(n,n)$ with $n\geq 4$.

\noindent
(i) If
$\mathfrak{a}_0$ is a torus of $\mathfrak{a}$, then all maximal
among open $\mathfrak{a}_0$-invariant subalgebras of $L$ 
are, up to conjugation,
the graded subalgebra of type $(1,\dots,1|0,\dots,0)$ and
 the graded subalgebras of type 
$(1,\dots,1,2,\dots,2|1,\dots,1,0,\dots,0)$ with $k$ 0's and $k$ 2's, for some
$k=0,\dots, n-2$.

\noindent
$(ii)$ If $\mathfrak{a}_0=\C\xi_1\dots\xi_n\rtimes \mathfrak{t}$,
where $\mathfrak{t}$ is a torus of $\mathfrak{a}$,
then all maximal among open $\mathfrak{a}_0$-invariant subalgebras of $L$ 
are, up to conjugation,
the graded subalgebras
of type $(1,\dots,1,2,\dots,2|1,\dots,1, 0,\dots,0)$,
with
$k$ 2's and $k$ zeros, for $k=0,\dots, n-2$.
%
\end{theorem}
{\bf Proof.} One uses Remark \ref{SHO'} and the same arguments  as in
the proof of
Theorem \ref{derS(1,n)}.
%
 \hfill$\Box$


\begin{theorem}\label{a-invSHO(3,3)} Let $L=SHO(3,3)$
and let  let
$\mathfrak{b}=\C e+\C h\subset \mathfrak{a}\cong gl_2$.

\noindent
(i) If $\mathfrak{a}_0$ is a 
one-dimensional subalgebra of $\mathfrak{a}$ or a two-dimensional torus of $\mathfrak{a}$,
then all maximal among open $\mathfrak{a}_0$-invariant subalgebras of $SHO(3,3)$ are, up to conjugation, the subalgebras of type $(1,1,1|1,1,1)$ and
$(1,1,2|1,1,0)$.

\noindent
(ii) If $\mathfrak{a}_0=\C e\rtimes\mathfrak{t}$, 
where $\mathfrak{t}$ is a torus of $\mathfrak{a}$, then
 the graded
subalgebra of type $(1,1,1|1,1,1)$  is, up to
conjugation, the only 
maximal among open $\mathfrak{a}_0$-invariant 
subalgebras of $SHO(3,3)$, which is not
invariant with respect to $\mathfrak{a}$.

\noindent
(iii) If ~$\mathfrak{a}_0=sl_2$ or
$\mathfrak{a}_0=\mathfrak{a}$, then
all maximal among open  $\mathfrak{a}_0$-invariant subalgebras
of $SHO(3,3)$ are, up to conjugation, the subalgebras of type $(1,1,2|1,1,0)$ and $(2,2,2|1,1,1)$.
\end{theorem}
{\bf Proof.} By Remark \ref{SHO'}, the proof of $(i)$ is the same
as the 
proof of $(i)$ and
$(ii)$ in
Theorem \ref{derS(1,n)}. Recall that the graded
subalgebras of principal and subprincipal type of $SHO(3,3)$ are
conjugate.

Now,  using Remark \ref{SHO(3,3)}, one verifies that
the graded subalgebras of $SHO(3,3)$ of type $(1,1,2|1,1,0)$ and  
$(2,2,2|1,1,1)$ are invariant with respect to $\mathfrak{a}$
(see also \cite[Example 5.5]{K}, \cite[Remark 4.4.1]{CK}).
On the other hand
the maximal graded subalgebra $L_0$ of $SHO(3,3)$  of type $(1,1,1|1,1,1)$
is 
invariant with respect to the action of  $h$, $e$ and $E$, 
but it
 is not $\mathfrak{a}$-invariant, indeed: $\xi_i\xi_j\in L_0$, $x_k\notin L_0$ and
  $f(\xi_i\xi_j)=\pm\frac{4}{3} x_k$ with $k\neq i,j$. 

%
Let $S_0$ be a maximal among open $\mathfrak{b}$-invariant
subalgebras of $SHO(3,3)$, then $S_0+\mathbb{C}\xi_1\xi_2\xi_3+\mathbb{C}
\sum_{i=1}^3 x_i\xi_i$ is a fundamental subalgebra of
$CSHO'(3,3)$, hence
it is contained in a maximal among fundamental subalgebras of
$CSHO'(3,3)$ containing $\xi_1\xi_2\xi_3$ and
$\sum_{i=1}^3 x_i\xi_i$. 
It follows, by Remark \ref{SHO'}, 
that $S_0$ is conjugate
either to
the graded subalgebra of type $(1,1,1|1,1,1)$ or  to the subalgebra
of type $(1,1,2|1,1,0)$. A similar argument holds if
$S_0$ is maximal among open $\mathfrak{a}_0$-invariant subalgebras,
with $\mathfrak{a}_0=\C e+\mathfrak{t}$ where $\mathfrak{t}$
is a one-dimensional torus of $\mathfrak{a}$.
Likewise, if $S_0$ is a maximal among open $\mathfrak{b}+\C E$-invariant
subalgebras of $SHO(3,3)$, then $S_0+
\mathbb{C}\xi_1\xi_2\xi_3+\mathbb{C}
\sum_{i=1}^3 x_i\xi_i+\C E$ is a fundamental maximal subalgebra of
$CSHO'(3,3)$, hence, by Remark \ref{SHO'}, 
it is conjugate
either to
the graded subalgebra of type $(1,1,1|1,1,1)$ or  to the subalgebra
of type $(1,1,2|1,1,0)$.


Finally, let $\mathfrak{a}_0=sl_2$
or $\mathfrak{a}_0=\mathfrak{a}$, and let $S'$ be   a maximal among open $\mathfrak{a}_0$-invariant subalgebras
of $SHO(3,3)$. Then $S'$ is  $\mathfrak{b}$-invariant,
hence $S'+\mathfrak{b}$ is contained in a maximal among
fundamental subalgebras of $CSHO'(3,3)$ containing $\mathfrak{b}$. It
follows that $S'$ is contained either in a conjugate of
the subalgebra of type $(1,1,2|1,1,0)$, thus coincides with it by maximality,
or 
in a conjugate $S_U$ of the subalgebra of type $(1,1,1|1,1,1)$. 
As we noticed in Remark \ref{inf.manySHO}, $S_U$ is conjugate to the subalgebra
of principal type by an automorphism $\varphi=\exp(ad~ a)$ 
for some $a\in\mathfrak{a}$. Since $S'$ is $\mathfrak{a}$-invariant,
$\varphi(S')=S'$, therefore
$S'$ is contained in the intersection of $S_U$ with the graded subalgebra
of principal type, i.e., it is contained in
the graded subalgebra of $SHO(3,3)$ of type $(2,2,2|1,1,1)$,
 thus
coincides with it by maximality.
\hfill$\Box$

\section{Maximal open subalgebras of $\boldsymbol{H(2k,n)}$}
Let $p_1, \dots, p_k, q_1,\dots,
q_k$ be $2k>0$ even indeterminates and $\xi_1,\dots, \xi_n$ be $n$ odd
indeterminates. Consider the differential
form $\omega=2\sum_{i=1}^kdp_i\wedge dq_i+\sum_{i=1}^nd\xi_id\xi_{n-i+1}$. The
Hamiltonian superalgebra $H(2k,n)$ is the Lie superalgebra defined as
follows \cite{K2}:
$$H(2k,n)=\{X\in W(2k,n)~|~ X\omega=0\}.$$
Let us consider the Lie superalgebra $\Lambda(2k,n)$ with the
following bracket:
\begin{equation}[f,g]=\sum_{i=1}^k(\frac{\partial f}{\partial p_i}\frac{\partial
g}{\partial q_i}-\frac{\partial f}{\partial q_i}\frac{\partial
g}{\partial p_i})-(-1)^{p(f)}\sum_{i=1}^n\frac{\partial f}{\partial \xi_i}
\frac{\partial g}{\partial \xi_{n-i+1}}.
\label{bracket}
\end{equation}
Then the map $$\Lambda(2k,n)\longrightarrow H(2k,n)$$
$$f\longmapsto \sum_{i=1}^k(\frac{\partial f}{\partial p_i}\frac{\partial
}{\partial q_i}-\frac{\partial f}{\partial q_i}\frac{\partial
}{\partial p_i})-(-1)^{p(f)}\sum_{i=1}^n\frac{\partial f}{\partial \xi_i}
\frac{\partial }{\partial \xi_{n-i+1}}$$ defines a surjective homomorphism
whose kernel consists of constant functions (cf.\ \cite[\S 1.2]{CK}). We
will therefore identify $H(2k,n)$ with $\Lambda(2k,n)/\C 1$ with
bracket (\ref{bracket}). Let us
fix the maximal torus $T=\langle p_iq_i,
\xi_j\xi_{n-j+1} ~|~ i=1,\dots, k, ~j=1,\dots, [n/2]\rangle$ of
$H(2k,n)$.
\begin{remark}\label{gradingsofH}\em
The $\Z$-grading of type $(a_1,\dots, a_{2k}|b_1, \dots, b_n)$ of
$W(2k,n)$ induces a grading of $H(2k,n)$ if and only if the
differential form $\omega$ is homogeneous in this grading (cf.\ \cite{K5}).

The $\Z$-grading of type  $(1,\dots,1|2,\dots,
2,1,\dots,1,0,\dots, 0)$ of
$W(2k,n)$, with
$t$ 2's and $t$ zeros,  induces an irreducible grading on
$H(2k,n)$ for every $t$ such that $0\leq t\leq
\left[\frac{n}{2}\right]$, where
$\g_0\cong
spo(2k,n-2t)\otimes\Lambda(t)+W(0,t)$ and $\g_{-1}\cong
\C^{2k|n-2t}\otimes\Lambda(t)$. One can check that
when $t>0$, $[\g_{-1}, \g_{-1}]\neq 0$ thus it coincides with $\g_{-2}$
by Remark \ref{gd}. Besides, property $(iii)^\prime$ of Proposition 
\ref{basic}$(b)$
is satisfied. Therefore the subalgebras $\prod_{j\geq 0}H(2k,n)_j$ of $H(2k,n)$
corresponding to the gradings of type $(1,\dots,1|2,\dots,
2,1,\dots,1,0,\dots, 0)$, with
$t$ 2's and $t$ zeros, for $0\leq t\leq
\left[\frac{n}{2}\right]$, are maximal open subalgebras of $H(2k,n)$.

The $\Z$-grading of $H(2k,n)$ induced by the principal grading of $W(2k,n)$ is also
called {\em principal}. The $\Z$-grading induced on $H(2k,2h)$ by
the $\Z$-grading of type $(1,\dots,1|2,\dots,2,0,\dots,0)$, with $h$ zeros,
is called {\em subprincipal}.
\end{remark}
\begin{remark}\em  The $\Z$-gradings of $H(2k,2h)$ of type 
$(1,\dots,1|2,\dots,2,$ $0,\dots,0)$ and
$(1,\dots,1|2,\dots,2,0,2,0,\dots,0)$, with $h$ zeros,  
are not conjugate by an element of $G$, but are conjugate by an outer automorphism.
\end{remark}

\begin{example}\label{boston}\em Let us identify $L=H(2k,n)$ with
  $\Lambda(2k,n)/\C 1$. Let $V$ be the
 $n$-dimensional odd
vector space spanned by $\xi_1,\dots, \xi_n$, with the bilinear form
 $(\xi_i, \xi_j)=\delta_{i, n-j+1}$. Let us fix a subspace $U$ of $V$
and let us repeat the same construction as in Example \ref{step2}.

 We define 
a valuation on $\Lambda(2k,n)/\C 1$ with values in $\Z_+$ by letting
$$\nu(p_i)=\nu(q_i)=1;$$
$$\nu(x)=2 ~\mbox{for}~ x\in U^0;~~\nu(x)=0 ~\mbox{for}~ x\in (U^0)^\prime;$$
$$\nu(x)=1 ~\mbox{for}~ x\in U^1;
~~\nu(x)=0 ~\mbox{for}~ x\in (U^1)^\prime.$$

Consider the following subspaces 
of $L$:
$$L_j(U)=\{x\in\Lambda(2k,n)/\C 1 ~|~ \nu(x)\geq j+2\}+\Lambda((U^1)^\prime)/\C 1
  ~~\mbox{for}~ j\leq 0;$$
$$L_j(U)=\{x\in\Lambda(2k,n)/\C 1 ~|~ \nu(x)\geq j+2\} ~~\mbox{for}~
j>0.$$
These subsets define, in fact, a filtration of $H(2k,n)$ for every subspace
$U$ of $V$, as one can verify using the definition of bracket (\ref{bracket}).
Notice that this filtration has depth 1 if and only if $U$ is
non-degenerate, including $U=0$.

Let us denote by $s$ the dimension of $U$ and by $s_i$ the dimension
of $U^i$ for $i=0,1$. Then $\overline{Gr L}\cong H(2k,n-r_1)\otimes\Lambda(r_1)+H(0,r_1)$
with respect to the grading of type
$(1,\dots,1|2,\dots,2,
1,\dots,1,0,\dots,0)$ of $H(2k,n-r_1)$,
with
$s_0$ 2's and $s_0$ zeros, where $r_1=n-2s_0-s_1=\dim(U^1)^\prime$,
and $\deg(\tau)=0$ for every $\tau\in\Lambda(r_1)$. 
This is an irreducible grading of
$H(2k,n-r_1)$ for every choice of $U$ (cf.\ Remark \ref{gradingsofH}), and, by Corollary \ref{cor}, $L_0(U)$ is a maximal open
subalgebra of $L$. 

Let us consider the
standard ideal $I_{\cal U}=(p_1,\dots, p_k, q_1,\dots, q_k , U)$ of $\Lambda(2k,n)$. Notice that $I_{\cal U}=\{x\in\Lambda(2k,n)/\mathbb{C}1~|~
\nu(x)\geq 1\}$. It follows that $L_0(U)$ stabilizes $I_{\cal U}$
hence, due to its maximality, $L_0(U)$ is the standard subalgebra
of $H(2k,n)$ corresponding to the ideal $I_{\cal U}$.
\end{example}

\begin{remark}\label{gradedofH}\em $L_0(U)$ is a maximal graded
  subalgebra of $L$ if and only if $U$ is a coisotropic subspace of $V$.
\end{remark}

\begin{remark}\label{regularofH}\em If $U$ is conjugate to a subspace of $V$ spanned by
  $\xi_{i_1}, \dots, \xi_{i_t}$ for some $i_1, \dots, i_t$, then $U$
  is stable under the action of the maximal torus $T$. It follows that
  in this case
  $L_0(U)$ is regular. If $n$ is odd, then any subspace of
  $V$ is conjugate to $\langle \xi_{i_1}, \dots, \xi_{i_t}\rangle$
  for some $i_1, \dots, i_t$, and the same holds  when $n$ is even for
  any subspace of $V$ whose non-degenerate part has even dimension.
\end{remark}

\begin{remark}\label{notregular}\em If $n$ is even and $s_1$ is odd then any maximal torus
  of $L$ has dimension $k+\frac{n}{2}$ and any maximal torus of $\overline{Gr L}$
  has dimension $k+\frac{n}{2}-1$. It follows that, under these
  hypotheses, $L_0(U)$ is not a regular subalgebra of $L$. For example
  any one-dimensional non-degenerate subspace $U$ of $V$ gives
  rise to a maximal open subalgebra $L_0(U)$ of $H(2k, 2t)$ which is not regular.
\end{remark}

\begin{lemma}\label{H1} Let us consider an ideal $J=(h_1,...,h_r)$ of
$\Lambda(0,n)$.
Suppose that 
$h_1=\eta_1+F$ and
$h_2=\eta^{\prime}_1+G$
where $\eta_1, \eta^{\prime}_1$  are non-degenerately paired  elements
of $V$ and $F, G$ contain no constant and linear terms. Then $J$ is
conjugate
 to an ideal $K=(\eta_1, \eta^{\prime}_1, f_1,...,f_{r-2})$
for some functions $f_i\in\Lambda(U)$ where $U$ is the orthogonal 
complement of $\langle \eta_1, \eta^{\prime}_1\rangle$ in $V$.
\end{lemma}
{\bf Proof.} Up to multiplying $h_1$  by some invertible function
 we can assume that $F$ does not depend on $\eta_1$, i.e., 
$\eta_1+F=\eta_1+f_1\eta^{\prime}_1+f_2$ where $f_1, f_2$ lie in $\Lambda(U)$.
Also, we can assume that $G$ lies in $\Lambda(U_1)$ where
$U_1=\langle U, \eta_1\rangle$.
Notice that $f_1\eta^{\prime}_1+f_1G$ lies in $J$, therefore
$J=(\eta_1+f_2-f_1G, \eta^{\prime}_1+G, h_3,\dots, h_r)$ where
$f_2-f_1G\in\Lambda(U_1)$. Therefore, up to multiplying
$\eta_1+f_2-f_1G$ by an invertible function, we can
write $J=(\eta_1+F', \eta^{\prime}_1+G, h_3,\dots, h_r)$ where
$F'\in\Lambda(U)$.

Now the automorphism $\exp ad(\eta^{\prime}_1F')$
maps $J$ to
the ideal $J'=(\eta_1, \eta^{\prime}_1+H, h'_3,...,h'_r)$ where the
$h'_i$'s lie in $\Lambda(U_1)$ and 
$H$ lies in $\Lambda(U)$. Then, similarly as above, 
the automorphism $\exp ad(\eta_1H)$, maps $J'$ to the ideal
$K=(\eta_1, \eta^{\prime}_1, f_1,...,f_{r-2})$, since
$H\in\Lambda(U)$. Since $\eta_1, \eta^{\prime}_1$ lie in $K$,
we can assume $f_1,...,f_{r-2}\in\Lambda(U)$.
\hfill$\Box$

\begin{remark}\label{H1+}\em Notice that if $\eta_1\in V$ is non-degenerately paired 
with itself, one can prove, arguing as in the proof of Lemma \ref{H1},
that if the ideal $J$ contains an element
of the form $\eta_1+F$, then, up to
automorphisms, $J=(\eta_1,f_1,...,f_{r-1})$
where the $f_i$'s lie in $\Lambda(U)$, $U$ being the orthogonal
complement of $\langle \eta_1\rangle$ in $V$.
\end{remark}

\begin{theorem}\label{standardforH} Let $L_0$ be a  maximal open subalgebra of $L=H(2k,n)$.
Then $L_0$ is conjugate to a standard subalgebra of $L$.
\end{theorem}
{\bf Proof.} By Remark \ref{rome} $L_0$ stabilizes an ideal of the form
$$J=(p_1+f_1, q_1+h_1, \dots, p_k+f_k, q_k+h_k,
 \eta_1+g_1,\eta_2+g_2, \dots, \eta_r+g_r)$$
for some
linear functions $\eta_j$ in odd indeterminates and even functions
$f_i, h_i$ and odd functions $g_j$ without constant and linear terms,
 and this ideal is maximal among the $L_0$-invariant
ideals of $\Lambda(2k,n)$. As in Lemma \ref{victorshelp} we can assume
$f_i, h_i$ and $g_j$ in $\Lambda(0,n)$.

Note that the automorphism
$\exp(ad(q_1f_1))$ maps
 $J$ to  
$J_1=(p_1, q_1+h'_1, p_2+f'_2, q_2+h'_2, \dots, p_k+f'_k, q_k+h'_k, \eta_1+g'_1, \eta_2+g'_2, \dots, \eta_r+g'_r)$.
As above we can assume $h'_1$  independent of all even
variables.
It follows that the automorphism $\exp ad(-p_1h'_1)$ maps $J_1$
to $J_2=(p_1, q_1, p_2+f^{\prime\prime}_2, q_2+h^{\prime\prime}_2, \dots, p_k+f^{\prime\prime}_k, q_k+h^{\prime\prime}_k, \eta_1+g^{\prime\prime}_1, \eta_2+g^{\prime\prime}_2, \dots, \eta_r+g^{\prime\prime}_r)$.
The same arguments applied to all generators $p_i+f^{\prime\prime}_i$ and
$q_j+h^{\prime\prime}_j$ show that $J$ is in fact conjugate to the ideal
$$I=(p_1, p_2,\dots, p_k, q_1, \dots, q_k, \eta_1+\ell_1, \eta_2+\ell_2, \dots, \eta_r+\ell_r)$$
where 
$\eta_1, \dots, \eta_r$ are linearly independent vectors in $V$ and 
$\ell_1, \dots, \ell_r$ are functions in $\Lambda(0,n)$
without constant and  linear terms. 
Since from now on we shall work only
with odd indeterminates, with an abuse of notation we shall
simply write
$$I=(\eta_1+\ell_1,\eta_2+\ell_2, \dots, \eta_r+\ell_r).$$

Let $U=\langle \eta_1,\dots, \eta_r\rangle\subset V$. Let 
$U^0=\langle \nu_1, \dots, \nu_s\rangle$ be
the kernel of the restriction of the bilinear form $(\cdot, \cdot)$
to $U$. Then, as in Example \ref{boston}, $U=U^0\oplus U^1$ where $U^1$ is a maximal
subspace of $U$ with non-degenerate metric.
Then, by Lemma \ref{H1} and Remark \ref{H1+}, $I=(U^1, \nu_1+\ell_1,\dots, \nu_s+\ell_s)$
where  
$\ell_1, \dots, \ell_s\in\Lambda((U^1)^\perp)$. 
In particular, $(U^1)^\perp$ contains $U^0$
and a subspace $(U^0)'$ non-degenerately paired with $U^0$. Let
$(U^0)'=\langle \nu'_1, \dots, \nu'_s\rangle$ with $(\nu_i, \nu'_j)=\delta_{i,j}$.

\medskip
Now, if $\ell_i=0$ for every $i=1,\dots, s$, then $I$ is standard.
Suppose that at least one of the $\ell_j$'s is not zero, i.e.,  
$$I=(U^1,\nu_1,\dots, \nu_{k-1},
\nu_k+\ell_k,\dots, \nu_s+\ell_s)$$ with 
$k=\min\{i=1, \dots, s~|~\ell_i\neq 0\}$.

Denote by
$I'$ the ideal $I'=(\nu_1,\dots, \nu_{k-1},
\nu_k+\ell_k,\dots, \nu_s+\ell_s)\subset I$.
Then, each function $f$ in $L_0$ (thus stabilizing $I$)
stabilizes the ideal
$K=(I, [I',I'])$. Indeed, for every $g, h\in I'$
we have:
$$[f,[g,h]]=[[f,g],h]\pm[g,[f,h]]\in[I, I']$$
and $[I, I']\subset K$ since every generator of $I'$ is orthogonal
to $U^1$.
Notice that $K$ is generated by the generators of $I$ and by
the 
brackets between every pair of generators of $I'$.
Therefore $K$ is a proper ideal of
$\Lambda(0,n)$ since among its generators there is
no invertible element. By the maximality of
$I$ among the ideals stabilized by $L_0$ we have 
$I=K$.

We first show that the function $\ell_k$ can be made 
independent of  $\nu'_1,\dots, \nu'_{k-1}$.
Indeed, let
$\nu_k+\ell_k=\nu_k+\nu'_1\phi_1+\phi_2$ where 
$\phi_1, \phi_2$ do not depend on $\nu'_1$.
Then $\phi_1=[\nu_1, \nu_k+\ell_k]\in [I^{\prime}, I^{\prime}]\subset
 K=I$, thus 
$I=(U^1,\nu_1,\dots, \nu_{k-1}, \nu_k+\phi_2, \nu_{k+1}+\ell_{k+1}, \dots, \nu_s+\ell_s)$,
where $\phi_2\in \Lambda((U^1)^\perp)$ does not
depend on 
$\nu'_1$. Arguing in the same way with the
variables $\nu'_2, \dots, \nu'_{k-1}$ we get
$$I=(U^1,\nu_1,\dots, \nu_{k-1}, \nu_k+\phi, \nu_{k+1}+\ell_{k+1}, \dots, \nu_s+\ell_s)$$
where $\phi$ does not depend on $\nu'_1,\dots, \nu'_{k-1}$.

Besides, multiplying $\nu_k+\phi$ by an invertible function, we can
assume that $\phi$ does not depend on $\nu_k$.
Now we can write $\phi=\nu'_k\psi_1+\psi_2$ with $\psi_1, \psi_2$
not depending on $\nu'_1,\dots, \nu'_{k}$.
Therefore, 
applying the automorphism $\exp(ad(\nu'_k\psi_2))$ to $I$, we can
assume $\psi_2=0$.  Then $\psi_1=1/2\{\nu_k+\phi, \nu_k+\phi\}\in 
[I^{\prime},I^{\prime}]\subset K=I$.
Therefore 
$$I=(U^1, \nu_1,\dots, \nu_{k-1}, \nu_k,
\nu_{k+1}+\ell_{k+1},\dots, \nu_r+\ell_r).$$ 
Arguing as above for $\ell_{k+1}, \dots, \ell_r$, we end up with
a standard ideal.
\hfill$\Box$

\begin{theorem}\label{H(2k,n)} All  maximal open subalgebras of
$L=H(2k,n)$ are, up to conjugation, the following:

\medskip
\noindent
$(a)$ if $n=2h+1$:
\begin{enumerate}
\item[$(i)$] the $\Z$-graded subalgebras of type $(1,\dots,
1|2,\dots,2, 1,\dots, 1, 0,\dots,0)$
with
$t$ 2's and $t$ zeros for
$0\leq t\leq h$;
\item[$(ii)$] the regular (non-graded) subalgebras $L_0(U)$ constructed in Example \ref{boston} 
  where $U$
is not coisotropic;
\end{enumerate}

\noindent
$(b)$ if $n=2h$:
\begin{enumerate}
\item[$(i)$] the $\Z$-graded subalgebras of type $(1,\dots,
1|2,\dots,2, 1,\dots, 1, 0,\dots,0)$
with
$t$ 2's and $t$ zeros for
$0\leq t\leq h$, and the
$\Z$-graded subalgebra of type 
$(1,\dots,1|2,\dots,2,$ $0,2,0,\dots,0)$ with $h$ zeros;
\item[$(ii)$] the regular (non-graded) subalgebras $L_0(U)$ constructed in Example \ref{boston} 
  where $U$ is not coisotropic and $\dim U^1$ is even;
\item[$(iii)$] the non-regular subalgebras $L_0(U)$ 
constructed in Example \ref{boston}, where $\dim U^1$ is odd.
\end{enumerate}
\end{theorem}
{\bf Proof.} By Theorem \ref{standardforH}, every maximal open subalgebra
of $L$ is conjugate to the standard subalgebra of $L$ stabilizing
the ideal $I_{\cal U}$ of $\Lambda(2k,n)$, for some subspace
${\cal U}=\langle p_1, \dots, p_k, q_1, \dots, q_k, U\rangle$ of
$\sum_{i=1}^k(\mathbb{C}p_i+\mathbb{C}q_i)+\sum_{j=1}^n\mathbb{C}\xi_j$,
where $U$ is a subspace of $V$. Then the statement follows from
Example \ref{boston} and Remarks \ref{gradingsofH}, \ref{gradedofH}, 
\ref{regularofH} and \ref{notregular}.\hfill$\Box$

\begin{corollary} All irreducible $\Z$-gradings of $H(2k,n)$ 
are, up to conjugation, the $\Z$-gradings of type $(1,\dots,
1|2,\dots,2, 1,\dots, 1, 0,\dots,0)$
with
$t$ 2's and $t$ zeros, 
for $t=0, \dots, [\frac{n}{2}]$, and the $\Z$-grading of type
$(1, \dots,1|2,\dots,2,0,2,0,\dots,0)$ with $n/2$ zeros if $n$ is even.
\end{corollary}

We recall that $Der H(2k,n)=H(2k,n)+\mathbb{C}E$ where
$E=\sum_{i=1}^k(p_i\frac{\partial}{\partial p_i}+
q_i\frac{\partial}{\partial q_i})+\sum_{j=1}^n\xi_j\frac{\partial}
{\partial \xi_j}$ is the Euler operator (cf.\ Proposition \ref{derivations}).
We now aim to classify all fundamental maximal subalgebras of 
$Der H(2k,n)$.
\begin{remark}\label{DerH1}\em All members of the filtration
$$H(2k,n)=L_{-d}(U)\supset\dots\supset L_0(U)\supset\dots$$ of $H(2k,n)$, constructed
in Example \ref{boston}, are invariant with
respect to the Euler operator, for every choice of
the subspace $U$. 
It follows that we can construct
a filtration $$Der L=L'_{-d}(U)\supset\dots\supset L'_0(U)\supset\dots$$
 of $Der L$ by setting $L'_k(U)=L_k(U)$ for every $k\neq 0$,
and $L'_0(U)=L_0(U)+\mathbb{C}E$. Then the completion
of the graded Lie superalgebra 
associated to this filtration is isomorphic to
$H(2k,n-r_1)\otimes\Lambda(r_1)+H(0,r_1)+\mathbb{C}(E_1+E_2)$,
with respect to the grading of type
$(1,\dots,1|2,\dots,2,
1,\dots,1,0,\dots,0)$ of $H(2k,n-r_1)$,
with
$s_0$ 2's and $s_0$ zeros,
where $s_0$ and $r_1$ are defined as in Example \ref{boston}, and
where  $E_1$ and $E_2$ are  the Euler operators of  $H(2k,n-r_1)$ and
 $H(0,r_1)$, respectively.
It follows that $L'_0(U)$ is a fundamental maximal subalgebra of $Der L$.
By Remark \ref{boston}, this is, in fact, the standard subalgebra of $Der L$ stabilizing
the ideal  $I_{\cal U}=(p_1,\dots, p_k, q_1,\dots, q_k , U)$.
\end{remark}
\begin{remark}\label{DerH2}\em The proof of Theorem \ref{standardforH}
works verbatim if we replace $L=H(2k,n)$ with $Der L=H(2k,n)+\mathbb{C}E$.
Therefore, every fundamental maximal subalgebra of $Der L$ is
conjugate  to a standard subalgebra.
\end{remark}
\begin{theorem} Let $L=H(2k,n)$. Then all maximal among
$E$-invariant subalgebras
of $L$ are, up to conjugation, the 
maximal open subalgebras of $L$ listed in Theorem \ref{H(2k,n)}.
\end{theorem}
{\bf Proof.} By Remark \ref{DerH2} every fundamental maximal subalgebra
of $Der L$ is conjugate to a standard subalgebra. Therefore,
by Remark \ref{DerH1},
all maximal fundamental subalgebras of $Der L$ are, up to
conjugation, the subalgebras $L_0+\C E$, where $L_0$ is one of the
maximal open subalgebras of $L$ listed in Theorem \ref{H(2k,n)}.
Let $S_0$ be a maximal among open $E$-invariant subalgebras of $L$.
Then $S_0+\C E$ is a fundamental maximal subalgebra of $Der L$. Hence
the thesis. \hfill$\Box$

\section{Maximal open subalgebras of $\boldsymbol{KO(n,n+1)}$ and
$\boldsymbol{SKO(n,n+1;\beta)}$} 
Let $x_1,\dots, x_n$ be $n$ even indeterminates and $\xi_1,\dots,
\xi_n, \xi_{n+1}=\tau$ be $n+1$ odd indeterminates. Consider
the differential form
$\Omega=d\tau+\sum_{i=1}^n(\xi_idx_i+x_id\xi_i)$. The odd contact
superalgebra is defined as follows (\cite{ALS}):
$$KO(n,n+1)=\{X\in W(n,n+1) ~|~X\Omega=f\Omega, ~f\in\Lambda(n,n+1)\}.$$ It is a simple Lie superalgebra
for every $n\geq 1$.

Define  the following bracket on $\Lambda(n,n+1)$ (cf.\ \cite[\S
1.4]{CK}):
\begin{equation}
[f,g]=(2-E)f\frac{\partial g}{\partial
\tau}+(-1)^{p(f)}\frac{\partial f}{\partial
\tau}(2-E)g-\sum_{i=1}^n(\frac{\partial f}{\partial
x_i}\frac{\partial g}{\partial \xi_i}+(-1)^{p(f)}\frac{\partial
f}{\partial \xi_i}\frac{\partial g}{\partial x_i})
\label{ko}
\end{equation}
where $E=\sum_{i=1}^n(x_i\frac{\partial}{\partial
x_i}+\xi_i\frac{\partial}{\partial \xi_i})$ is the Euler operator. Then
the map\break $\rho: \Lambda(n,n+1)\longrightarrow KO(n,n+1)$,
$$\rho(f)=X_f:=(2-E)f\frac{\partial}{\partial
\tau}-(-1)^{p(f)}\frac{\partial f}{\partial
\tau}E-\sum_{i=1}^n(\frac{\partial f}{\partial
x_i}\frac{\partial}{\partial \xi_i}+(-1)^{p(f)}\frac{\partial
f}{\partial \xi_i}\frac{\partial}{\partial x_i})$$ is an isomorphism
between $KO(n,n+1)$ and $\Lambda(n,n+1)$ with reversed parity. We
will therefore identify $KO(n,n+1)$ with $\Lambda(n,n+1)$ with
reversed parity and fix its maximal torus $T=\langle \tau, x_i\xi_i
~|~i=1,\dots, n
\rangle$. 
\begin{remark}\label{Leibniz}\em Bracket (\ref{ko}) satisfies the
following rule:
$$[f,gh]=[f,g]h+(-1)^{p(X_f)p(g)}g[f,h]-2(-1)^{p(f)}\frac{\partial f}{\partial \tau}gh.$$
It follows, in particular, that an ideal $I=(f_1, \dots, f_r)$ of
$\Lambda(n, n+1)$ is stabilized by a function $f$ in $KO(n, n+1)$ if
and only if $[f,f_i]$ lies in $I$ for every $i=1,\dots, r$.

Besides, if $f$ is an odd function independent of $\tau$, then $\varphi=\exp ad(f)$ is an
automorphism of $\Lambda(n, n+1)$ with respect to both the Lie bracket
and the usual product of polynomials.
It follows that a subalgebra $L_0$ of $KO(n, n+1)$ stabilizes an ideal
$I=(f_1,\dots, f_r)$ of $\Lambda(n, n+1)$ if and only if the subalgebra $\varphi(L_0)$ stabilizes
the ideal $J=(\varphi(f_1), \dots, \varphi(f_r))$.
\end{remark}

For $\beta\in\C$ let $div_{\beta}:=\Delta+(E-n\beta)\frac{\partial}{\partial\tau}$, where $\Delta=\sum_{i=1}^n
\frac{\partial^2}{\partial x_i\partial\xi_i}$ is the odd Laplace operator, and let 
$$SKO^{\prime}(n,n+1;\beta)=\{f\in\Lambda(n,n+1)~|~ div_{\beta}(f)=0\}=:
\Lambda^{\beta}(n,n+1)$$
(cf.\ \cite{Ko}, \cite[Example 4.9]{K}, \cite[\S 1.4]{CK}).

\begin{remark}\label{formula}\em
If $f,g\in\Lambda(n,n+1)$, then:
$$div_{\beta}([f,g])=X_f(div_{\beta}(g))-
(-1)^{p(X_f)p(X_g)}X_g(div_{\beta}(f)).$$
It follows that the function $div_{\beta}: KO(n,n+1)\rightarrow \Lambda(n,n+1)$
is a divergence (see Definition \ref{divergence}).
Therefore $SKO'(n,n+1;\beta)$ is a subalgebra of the Lie superalgebra $\Lambda(n,n+1)$ with bracket (\ref{ko}). According to Remark
\ref{annihilating},
$$SKO'(n, n+1; \beta)=S'KO(n,n+1)=\{X\in KO(n,n+1)~|~ X \omega_{\beta}=0\}$$
where $\omega_{\beta}$ is the volume form attached to the divergence
 $div_{\beta}$.
\end{remark}
Let $SKO(n,n+1;\beta)$ denote the derived algebra of 
$SKO^{\prime}(n,$ $n+1;\beta)$. Then
$SKO(n,n+1;\beta)$ is simple for $n\geq 2$ and coincides with
$SKO^{\prime}(n,n+1;\beta)$ unless $\beta=1$ or
$\beta=\frac{n-2}{n}$. The Lie superalgebra $SKO(n,n+1;1)$ (resp.\
$SKO(n,n+1;(n-2)/n)$) consists of the elements of $SKO^\prime(n,n+1;1)$ (resp.\
$SKO^\prime(n,n+1;(n-2)/n)$) not containing the monomial $\tau\xi_1\dots
\xi_n$ (resp.\ $\xi_1\dots\xi_n$). 

Since the Lie superalgebra $KO(1,2)$ is isomorphic
to the Lie superalgebra $W(1,1)$ (cf.\ \cite[Remark 6.6]{K}), and 
since $SKO(n,n+1;\beta)$ is simple for $n\geq 2$,
when talking about $KO(n,n+1)$ and  $SKO(n,n+1;\beta)$
 we shall assume $n\geq 2$. 

\begin{remark}\label{gradingsofKO}\em
The $\Z$-grading of type $(1,\dots,
1|0,\dots, 0,1)$ of $W(n,n+1)$ induces
on
$KO(n,$ $n+1)$ (resp.\ $SKO(n,n+1;\beta)$) a grading of
depth 1 which is irreducible by Remark \ref{gd}. This grading is called
the
{\em subprincipal} grading of $KO(n,$ $n+1)$ (resp.\ $SKO(n,n+1;\beta)$). 

The $\Z$-grading of type $(1,\dots, 1,2,\dots,
2|1,\dots, 1,0,\dots, 0,2)$
of $W(n,n+1)$, with
$t+1$ 2's and $t$ zeros, induces, for every
$t=0, \dots, n-2$, an irreducible grading on $\g=KO(n,n+1)$ (resp.\ $SKO(n,n+1;\beta)$ for $(t,\beta)\neq (n-2, (n-2)/n)$) where 
$\g_0$ is isomorphic to the Lie superalgebra
$c\tilde{P}(n-t)\otimes\Lambda(t)+W(0,t)$ (resp.\ $\tilde{P}(n-t)\otimes\Lambda(t)+W(0,t)$), $\g_{-1}$ is isomorphic to
$\C^{n-t|n-t}\otimes\Lambda(t)$ and $\g_{-2}$ is isomorphic to
$\C\otimes\Lambda(t)$. When $\g=SKO(n,n+1;(n-2)/n)$ and $t=n-2$,
$\g_0$ does not contain the element $\xi_1\dots\xi_n$, and the
grading is irreducible if and only if $n>2$. These gradings satisfy the hypotheses of Proposition 
\ref{basic}$(b)$, therefore the corresponding graded subalgebras of $KO(n,n+1)$ and 
$SKO(n,n+1;\beta)$ are maximal. 

The grading of type $(1,\dots,1|1,\dots,1,2)$
is called the {\em principal} grading of $KO(n,$ $n+1)$ (resp.\ $SKO(n,n+1;\beta)$).
\end{remark}

\begin{remark}\em\label{irrenotKO(n,n+1)}
The $\Z$-grading of type $(1,
2,\dots, 2|1,0,\dots, 0,2)$
of $W(n,n+1)$ induces  on $KO(n,n+1)$ (resp.\ $SKO(n,n+1;\beta)$)
a grading
which is not irreducible. In fact the corresponding subalgebra
$\prod_{j\geq 0}\g_j$ of $KO(n,n+1)$ (resp.\ $SKO(n,n+1;\beta)$) is contained in the subalgebra of type $(1,\dots,
1|0,\dots, 0,1)$. 
\end{remark}
\begin{remark}\label{weightsforKO}\em
The subspaces  $\mathbb{C}1$, $\C x_i$, $\mathbb{C}\xi_{i_1}\dots\xi_{i_h}$,
$\mathbb{C}x_k\xi_{j_1}\dots\xi_{j_h}$, 
 with $k\neq j_1,\dots, j_h$, 
$\C\xi_{i_1}\dots\xi_{i_h}\otimes T$,
and $\C x_k\xi_{j_1}\dots\xi_{j_h}\otimes T$ with 
$k\neq j_1,\dots ,j_h$, are $T$-weight spaces of $KO(n,n+1)$.
\end{remark}
\begin{remark}\label{stabKO}\em Let $L=KO(n,n+1)$. Then
the graded subalgebra $L_k$ of $L$ of type $(1,\dots,1, 2,\dots, 2$ 
$|1,\dots, 1, 0, \dots, 0,2)$ with $n-k+1$ 2's and $n-k$ zeros, is,
for every $k=1,\dots, n$, the standard
subalgebra $L_U$ of $L$ stabilizing the ideal 
$I_U=(x_1,\dots, x_n, \xi_1, \dots,
\xi_k, \tau)$. Indeed, for every $k$, $L_k\subset L_U$
since $L_k$
is contained in the graded subalgebra of $W(n,n+1)$ of type $(1,\dots,1|1,\dots,1,0,\dots,0,1)$ with $n-k$ zeros, which stabilizes $I_U$ (cf.\ the proof of Theorem
\ref{W(m,n)}). Since, for every $k\neq 1$, $L_k$ is a maximal  
subalgebra of $L$ (cf.\ Remark \ref{gradingsofKO}), 
$L_k=L_U$.

Now suppose $k=1$. Notice that $L_U$ contains the standard
torus $T$ of $KO(n,n+1)$, hence it is regular and decomposes
into the direct product of $T$-weight spaces. 
The subspace $S=<1, x_1, \xi_1>\otimes\Lambda(\xi_2, \dots, \xi_n)$
is a $T$-invariant complementary subspace of the subalgebra $L_1$ and,
according to Remark \ref{weightsforKO},
the subspaces $\mathbb{C}1$, $\mathbb{C}\xi_{j_1}\dots\xi_{j_h}$,
$\mathbb{C}x_1\xi_{i_1}\dots\xi_{i_h}$, 
 with $1\neq i_1\neq\dots\neq i_h$, are one-dimensional 
$T$-weight spaces.
Therefore, in order to prove that $L_U\subset L_1$, it is sufficient to show 
that no element of $S$ lies in $L_U$. Notice that
$L_U$ contains the elements $x_2, \dots, x_n$ but it does not contain
the elements $1, x_1, \xi_i$ for any $i=1, \dots, n$. 
Since
$[x_{i_1}, \xi_{i_1}\dots \xi_{i_h}]=-\xi_{i_2}\dots\xi_{i_h}$ and
$[x_{i_1}, x_j\xi_{i_1}\dots\xi_{i_h}]=-x_j\xi_{i_2}\dots\xi_{i_h}$, it follows
 that 
$S$ is contained in the $T$-invariant complementary subspace of $L_U$, 
therefore
$L_U\subset L_1$, hence $L_U=L_1$.

Likewise, the graded subalgebra of $L$ of type $(1,\dots,1|0,\dots,0,1)$
is the standard subalgebra of $L$ stabilizing the ideal $(x_1, \dots, x_n, \tau)$.
\end{remark}

\begin{example}\label{ex1KO}\em Throughout this example 
we shall identify
$KO(n,n+1)$ with $\Lambda(n,n+1)$ as described at the beginning of this paragraph. On $\Lambda(n,n+1)$ we define a valuation  with values in $\Z_+$ 
by setting:
$$\nu(x_i)=1, ~~\nu(\xi_i)=0, ~~i=1,\dots, n, ~~\nu(\tau)=0,$$
(see Remark \ref{valuation1}).
Consider the following subspaces 
of $KO(n,n+1)$:
$$L_i=\{f\in\Lambda(n,n+1)~|~\nu(f)\geq i+1\}+\Lambda(\tau) ~~\mbox{for}~ i=-1, 0;$$
$$L_i=\{f\in\Lambda(n,n+1)~|~\nu(f)\geq i+1\} ~~\mbox{for}~ i> 0.$$
Using commutation rules (\ref{ko}) one can check that the subspaces $L_i$
define in fact a filtration of $KO(n,n+1)$ of depth 1 whose associated
graded superalgebra $Gr L$ has the following structure:
$$Gr_{-1}L=\Lambda(\xi_1,\dots, \xi_n, \tau)/\Lambda(\tau);$$
$$Gr_0L=\langle x_i\rangle\otimes\Lambda(\xi_1,\dots, \xi_n, \tau)+\Lambda(\tau);$$
$$Gr_jL=\langle f\in\C[[x_1,\dots, x_n]]~|~\deg(f)=j+1\rangle\otimes\Lambda
(\xi_1,\dots, \xi_n, \tau) ~~\mbox{for}~ j\geq 1.$$
It follows that $\overline{Gr L}\cong HO(n,n)\otimes\Lambda(\eta)+\C\frac{\partial}
{\partial\eta}+\C(E-2+2\eta\frac{\partial}
{\partial\eta})$ 
with respect to the grading of type $(1,\dots, 1|0,\dots, 0)$ of $HO(n,n)$,
where $E=\sum_{i=1}^n(x_i\frac{\partial}{\partial x_i}+\xi_i\frac{\partial}{\partial \xi_i})$, and $\deg(\eta)=0$.
Since this grading is irreducible (cf.\ Remark \ref{gradingsofHO}) and satisfies property $(iii)^{\prime}$ of Proposition \ref{basic}$(b)$,  $L_0$ is a maximal subalgebra
of $KO(n,n+1)$ by Corollary \ref{cor}.

Note that the subalgebra $L_0$ stabilizes the ideal $I_U=(x_1, \dots, x_n)$ 
of $\Lambda(n,n)$, hence, due to its maximality, $L_0$ is the
standard subalgebra of $KO(n,n+1)$ corresponding to the ideal $I_U$.
\end{example}
\begin{example}\label{ex2KO}\em Throughout this example 
we shall   identify
$KO(n,n+1)$ with $\Lambda(n,n+1)$ as above. Let us fix an integer $t$ such 
that $1\leq t\leq n$ and let us define a valuation on $\Lambda(n,n+1)$ by  setting:
$$\nu(x_i)=1, ~~\nu(\xi_i)=1, ~\mbox{for}~ i=1,\dots,t;$$
$$\nu(\tau)=0, ~~\nu(x_i)=2, ~~\nu(\xi_i)=0, ~\mbox{for}~ i=t+1,\dots,n.$$
Consider the following subspaces 
of $KO(n,n+1)$:
$$L_i(t)=\{f\in\Lambda(n,n+1)~|~\nu(f)\geq i+2\}+\Lambda(\tau) ~~\mbox{for}~ i\leq  0;$$
$$L_i(t)=\{f\in\Lambda(n,n+1)~|~\nu(f)\geq i+2\} ~~\mbox{for}~ i> 0.$$
Using commutation rules (\ref{ko}) one verifies that the subspaces $L_i(t)$
define in fact a filtration of $KO(n,n+1)$. This filtration has depth 1 if $t=n$, otherwise it has depth 2. 
We have: $\overline{Gr L}\cong HO(n,n)\otimes\Lambda(\eta)+\C\frac{\partial}
{\partial\eta}+\C(E-2+2\eta\frac{\partial}
{\partial\eta})$ 
with respect to the grading of type $(1,\dots,1,2,\dots, 2|1,\dots, 1,0,\dots,0)$ of $HO(n,n)$, with
$n-t$ 2's and $n-t$ zeros, and $\deg(\eta)=0$.
Since this grading is irreducible for every $t=2,\dots, n$
(cf.\ Remark \ref{gradingsofHO}) and satisfies property $(iii)^{\prime}$ of Proposition \ref{basic}$(b)$,
by Corollary \ref{cor}  $L_0(t)$ is a maximal (regular) subalgebra
of $KO(n,n+1)$ for every $t=2,\dots, n$. On the contrary, the subalgebra
$L_0(1)$ is contained in the subalgebra $L_0$ of $KO(n,n+1)$
constructed in Example \ref{ex1KO}, hence it is not maximal.

Note that $L_0(t)$  is contained in the graded
subalgebra of $W(n,n+1)$ of type $(1,\dots, 1|1,\dots,1,0,\dots,0)$
with $n+t$ $1$'s, therefore it stabilizes the ideal $I_{\cal{U}}=(x_1, \dots, x_n, 
\xi_1, \dots, \xi_t)$ of $\Lambda(n,n+1)$. It follows that, for every $t=2, \dots, n$,
$L_0(t)$ is the standard subalgebra of $KO(n,n+1)$ corresponding
to the ideal $I_{\cal{U}}$, due to its maximality.

Likewise, $L_0(1)$ is the standard subalgebra $L_{{\cal{U}}_1}$ of
$KO(n,n+1)$ stabilizing the ideal $I_{{\cal{U}}_1}=(x_1, \dots, x_n, \xi_1)$.
Indeed,  $L_{{\cal{U}}_1}$ contains the standard torus $T$ of $KO(n,n+1)$,
hence it is regular and decomposes into the direct product of
$T$-weight spaces. By definition $L_{{\cal{U}}_1}$
contains the elements $1, x_2, \dots, x_n$ and does not contain
the elements $x_1$ and $\xi_j$ for any $j=1, \dots, n$.
 Notice that
$Gr_{<0}L:=L_{-2}(1)/L_0(1)=
(\langle 1, x_1, \xi_1\rangle \otimes\Lambda(\xi_2, \dots, \xi_n, \tau))/
\Lambda(\tau)$.
Then the same arguments as in Remark \ref{stabKO} show that
no element in $\langle 1, x_1, \xi_1\rangle \otimes\Lambda(\xi_2, \dots, \xi_n))/\C 1$ lies in $L_{{\cal{U}}_1}$. 
Now suppose that an element of the form $\xi_{i_1}\dots\xi_{i_k}\tau+\varphi$
lies in $L_{{\cal{U}}_1}$  for some $\varphi\in\Lambda(n,n)$, where by
$\Lambda(n,n)$ we mean the subalgebra of $\Lambda(n,n+1)$ generated by all
even indeterminates and by the odd indeterminates except $\tau$.
Then $L_{{\cal{U}}_1}$ contains the element $[1, \xi_{i_1}\dots\xi_{i_k}\tau+\varphi]=
\pm 2 \xi_{i_1}\dots\xi_{i_k}$ and this is a contradiction.
Therefore $L_{{\cal{U}}_1}$ cannot contain any element of the form
$\xi_{i_1}\dots\xi_{i_k}\tau+\varphi$
for any function $\varphi\in\Lambda(n,n)$ and, similarly, it cannot
contain any element of the form $x_1\xi_{i_1}\dots\xi_{i_k}\tau+\varphi$
for any $i_1\neq\dots\neq i_k\neq 1$ and any function $\varphi\in\Lambda(n,n)$.
By Remark \ref{weightsforKO}
 it follows
that $L_{{\cal{U}}_1}$ is contained in $L_0(1)$,
hence $L_{{\cal{U}}_1}=L_0(1)$.  
\end{example}

\begin{remark}\label{coord1}\em Let $1\leq i<j\leq n$. Then the change
of indeterminates that leaves $\tau$ invariant and exchanges $x_i$ with $x_j$ and $\xi_i$ with $\xi_j$,
preserves the form $\Omega$. 
\end{remark}
\begin{remark}\label{coord2}\em
Let $\eta=\alpha_{i_1}\xi_{i_1}+\dots+\alpha_{i_k}\xi_{i_k}$
for some $k\leq n$, with $\alpha_{i_j}\in\mathbb{C}$, $\alpha_{i_j}\neq 0$.
According to Remark \ref{coord1}, up to changes of variables, we can assume 
$\eta=\alpha_{1}\xi_{1}+\dots+\alpha_{k}\xi_{k}$ with
$\alpha_i\neq 0$ for $i=1,\dots,k$.
Then the following change of indeterminates preserves the form
$\Omega$:
$$\begin{array}{ll}
\tau'=\tau & \\
x'_1=\frac{1}{\alpha_1}x_1 & \xi'_1=\eta\\
x'_2=x_2-\frac{\alpha_2}{\alpha_1}x_1 & \xi'_2=\xi_2\\
\vdots & \vdots\\
x'_k=x_k-\frac{\alpha_k}{\alpha_1}x_1 & \xi'_k=\xi_k\\
x'_i=x_i & \xi'_i=\xi_i ~~\forall i>k.
\end{array}$$
\end{remark}
\begin{theorem}\label{standardforKO} Let $L_0$ be a maximal open subalgebra of $L=KO(n,n+1)$.
Then $L_0$ is conjugate to a standard subalgebra of $L$.
\end{theorem}
{\bf Proof.} By Remark \ref{rome} $L_0$ stabilizes an ideal of the form
$$J=(x_1+f_1, \dots, x_n+f_n,
 \eta_1+g_1, \dots, \eta_s+g_s)$$
for some
linear functions $\eta_j$ in odd indeterminates, and even functions
$f_i$ and odd functions $g_j$ without constant and linear terms,
 and $J$ is maximal among the $L_0$-invariant
ideals of $\Lambda(n,n+1)$.

We distinguish the following two cases:

\medskip

\noindent
Case 1: $\eta_i$ lies in $\Lambda(\xi_1,\dots, \xi_n)$ for every
$i=1,\dots, s$.
By Remark \ref{coord2}, up to changes of indeterminates, we have:
$$J=(x_1+F_1, \dots, x_n+F_n,
 \xi_1+G_1,\dots, \xi_{s}+G_{s})$$
for some even functions
$F_i$ and odd functions $G_j$ without constant and linear terms,
where the functions $G_j$'s  are independent of 
$\xi_1, \dots, \xi_{s}$, for every $j=1,\dots, s$
and where, since the ideal $J$
is closed, 
we can assume the functions $F_i$ and $G_i$ independent of
all even indeterminates, i.e., $F_i, G_i\in\Lambda(0,n+1)$
(cf.\ Lemma \ref{victorshelp}).

Suppose that $x_1+F_1=x_1+\xi_1F'_1+F^{\prime\prime}_1$ with
$F'_1$ and $F^{\prime\prime}_1$ independent of $\xi_1$.
Then we can replace $x_1+F_1$ by $x_1+F_1-(\xi_1+G_1)F'_1=x_1+H_1$
with $H_1$ independent of $\xi_1$. Similarly we can make
every function 
$F_i$ independent of $\xi_j$ for every $j=1,\dots,s$.

Now suppose $x_1+F_1=x_1+\tau\varphi_0+\varphi_1$ with $\varphi_0$ and $\varphi_1$
independent of $\tau$. Notice that, although the map $ad(\tau\xi_1\varphi_0)$ is
not a derivation of $\Lambda(n,n+1)$ with respect to the usual product, the map 
$\psi:=ad(\tau\xi_1\varphi_0)+2\xi_1\varphi_0id$
is a derivation, as one can verify using Remark \ref{Leibniz}. 
Thus $\exp(\psi)$ is an automorphism of $\Lambda(n,n+1)$
with respect both to bracket (\ref{ko}) and to the usual product.
Notice that  $\exp(\psi)(x_1+F_1)=x_1+\Phi_1$ 
 for some function $\Phi_1$ independent of $\tau$.
Thus, up to automorphisms, we can assume $F_1$ and, similarly, every function $F_i$,  for every $i=1, \dots, s$, independent of $\tau$.
As a consequence, the map
$\exp(ad(-\xi_1F_1))$ is an automorphism of $\Lambda(n,n+1)$, mapping $J$
to the ideal 
$$I=(x_1, x_2+F'_2, \dots, x_n+F'_n,
 \xi_1+G_1, \dots, \xi_s+G_s).$$
Arguing in the same way for every function $F'_j$ with
$1\leq j\leq s$, we have, up
to automorphisms,
$$I=(x_1, \dots, x_s, x_{s+1}+h_{s+1}, \dots, x_n+h_n,
 \xi_1+G_1, \dots, \xi_{s}+G_s)$$
for some functions $h_{s+1}, \dots, h_n\in \Lambda(\xi_{s+1},\dots, \xi_n, \tau)$.

Suppose $G_1=\tau \rho_0+\rho_1$ with $\rho_0$, $\rho_1$ independent of $\tau$.
Then $\exp(ad(x_1\tau \rho_0)+2x_1 \rho_0 id)$ is an automorphism of $\Lambda(n,n+1)$
 mapping the ideal $I$ to
$$I'=(x_1, \dots, x_s, x_{s+1}+h_{s+1}, \dots, x_n+h_n,
 \xi_1+\rho_1, \xi_2+G'_2, \dots, \xi_{s}+G'_s)$$
where $\rho_1$ is independent of $\tau$.
Arguing in the same way for every function $G_j$ we can assume, up
to automorphisms, that 
$$I=(x_1, \dots, x_s, x_{s+1}+h_{s+1}, \dots, x_n+h_n,
 \xi_1+\rho_1, \dots, \xi_{s}+\rho_s)$$
where
$\rho_j$ lies in $\Lambda(\xi_{s+1},\dots, \xi_n)$
for every $j$. 
It follows that the map $\exp(ad(-x_1\rho_1))$ is an automorphism of $L$
sending the ideal $I$ to the ideal
$$Y=(x_1, \dots, x_s, x_{s+1}+h'_{s+1}, \dots, x_n+h'_n,
 \xi_1, \xi_2+\rho'_2, \dots, \xi_{s}+\rho'_s),$$
for some functions $h'_i\in\Lambda(\xi_{s+1}, \dots, \xi_n, \tau)$,
$\rho'_j\in\Lambda(\xi_{s+1}, \dots, \xi_n)$. Analogous automorphisms
yield to the ideal
$$Y'=(x_1, \dots, x_s, x_{s+1}+h''_{s+1}, \dots, x_n+h''_n,
 \xi_1, \dots, \xi_{s}),$$
for some functions $h''_i\in\Lambda(\xi_{s+1}, \dots, \xi_n, \tau)$.

Let $h''_{s+1}=\xi_{s+1}\psi_1+\psi_2$ for some $\psi_1, \psi_2$ independent
of $\xi_{s+1}$. By the same argument as above we can assume $\psi_2$
independent of $\tau$ and, applying the automorphism $\exp(ad(\xi_{s+1}\psi_2))$,
we can assume $\psi_2=0$.
Now the proof can be concluded as in the case of the Lie superalgebra
$HO(n,n)$ (cf.\ Theorem \ref{standardforHO}). Namely, 
let  $Y''=(x_{s+1}+h_{s+1}, \dots,  x_n+h_n)\subset Y'$.
Then, each function $f$ in $L_0$ (thus stabilizing $Y$)
stabilizes the ideal
$K=(Y', [Y'',Y''])$, i.e., the ideal
generated by the generators of $Y'$ and 
by
the commutators between every pair of generators of $Y''$. Indeed, for every $g, h\in Y''$
we have:
$$[f,[g,h]]=[[f,g],h]\pm[g,[f,h]]\in[Y', Y'']$$
and $[Y', Y'']\subset K$.
Notice that $K$ is a proper ideal of
$\Lambda(2k+1,n)$ since among its generators there is
no invertible element. By the maximality of
$J$ among the ideals stabilized by $L_0$ we have 
$Y'=K$.

Now $1/2[x_{s+1}+\xi_{s+1}\psi_1, x_{s+1}+\xi_{s+1}\psi_1]=-\psi_1+\xi_{s+1}\tilde{\varphi}\in
[Y'',Y'']\subset K=Y'$, therefore $(\psi_1-\xi_{s+1}\tilde{\varphi})\xi_{s+1}=
\psi_1\xi_{s+1}$ lies in $Y'$. It follows that
$$Y'=(x_1, \dots, x_s, x_{s+1}, x_{s+2}+h'_{s+2}, \dots,  x_n+h'_n,
 \xi_1, \xi_2,\dots, \xi_s).$$
Arguing in the same way for every function $h'_j$, we end up with a standard
ideal.

\medskip

Case 2: there exists one $i$ such that
$\eta_i=\tau+\eta$ with $\eta\in\Lambda(\xi_1,\dots, \xi_n)$, i.e.,
up to changes of indeterminates,
$$J=(x_1+f_1, \dots, x_n+f_n,
 \xi_1+g_1, \dots, \xi_{s-1}+g_{s-1}, \tau+\eta_{s}+g_{s})$$
for some
linear function $\eta_s$ in $\Lambda(\xi_1,\dots, \xi_n)$, and even functions
$f_i$ and odd functions $g_j$ without constant and linear terms.
We can assume $f_i, g_j$ and $\eta_s$
in $\Lambda(\xi_s,\dots, \xi_n)$.

Besides, arguing similarly as above and as in the proof of
Theorem \ref{standardforHO}, one shows that, up
to automorhisms,
$$J=(x_1, \dots,x_{s-1}, x_{s}+h_{s}, \dots, x_n+h_n,
 \xi_1, \dots, \xi_{s-1}, \tau+\eta_{s}+H)$$
for some functions $h_i, \eta_s, H\in\Lambda(\xi_{s}, \dots, \xi_{n})$.

(i) Suppose $\eta_s=0$. 
Denote by
$J'$ the ideal $J'=(x_{s}+h_{s}, \dots,  x_n+h_n, \tau+H)\subset J$.
Then, each function $f$ in $L_0$ (thus stabilizing $J$)
stabilizes the ideal
$K=(J, [J',J'])$, i.e., the ideal
generated by the generators of $J$ and 
by
the commutators between every pair of generators of $J'$. Indeed, for every $g, h\in J'$
we have:
$$[f,[g,h]]=[[f,g],h]\pm[g,[f,h]]\in[J, J']$$
and $[J, J']\subset K$.
Notice that $K$ is a proper ideal of
$\Lambda(n, n+1)$ since among its generators there is
no invertible element. By the maximality of
$J$ among the ideals stabilized by $L_0$ we have 
$J=K$.

Suppose that $h_{s}=\xi_{s}\psi_1+\psi_2$ with 
$\psi_1$ and $\psi_2$ independent of $\xi_{s}$. Then
applying the automorphism $\exp(ad(-\xi_{s}\psi_2))$ we can assume
$$J=(x_1, \dots, x_{s-1}, x_{s}+\xi_{s}\psi_1, x_{s+1}+h'_{s+1}, \dots,  x_n+h'_n,
 \xi_1, \xi_2,\dots, \xi_{s-1}, \tau+H').$$
Now $\psi_1=-1/2[x_{s+1}+\xi_{s+1}\psi_1, x_{s+1}+\xi_{s+1}\psi_1]\in
[J',J']\subset K=J$, therefore
$$J=(x_1, \dots, x_{s-1}, x_{s}, x_{s+1}+h'_{s+1}, \dots,  x_n+h'_n,
 \xi_1, \xi_2,\dots, \xi_{s-1}, \tau+H').$$
Repeating a  similar argument for every function $h'_j$ and, finally,
for the function $H'$, we end up 
with the standard ideal $J=(x_1, \dots, x_n,
 \xi_1, \xi_2,\dots, \xi_{s-1}, \tau).$

(ii) If $\eta_s\neq 0$, by Remark \ref{coord2}, we can assume:
$$J=(x_1, \dots,x_{s-1}, x_{s}+h_{s}, \dots, x_n+h_n,
 \xi_1, \dots, \xi_{s-1}, \tau+\xi_{s}+H).$$
Thus the automorphism $\exp(ad(x_s\tau)+2x_sid)$ maps $J$ to the
ideal 
$$J'=(x_1, \dots,x_{s-1}, x_{s}+h'_{s}, \dots, x_n+h'_n,
 \xi_1, \dots, \xi_{s-1}, \xi_{s}+\tau\rho+H')$$
where $H'$ is independent of  $\tau$ and $\deg(\rho)\geq 1$.
In the limit, since $J'$ is closed, we get the ideal
$$J^{\prime\prime}=(x_1, \dots,x_{s-1}, x_{s}+h_{s}, \dots, x_n+h_n,
 \xi_1, \dots, \xi_{s-1}, \xi_{s}+M)$$
where $M$ is independent of $\tau$.
 We thus proceed as in case 1.
 \hfill$\Box$

\begin{theorem}\label{KO} All maximal open subalgebras of
$L=KO(n,n+1)$ are, up to conjugation, the following:

\begin{itemize}
\item[$(i)$] the graded subalgebra of type $(1,\dots,
1|0,\dots, 0,1)$;
\item[$(ii)$] the graded subalgebras of type $(1,\dots, 1,2,\dots,
2|1,\dots, 1,0,\dots, 0,2)$ with
$n-t+1$ 2's and $n-t$ zeros,
for $t=2, \dots, n$;
\item[$(iii)$] the non-graded subalgebra $L_0$ described in Example \ref{ex1KO} and the
non-graded subalgebras $L_0(t)$ described in Example \ref{ex2KO} for $t=2,\dots, n$.
\end{itemize}
\end{theorem}
{\bf Proof.} Let $L_0$ be
a maximal open subalgebra of $L$. By Theorem \ref{standardforKO},
$L_0$ is, up to conjugation, the standard subalgebra of $L$ stabilizing
either the ideal $I_{\cal{U}}=(x_1, \dots, x_n, \xi_1,
\dots, \xi_s)$ for some $s=0, \dots,n$, or the ideal $I_{\cal{U}'}=
(x_1, \dots, x_n, \xi_1,$ $\dots,\xi_t, \tau)$ for some $t=0, \dots, n$.
The statement then follows using Remarks \ref{stabKO}, \ref{gradingsofKO}, 
\ref{irrenotKO(n,n+1)}, and Examples  \ref{ex1KO}, \ref{ex2KO}.
\hfill $\Box$

\begin{corollary} All irreducible $\Z$-gradings of  $KO(n,n+1)$ are, up to conjugation, 
 the grading of
type $(1,\dots,
1|0,\dots, 0,1)$ and
the gradings of type $(1,\dots, 1,$ $2,\dots,
2|1,\dots, 1,0,\dots, 0,2)$ with
$t+1$ 2's and $t$ zeros, 
for $t=0, \dots, n-2$.
%
\end{corollary} 

We shall now focus on the Lie superalgebra $SKO(n,n+1;\beta)$ introduced
at the beginning of this section.

\begin{remark}\label{irrsko2}\em 
The $\Z$-grading  of type $(1,\dots, 1|-1, \dots,-1,0)$ of $W(n,n+1)$
induces on $\g=SKO(n,n+1;\beta)$
a $\Z$-grading  $\g=\prod_j \g_j$, where
$$\g_{-1}=\{f\in\oplus_{h=0}^{n-2}\langle x_{i_1}\dots x_{i_h}\xi_{j_1}\dots
\xi_{j_{h+1}}\rangle\otimes\langle 1, \tau, \Phi\rangle ~|~  
div_{\beta}(f)=0\}$$ and 
$$\g_0=\{f\in\oplus_{h=0}^{n-1}\langle x_{i_1}\dots x_{i_h}\xi_{j_1}\dots
\xi_{j_{h}}\rangle\otimes\langle 1, \tau, \Phi\rangle ~|~ 
div_{\beta}(f)=0\}+\langle 1, \tau+\beta\Phi\rangle$$
where $\Phi=\sum_{i=1}^nx_i\xi_i$.
One can check that 
$S=\{f\in\oplus_{h=1}^{n-2}\langle x_{i_1}\dots x_{i_h}\xi_{j_1}\dots
\xi_{j_{h+1}}
\rangle$ $\otimes \langle 1, \tau, \Phi\rangle ~|~  div_{\beta}(f)=0\}$,
is a $\g_0$-stable subspace of $\g_{-1}$. Notice that $S=0$ if
and only if $n=2$. It follows that
for $n>2$  the grading of type
$(1,\dots, 1|$ $-1,\dots,-1,0)$ induces on $SKO(n,n+1,\beta)$ a grading
which is not irreducible.

Now suppose $n=2$ and $\beta\neq 0$. Then
the $\Z$-grading of type $(1,1|-1,-1,0)$  has depth $2$. 
One has:  $\g_0\cong sl_2\otimes\Lambda(1)+W(0,1)$,
$\g_{-1}\cong\C^2\otimes\Lambda(1)$, where $\C^2$ is the standard
$sl_2$-module, and
$\g_{-2}=\C\xi_1\xi_2=[\g_{-1}, \g_{-1}]$. It follows
that the grading of type $(1,1|-1,-1,0)$ of $\g=SKO(2,3;\beta)$  is irreducible. The $\g_0$-module $\g_1$ consists of the elements $f\in\langle
x_i, x_ix_j\xi_k\rangle\otimes\langle 1, \tau, \Phi\rangle$ such that 
$div_{\beta}(f)=0$.
Notice that $\g_1$ is not irreducible: 
it has an irreducible $\mathfrak{g}_0$-submodule
$S\cong S^3(\C^2)\otimes\Lambda(1)$ and $\g_1/S\cong \C^2\otimes\Lambda(1)$. 
Besides, for every $j>1$, $\g_j=\g_1^j$. 
 One can check that
property $(iii)^\prime$ of Proposition \ref{basic}$(b)$ is satisfied,
hence $\prod_{j\geq 0}\g_{j}$ is a maximal subalgebra of $\g$.

Finally, if $n=2$ and $\beta=0$ the grading of type $(1,1|-1,-1,0)$
has depth 1, hence it is irreducible by Remark \ref{gd}.
\end{remark}
\begin{remark}\label{outerSKO(2,3;1)}\em When $\beta\neq 0, -1$
the even part of the Lie superalgebra \break $SKO(2,3;\beta)$ is 
isomorphic to $W(2,0)$ and its odd part is isomorphic to \break
$\Omega^0(2)^{-\frac{1}{\beta+1}}\oplus \Omega^0(2)^{-\frac{\beta}{\beta+1}}$
(cf.\ Definition \ref{twisted}).
It follows that, when $\beta=1$,  $SKO(2,3;\beta)_{\bar{1}}$ is the
direct sum of two irreducible $SKO(2,3;\beta)_{\bar{0}}$-submodules
each of which is isomorphic to $\Omega^0(2)^{-1/2}$ (cf.\ \cite[Proposition 5.3.4]{CK}). Therefore
if $S=SKO(2,3;1)$, then $Der S=S+\mathfrak{a}$ with
$\mathfrak{a}\cong sl_2$ (cf.\ \cite[Proposition 6.1]{K},
Proposition \ref{derivations}).
Let $e,h,f$ be the standard basis of $\mathfrak{a}$ where
$e=ad(\xi_1\xi_2\tau)$ and 
$h=ad(\sum_{i=1}^2x_i\xi_i)$. We will denote by
$\mathfrak{b}$ the
subalgebra of $\mathfrak{a}$ spanned by $e$ and $h$.
\end{remark}
\begin{remark}\label{inf.manySKO}\em Let $S=SKO(2,3;1)$.
Let us denote by $S_0$ the intersection
between the graded subalgebras of $S$ of type $(1,1|1,1,2)$ and 
$(1,1|-1,-1,0)$, and  let $S=S_{-2}\supset S_{-1}\supset S_0\supset \dots$
be the  Weisfeiler filtration associated to $S_0$, where
$S_{-1}=\langle 1, x_i, \xi_1\xi_2, \xi_i(\tau+\Phi) | i=1,2\rangle+S_0$.
  Then $Gr S$ is a graded
Lie superalgebra of depth 2 where $Gr_0S\cong S(0,2)+\mathbb{C}E$
and $Gr_{-1}S$
 is isomorphic, as a $Gr_0S$-module, to the direct sum of two
copies of $\Lambda(2)/\mathbb{C}1$. Let $V$ be the subspace of
$Gr_{-1}S$ spanned by the elements $\xi_1\xi_2$ and $\xi_i(\tau+\Phi)$ for
$i=1,2$. Then $V$ is a $Gr_0S$-submodule of $Gr_{-1}S$ and
$\prod_{j\geq 0}Gr_j S+V$ is the graded subalgebra of $S$ of type $(1,1|1,1,2)$.
Likewise, for every $\gamma\in\C$,  the subspace
$V_{\gamma}=\langle 1+\gamma\xi_1\xi_2, -2x_1+\gamma\xi_2(\tau+\Phi),
2x_2+\gamma\xi_1(\tau+\Phi)\rangle$ is a $Gr_0S$-submodule of $Gr_{-1}S$ and
$\prod_{j\geq 0}Gr_j S+V_0$ is the graded subalgebra of $S$ of type $(1,1|-1,-1,0)$.
Notice that, for every $\gamma\neq 0$, the automorphism
$\exp(\frac{\gamma}{2} e)$ maps $V_{\gamma}$ to $V_0$. It
follows that every subalgebra $S_{\gamma}:=\prod_{j\geq 0}Gr_j S+V_{\gamma}$, with
$\gamma\in\C$, is conjugate to the maximal subalgebra of type
$(1,1|-1,-1,0)$.
On the other hand,
the grading of type $(1,1|-1,-1,0)$ is 
conjugate to the grading of type $(1,1|1,1,2)$ by the
automorphism $\exp(e)\exp(-f)\exp(e)$. 
Therefore the maximal subalgebras of $S$ of type $(1,1|-1,-1,0)$
and $(1,1|1,1,2)$ lie in the same $G$-orbit. This orbit consists of
the  subalgebras $S_{\gamma}$, with 
$\gamma\in\C$, and of the subalgebra of principal type, and
the intersection of any pair of subalgebras in this orbit is
the subalgebra $\prod_{j\geq 0}Gr_j S$. Notice that $\prod_{j\geq 0}Gr_j S$ is
contained also in the (maximal) subalgebra of type $(1,1|0,0,1)$.  
\end{remark}
\begin{remark}{\em 
If $\beta\neq -1$, then
the subalgebra of $SKO^{\prime}(n,n+1;\beta)$ consisting of the elements
$f\in \langle P\xi_k, Q\tau
~|~ P, Q\in\C[[x_1, \dots, x_n]]\rangle$ such that $div_{\beta}(f)=0$ is 
isomorphic to $W_n$.}
\end{remark}

\begin{remark}\label{outerSKO(2,3;0)}\em 
The even part of the Lie superalgebra $SKO(2,3;0)$ is 
isomorphic to $W(2,0)$ and its odd part is isomorphic to
$\Omega^0(2)^{-1}\oplus \Omega^0(2)/\C 1$.
The outer derivation $D=ad(\xi_1\xi_2)$ of $SKO(2,3;0)$ can then
be described as follows. Let $p:\Omega^0(2)\rightarrow \Omega^0(2)/\C 1$
be the natural projection. Then:
$$D(X)=p(div(X)) ~~\mbox{if}~~ X\in W(2,0);$$
$$D(f)=df ~~\mbox{if}~~ f\in \Omega^0(2)^{-1};$$
$$D(f)=0 ~~\mbox{if}~~ f\in \Omega^0(2)/\C 1.$$
The image of $D$ is thus given by $(\Omega^1(2)_{closed})^{-1}+\Omega^0(2)/\C 1$
where $(\Omega^1(2)_{closed})^{-1}$
can be identified with $S(2,0)$ via contraction with the
volume form $dx_1\wedge dx_2$.
\end{remark}

\begin{remark}\em
Let us describe the structure of the
Lie superalgebra  $SKO(2,3;$ $-1)$.
Its even part is not simple: 
it has a commutative ideal 
consisting of elements in $\Omega^0(2)(\tau-\Phi)$.
We have:
$$SKO(2,3;-1)_{\bar{0}}\cong\Omega^0(2)\rtimes S(2,0);
~~~~SKO(2,3;-1)_{\bar{1}}\cong\Omega^0(2)+\Omega^0(2).$$
Here $S(2,0)$ acts on each odd copy of $\Omega^0(2)$ in the natural way,
and the even functions in $\Omega^0(2)$ act by multiplication on
one copy and by $-$multiplication on the other.
\end{remark}

\begin{example}\label{ex1sko}\em 
Throughout this example we shall consider
the Lie superalgebra $S'=SKO'(n,n+1;\beta)$ for $n>2$, and we shall
identify it with $\Lambda^{\beta}(n,n+1)$ as explained at the
beginning of this section.

Notice that $\Lambda^{\beta}(n,n+1)\subset
\Lambda^{\Delta}(n,n)\otimes\langle
1, \tau, \Phi\rangle$, where  
$\Lambda^{\Delta}(n,n)=\{f\in\Lambda(n,n)~|~ \Delta(f)=0\}$
and 
$\Phi=\sum_{i=1}^n x_i\xi_i$.
We define a valuation
$\nu$ on $\Lambda^{\Delta}(n,n)\otimes\langle
1, \tau, \Phi\rangle$ (hence on $\Lambda^{\beta}(n,n+1)$) 
by setting:
$$\nu(1)=\nu(\tau)=\nu(\Phi)=0$$
$$\nu(x_i)=1 ~~\forall~i=1,\dots,n, ~~\nu(\xi_{i_1}\dots\xi_{i_k})=0 ~\forall~k<n, ~~\nu(\xi_1...\xi_n)=-1$$
and we extend it  on $\Lambda(0,n)$ by property
$b)$ in Remark \ref{valuation1}, on $\C[[x_1, \dots, x_n]]$ by
properties $a)$ and $b)$ in Remark \ref{valuation1}, and finally  on
 $\Lambda^{\Delta}(n,n)\otimes\langle
1, \tau, \Phi\rangle$ 
by setting
$\nu(\sum_iP_i(x)Q_i(\xi)\eta_i)=\min_i(\nu(P_i(x))+\nu(Q_i(\xi)))$ where
$P_i(x)\in\C[[x_1, \dots, x_n]]$, $Q_i(\xi)\in
\Lambda(0,n)$ and
$\eta_i\in\langle
1, \tau, \Phi\rangle$.

Then the following subspaces define a filtration of $SKO'(n,n+1;\beta)$:
$$S'_j=\{f\in \Lambda^{\beta}(n,n+1) ~|~ \nu(f)\geq j+1\}+\langle 1, 
\tau+\beta\Phi\rangle ~~\mbox{if}~
j\leq 0;$$
$$S'_j=\{f\in \Lambda^{\beta}(n,n+1) ~|~ \nu(f)\geq j+1\} ~~\mbox{if}~ j>0.$$

This filtration has depth $2$, with
$Gr_{-2}S'=\langle\xi_1\dots\xi_n\rangle$ if $\beta\neq 1$ and
$Gr_{-2}S'=\langle\xi_1\dots\xi_n,\xi_1\dots\xi_n\tau \rangle$ if
$\beta=1$. In fact
$Gr_{-2}S'$ is an ideal of $Gr S'$, since for any $g\in Gr_j \,S'$, $j\geq 1$,
and any $f\in Gr_{-2}S'$, 
$\nu([f, g])=\nu(g)-1$, hence $[f, g]$ lies 
in $S'_{j-1}$, i.e., $[f, g]=0$ in $Gr S'$.
We have:
$$\overline{Gr S'}/Gr_{-2}S'\cong SHO(n,n)\otimes\Lambda(\eta)+
\C\frac{\partial}{\partial\eta}+\C(E-2-\beta \,ad(\Phi)+2\eta\frac{\partial}{\partial\eta})$$
with respect to the grading of type $(1,\dots,1|0,\dots,0)$ on 
$SHO(n,n)$ and $\deg(\eta)=0$.
$\prod_{j\geq 0}Gr_j S'$  is thus not a maximal subalgebra
of $\overline{GrS'}$ since it is contained in $\prod_{j\geq 0}Gr_j S'+Gr_{-2}S'$.
Nevertheless, note that, for every $\beta$, $S'_0$ is contained in
$S=SKO(n,n+1;\beta)$, and $S'_0+\C\xi_1\dots\xi_n$ generates
the whole $S$.
It follows that, for every $\beta\neq 1, (n-2)/n$,
since $S=S'$, $S'_0$ is a maximal open subalgebra of $S$.
If $\beta=1$ or $\beta=(n-2)/n$, $S'_0$ is not a 
maximal subalgebra of $S'$ but it is a maximal subalgebra of $S$. 



Finally, for every $\beta$,
$S'_0=SKO(n,n+1;\beta)\cap L_0$ where $L_0$ is the standard subalgebra
of $KO(n,n+1)$ constructed in Example \ref{ex1KO}. It follows
that $S'_0$ is the standard subalgebra of $SKO(n,n+1;\beta)$ stabilizing
the ideal $I_U=(x_1, \dots, x_n)$.
\end{example}
\begin{example}\label{ex2sko}\em Let $t$ be an integer
such that $1\leq t\leq n$ and let us consider the valuation
$\nu$ on $\Lambda(n,n+1)$ defined in Example \ref{ex2KO}.
Consider the following subspaces of $S'=SKO'(n,n+1;\beta)$:
$$S'_i(t)=\{f\in\Lambda^{\beta}(n,n+1) ~|~ \nu(f)\geq i+2\}+\langle 1, \tau+\beta\Phi\rangle
~~\mbox{if}~ i\leq 0,$$
$$S'_i(t)=\{f\in\Lambda^{\beta}(n,n+1) ~|~ \nu(f)\geq i+2\} 
~~\mbox{if}~ i>0.$$
By  commutation rules (\ref{ko}),  these subspaces
define in fact a filtration of \break $SKO'(n,n+1;\beta)$, having
depth $2$ if $t\neq n$ and  depth $1$ if $t=n$.
Then, if $\beta\neq 1$,
$$\overline{GrS'}\cong SHO(n,n)\otimes\Lambda(\eta)+\mathbb{C}\xi_1\dots\xi_n+\C\frac{\partial}{\partial\eta}+\C(E-2-\beta \,ad(\Phi)+2\eta\frac{\partial}{\partial\eta})$$
and, if $\beta=1$,
$$\overline{GrS'}\cong SHO'(n,n)\otimes\Lambda(\eta)+\C\frac{\partial}{\partial\eta}+\C(E-2-\beta \,ad(\Phi)+2\eta\frac{\partial}{\partial\eta}),$$
with respect to the grading of type $(1,...,1,2,...,2|1,...,1,0,...,0)$ of $SHO^\prime(n,n)$,
with
$n-t$ 2's and $n-t$ zeros, and $\deg(\eta)=0$.
When $n>2$ these gradings are irreducible
for every $t=2,...,n$ (cf.\ Remark \ref{SHO'}) and
satisfy property $(iii)^{\prime}$ of Proposition \ref{basic}$(b)$. Therefore,
by Corollary \ref{cor}, when $n>2$, $S'_0(t)$ is a
maximal subalgebra of $SKO'(n,n+1;\beta)$ for every $t=2,...,n$.

Let $S=SKO(n,n+1;\beta)$ and let $S_j(t):=S'_j(t)\cap S$.
If $\beta\neq 1,(n-2)/n$, then $S=S'$, hence $S_0(t)$ is, for every
$t=2, \dots, n$, a maximal open subalgebra of $S$.
If $\beta=(n-2)/n$ or $\beta=1$,
then the subspaces $S_j(t)$ define a filtration of $S$ such that:
$$\overline{GrS}\cong SHO(n,n)\otimes\Lambda(\eta)+\C\frac{\partial}{\partial\eta}+\C(E-2-\beta \,ad(\Phi)+2\eta\frac{\partial}{\partial\eta})$$
or
$$\overline{GrS}\cong SHO(n,n)\otimes\Lambda(\eta)+\mathbb{C}\xi_1\dots\xi_n+\C\frac{\partial}{\partial\eta}+\C(E-2-\beta \,ad(\Phi)+2\eta\frac{\partial}{\partial\eta}),$$
respectively, with respect to the grading of type $(1,...,1,2,...,2|1,...,1,0,...,0)$ of $SHO(n,n)$,
with
$n-t$ 2's and $n-t$ zeros, and $\deg(\eta)=0$. It follows that, if $n>2$,
then
$S_0(t)$ is a maximal open subalgebra of $S$, for every $t=2,\dots,n$
(cf.\ Remark \ref{gradingsofHO}).

Notice that the grading of principal type of $W(n,n)$ induces an irreducible
grading on $SHO'(n,n)$ also for $n=2$,
but it induces on $SHO(2,2)$ a grading
which is not irreducible. 
It follows that $S'_0(2)$ is a fundamental maximal subalgebra
of $SKO'(2,3;\beta)$ for every $\beta$, but $S_0(2)$ is a maximal subalgebra
of  $SKO(2,3;\beta)$ if and only if 
$\beta\neq 0$. When $\beta=0$ the subalgebra $S_0(2)$ of
$SKO(2,3;0)$ is indeed contained in the graded subalgebra
of type $(1,1|$ $-1,-1,0)$.

Finally, note that $S_0(t)=L_0(t)\cap SKO(n,n+1;\beta)$ where $L_0(t)$ is
the subalgebra of $KO(n,n+1)$ constructed in Example \ref{ex2KO}.
It follows that $S_0(t)$ 
stabilizes the ideal $I_{\cal{U}}=(x_1, \dots, x_n, 
\xi_1, \dots, \xi_t)$ of $\Lambda(n,n+1)$.
\end{example}

\begin{remark}\label{weightsforSKO}\em
Let $S=SKO(n,n+1;\beta)$ and consider its grading of
principal type: $S=\prod_{j\geq -2}S_j$.
Then $\tau$ acts on $S_j$ by multiplication by $j$. 
By Remark \ref{weightsforKO}, $(\C\xi_{i_1}\dots\xi_{i_h}\otimes T)
\cap S_{h}$,
and $(\C x_k\xi_{j_1}\dots\xi_{j_h}\otimes T)\cap S_{h+1}$ with 
$k\neq j_1,\dots ,j_h$, are $T$-weight spaces of $SKO(n,n+1;\beta)$.
\end{remark}

\begin{remark}\label{SKO'}\em 
The same arguments as in the proof of Theorem \ref{newS}
 show, by Remark \ref{formula}, that
every maximal open subalgebra of $SKO(n,n+1;\beta)$,
 $SKO'(n,n+1;\beta)$ and
$CSKO'(n,n+1;\beta)=SKO'(n,n+1;\beta)+\C\Phi$ is regular.

\begin{theorem}\label{SKO} Let $S=SKO(n,n+1;\beta)$. Then all  maximal open 
subalgebras of $S$ are, up to conjugation, the following:
\begin{itemize}
\item[$(a)$] if $n=2$ and $\beta\neq 0, 1$, 
\begin{itemize}
\item[$(i)$] 
the graded subalgebras of type $(1,1|0,0,1)$, $(1,1|1,1,2)$
and \break $(1,1|-1,-1,0)$;
\item[$(ii)$] the non-graded subalgebra $S_0(2)$ constructed in Example \ref{ex2sko};
\end{itemize}
\item[$(b)$]
if $n=2$ and $\beta=1$,
\begin{itemize}
\item[$(i)$] 
the graded subalgebras of type $(1,1|0,0,1)$, $(1,1|1,1,2)$;
\item[$(ii)$] the non-graded subalgebra $S_0(2)$ constructed in Example \ref{ex2sko};
\end{itemize}
\item[$(c)$]
if $n=2$ and $\beta=0$,
\begin{itemize}
\item[$(i)$] 
the graded subalgebras of type $(1,1|0,0,1)$ and $(1,1|-1,-1,0)$;
\end{itemize}
\item[$(d)$]
if $n>2$,
\begin{itemize}
\item[$(i)$] the graded subalgebra of type $(1,\dots,
1|0,\dots, 0,1)$ and the graded subalgebras of type $(1,\dots, 1,2,\dots,
2|1,\dots, 1,0,\dots, 0,2)$ with
$n-t+1$ 2's and $n-t$ zeros,
for $t=2, \dots, n$;
\item[$(ii)$]
the non-graded subalgebra $S'_0$ described in Example \ref{ex1sko} and the
non-graded subalgebras $S_0(t)$ described in Example \ref{ex2sko} for $t=2,\dots, n$.
\end{itemize}
\end{itemize}
\end{theorem}
{\bf Proof.} Let $L_0$ be a maximal open subalgebra of $S$. By Remark 
\ref{SKO'}, $L_0$ is regular.
Therefore, by Remark \ref{standard} and Proposition \ref{derivations},
 we can assume that
$L_0$ is invariant with respect
to the standard torus $T$ of $KO(n,n+1)$. 
It follows that
$L_0$ decomposes into the direct product of weight spaces with respect to $T$.
Notice that
$\mathbb{C}1$, $\mathbb{C}x_i$,  $\mathbb{C}\xi_{i_1}\dots \xi_{i_h}$,
$\mathbb{C}x_j\xi_{i_1}\dots \xi_{i_h}$, with $j\neq i_1\neq \dots \neq i_h$,
are one-dimensional $T$-weight spaces
(see Remark
\ref{weightsforKO}). 
Besides, note that
the elements $\xi_i$ cannot
lie in $L_0$ since the corresponding vector fields
$\rho(\xi_i)=\xi_i\frac{\partial}{\partial\tau}+\frac{\partial}{\partial x_i}$ 
are not exponentiable.  

Let us first assume $n=2$. We distinguish two cases:

CASE I: $1$ does not lie in $L_0$. We may assume that one of the
following possibilities occurs:

\noindent
1) No $x_i$ lies in $L_0$. Then the $T$-invariant complement of $L_0$ contains
the $T$-invariant complement of the maximal graded subalgebra
of $S$ of type $(1,1|1,1,2)$, hence $L_0$ is contained in the graded
subalgebra of principal type. If $\beta=0$ then the subalgebra of principal type 
is not maximal therefore this contradicts the maximality of $L_0$.
 If $\beta\neq 0$ then
$L_0$ coincides with the graded subalgebra of type 
$(1,1|1,1,2)$ by maximality. 

\noindent
2) The elements $x_{1}, x_2$ lie in $L_0$.
Then the $T$-invariant complement of $L_0$ contains 
 the $T$-invariant complement of the maximal graded subalgebra of $L$ of
type $(1, 1|0, 0,1)$. Since $L_0$ is maximal it
coincides with this graded subalgebra.

Notice that if $L_0$ contains $x_2$ (resp.\ $x_1$)
 then, due to its maximality, it contains also $x_1$ (resp.\ $x_2$). Indeed,
any open regular subalgebra of $S$ containing $x_2$ and not
containing $1$ and $x_1$ (resp.\ containing $x_1$ and not containing $1$ and $x_2$)  is
contained in the subalgebra of type $(1,
2|1,0,2)$ (resp.\ $(2,1|0,1,2)$) which is not maximal by Remark
\ref{irrenotKO(n,n+1)}.

\medskip

CASE II: $1$ lies in $L_0$. Since the elements $\xi_i$'s
do not lie in $L_0$, the elements $\xi_i\tau+\varphi$ cannot
lie in $L_0$ for any $\varphi\in\Lambda(2,2)$, where
by $\Lambda(n,n)$ we mean the subalgebra of $\Lambda(n,n+1)$
generated by all even indeterminates and all odd indeterminates
except $\tau$. Indeed, by commutation rules
(\ref{ko}), we have: $[1, \xi_i\tau+\varphi]=-2\xi_i$.
Note that if $\beta=0$ then the grading of type $(1,1|-1,-1,0)$ has depth 1
with $-1$-st graded component spanned by the elements $\xi_i$ and
$\xi_i(\tau-\Phi)$ for $i=1,2$. It follows that if $\beta=0$, 
then $L_0$ is contained in
the graded subalgebra of $SKO(2,3;0)$ of type $(1,1|-1,-1,0)$, thus
coincides with it, due to its maximality.

Now suppose $\beta\neq 0$. 
Since, for every $i$, $\C x_i$ is a one-dimensional weight space
of $SKO(n,n+1;\beta)$, we may assume that one of the following situations
holds:

\noindent
1) No $x_i$  lies in $L_0$. Then the same arguments
as in Example \ref{ex2KO} show that $L_0$ coincides with the subalgebra $S_0(2)$
constructed in Example \ref{ex2sko};

\noindent
2) $x_1, x_2$ lie in $L_0$.  Then 
$L_0$ is contained
in  the graded subalgebra of $S$ 
of type $(1,1|$ $-1,-1,0)$. Since $L_0$ is maximal the two subalgebras coincide.

Notice that if $1, x_2$ lie in $L_0$, by the maximality
of $L_0$, also $x_1\in L_0$. Indeed, any open regular subalgebra of $S$ containing the elements
$1, x_2$ and not containing $x_1$ is contained in the 
maximal subalgebra of type $(1,1|-1,-1,0)$.

Finally, as we pointed out in Remark \ref{outerSKO(2,3;1)}, when $\beta=1$,
the subalgebras of type $(1,1|-1,-1,0)$ and $(1,1|1,1,2)$ are conjugate
by an element of $G$.

\medskip   
Let us now suppose $n>2$.
We distinguish two cases:

CASE I: $1$ does not lie in $L_0$. We may assume that one of the
following possibilities occurs:

\noindent
1) No $x_i$ lies in $L_0$. Then the $T$-invariant complement of $L_0$ contains
the $T$-invariant complement of the maximal graded subalgebra
of $S$ of type $(1,\dots, 1|$ $1,\dots,1,2)$. By the maximality of $L_0$
it follows that $L_0$ coincides with the graded subalgebra of type 
$(1,\dots, 1|1,\dots,1,2)$;

\noindent
2) the elements $x_{t+1}, \dots, x_n$ lie in $L_0$ for some $t=2,\dots,
n-1$, and the elements $x_{1}, \dots,
x_t$ do not. It follows, using commutation rules (\ref{ko}),
 that the $T$-invariant complement of $L_0$ contains 
the $T$-invariant complement of the maximal graded subalgebra of $L$ of
type $(1,\dots, 1,2,\dots,
2|1,\dots, 1,0,\dots, 0,2)$ with
$n-t+1$ 2's and $n-t$ zeros. Since $L_0$ is maximal it
coincides with this graded subalgebra;

\noindent
3) the elements $x_i$ lie in $L_0$ for every $i$. Then the $T$-invariant complement
of $L_0$ contains 
the $T$-invariant complement of the maximal graded
subalgebra of $S$ of type $(1,\dots, 1|0,\dots, 0,1)$. It follows that
$L_0$ coincides with this subalgebra.

Notice that if $L_0$ contains the elements $x_2, \dots,
x_{n}$ then, due to its maximality, it contains also $x_1$. Indeed,
any  regular subalgebra of $S$ containing $x_2, \dots, x_{n}$ and not 
containing $1$ and $x_1$ is
contained in the subalgebra of type $(1,
2,\dots, 2|1,0,\dots, 0,2)$ which is not maximal by Remark
\ref{irrenotKO(n,n+1)}.

\medskip

CASE II: $1$ lies in $L_0$. We may assume that one of the following situations
holds:

\noindent
1) For some $t=2, \dots, n$ the elements $x_{1},
\dots, x_t$ do not lie in $L_0$ and  $x_{t+1}, \dots, x_n$ do. Then
the same arguments as in Example \ref{ex2KO} show that $L_0$ coincides with the subalgebra $S_0(t)$
constructed in Example \ref{ex2sko};

\noindent
2) $x_1, \dots, x_n$ lie in $L_0$. Then the same arguments
as in Example \ref{ex1KO} show that $L_0$
is contained in
the subalgebra $S'_0$ of $S$ constructed in
Example \ref{ex1sko}. 
Since $L_0$ is maximal the two subalgebras coincide.

Notice that if $1, x_2, \dots, x_{n}$ lie in $L_0$, by the maximality
of $L_0$, also $x_1\in L_0$. Indeed any open regular subalgebra of $S$ containing the elements
$1, x_2, \dots, x_{n}$ and not containing $x_1$ is contained in the 
maximal subalgebra constructed in Example \ref{ex1sko}.   
\hfill$\Box$

\begin{corollary}\label{SKO-graded} All irreducible
$\Z$-gradings of $SKO(n,n+1;\beta)$ are, up to conjugation,
the following:

\noindent
$(i)$ the gradings of type $(1,1|0,0,1)$, $(1,1|1,1,2)$
and $(1,1|-1,-1,0)$, if $n=2$, $\beta\neq 0, 1$;

\noindent
$(ii)$  the gradings of type $(1,1|0,0,1)$ and $(1,1|1,1,2)$
if $n=2$, $\beta=1$;

\noindent
$(iii)$ the gradings of type $(1,1|0,0,1)$, $(1,1|-1,-1,0)$
if $n=2$,
$\beta=0$;

\noindent
$(iv)$ the gradings of type   $(1,\dots, 1,2,\dots,
2|1,\dots, 1,0,\dots, 0,2)$ with
$t+1$ 2's and $t$ zeros, 
for $t=0, \dots, n-2$ and $(1,\dots,
1|0,\dots, 0,1)$, if $n>2$.
\end{corollary}

\end{remark}

We recall that if $S=SKO(n,n+1;\beta)$ with
$n\geq 2$ and $\beta\neq 1, (n-2)/n$, then
$Der S=S+\C\Phi$ with $\Phi=\sum_{i=1}^n x_i\xi_i$;
if $S=SKO(n,n+1;(n-2)/n)$ with $n\geq 2$, then
$Der S=S+\C\Phi+\C\xi_1\dots\xi_n$; if
$S=SKO(n,n+1;1)$ with $n>2$, then $DerS=S+\C\Phi+\C\xi_1\dots\xi_n\tau$;
finally, if $S=SKO(2,3;1)$ then $Der S=S+sl_2$
 (cf.\ Proposition \ref{derivations}, Remark \ref{outerSKO(2,3;1)}). 

\begin{theorem} 
Let $S=SKO(n,n+1;\beta)$ with $n\geq 2$ and $\beta\neq 1,
(n-2)/n$, so that $SKO(n,n+1;\beta)=SKO'(n,n+1;\beta)$
and $Der S=CSKO'(n,n+1;\beta)$. Then all
maximal among open $\Phi$-invariant subalgebras of $S$ 
are, up to conjugation,
the subalgebras of $S$
listed in Theorem \ref{SKO} $(a)$ and $(d)$.
\end{theorem}
{\bf Proof.} Let $L_0$ be a
maximal among open $\Phi$-invariant subalgebras of $S$.
Then $L_0+\C\Phi$ is a maximal open subalgebra of
$CSKO'(n,n+1;\beta)$, hence it is regular by
Remark \ref{SKO'}.
Then one uses the same arguments as in the proof
of Theorem \ref{SKO}. \hfill$\Box$

\bigskip

We shall now classify the open
subalgebras of  $S=SKO(n,n+1;(n-2)/n)$ and
$S=SKO(n,n+1;1)$, which are maximal
among the $\mathfrak{a}_0$-invariant subalgebras of $S$, for every
subalgebra $\mathfrak{a}_0$ 
of $\mathfrak{a}$. 

\begin{remark}\label{mfSKO'}\em 
By Remark \ref{SKO'} every maximal open subalgebra
of $SKO'(n,n+1;\beta)$ or $CSKO'(n,n+1;\beta)$ is regular. 
Therefore the same
arguments as in the proof of Theorem \ref{SKO}
show that all fundamental among maximal subalgebras of $SKO'(n,n+1;(n-2)/n)$
(resp.\ $CSKO'(n,n+1;(n-2)/n)$), with $n>2$,
are, up to conjugation, 
the graded subalgebras of type 
$(1,\dots,1,2,\dots,2|1,\dots,1,$ $0,\dots,0,2)$,
with
$n-t+1$ 2's and $n-t$ zeros, 
and the non-graded subalgebras $S'_0(t)$
(resp.\ $S'_0(t)+\C\Phi$) constructed in
Example \ref{ex2sko},
for $t=2,\dots, n$.
Indeed, the graded subalgebra of $SKO'(n,n+1;(n-2)/n)$ 
(resp.\ $CSKO'(n,n+1;(n-2)/n)$) of 
type $(1,\dots,1|0,\dots,0,1)$ and the subalgebra $S'_0$
constructed in Example \ref{ex1sko},
are not  maximal, since they are  contained 
in $SKO(n,n+1;(n-2)/n)$ (resp.\ $SKO(n,n+1;(n-2)/n)+\C \Phi$). 
By the same arguments,
all maximal among fundamental subalgebras of $SKO'(n,n+1;(n-2)/n)$ 
and $CSKO'(n,n+1;(n-2)/n)$, for $n>2$,
are, up to conjugation, the graded subalgebra of subprincipal type,
the graded subalgebras of type 
$(1,\dots,1,2,\dots,2|1,\dots,1,$ $0,\dots,0,2)$,
with
$n-t+1$ 2's and $n-t$ zeros,
the non-graded subalgebras $S'_0(t)$ constructed in
Example \ref{ex2sko}, for $t=2,\dots, n$,
and, the subalgebra $S'_0$
constructed in Example \ref{ex1sko}.

Likewise, all fundamental among maximal subalgebras of $SKO'(2,3;0)$
(resp.\break $CSKO'(2,3;0)$) are, up to conjugation,
 the graded subalgebra of type $(1,1|1,1,2)$
and the subalgebra $S'_0(2)$ (resp.\ $S'_0(2)+\C\Phi$).
All maximal among fundamental subalgebras of $SKO'(2,3;0)$ (resp.
$CSKO'(2,3;0)$) are the graded subalgebras of type
$(1,1|1,1,2)$, $(1,1|0,0,1)$, $(1,1|-1,-1,0)$ and
the non-graded subalgebra $S'_0(2)$ (resp.\ $S'_0(2)+\C\Phi$).
\end{remark}

\begin{theorem}\label{derSKO} Let $S=SKO(n,n+1;(n-2)/n)$ with $n\geq 2$. 

\noindent
(i) All maximal among open $\Phi$-invariant subalgebras of $S$
are, up to conjugation,
 the maximal open subalgebras listed in Theorem \ref{SKO} (c) and
(d).

\noindent
(ii) If $\mathfrak{a}_0=\C\xi_1\dots\xi_n$ or $\mathfrak{a}_0=
\mathfrak{a}$, then
all  maximal among $\mathfrak{a}_0$-invariant open subalgebras
of $S$ are, up to conjugation,
the graded subalgebras of type 
$(1,\dots,1,$ $2,\dots,2|1,\dots,1,0,\dots,0,2)$,
with
$n-t+1$ 2's and $n-t$ zeros, 
and the non-graded subalgebras $S_0(t)$ constructed in
Example \ref{ex2sko}, for $t=2,\dots, n$.
%
\end{theorem}
{\bf Proof.} One uses Remark \ref{mfSKO'} and the same arguments as in
the proof of Theorem \ref{derS(1,n)}.
%
%
\hfill$\Box$

\begin{remark}\label{mfSKO(n,n+1;1)}\em
By Remark \ref{SKO'} every maximal open subalgebra
of $SKO'(n,n+1;\beta)$ or $CSKO'(n,n+1;\beta)$, 
for every $n\geq 2$, is regular. 
Therefore the same
arguments as in the proof of Theorem \ref{SKO}
show that all fundamental among maximal subalgebras of
$SKO'(n,n+1;1)$ (resp.\ $CSKO'(n,n+1;1)$) are, up to conjugation,
the graded subalgebras of type $(1,\dots,
1|0,\dots, 0,1)$ and $(1,\dots, 1,2,\dots,
2|$ $1,\dots, 1,0,\dots, 0,2)$ with
$n-t+1$ 2's and $n-t$ zeros,
and the
non-graded subalgebras $S'_0(t)$ (resp.\ $S'_0(t)+\C\Phi$) constructed
 in Example \ref{ex2sko}, for $t=2,\dots, n$.
By the same arguments, all maximal among fundamental subalgebras
of $SKO'(n,n+1;1)$ (resp.\ $CSKO'(n,n+1;1)$) are, up
to conjugation, all the subalgebras 
listed above and the subalgebra $S'_0$ constructed in Example
\ref{ex1sko}, if $n>2$, or the graded subalgebra of type $(1,1|-1,-1,0)$
if $n=2$. Note that the subalgebras of $SKO'(2,3;1)$ or $CSKO'(2,3;1)$
of type $(1,1|1,1,2)$ and
$(1,1|-1,-1,0)$ are not conjugate.
\end{remark}

\begin{theorem}\label{SKO(n,n+1;1)} Let $S=SKO(n,n+1;1)$ with $n>2$.

\noindent
(i) All maximal among open $\Phi$-invariant subalgebras of $S$
are, up to conjugation,
 the maximal open subalgebras listed in Theorem \ref{SKO} 
(d).

\noindent
(ii) If $\mathfrak{a}_0=\C\xi_1\dots\xi_n\tau$ or $\mathfrak{a}_0=
\mathfrak{a}$, then
all  maximal among $\mathfrak{a}_0$-invariant open subalgebras
of $S$ are, up to conjugation,
the graded subalgebras of type $(1, \dots,1|$ $0,\dots,0)$ and
$(1,\dots,1,2,\dots,2|$ $1,\dots,1,0,\dots,0,2)$,
with
$n-t+1$ 2's and $n-t$ zeros, 
and the non-graded subalgebras $S_0(t)$ constructed in
Example \ref{ex2sko}, for $t=2,\dots, n$.
\end{theorem}
{\bf Proof.} One uses Remark \ref{mfSKO(n,n+1;1)} and the same arguments as in
the proof of Theorem \ref{derS(1,n)}. \hfill$\Box$

\begin{theorem}\label{sl_2-inv} Let $S=SKO(2,3;1)$ and
let $\mathfrak{b}=\C e+\C h\subset\mathfrak{a}\cong sl_2$.

\noindent
(i) If $\mathfrak{a}_0$ is a one-dimensional subalgebra of $\mathfrak{a}$,
then all maximal among open $\mathfrak{a}_0$-invariant subalgebras of $S$
are, up to conjugation, the maximal subalgebras listed in
Theorem \ref{SKO} (b).

\noindent
(ii) The graded
subalgebra of type $(1,1|1,1,2)$  is, up to
conjugation, the only 
maximal among open $\mathfrak{b}$-invariant 
subalgebras of $S$, which is not
invariant with respect to $\mathfrak{a}$.

\noindent
(iii) All maximal among open $\mathfrak{a}$-invariant subalgebras
of $S$ are, up to conjugation,
the graded subalgebra of type $(1,1|0,0,1)$ and
the non-graded subalgebra $S_0(2)$ constructed in Example
\ref{ex2sko}.
\end{theorem}
{\bf Proof.} By Remark \ref{mfSKO(n,n+1;1)}, the proof of
$(i)$ is the same as the proof of $(i)$ and
$(ii)$ in
Theorem \ref{derS(1,n)}. Recall that the graded
subalgebras of type $(1,1|1,1,2)$ and $(1,1|-1,-1,0)$ are
conjugate.

Now, using \cite[Proposition 5.3.4]{CK}
one can check that the maximal graded subalgebra of $SKO(2,3;1)$ of
type $(1,1|0,0,1)$ and the subalgebra $S_0(2)$ constructed in Example
\ref{ex2sko} are invariant with respect to
$\mathfrak{a}$. On the other hand,  the maximal
subalgebra $L_0$ of $S$ of type $(1,1|1,1,2)$  is
invariant with respect to $\mathfrak{b}$
but it is
not $\mathfrak{a}$-invariant. Indeed $L_0$ contains $\xi_1\xi_2$,
it does not contain $1$, but
$f(\xi_1\xi_2)=1$.
Let $M_0$ be a maximal among open $\mathfrak{b}$-invariant
subalgebras of $SKO(2,3;1)$, then
$M_0+\mathbb{C}\xi_1\xi_2\tau+\mathbb{C}\Phi$
 is a fundamental maximal subalgebra of
$CSKO'(2,3;1)$ containing $\xi_1\xi_2\tau$
and $\Phi$, hence,
by Remark \ref{mfSKO(n,n+1;1)},  
$M_0$ is conjugate
to
the graded subalgebra of type $(1,1|1,1,2)$, or  to the subalgebra
of type $(1,1|0,0,1)$, or to the subalgebra $S_0(2)$.

Now suppose that $\tilde{S}$ is  a maximal among open $\mathfrak{a}$-invariant subalgebras
of $SKO(2,3;1)$. Then $\tilde{S}$ is   $\mathfrak{b}$-invariant, hence
it is conjugate either to the graded subalgebra of type $(1,1|0,0,1)$, or
to the subalgebra $S_0(2)$ constructed in Example
\ref{ex2sko}. Indeed, otherwise, $\tilde{S}$ is contained
either in the subalgebra of type $(1,1|1,1,2)$ or in
a conjugate $S_{\gamma}$ of it (see  Remark \ref{inf.manySKO}).
Since $\tilde{S}$ is $\mathfrak{a}$-invariant, it is
invariant with respect to all outer automorphisms of
$S$, hence it is contained
in the intersection of all the subalgebras in the orbit
of the subalgebra of principal type. It follows, by
Remark \ref{inf.manySKO}, that $\tilde{S}$ is
contained  in the subalgebra of type $(1,1|0,0,1)$.
This contradicts the maximality of $\tilde{S}$ among
$\mathfrak{a}$-invariant subalgebras.
\hfill$\Box$
 
\section{Maximal open subalgebras of $\boldsymbol{SHO^{\sim}(n,n)}$ and
  $\boldsymbol{SKO^\sim(n,n+1)}$}
{\bf {\em The Lie superalgebra $\boldsymbol{SHO^{\sim}(n,n)}$}.}
Let $n$ be even. The Lie superalgebra $SHO^\sim(n,n)$
is the subalgebra of $HO(n,n)$ defined as follows:
$$SHO^\sim(n,n)=\{X\in HO(n,n)~|~ X(F\omega)=0\}$$
where $\omega$ is the volume form associated to the usual
 divergence and $F=1-2\xi_1\dots\xi_n$.
By Remark \ref{tilde}, $SHO^\sim(n,n)$ consists of vector fields $X$ in
$HO(n,n)$ such that $div_F(X)=0$ or, equivalently, by Remark \ref{usual},
such that $div(FX)=0$.

Using the isomorphism between $HO(n,n)$ and $\Lambda(n,n)/\C 1$ with
the Buttin bracket, it is possible to realize $SHO^\sim(n,n)$ as
follows (cf.\ \cite[\S 2]{Gafa}):
$$SHO^{\sim}(n,n)=((1+\xi_1\dots\xi_n)\Lambda^{\Delta}(n,n))/\C 1$$
where $\Lambda^{\Delta}(n,n)=\{f\in\Lambda(n,n)~|~\Delta(f)=0\}$
and $\Delta$ is the odd Laplacian. 
Equivalently, $SHO^{\sim}(n,n)$
 can be
identified with the space $\Lambda(n,n)^{\Delta}/\C 1$ 
with 
the following deformed bracket (\cite[\S 5]{ChK}): 
\begin{eqnarray}
{[f,g]~}&=&{[\xi_1\dots\xi_n, fg]_{ho} ~~\mbox{if}~ ~f,g\in\C [[x_1, \dots, x_n]],}\nonumber\\
{[x_i, \xi_j]}&=&{\delta_{ij}\xi_1\dots\xi_n,}\\
{[f,g]~}&=&{[f,g]_{ho}  ~~\mbox{otherwise},\nonumber}
\label{deformedSHO}
\end{eqnarray}

\noindent
where $[\cdot,\cdot]_{ho}$ denotes the bracket in $HO(n,n)$.

\medskip

The superalgebra $SHO^\sim(n,n)$
is simple for $n\geq 2$ ($n$ even) \cite[Example 6.2]{K}. Since, as we recalled in the introduction,
$SHO^\sim(2,2)\cong H(2,1)$, when dealing with 
$SHO^\sim(n,n)$ we will assume $n>2$.
\begin{remark}\label{gradingsofSHOsim}\em
A $\Z$-grading of $W(n,n)$ induces a $\Z$-grading  on $SHO^\sim(n,n)$
if and only if $\deg x_i+\deg\xi_i=\mbox{const.}$ and
$\sum_{i=1}^n\deg\xi_i=0$. In particular
the $\Z$-grading of type $(1,\dots, 1|0,\dots,0)$ induces on  $SHO^\sim(n,n)$ 
a grading of depth $1$ which is irreducible by Remark \ref{gd}.
\end{remark}

 
In what follows we will identify $SHO^\sim(n,n)$ with 
$\Lambda(n,n)^{\Delta}/\C 1$ with bracket (5.1) and fix its
maximal torus $T=\langle
x_i\xi_i-x_{i+1}\xi_{i+1}~|~ i=1,\dots, n-1\rangle$.

%
\begin{example}\label{9th}\em  On $\Lambda(n,n)$, for any fixed  integer $t$
such that $1\leq t\leq n$,
  let us define  the following valuation $\nu$: 
$$\nu(x_i)=1,  ~~\nu(\xi_i)=1 ~~\mbox{for}~
  i=1,\dots,t;$$
$$\nu(x_i)=2, ~~ \nu(\xi_i)=0 ~~\mbox{for}~ i=t+1,\dots,n.$$
 Let us define
  the following filtration of $L=SHO^\sim(n,n)$:
$$L_j(t)=\{x\in \Lambda^{\Delta}(n,n)/\C 1 ~|~ \nu(x)\geq j+2\}$$
Then $\overline{Gr L}\cong SHO^{\prime}(n,n)$  with respect to the
  $\Z$-grading of type  $(1,\dots, 1,
  2,\dots,2|$ $1,\dots,1,0,\dots,0)$
with
$n-t$ 2's and $n-t$ zeros. Since this grading is
irreducible for every $t=2,\dots, n$ (cf.\ Remarks
\ref{gradingsofHO}, \ref{SHO'}), it follows, using
  Corollary \ref{cor}, that $L_0(t)$ is a maximal  regular subalgebra
  of $L$ for every $t=2,\dots, n$.
\end{example}

\begin{remark}\label{weightsforSHOsim}\em 
Let 
$\Sigma_0:=\langle x_{i_1}\dots x_{i_k}\xi_{i_1}\dots \xi_{i_k}
~|~ k=1, \dots,n\rangle$. All elements of $SHO^\sim(n,n)$
lying in $\Sigma_0$ have $T$-weights equal to zero.
 
Let $i_1\neq \dots \neq i_h$ and 
$\{i_1, \dots, i_h, j_1, \dots, j_{n-h}\}=\{1, \dots, n\}$. Then  
$\{f\in \langle \xi_{i_1}\dots\xi_{i_h},$ $x_{j_1}\dots x_{j_{n-h}}\rangle\otimes
\Sigma_0 ~|~ \Delta(f)=0\}$ is a weight space with respect to $T$. 
Likewise, if $i_1\neq \dots \neq i_h\neq j$,  then
$\{f\in\langle x_j\xi_{i_1}\dots\xi_{i_h}, x_jx_{j_1}\dots x_{j_{n-h}}
\rangle\otimes
\Sigma_0 ~|~ \Delta(f)=0\}
$ is a weight space with respect to $T$.
\end{remark}

\begin{theorem}\label{SHOsim} Let $L=SHO^\sim(n,n)$ with $n>2$ even.
All maximal open subalgebras 
  of $L$ are, up to conjugation, the following:
\begin{itemize}
\item[$(i)$] the graded subalgebra of type $(1,\dots,1|0,\dots,0)$;
\item[$(ii)$] the non-graded subalgebras $L_0(t)$ constructed in Example \ref{9th},
for $t=2,\dots, n$.
\end{itemize} 
\end{theorem}
{\bf Proof.} Let $L_0$ be a   maximal open
subalgebra of $L$. The same argument as in the proof of Theorem \ref{newS} shows that $L_0$ is regular
hence we can assume, by Remark \ref{standard}, that it is invariant 
with respect to the torus $T$. It follows that
$L_0$ decomposes into the direct product of $T$-weight spaces.
Note that  the elements $\sum_j\alpha_j\xi_j+f$ cannot
lie in $L_0$ for any non zero linear combination $\sum_j\alpha_j\xi_j$
and any odd function
$f\in\Lambda^{\Delta}(n,n)/\C 1$ with no linear terms,  since the elements $\xi_j$
are not exponentiable. 
We may therefore assume that one of the
following situations occurs:

\noindent
1) the elements $x_i+\varphi_i$
 lie in $L_0$ for some elements $\varphi_i$ 
with no linear terms,
for every $i=1, \dots,n$.
Then the elements $\xi_i\xi_j+\psi$
do not lie in $L_0$ for any $\psi$ in the $T$-weight space
of $\xi_i\xi_j$, $\psi\notin\C\xi_i\xi_j$, since, for such a $\psi$,
by  Remark \ref{weightsforSHOsim},
$[x_i+\varphi_i, \xi_i\xi_j+\psi]=\xi_j+\eta$
for some function $\eta\in\Lambda^{\Delta}(n,n)/\C 1$ without
linear terms. It follows
that $L_0$ does not contain any element $\xi_i\xi_j+\psi$
for any $\psi\notin\C\xi_i\xi_j$.
The same argument shows, by induction on $k=1, \dots,n$, 
that $L_0$ does not contain the elements
$\xi_{i_1}\dots\xi_{i_k}+\psi_k$ for any
function $\psi_k\notin\C\xi_{i_1}\dots\xi_{i_k}$, 
for any
$k=1, \dots,n$. 
$L_0$ is therefore contained in the maximal graded subalgebra
   of $L$ of type $(1,\dots,1|0,\dots,0)$, hence coincides with it
since it is maximal;

\noindent
2) there exists some $t=2, \dots, n$ such that the elements
$x_1+\varphi_1, \dots, x_t+\varphi_t$ do not  lie in $L_0$ for any
functions $\varphi_1, \dots, \varphi_t$ without
linear terms,
and $x_{t+1}+\varphi_{t+1},
\dots, x_n+\varphi_n$ lie in $L_0$ for
some functions $\varphi_{t+1}, \dots, \varphi_n$
with no linear terms.
Then arguing as in 1) and using Remark \ref{weightsforSHOsim},
one shows that
$L_0$ is contained in the subalgebra $L_0(t)$
constructed in Example \ref{9th}. Thus $L_0=L_0(t)$ due to the maximality
of $L_0$.

Notice that if $x_2+\varphi_2,\dots, x_{n}+\varphi_n$ lie in $L_0$
for some functions $\varphi_2, \dots, \varphi_n$
with no linear terms,
then also $x_1+\varphi_1$ lies in $L_0$ for some
$\varphi_1\in
\Lambda^{\Delta}(n,n)/\C 1$ with no linear terms.
Indeed, any open $T$-invariant subalgebra of $L$ containing
$x_2+\varphi_2,\dots, x_{n}+\varphi_n$
and not containing $x_1+\varphi$ for any function
$\varphi\in
\Lambda^{\Delta}(n,n)/\C 1$ with no linear terms,
 is properly contained in the maximal graded subalgebra
of type $(1,\dots,1|0,\dots,0)$, hence it is not maximal.
\hfill$\Box$

\begin{corollary} The Lie superalgebra $SHO^\sim(n,n)$ 
 has, up to conjugation, only one irreducible $\Z$-grading:
the grading of type $(1,\dots,1|0,\dots,0)$.
\end{corollary}

\noindent
{\bf {\em The Lie superalgebra $\boldsymbol{SKO^{\sim}(n,n+1)}$}.}
Let $n$ be odd. The Lie superalgebra $SKO^\sim(n,n+1)$ is the
subalgebra of $KO(n,n+1)$ defined as follows:
$$SKO^\sim(n,n+1)=\{X\in KO(n,n+1)~|~X(F\omega_{\beta})=0\}$$
where $\omega_{\beta}$ is the volume form attached to
the divergence $div_{\beta}$ for $\beta=(n+2)/n$
 and $F=1+\xi_1\dots\xi_n\tau$.
By Remark \ref{tilde}, $SKO^\sim(n,n+1)$ consists of vector fields $X$
in $KO(n,n+1)$ such that $X(F)F^{-1}+div_{\beta}(X)=0$, where
$\beta=(n+2)/n$.

Using the isomorphism between $KO(n,n+1)$ and $\Lambda(n,n+1)$
with bracket (\ref{ko}), it is possible to realize $SKO^\sim(n,n+1)$
as follows  (cf.\ \cite[\S 2]{Gafa}):
$$SKO^{\sim}(n,n+1)=(1+\xi_1\dots\xi_n\tau)\Lambda^{\Delta'}(n,n+1)$$
where
$\Lambda^{\Delta'}(n,n+1)=\{f\in\Lambda(n,n+1)~|~\Delta'(f)=0\}$ and
 $\Delta':=div_{(n+2)/n}=
\Delta+(E-(n+2))\frac{\partial}{\partial \tau}$.
Equivalently, $SKO^\sim(n,n+1)$  can be
identified with  the space $\Lambda(n,n+1)^{\Delta'}$ 
with the following deformed bracket: 
\begin{equation}
[f,g]=[f,g]_{ko}+\alpha(fg)
\label{deformedsko}
\end{equation}
where $[\cdot,\cdot]_{ko}$ denotes the bracket in the Lie superalgebra $KO(n,n+1)$ and
$\alpha(b)=[\xi_1\dots\xi_n\tau, b]_{ko}-2b\xi_1\dots\xi_n$ if $b$ is a monomial in the
$x_i$, and $\alpha(b)=0$ for all other monomials
(\cite{ChK}, \cite[Example 6.3]{K}). The superalgebra $SKO^\sim(n,n+1)$
is simple for $n\geq 3$ ($n$ odd).

\begin{remark}\em If $F=1+\xi_1\dots\xi_n\tau$ and
 $\beta\neq (n+2)/n$, then 
$\{X\in KO(n,n+1)~|~X(F\omega_{\beta})=0\}=\{X\in SKO(n,n+1)~|~X(F)=0\}$.
In particular this is a proper subalgebra of $KO(n,n+1)$ which
is not transitive.
\end{remark}

In what follows we will identify $SKO^\sim(n,n+1)$
with $\Lambda(n,n+1)^{\Delta'}$ with bracket (\ref{deformedsko}). 
Let us fix the torus $T=\langle \tau+\frac{n+2}{n}\Phi, ~x_i\xi_i-x_{i+1}\xi_{i+1} ~|~
i=1, \dots, n-1\rangle$, where $\Phi=\sum_{i=1}^n x_i\xi_i$.
\begin{example}\label{10thII}\em 
  Let us define  the following valuation $\nu$ on $\Lambda(n,n+1)$: 
$$\nu(x_i)=1,~~\nu(\xi_i)=0,  ~~\nu(\tau)=1$$
and let us consider the following filtration of $L=SKO^\sim(n,n+1)$:
$$L_j=\{x\in \Lambda^{\Delta'}(n,n+1) ~|~ \nu(x)\geq j+1\}.$$
then $\overline{Gr
  L}\cong SKO^{\prime}(n,n+1; \frac{n+2}{n})$ with respect to the
  $\Z$-grading of type  $(1,\dots,1|$ $0,\dots,0,1)$.  It follows, using
  Corollaries \ref{cor} and \ref{SKO-graded}, that $L_0$ is a maximal  open subalgebra
  of $L$.
\end{example}
\begin{example}\label{10thI}\em
On   $\Lambda(n,n+1)$, for any fixed  integer $t$, $1\leq t\leq n$,
  let us define  the following valuation $\nu$: 
$$\nu(x_i)=1,
 ~~\nu(\xi_i)=1 ~~\mbox{for}~
  i=1,\dots,t;$$ 
 $$\nu(x_i)=2,
~~\nu(\xi_i)=0 ~~\mbox{for}~ i=t+1,\dots,n; ~~\nu(\tau)=2;$$
where  by $\tau$ we denoted
  the $n+1$-th odd indeterminate of $\Lambda(n,n+1)$. Let us define
  the following filtration of $L=SKO^\sim(n,n+1)$:
$$L_j(t)=\{x\in \Lambda^{\Delta'}(n,n+1) ~|~ \nu(x)\geq j+2\}.$$ 
Then $\overline{Gr
  L}\cong SKO^{\prime}(n,n+1; \frac{n+2}{n})$ with respect to the
  $\Z$-grading of type  $(1,\dots, 1,$
$2,\dots,2|1,\dots,1,0,\dots,0,2)$ with
$n-t+1$ 2's and $n-t$ zeros. It follows, using
  Corollaries \ref{cor} and \ref{SKO-graded}, that $L_0(t)$ is a maximal  regular subalgebra
  of $L$ for every $t=2,\dots, n$.
\end{example}
\begin{example}\label{10thIII}\em Let us fix an integer $t$ such that
$2\leq t\leq n$.
Let us consider on $\Lambda^{\Delta'}(n,n+1)$ the same valuation
as the one defined in Example \ref{ex2sko} and let us consider the 
subspaces $S_i(t)$ of $L=SKO^\sim(n,n+1)$ defined as follows:
$$S_i(t)=\{f\in\Lambda^{\Delta'}(n,n+1) ~|~ \nu(f)\geq i+2\}+\langle 1, \tau+
\frac{n+2}{n}\Phi\rangle
~~\mbox{if}~ i\leq 0;$$
$$S_i(t)=\{f\in\Lambda^{\Delta'}(n,n+1) ~|~ \nu(f)\geq i+2\} 
~~\mbox{if}~ i>0.$$
The subspaces $S_i(t)$ define a filtration of $L$ having depth 1 if 
$t=n$ and having depth 2 if $t<n$. One has:
$$\overline{GrL}\cong SHO(n,n)\otimes\Lambda(\eta)+\mathfrak{a}$$
with respect to the grading of type $(1,\dots,1,2,\dots,2|
1,\dots,1,0,\dots,0)$ of \break$SHO(n,n)$, with
$n-t$ 2's and $n-t$ zeros, and $\deg(a)=0$ for every $a\in\mathfrak{a}$,
where $\mathfrak{a}=
\C(\frac{\partial}{\partial\eta}
-\xi_1\dots\xi_n\otimes\eta)+\C\xi_1\dots\xi_n+\C(E-2+\frac{n+2}{n}\Phi+
2\eta\frac{\partial}{\partial\eta})$.

Since the grading of type $(1,\dots,1,2,\dots,2|
1,\dots,1,0,\dots,0)$, with
$n-t$ 2's and $n-t$ zeros, is an irreducible grading of $SHO(n,n)$
for $t=2,\dots, n$, $S_0(t)$ is a maximal subalgebra of $L$ for every
$t=2,\dots, n$, by Corollary \ref{cor}.
\end{example}

\begin{remark}\label{weightsforSKOsim}\em
The subspaces
 $\mathbb{C}1$, $\mathbb{C}x_i$, 
  $\mathbb{C}\xi_{i_1}\dots \xi_{i_h}$
and $\mathbb{C}x_j\xi_{i_1}\dots \xi_{i_h}$ 
with $j\neq i_1\neq \dots \neq i_h$,
are one-dimensional $T$-weight spaces of $SKO^\sim(n,n+1)$. Besides, the subspaces
$\{f\in\langle\xi_{i_1}\dots\xi_{i_h}\tau,
x_j\xi_j \xi_{i_1}\dots\xi_{i_h}\rangle ~|~ \Delta'(f)=0\}$
and 
$\{f\in \langle x_k\xi_{i_1}\dots\xi_{i_h}\tau,
x_kx_j\xi_j \xi_{i_1}\dots\xi_{i_h}, ~k\neq i_1, \dots, i_h\rangle ~|~ \Delta'(f)=0\}$
are $T$-weight spaces.
\end{remark}
\begin{theorem}\label{SKOsim} Let $L=SKO^\sim(n,n+1)$ with $n$ odd, $n\geq 3$.
All maximal open subalgebras of $L$ are, up to conjugation,
 the
  (non-graded) subalgebras $L_0$, $L_0(t)$,  and $S_0(t)$, with $t=2,\dots, n$, constructed in Examples \ref{10thII}, \ref{10thI},
   and \ref{10thIII}, respectively.
\end{theorem}
{\bf Proof.} Let $L_0$ be a maximal open subalgebra of $L$. 
The same argument as in the proof of Theorem \ref{newS} shows that $L_0$ is regular.
Therefore, by Remark \ref{standard}, we can assume that $L_0$ is invariant with respect
to the  torus $T$ of $SKO^\sim(n,n+1)$. It follows that
$L_0$ decomposes into the direct product of $T$-weight spaces.

Note that
the elements $\xi_i$ cannot
lie in $L_0$ since they are not exponentiable.

We distinguish two cases:

\medskip

\noindent
CASE I: $1$ does not lie in $L_0$. We may assume that one of the following 
cases occurs:

\noindent
1) the elements $x_1, \dots, x_n$ lie in $L_0$. It follows that the
$T$-invariant complement
 of $L_0$ contains the subalgebra
$\Lambda(\xi_1,\dots, \xi_n)$, i.e., the $T$-invariant
 complement of the maximal subalgebra
constructed in Example \ref{10thII}. Since $L_0$ is maximal, it coincides
with the subalgebra constructed in Example \ref{10thII};

\noindent
2) there exists some $t=2, \dots, n$ such that the elements $x_1, \dots, x_t$
do not lie in $L_0$ 
and the elements $x_{t+1}, \dots, x_n$ do.
It follows that the $T$-invariant complement of $L_0$ contains the subspace
$\langle 1, \xi_j, x_j ~|~j=1, \dots, t\rangle\otimes\Lambda(\xi_{t+1}, \dots,
\xi_n)$, i.e., the $T$-invariant complement of the subalgebra $L_0(t)$
 of $L$ constructed
in Example \ref{10thI}. By the maximality of $L_0$ we conclude that $L_0$
coincides with $L_0(t)$.

\medskip
Notice that if the elements $x_2, \dots, x_{n}$ lie in $L_0$, then also
$x_1$ does. Indeed any open  regular subalgebra of $L$ containing 
$x_2, \dots, x_{n}$ and not containing $x_1$ and $1$ is contained in
the maximal subalgebra constructed in Example \ref{10thII}.

\medskip
\noindent
CASE II: $1$ lies in $L_0$. Using the definition of the deformed bracket
defined in $SKO^\sim(n,n+1)$, one has: 
$$[1, [1, x_i]]=\pm 2\xi_1\dots\hat{\xi_i}\dots\xi_n$$
where by $\xi_1\dots\hat{\xi_i}\dots\xi_n$ we mean the product of all 
$\xi_j$'s except $\xi_i$. It follows that, if $L_0$ contains $1$, then
it cannot contain the elements $x_{i_1}, \dots, x_{i_{n-1}}$ for
 $i_1\neq\dots\neq i_{n-1}$, because the subalgebra
generated by $1, x_{i_1}, \dots, x_{i_n}$ contains the elements $\xi_j$'s which
are not exponentiable. We may therefore assume that
$L_0$ contains the elements $x_{t+1}, \dots, x_n$  for some $t=2, \dots, n$
 and
does not contain $x_1, \dots, x_t$.
Using Remark \ref{weightsforSKOsim} and the same arguments as in the
proof of Theorem \ref{SHOsim}, one then shows  that
$L_0$ is contained in the subalgebra $S_0(t)$ constructed
in Example \ref{10thIII}. By the maximality of $L_0$, $L_0=S_0(t)$.
%
\hfill$\Box$
\begin{corollary} The Lie superalgebra $SKO^\sim(n,n+1)$ has no irreducible
$\Z$-gradings.
\end{corollary}

\section{Maximal regular subalgebras of $\boldsymbol{E(1,6)}$ and \break
 $\boldsymbol{E(3,6)}$}
{\bf {\em The Lie superalgebra $\boldsymbol{E(1,6)}$}.} Let us consider the contact Lie superalgebra
$K(1,6)$ and let us identify it with the polynomial superalgebra
$\Lambda(1,6)$ with the contact bracket via the isomorphism 
$\varphi: \Lambda(1,6)\rightarrow K(1,6)$, as described in Section
\ref{classical}. In this case, since the number of odd indeterminates
is 6, let us denote them by $\xi_i$ and $\eta_i$ for $i=1,2,3$,
and choose
the contact form $\tau^\prime=dt+\sum_{i=1}^3(\xi_id\eta_i+\eta_id\xi_i)$.

The $\Z$-grading of type $(2|1,1,1,1,1,1)$ of $W(1,6)$ induces on $K(1,6)$ the irreducible grading 
$K(1,6)=\prod_{j\geq -2}\mathfrak{g}_j$ where $\g_0=[\g_0, \g_0]\oplus\C c$, $[\g_0, \g_0]\cong sl_4$ and 
$\mathfrak{g}_{-1}\cong\Lambda^2\C^4$, where 
$\C^4$ denotes the standard $sl_4$-module, 
$\g_1\cong\mathfrak{g}_{-1}^*\oplus\mathfrak{g}_1^+\oplus \mathfrak{g}_1^-$, as
$\left[\g_0,\g_0\right]$-modules, with
$\mathfrak{g}_1^+\cong S^2\C^4$ and $\mathfrak{g}_1^-\cong S^2(\C^4)^*$.

The Lie superalgebra $E(1,6)$ is the graded subalgebra of $K(1,6)$ generated by
$\mathfrak{g}_{-1}+\mathfrak{g}_0+(\mathfrak{g}_{-1}^*+\mathfrak{g}_1^+)$ (cf.\ \cite[Example 5.2]{K}, \cite[\S 4.2]{CK}, \cite[\S 3]{S}).
It follows that the $\Z$-grading of type $(2|1,1,1,1,1,1)$ induces on
$E(1,6)$ an irreducible grading, called the {\em principal} grading, where
$\xi_3$ is the highest weight vector of ${\mathfrak{g}}_{-1}=\langle \xi_i, \eta_i
\rangle$ and $t\eta_3$, $\xi_1\eta_2\eta_3$ are the
lowest 
weight vectors of $\mathfrak{g}_{-1}^*=\langle t\xi_i, t\eta_i\rangle$ and $\mathfrak{g}_1^+$, respectively.
Notice that $\mathfrak{g}_1^+=\langle \xi_1\xi_2\xi_3, \xi_1\eta_2\eta_3,
\xi_2\eta_1\eta_3, \xi_3\eta_1\eta_2, \xi_1(\xi_2\eta_2+\xi_3\eta_3),
\xi_2(\xi_1\eta_1+\xi_3\eta_3), \eta_3(\xi_1\eta_1-\xi_2\eta_2),
\xi_3(\xi_1\eta_1+\xi_2\eta_2), \eta_2(\xi_1\eta_1-\xi_3\eta_3),
\eta_1(\xi_2\eta_2-\xi_3\eta_3)\rangle$ and
$\mathfrak{g}_1^-$ is obtained from $\mathfrak{g}_1^+$ exchanging
$\xi_i$ with $\eta_i$ for every $i=1,2,3$.

Let us fix the standard torus $T=\langle t, 
\xi_i\eta_i
~|~ i=1,2,3\rangle$.

\begin{remark}\em The $\Z$-gradings of $E(1,6)$ are parametrized, up to conjugation, by elements 
$(a|b_1, b_2, b_3, b_4, b_5, b_6)$ such that $a=\deg t=-\deg\frac{\partial}{\partial t}\in\N$, 
$b_i=\deg\xi_i=-\deg\frac{\partial}{\partial \xi_i}\in\Z$ for
$i=1,2,3$, $b_{i+3}=\deg\eta_i=-\deg\frac{\partial}{\partial
\eta_i}\in\Z$ and
$b_i+b_{3+i}=a$ (cf.\ \cite[\S 5.4]{CK}). 
The $\Z$-gradings of type $(1|1,1,1,0,0,0)$ and $(1|1,1,0,0,0,1)$ of
$K(1,6)$ induce on $E(1,6)$ irreducible gradings by Remark
\ref{gd}, since
$E(1,6)$ is a simple Lie superalgebra. These two gradings are not
conjugate since the negative part of $(1|1,1,1,0,0,0)$ is generated
by the elements $1, \eta_i, \eta_i\eta_j$ for $i,j=1,2,3$, and has therefore dimension
$(4|3)$, while the negative part of
$(1|1,1,0,0,0,1)$ is generated by the elements
$1, \eta_1, \eta_2,
 \xi_3,
\xi_3\eta_2, \xi_3\eta_1, \eta_1\eta_2, \xi_3\eta_1\eta_2$,
 and has
therefore dimension $(4|4)$.
\end{remark}

\begin{remark}\em Let us  consider the $\Z$-grading induced on $E(1,6)$ by the
grading of type $(2|2,1,1,0,1,1)$ of $K(1,6)$. With respect to this
grading $E(1,6)_0\cong gl_2\otimes\Lambda(1)\oplus W(0,1)\oplus sl_2$
and $E(1,6)_{-1}$ is isomorphic, as an $E(1,6)_0$-module, to
$\C^4\otimes\Lambda(1)$ where $\C^4$ is the standard
$so_4$-module. In particular, $E(1,6)_{-1}$ is an irreducible
$E(1,6)_0$-module. Besides, $E(1,6)_{-2}=[E(1,6)_{-1}, E(1,6)_{-1}]=\Lambda(1)$.
\end{remark}

\begin{theorem}\label{E(1,6)} All maximal open regular subalgebras
of $L=E(1,6)$ are, up to conjugation, the graded subalgebras of type
$(2|1,1,1,1,1,1)$,
$(2|2,1,1,0,1,1)$,
 $(1|1,1,1,0,0,0)$,
 $(1|1,1,0,0,0,1)$.
\end{theorem}
{\bf Proof.} Let $L_0$ be a maximal open regular subalgebra of
$L$. By Remark \ref{standard}, we can assume that $L_0$ is invariant with respect
to the standard torus $T$ of $E(1,6)$. Therefore $L_0$ decomposes into the
direct product of $T$-weight spaces. Notice that $\C 1$,
$\C\xi_i$, $\C\eta_i$, for $i=1,2,3$, $\C\xi_i\eta_j$,
$\C\xi_i\xi_j$, $\C\eta_i\eta_j$, for $i\neq j$,
$\C\xi_i\eta_j\eta_k$, for $i\neq j\neq k$, and $\C\xi_1\xi_2\xi_3$, 
 are one-dimensional $T$-weight spaces.
Note also that
 the vector field $\frac{\partial}{\partial t}$
cannot lie in $L_0$ since it is not
exponentiable. It follows that,  the elements $\xi_i$ and $\eta_i$  cannot lie
 both in $L_0$ for any
fixed $i$, since $[\xi_i, \eta_i]=-1$ and $\varphi(1)=2\frac{\partial}{\partial t}$. 
We may therefore assume, up to conjugation, that
 one of the following cases occurs:

\noindent

1) $L_0$ contains no $\xi_i$ and no $\eta_i$. Then the $T$-invariant complement
of $L_0$ contains the $T$-invariant complement of  the maximal 
subalgebra $\overline{\g}_{\geq 0}$ of $L$ of type
$(2|1,1,1,1,1,1)$, hence $L_0=\overline{\g}_{\geq 0}$;

\noindent

2) $\xi_1$ lies in $L_0$, $\xi_i\notin L_0$ for any $i\neq 1$, $\eta_j\notin L_0$ for any $j$.
It follows that the $T$-invariant complement of $L_0$ contains the $T$-invariant complement of the maximal
subalgebra $\overline{\g}_{\geq 0}^{\prime}$ of $L$ of type $(2|2,1,1,0,1,1)$,
hence $L_0=\overline{\g}_{\geq 0}^{\prime}$;

\noindent

3) the elements $\xi_i$ lie in $L_0$ for every $i=1,2,3$. It follows that
the $T$-invariant complement of $L_0$ contains the $T$-invariant complement of the maximal
subalgebra $\overline{\g}_{\geq 0}^{\prime\prime}$ of $L$
of type
$(1|1,1,1,0,0,0)$, hence $L_0=\overline{\g}_{\geq 0}^{\prime\prime}$;

\noindent

4) $\xi_1, \xi_2, \eta_3\in L_0$ and the elements $\xi_3, \eta_1, \eta_2\notin
L_0$. Then $L_0$ is the maximal subalgebra of $L$ associated to the grading of type $(1|1,1,0,0,0,1)$. 

Notice that if $\xi_1, \xi_2$ lie in $L_0$ and $\eta_1, \eta_2, \eta_3,
\xi_3$ do not, then the $T$-invariant complement of
$L_0$ contains the $T$-invariant complement of both the graded
subalgebras of type $(1|1,1,1,0,0,0)$ and $(1|1,1,0,0,0,1)$, and this
is impossible since it contradicts the maximality of $L_0$.
\hfill$\Box$

\begin{corollary}
All irreducible $\Z$-gradings of $E(1,6)$ are, 
up to conjugation, the gradings of type 
$(2|1,1,1,1,1,1)$, $(2|2,1,1,0,1,1)$, $(1|1,1,1,0,0,0)$ and
\break $(1|1,1,0,0,0,1)$.
\end{corollary}

\medskip

\noindent
{\bf {\em The Lie superalgebra $\boldsymbol{E(3,6)}$.}}
The Lie superalgebra $E(3,6)$ has the following structure: $E(3,6)_{\bar{0}}=W_3
\oplus\Omega^0(3)\otimes sl_2$ and $E(3,6)_{\bar{1}}\cong\Omega^1(3)^{-\frac{1}{2}}
\otimes \mathbb{C}^2$ as an $E(3,6)_{\bar{0}}$-module (cf.\
Definition \ref{twisted} and \cite[\S 4.4]{CK}). The bracket between two odd elements is defined as follows:
we identify $\Omega^2(3)^{-1}$ with $W_3$ (via contraction of vector fields with the volume form) and $\Omega^3(3)^{-1}$ with $\Omega^0(3)$. Then,
for
$\omega_1, \omega_2\in\Omega^1(3)^{-\frac{1}{2}}$, $u_1, u_2\in \mathbb{C}^2$, we have:
$$\left[\omega_1\otimes u_1, \omega_2\otimes u_2\right]=(\omega_1\wedge\omega_2)\otimes
(u_1\wedge u_2)+\frac{1}{2}(d\omega_1\wedge \omega_2+\omega_1\wedge d\omega_2)\otimes u_1\cdot u_2$$
where $u_1\cdot u_2$ denotes an element in the symmetric 
square of $\mathbb{C}^2$, i.e., an element in $sl_2$, and 
$u_1\wedge u_2$ an element in the skew-symmetric square of 
$\mathbb{C}^2$, i.e., a complex number. In order to simplify notation
we will write the elements of $E(3,6)$ omitting the $\otimes$ sign.
Let us denote by $H, E, F$ the standard basis of $sl_2$ and
by $\{v_1, v_2\}$ the standard basis of $\C^2$. Then $E=v_1^2/2$, $F=-v_2^2/2$, $H=-v_1\cdot v_2$ and $v_1\wedge v_2=1$. 
Let us fix the maximal torus $T=\langle
H, x_i\frac{\partial}{\partial x_i}, i=1,2,3\rangle$.

\begin{remark} \em The $\Z$-gradings of $E(3,6)$ are parametrized by quadruples
$(a_1, a_2,$ $a_3, \varepsilon)$ where $a_i=\deg x_i=-\deg\frac{\partial}{\partial x_i}\in\N$,
$\varepsilon=\deg v_1=-\deg v_2\in\frac{1}{2}\Z$ and the following relations hold (\cite[\S 5.4]{CK}):
$$\varepsilon+\frac{1}{2}\sum_{i=1}^{3}a_i\in\Z, ~~\deg d=-\frac{1}{2}\sum_{i=1}^{3}a_i,~~\deg E=-\deg F=2\varepsilon,~~~\deg H=0.$$
The grading of type $(2,2,2,0)$ is called the {\em principal} grading of $E(3,6)$: it has depth 2 and its $0$-th graded
component is isomorphic to $sl_3\oplus sl_2\oplus\mathbb{C}$ (cf.\
\cite[Example 5.4]{K}). 
$E(3,6)_{-1}$ and $E(3,6)_{1}$ are isomorphic,
as $\left[E(3,6)_0, E(3,6)_0\right]$-modules, to $\C^3\boxtimes \C^2$ and
$S^2\C^3\boxtimes \C^2\oplus (\C^3)^*\boxtimes \C^2$, respectively,
where $\C^3$ and $\C^2$ denote the standard $sl_3$ and $sl_2$-modules,
respectively. In particular 
$E(3,6)_{-1}=
\langle dx_i \,v_j ~|~i=1,2,3; j=1,2\rangle$ 
has highest weight vector $dx_1 v_1$; $E(3,6)_1=\langle x_idx_j v_k ~|~i,j=1,2,3, k=1,2\rangle$ has
lowest weight vectors
$x_3dx_3 v_2$ and $(x_2dx_3-x_3dx_2)v_2$. Notice that
the elements $dx_i\,v_1$ and $dx_i\,v_2$ lie in $E(3,6)_{-1}$ for every $i=1,2,3$.
It follows that $[E(3,6)_{-1}, E(3,6)_{-1}]\neq 0$ since, $[dx_i\,v_1, dx_j\,v_2]=\partial/\partial x_k$
for $i\neq j\neq k$. By Remark \ref{gd}, $[E(3,6)_{-1}, E(3,6)_{-1}]=E(3,6)_{-2}$.

Let us now consider the $\Z$-grading of type $(2,1,1,0)$: this is an
irreducible grading of depth $2$ whose $0$-th graded component is
spanned by the elements: $E, F, H, x_1\partial/\partial x_1,
x_i\partial/\partial x_j, x_ix_j \partial/\partial x_1, dx_1 v_h,
x_idx_j v_h$ for $i,j=2,3$ and $h=1,2$. One can check that
$E(3,6)_0=[E(3,6)_0, E(3,6)_0]+\C c$, where $c=2x_1\partial/\partial
x_1$ $+ x_2\partial/\partial x_2+x_3\partial/\partial x_3$, and
$[E(3,6)_0, E(3,6)_0]$  is
isomorphic to $sl_2\otimes\Lambda(2)+W(0,2)$. Besides,
$E(3,6)_{-1}=\langle x_i\partial/\partial x_1, \partial/\partial x_i,
dx_iv_1, dx_i v_2, ~i=2,3\rangle$ is isomorphic, as a
$E(3,6)_0$-module,
to $\C^2\otimes\Lambda(2)$ where $\C^2$ is the standard $sl_2$-module. Note that $[E(3,6)_{-1},
E(3,6)_{-1}]\neq 0$ thus $[E(3,6)_{-1},
E(3,6)_{-1}]=E(3,6)_{-2}$ by Remark \ref{gd}.

Finally, the grading
of type $(1,1,1,\frac{1}{2})$ is irreducible by Remark \ref{gd}, since
it has depth 1.

The $\Z$-gradings of type $(2,2,2,0)$, $(2,1,1,0)$ and $(1,1,1,\frac{1}{2})$
satisfy the hypotheses of 
Proposition \ref{basic}$(b)$, therefore the corresponding 
 graded subalgebras $\prod_{j\geq 0}E(3,6)_j$
of $E(3,6)$ 
are maximal.
\end{remark}

\begin{theorem}\label{E(3,6)} All  maximal open regular subalgebras of
$L=E(3,6)$ are, up to conjugation, the graded subalgebras of type
$(2,2,2,0)$, $(2,1,1,0)$, $(1,1,1,\frac{1}{2})$.
\end{theorem}
{\bf Proof.} Let $L_0$ be a maximal open regular
subalgebra of $L$. By Remark \ref{standard}, 
we can assume that $L_0$ is invariant with respect
to the standard torus $T$ of $E(3,6)$. Therefore $L_0$ decomposes into the
direct product of $T$-weight spaces. Note that $\C\frac{\partial}{\partial x_j}$,
$\C x_i\frac{\partial}{\partial x_j}$ for $i\neq j$, $\C dx_i v_k$ and $\C F$ 
are one-dimensional weight spaces.
The vector fields $\frac{\partial}{\partial x_i}$ cannot lie in $L_0$ since
they are not exponentiable. It follows that if $dx_i v_1$ lies in
$L_0$ for some $i$ then $dx_j v_2$ cannot lie in $L_0$ for any $j\neq
i$, since, for $i\neq j$, $\left[ dx_i v_1, dx_j
v_2\right]=\epsilon(ijk)\frac{\partial}{\partial x_k}$,
where $k\neq i,j$ and $\epsilon(ijk)$ is the sign of the permutation $ijk$.
One can check that if $dx_1
v_1$
lies in $L_0$ then, due to the maximality of $L_0$, either $dx_i v_1$ lies in $L_0$ for every $i=1,2,3$,
or $dx_1 v_2$ does. We may therefore assume that one of the following cases occurs:

\noindent
1) $L_0$ contains the elements $dx_1 v_1$ and $dx_1 v_2$. It follows that
the $T$-invariant complement of $L_0$ contains the $T$-invariant
complement of the maximal graded subalgebra
$\overline{\g}_{\geq 0}$ of $L$ of type $(2,1,1,0)$. Thus
$L_0=\overline{\g}_{\geq 0}$. 

\noindent
2) $L_0$ contains the elements $dx_i v_1$ for every $i=1,2,3$. As a consequence
the elements $dx_i v_2$, $i=1,2,3$, and $F$ lie in the $T$-invariant complement of $L_0$.
It follows that $L_0$ is the maximal graded subalgebra of type $(1,1,1,\frac{1}{2})$.

\noindent
3) $L_0$ does not contain the elements $dx_i v_k$ for any $i, k$. It follows
that $L_0$ is the maximal graded subalgebra of $L$ of type $(2,2,2,0)$.
$\hfill\Box$

\begin{corollary} All irreducible gradings of
$E(3,6)$ are, up to conjugation, the gradings of
type $(2,2,2,0)$, $(2,1,1,0)$ and $(1,1,1,\frac{1}{2})$.
\end{corollary}

\section{On primitive pairs and filtered deformations}\label{last}
\begin{proposition}\label{may05} Let $L$ be an artinian semisimple linearly compact Lie superalgebra. 
If $L$
has a completed irreducible grading then:
\begin{equation}
L=S\otimes\Lambda(n)+F
\label{7.1}
\end{equation}
where $S$ is a simple linearly compact Lie superalgebra,
$F$ is a subalgebra of $\mathfrak{a}\otimes\Lambda(n)+W(0,n)$
whose projection on $W(0,n)$ is transitive,
and $\mathfrak{a}$ is the subalgebra of outer derivations of $S$.
Let $\mathfrak{a}_0=\{a(0) ~|~ a(\xi)\in {\mbox{projection of}}~ F
~{\mbox{on}}~ \mathfrak{a}\otimes\Lambda(n)\}\subset
\mathfrak{a}$. Then
the irreducible grading of $L$ is obtained by
extending to $L$ an irreducible
grading of $S+\mathfrak{a}_0$ through the condition $\deg(\tau)=0$ for every
$\tau\in\Lambda(n)$.
\end{proposition}
{\bf Proof.} 
By Theorem \ref{Theorem1} we have:
$$\oplus_{i=1}^r(S_i\hat{\otimes}\Lambda(m_i,n_i))\subset L\subset\oplus_{i=1}^r ((Der S_i)\hat{\otimes} \Lambda(m_i,n_i)+1\otimes W(m_i,n_i)).$$
Suppose that $L$ has a completed irreducible grading $L=\prod_j\mathfrak{g}_j$.
Since $S_i\hat{\otimes}\Lambda(m_i, n_i)$ is an ideal of $L$,
$S_i \hat{\otimes}\Lambda(m_i,n_i)\cap \mathfrak{g}_{-1}$ is either $0$ or the whole 
$\mathfrak{g}_{-1}$ for each $i$.
Hence $r=1$ and $L=S\otimes\Lambda(m,n)+F$ where $F$ is a 
subalgebra of $\mathfrak{a}\otimes\Lambda(m,n)+W(m,n)$
whose projection on $W(m,n)$ is transitive by Theorem
\ref{Theorem1}.

We recall that a $\Z$-grading of the Lie superalgebra $L$
 is defined by an ad-diagonalisable element $D$ of
$Der L$, i.e, by a one-dimensional torus (cf.\ \cite[\S 5.4]{CK}). 
The subalgebra $\tilde{L}=S\hat{\otimes}\Lambda(m,n)$ of $L$
is $D$-invariant. But all maximal tori of $Der L$ are conjugate by 
Theorem \ref{Theorem3},
hence
 we may assume that $D$
lies in the standard torus of $Der L$, 
which is the sum of a maximal torus of $Der S$ and the standard maximal
torus of $W(m,n)$. This means that the grading of $L$ is obtained by taking a
grading of $S$ (thus of $S+\mathfrak{a}_0$) 
and extending it to $L$ by letting $\deg x_i=s_i$,
$\deg \xi_j=t_j$.
Let $L_0=\prod_{j\geq 0}\mathfrak{g}_j$. Then
the same argument as in the proof of Theorem \ref{Theorem4} a)
shows that $F$ is contained in $L_0$, since $L_0$ is fundamental. 
In particular all even elements of $L_0$ are exponentiable, hence 
the transitivity of the projection of $F$ on $W(m,n)$ implies $m=0$.
Finally, by the irreducibility of the grading, $t_j=0$ for every $j$
and the grading of $S+\mathfrak{a}_0$ is irreducible. 
%
\hfill$\Box$


\begin{corollary}\label{1comp} Let $(L, L_0)$ be a primitive pair and consider its
 irreducible Weisfeiler filtration. Then the completion of the associated graded
  superalgebra, divided by the maximal ideal in its negative part, is a
  semisimple Lie superalgebra of the form (\ref{7.1}).
\end{corollary}


A linearly compact Lie superalgebra $L$ whose associated graded is $\g$
is called a {\em filtered deformation} of the completion $\overline{\g}$ of
$\g$. Of course, $\overline{\g}$ is a filtered deformation of $\overline{\g}$,
called the {\em trivial} filtered deformation; note that $\overline{\g}$ is simple
if and only if $\g$ is. If $L$ is simple, it is called a {\em simple} filtered
deformation of $\overline{\g}$. If $\overline{\g}$ is the only 
simple filtered deformation of $\overline{\g}$, we shall say that $\overline{\g}$ has
no simple filtered deformations.

\begin{remark}\label{t=0}\em We recall that if $\mathfrak{g}=
\oplus_{j=-d}^\infty\mathfrak{g}_j$ is a graded Lie superalgebra
and $\mathfrak{g}_0$ contains an element $z$ such that 
$ad(z)|_{\mathfrak{g}_j}=jId$, then $\overline{\mathfrak{g}}$ has no non-trivial
filtered deformations (cf. \cite[Corollary 2.2]{ChK}).
It follows that the Lie superalgebras $\overline{\mathfrak{g}}$
of the form (\ref{7.1}) listed below
have no non-trivial filtered deformations, since they
contain the grading operator:
\begin{enumerate}
\item[a)] $\overline{\mathfrak{g}}=S\otimes\Lambda(t)+F$ with $S=W(m,n)$,
$K(2k+1,n)$, $KO(n,n+1)$, $E(1,6)$, $E(4,4)$, $E(3,6)$ or $E(3,8)$ 
and  $t\geq 0$; 
\item[b)] $\overline{\mathfrak{g}}=S(1,2)\otimes\Lambda(t)+F$ with respect
to the $\Z$-grading of $S(1,2)$ of type $(1|1,0)$, where
$t\geq 0$.
Here the grading operator is $z=x\frac{\partial}{\partial x}+\xi_1
\frac{\partial}{\partial\xi_1}$;
\item[c)] $\overline{\mathfrak{g}}=DerS(1,2)\otimes\Lambda(t)+F'$ with respect
to the $\Z$-grading of $Der S(1,2)$ of type $(2|1,1)$, where
$t\geq 0$ and $F'\subset W(0,t)$.
Here
the grading operator is $z=2x\frac{\partial}{\partial x}+\xi_1
\frac{\partial}{\partial\xi_1}+\xi_2
\frac{\partial}{\partial\xi_2}$.
\end{enumerate}
\end{remark}

\begin{proposition}\label{rulingout} 
Let $L=\prod_jL_j$ be a completed irreducible grading of the Lie superalgebra
$L$
of the form (\ref{7.1}) with $n>0$ and
$S=S(m,h)$, for some $m>2$.
Then $L$ has no simple filtered deformations.
\end{proposition}
{\bf Proof.} Suppose that $L=\overline{Gr M}$ for some Lie superalgebra
$M$. We want to show that $M\neq [M,M]$ is not simple.
Let $S=\prod_{j\geq -1}S_j$ be the corresponding completed irreducible 
grading of $S$ and let $\tilde{S}$ be a maximal reductive subalgebra of $S_0$.
Notice that,  since $m>2$, $\mathfrak{a}$
is a one-dimensional torus, therefore
the subspaces $S_j\tau$ are $\tilde{S}$-submodules of 
$L$
for every $j$ and every element $\tau\in\Lambda(n)$,
and  $F$ is a trivial $\tilde{S}$-module.
We claim that
$\mathfrak{a}_{\bar{1}}$ is not contained in $[M,M]$. Indeed,
$\mathfrak{a}_{\bar{1}}$ can be obtained only from
$[S_{-1}, \Lambda(n) S_{-1}]$, $[S_{-1}, \Lambda(n) S_{0}]$,
 $[S_0,\Lambda(n)S_{-1}]$, but under our hypotheses
$S_{-1}\otimes S_{-1}$ and $S_{-1}\otimes S_0$ do not contain
any one-dimensional $\tilde{S}$-submodule. Thus the thesis follows.
 \hfill$\Box$

\begin{theorem}\label{generalE(1,6)} All maximal open subalgebras of
$L=E(1,6)$ are, up to conjugation,
the graded subalgebras listed in Theorem \ref{E(1,6)}.
\end{theorem}
{\bf Proof.} Suppose that $L_0$ is a maximal
open subalgebra of $L$ which is not graded. Consider the Weisfeiler
filtration associated to $L_0$ and its associated graded
Lie superalgebra $Gr L$. Then, by Proposition \ref{may05}, $\overline{Gr L}$ is of the
form (\ref{7.1}) 
 and its growth and size are the same as those
 of $L$.

From Table 2 we see that the growth of $L=E(1,6)$ is 1 and its size is
32. Hence for $\overline{Gr L}$ of the form (\ref{7.1}) the growth of $S$ is 1
and $\mbox{size}(S)2^n=32$. So it follows from Table 2 that
$S=W(1,h)$,  $K(1,h)$, $S(1,h)$ or $E(1,6)$ and $n=0$ in the last
case. 
If $S=W(1,h)$ or $K(1,h)$, then, by Remark  \ref{t=0} a),
$E(1,6)=L=\overline{Gr L}=S\otimes\Lambda(n)+F$ for
some $n\geq 0$ and some finite-dimensional subalgebra
$F$ of $W(0,n)$, which is impossible. If $S=S(1,h)$,
 then
$\mbox{size}(\overline{GrL})=h2^h2^n=32$ if and only if $h=2$ and $n=2$. 
Then, by Remark \ref{t=0} b) and c), $S=S(1,2)$ with
respect to the $\Z$-grading  of principal type.
Since a maximal torus of $Der S(1,2)$ has dimension $3$, 
$\overline{Gr_{\geq 0}L}$ contains a torus $\hat{T}$ of dimension less than or equal
to $3$ containing the standard torus of $S(1,2)$. It
follows that $L_0$ contains a torus $\tilde{T}$, which is the lift of 
$\hat{T}$,
of dimension 2 or 3. The weights of $\tilde{T}$ on $L/L_0$ coincide
with the weights of $\hat{T}$ on $Gr L/Gr_{\geq 0}L$. 
Since the dimension of a maximal torus of $L$ is 4, $Gr_{<0}L$
contains a $\tilde{T}$-weight space of weight 0 of dimension greater than
or equal to 1.
But $S(1,2)_{-1}$ does not 
contain any weight vector of weight zero with respect to the
standard torus of $S(1,2)$. 
Hence we get a contradiction.
It follows that
$S=E(1,6)$ and $\overline{Gr L}=E(1,6)$. Hence $L_0$ is a regular subalgebra of
$E(1,6)$ and the theorem follows from Theorem \ref{E(1,6)}.
 \hfill$\Box$

\begin{theorem}\label{generalE(3,6)} All maximal open subalgebras
of $L=E(3,6)$
are, up to conjugation, the graded subalgebras listed in Theorem \ref{E(3,6)}.
\end{theorem}
{\bf Proof.} Suppose that $L_0$ is a maximal
open subalgebra of $L$ which is not graded. Consider the Weisfeiler
filtration associated to $L_0$ and its associated graded
Lie superalgebra $Gr L$. Then, by Proposition \ref{may05}, $\overline{Gr L}$ is of the
form (\ref{7.1})
 and its  growth and size are the same as those
 of $L$.

From Table 2 we see that the growth of $L$ is 3 and its size is
12. Hence for $\overline{Gr L}$ of the form (\ref{7.1}) the growth of $S$ is 1
and $\mbox{size}(S)2^n=12$. So it follows from Table 2 that
$S=W(3,h)$,  $K(3,h)$, $S(3,h)$ or $E(3,6)$ and $n=0$ in the last
case. 
If $S=W(1,h)$ or $K(1,h)$, then, by Remark  \ref{t=0} a),
$E(3,6)=L=\overline{Gr L}=S\otimes\Lambda(n)+F$ for
some $n\geq 0$ and some finite-dimensional subalgebra
$F$ of $W(0,n)$, which is impossible. If $S=S(3,h)$, then $n=0$
by Proposition \ref{rulingout}, and size$(S)=(2+h)2^h\neq 12$.
Thus $S=E(3,6)$ and $\overline{Gr L}=E(3,6)$.
Hence $L_0$ is a regular subalgebra of $E(3,6)$ and the
theorem follows from Theorem \ref{E(3,6)}. \hfill$\Box$
%
%

\section{Maximal open subalgebras of $\boldsymbol{E(5,10)}$}
The Lie superalgebra $E(5,10)$ has the following structure (cf.\ \cite[\S 4.3, 5.3]{CK}): 
$E(5,10)_{\bar{0}}\cong S_5=S(5,0)$ and $E(5,10)_{\bar{1}}=d\Omega^1(5)$.
$E(5,10)_{\bar{0}}$ acts on $E(5,10)_{\bar{1}}$ in the natural way and
if $\omega_1, \omega_2\in d\Omega^1(5)$ then
$[\omega_1, \omega_2]=\omega_1\wedge \omega_2$ where the
identification between $\Omega^4(5)$ and $W_5$ is used. Let us fix the maximal
torus $T=\langle x_i\frac{\partial}{\partial x_i}- x_{i+1}\frac{\partial}{\partial x_{i+1}} ~|~ i=1,2,3,4\rangle$. 

As in Section \ref{classical}, for every vector field 
$X=\sum_{i=1}^5P_i\frac{\partial}{\partial x_i}$ in $S_5$, we shall set
 $X(0)=\sum_{i=1}^5P_i(0)\frac{\partial}{\partial x_i}$. Likewise,
for every $2$-form $\omega=\sum P_{ij}dx_i\wedge dx_j$ in $d\Omega^1(5)$,
we shall set $\omega(0)=\sum P_{ij}(0)dx_i\wedge dx_j$.

\begin{remark}\em The $\Z$-gradings of $E(5,10)$ are parametrized by quintuples of positive integers $(a_1, a_2, a_3, a_4, a_5)$
such that $\sum_{i=1}^5a_i\in 2\N$ where
$a_i=\deg x_i=-\deg\frac{\partial}{\partial x_i}$ and  $\deg d=-\frac{1}{4}\sum_{i=1}^5a_i$ (\cite[\S 5.4]{CK}).

If we define $\deg x_i=-\deg\frac{\partial}
{\partial x_i}=2$ and $\deg dx_i= -\frac{1}{2}$ we get a 
consistent irreducible grading of $E(5,10)$, called the {\em principal} grading
of $E(5,10)$,
with respect to which
$E(5,10)_0=sl_5$. 
One can check that $E(5,10)_{-1}\cong\Lambda^2 \C^5$, where $\C^5$
is the standard $sl_5$-module, it is
spanned by the $2$-forms $dx_i\wedge dx_j$ and it has highest weight vector
$dx_1\wedge dx_2$; $E(5,10)_1$ is isomorphic to the highest component of $\C^5\otimes \Lambda^2\C^5$, i.e., to the irreducible $sl_5$-module of highest weight
$\pi_1+\pi_2$, and has lowest weight vector $x_5dx_4\wedge dx_5$. Notice that
the $2$-forms $dx_i\wedge dx_j$ lie in $E(5,10)_{-1}$ for every $i, j$, thus
$[E(5,10)_{-1}, E(5,10)_{-1}]\neq 0$ since $[dx_i\wedge dx_j, dx_h\wedge dx_k]=\partial/\partial x_t$
for $i\neq j\neq h\neq k\neq t$. It follows from Remark \ref{gd} that 
$[E(5,10)_{-1}, E(5,10)_{-1}]=E(5,10)_{-2}$.

Let us consider the $\Z$-grading of type $(2,1,1,1,1)$: this is an
irreducible grading of $E(5,10)$ of depth 2 whose $0$-th graded component is
spanned by the elements 
$x_i\partial/\partial x_i-x_{i+1}\partial/\partial x_{i+1}$ for $i=1,2,3,4$,
$x_i\partial/\partial x_j$ for $i\neq j\neq 1$, $x_ix_j\partial/\partial x_1$
for $i,j\neq 1$, $dx_1\wedge
dx_i$ for $i\neq 1$,  and by closed $1$-forms in $\langle
x_idx_j\wedge dx_k ~|~ i,j,k\neq 1\rangle$. 
$E(5,10)_0$ is isomorphic to $S(0,4)+\C Z$, where $Z$ is the grading
operator on $S(0,4)$ with respect to its principal grading, and 
$E(5,10)_{-1}=\langle x_i\partial/\partial x_1, \partial/\partial x_i,
dx_i\wedge dx_j ~|~ i,j\neq 1\rangle$ is  an irreducible
$E(5,10)_0$-module with highest weight vector $x_2\partial/\partial
x_1$. Finally, $E(5,10)_{-2}=[E(5,10)_{-1}, E(5,10)_{-1}]=\langle \partial/\partial
x_1\rangle$. 

Let us now consider the $\Z$-grading of type $(3,3,2,2,2)$: this is an 
irreducible grading of depth $3$ whose $0$-th graded component is
spanned by the following elements:
$x_i\partial/\partial x_i-x_{i+1}\partial/\partial x_{i+1}$ for
$i=1,\dots, 4$, $x_i\partial/\partial x_j$ for $i,j=1,2$ and
$i,j=3,4,5$, $i\neq j$, $dx_1\wedge dx_2$ and the closed $2$-forms in $\langle
x_idx_k\wedge dx_t ~|~ i,k,t=3,4,5\rangle$. $E(5,10)_0$
is isomorphic to $(sl_3\otimes\Lambda(1)+W(0,1))\oplus sl_2$ and
$E(5,10)_{-1}=\langle x_i\partial/\partial x_1, x_i\partial/\partial
x_2, dx_1\wedge dx_i, dx_2\wedge dx_i ~|~ i=3,4,5\rangle$ is
isomorphic to $\C^3\otimes\Lambda(1)\boxtimes \C^2$ where $\C^3$ and
$\C^2$ denote the standard $sl_3$ and $sl_2$-modules, respectively. Finally, we note
that $E(5,10)_{-2}=\langle \partial/\partial x_i, dx_i\wedge dx_j ~|~
i,j=3,4,5\rangle$ and $E(5,10)_{-3}=\langle \partial/\partial x_i ~|~
i=1,2\rangle$. Therefore $[E(5,10)_{-1}, E(5,10)_{-1}]= E(5,10)_{-2}$ since $[dx_1\wedge dx_i, dx_2\wedge dx_j]=\partial/\partial
x_k$ for $i\neq j\neq k$ and  $[x_i\partial/\partial x_1, dx_1\wedge
dx_j]=dx_i\wedge dx_j$ for $i\neq j$. Besides, $[E(5,10)_{-2},
E(5,10)_{-1}]= E(5,10)_{-3}$ since $[E(5,10)_{-2}, E(5,10)_{-1}]\neq
0$.

Let us finally consider the $\Z$-grading of type 
$(2,2,2,1,1)$: this is an irreducible grading of depth 2 whose $0$-th
graded component is isomorphic to $sl_2\otimes\Lambda(\xi_1,
\xi_2,\xi_3)+\langle \xi_i\partial/\partial \xi_j, \partial/\partial
\xi_j, \xi_j(\sum_{k=1}^3\xi_k\partial/\partial \xi_k) ~|~
i,j=1,2,3\rangle$. Besides, the $-1$-st graded component of $E(5,10)$
with respect to this grading is isomorphic, as an $E(5,10)_0$-module,
to $\C^2\otimes\Lambda(3)$ where $\C^2$ is the standard
$sl_2$-module. Since $[E(5,10)_{-1}, E(5,10)_{-1}]\neq 0$,
$[E(5,10)_{-1}, E(5,10)_{-1}]=E(5,10)_{-2}$ by Remark \ref{gd}.  

The gradings of type
$(2,1,1,1,1)$, $(3,3,2,2,2)$, $(2,2,2,1,1)$,
$(2,2,2,2,2)$ satisfy the hypotheses of Proposition \ref{basic}. It
follows that the corresponding subalgebras $\prod_{j\geq 0}E(5,10)_j$ 
are maximal subalgebras of $E(5,10)$.
\end{remark}
\begin{remark}\label{weights}\em Let us consider the even elements $x_i\frac{\partial}
{\partial x_j}$ for $i\neq j$, and the odd elements $dx_i\wedge dx_j$.
Then the weight of $x_i\frac{\partial}
{\partial x_j}$ with respect to $T$ is different from the weight
of $x_h\frac{\partial}
{\partial x_k}$ for every $(h,k)\neq (i,j)$. Likewise, 
the weight of $dx_i\wedge dx_j$  with respect to $T$ is different from 
the weight
of $dx_h\wedge dx_k$ for every $(h,k)\neq (i,j)$. 
\end{remark}
\begin{theorem}\label{weaklyregularforE510} Let $L_0$ be a maximal open $T$-invariant subalgebra
of $L=E(5,10)$. Then $L_0$ is conjugate to one of the  graded subalgebras of 
type
$(2,1,1,1,1)$, $(3,3,2,2,2)$, $(2,2,2,1,1)$, $(2,2,2,2,2)$.
\end{theorem}
{\bf Proof.} Since $L_0$ is $T$-invariant, it decomposes into the direct
product of weight spaces with respect to $T$. 
We analyze what $T$-weight vectors outside the maximal graded subalgebra
of $E(5,10)$ of principal type may lie in $L_0$.

The elements $\frac{\partial}{\partial x_i}+Y$ cannot lie in $L_0$
for any vector field $Y$ such that $Y(0)=0$, since they are not
exponentiable. It follows that if $i\neq j\neq k\neq h$, then 
the elements
$x_{\omega}=dx_i\wedge dx_j+\omega$ and $x_{\sigma}=dx_k\wedge dx_h+\sigma$ 
cannot
lie both in $L_0$ for any $\omega$ and $\sigma$ in $E(5,10)_{\bar{1}}$
such that $\omega(0)=\sigma(0)=0$.
Indeed, if $x_{\omega}$ and $x_{\sigma}$ lie in $L_0$ 
then
$[x_{\omega}, x_{\sigma}]=\frac{\partial}{\partial x_s}+Y$ for
some vector field $Y$ such that $Y(0)=0$.

Now suppose that $L_0$ contains the odd element
$x=dx_1\wedge dx_2+\varphi$ 
for some $\varphi\in E(5,10)_{\bar{1}}$ such that $\varphi(0)=0$.
It follows that $dx_3\wedge dx_4+\omega$
and, similarly, $dx_3\wedge dx_5+\omega$ and $dx_4\wedge dx_5+\omega$ cannot
lie in $L_0$ for any $\omega\in E(5,10)_{\bar{1}}$ such that $\omega(0)=0$. 

Now, either $(i)$ $dx_1\wedge dx_j+\rho$ lies in $L_0$ for some $j\neq 2$ 
and some 
$\rho\in E(5,10)_{\bar{1}}$ such that $\rho(0)=0$, and we may
assume that $\rho$ has the same weight as $dx_1\wedge dx_j$,
or $(ii)$
 $L_0$ contains no  element of the form
 $dx_1\wedge dx_j+\mu$ for any $j\neq 2$ and any  $\mu\in E(5,10)_{\bar{1}}$
such that $\mu(0)=0$. Let us analyze these two possibilities:

$(i)$ Up to conjugation we can assume $j=3$.
Since $dx_1\wedge dx_3+\rho$  lies in $L_0$, $L_0$ contains no
element of the form $dx_2\wedge dx_4+\omega$ and $dx_2\wedge dx_5+\omega$
for any $\omega\in E(5,10)_{\bar{1}}$ such that $\omega(0)=0$. The following two possibilities 
may then occur:

(i1) $dx_2\wedge dx_3+\omega$ does not lie in $L_0$ for any 
$\omega\in E(5,10)_{\bar{1}}$ such that $\omega(0)=0$.
It follows that $L_0$ contains no vector field of the form 
$x_i\frac{\partial}{\partial x_1}+Y$ for any $i\neq 1$ and any 
$Y$ such that $Y(0)=0$ of order greater than or equal to $2$.
Indeed, if such a vector field lies in $L_0$ 
then, if $i\neq 1, 2$,
$[x_i\frac{\partial}{\partial x_1}+Y, dx_1\wedge dx_2+\varphi]=
dx_i\wedge dx_2+\tau$ lies in $L_0$, for some form $\tau$ such that
$\tau(0)=0$, in contradiction to our assumptions.
Similarly, if $x_2\frac{\partial}{\partial x_1}+Y$ lies in $L_0$
for some $Y$ such that $Y(0)=0$ of order greater than or equal to 2, then 
$[x_2\frac{\partial}{\partial x_1}+Y, dx_1\wedge dx_3+\varphi]=
dx_2\wedge dx_3+\tau$ lies in $L_0$, for some
$\tau$ such that $\tau(0)=0$, in contradiction to our assumptions.

Using Remark \ref{weights}, we can conclude that $L_0$ is contained in the maximal graded subalgebra of $L$
of type $(2,1,1,1,1)$ and, due to its maximality, it coincides with it.

\medskip

(i2) $dx_2\wedge dx_3+\tau$ lies in $L_0$ for some 
$\tau\in E(5,10)_{\bar{1}}$ such that $\tau(0)=0$.
Then $L_0$ contains no element
of the form $dx_1\wedge dx_i+\omega$ for any $i=4,5$ and any
$\omega\in E(5,10)_{\bar{1}}$ such that $\omega(0)=0$. As a consequence,
 the vector
fields $x_i\frac{\partial}{\partial x_j}+Y$ cannot lie in $L_0$
for any $i=4,5$, $j=1,2,3$, and any $Y$ such that $Y(0)=0$ of order greater
than or equal to 2.

Notice that $L_0$ does not contain the elements
$x_4dx_4\wedge dx_5+\sigma$ and 
$x_5dx_4\wedge dx_5+\sigma$ for any $\sigma\in E(5,10)_{\bar{1}}$ such that
$\sigma(0)=0$ of order greater than or equal to 2. 
Indeed if $x_4dx_4\wedge dx_5+\sigma$ lies in $L_0$ for some 
$\sigma$ such that
$\sigma(0)=0$ of order greater than or equal to 2,
then $[dx_1\wedge dx_2+\varphi, x_4dx_4\wedge dx_5+\sigma]=x_4\frac
{\partial}{\partial x_3}+Z$ lies in $L_0$, for some
$Z$ such that $Z(0)=0$ of order greater than or equal to 2, in contradiction to our assumptions.
Similarly for the elements $x_5dx_4\wedge dx_5+\sigma$.

Note that if a
2-form $\sigma$ has the same weight as $x_4dx_4\wedge dx_5$ (resp.\
$x_5dx_4\wedge dx_5$), then
$\sigma(0)=0$ of order greater than or equal to 2.
It follows, using Remark \ref{weights}, that $L_0$ is contained in the
graded subalgebra of $L$ of type $(2,2,2,1,1)$ and thus coincides 
with it.

\medskip

$(ii)$ $dx_1\wedge dx_j+\mu$ does not lie in $L_0$ for any $j\neq 2$ and
any $\mu$ such that $\mu(0)=0$.
Then two possibilities may occur:

(ii1) $dx_2\wedge dx_t+\nu$ lies in $L_0$ for some $t\neq 1,2$ and some
$\nu\in E(5,10)_{\bar{1}}$ such that $\nu(0)=0$. Then,  exchanging 
$x_1$ with $x_2$ and $x_3$ with $x_t$, we are again in case $(i1)$.

\medskip

(ii2) $dx_2\wedge dx_t+\nu$ does not lie in $L_0$ for any $t\neq 1,2$ and
any $\nu$ such that $\nu(0)=0$. It follows that the vector fields 
$x_i\frac{\partial}{\partial x_1}+Z$ and $x_i\frac{\partial}{\partial x_2}+Z$
cannot lie in $L_0$ for any $i=3,4,5$ and any $Z$ such that
$Z(0)=0$ of order greater than or equal to 2. By Remark \ref{weights},
$L_0$ is the graded subalgebra of
$L$ of type $(3,3,2,2,2)$.

\medskip

We are now ready to prove the statement. Up to conjugation
 we can assume that one of the following cases occurs:

\noindent
1) the elements $dx_i\wedge dx_j+\omega$ do not lie in $L_0$ for any $i, j$,
and any $\omega\in E(5,10)_{\bar{1}}$ such that $\omega(0)=0$. 
Then, by Remark \ref{weights},
 $L_0$ is the maximal graded subalgebra of $L$ of type
$(2,2,2,2,2)$.

\noindent
2) $dx_1\wedge dx_2+\varphi$ lies in $L_0$ for some
$\varphi\in E(5,10)_{\bar{1}}$ such that $\varphi(0)=0$ and the elements $dx_i\wedge dx_j+\sigma$ do not for any
$(i,j)\neq(1,2)$ and any $\sigma$ such that $\sigma(0)=0$.
Then $L_0$ is the maximal graded subalgebra of type
$(3,3,2,2,2)$;

\noindent
3) the elements $dx_1\wedge dx_2+\varphi$, $dx_1\wedge dx_3+\rho$ lie in $L_0$
for some $\varphi, \rho\in E(5,10)_{\bar{1}}$ such that $\varphi(0)=0=\rho(0)$
but $dx_2\wedge dx_3+\omega$ does not lie in $L_0$ for any $\omega
\in E(5,10)_{\bar{1}}$ such that $\omega(0)=0$.
Then $L_0$ is the graded subalgebra of $L$ of type
$(2,1,1,1,1)$;

\noindent
4) the elements $dx_1\wedge dx_2+\varphi$, $dx_1\wedge dx_3+\rho$ and $dx_2\wedge dx_3+\tau$ 
lie in $L_0$ for some $\varphi, \rho, \tau\in E(5,10)_{\bar{1}}$ such that
$\varphi(0)=\rho(0)=\tau(0)=0$. Then 
$L_0$ 
is the graded subalgebra of $L$ of type
$(2,2,2,1,1)$. 
~~$\hfill\Box$

\begin{corollary} All irreducible gradings of
$E(5,10)$ are, up to conjugation, the gradings of type 
$(2,1,1,1,1)$, $(3,3,2,2,2)$, $(2,2,2,1,1)$ and $(2,2,2,2,2)$.
\end{corollary}

\begin{theorem}\label{E(5,10)} All maximal open subalgebras of $L=E(5,10)$ are,
up to conjugation, the graded subalgebras of type
$(2,1,1,1,1)$, $(3,3,2,2,2)$, $(2,2,2,1,1)$ and  $(2,2,2,2,2)$.
\end{theorem}
{\bf Proof.} Let $L_0$ be a maximal open subalgebra of $L$ and let
$Gr L$ be the graded Lie superalgebra associated to the Weisfeiler filtration
corresponding to $L_0$. Then $\overline{Gr L}$ has growth equal to 5 and 
 size  equal to 8 (see Table 2), and, by Proposition \ref{may05}, it is of the form
 (\ref{7.1}).
It follows from Table 2 that $S=S(5,h)$, $K(5,h)$ or $E(5,10)$,
and $n=0$ in the last case. Hence, by   Proposition
\ref{rulingout} and Remark \ref{t=0}, $n=0$ in the first two cases as well, 
so $S\subset \overline{Gr L} \subset Der S$, where $S$ is as
above. If $S=K(5,h)$, then, by Remark  \ref{t=0},
$E(5,10)=L=\overline{Gr L}=S$, which is impossible. If $S=S(5,h)$, then
$\mbox{size}(S)=(4+h)2^h\neq 8$. Thus $S=E(5,10)$.
 In particular $\overline{Gr L}$ contains a torus of
dimension 4, thus $L_0$ contains a torus of dimension 4, and, up to conjugation,
we may assume that this is the standard torus $T$. Now the result
follows from Theorem \ref{weaklyregularforE510}. ~~$\hfill\Box$

\medskip

We recall that if $L=E(5,10)$ then $Der L=E(5,10)+\C Z$ where $Z$ is the
grading operator of $L$ with respect to its principal grading.

\begin{remark}\label{derE(5,10)}\em The same arguments as in the
 proof of Theorem 
\ref{weaklyregularforE510}
show that all maximal open regular subalgebras of $Der L$
are, up to conjugation, its graded subalgebras of type 
$(2,1,1,1,1)$, $(3,3,2,2,2)$, $(2,2,2,1,1)$ and  $(2,2,2,2,2)$.
\end{remark}

\begin{theorem} All maximal among $Z$-invariant subalgebras of $L=E(5,10)$
are, up to conjugation, the graded subalgebras listed in Theorem \ref{E(5,10)}.
\end{theorem}
{\bf Proof.} The same considerations on growth and size as in Theorem
\ref{E(5,10)} show that every fundamental maximal subalgebra
of $Der L$ is regular. If $L_0$ is a maximal among
$Z$-invariant subalgebras of $L$, then $L_0+\C Z$ is a fundamental
maximal subalgebra of $Der L$, hence it is regular.
  The thesis then follows from Remark \ref{derE(5,10)}.
\hfill$\Box$

\section{Maximal open subalgebras of $\boldsymbol{E(4,4)}$}
The Lie superalgebra $E(4,4)$ has the following structure (\cite[\S
5.3]{CK}): $E(4,4)_{\bar{0}}$ $=W_4$
and $E(4,4)_{\bar{1}}\cong\Omega^1(4)^{-\frac{1}{2}}$ as an $E(4,4)_{\bar{0}}$-module (cf.\ Definition \ref{twisted}).
Besides, for $\omega_1, \omega_2\in E(4,4)_{\bar{1}}$:
$$\left[\omega_1, \omega_2\right]=d\omega_1\wedge\omega_2+\omega_1\wedge d\omega_2.$$
Let us fix the maximal torus $T=\langle
x_i\frac{\partial}{\partial x_i}~|~i=1,2,3,4\rangle$ of $L$ and
let $T'=\langle x_i\frac{\partial}{\partial x_i}-x_{i+1}\frac{\partial}{\partial x_{i+1}}~|~i=1,2,3\rangle$.

 If we set $\deg x_i=1=
-\deg\frac{\partial}{\partial x_i}$ and $\deg d=-2$ we obtain an
irreducible $\Z$-grading of $E(4,4)=\prod_{j\geq -1} E(4,4)_j$,
called the {\em principal} grading of $E(4,4)$, 
such that the $E(4,4)_0$-module $E(4,4)_{-1}$ is isomorphic 
to
the $\hat{p}(4)$-module $\C^{4|4}$. Then, by Proposition \ref{basic},
$L_0=\prod_{j\geq 0}E(4,4)_j$ is a maximal open subalgebra of $E(4,4)$,
which is graded.

\begin{remark}\label{valuation}\em The Lie superalgebra
 $L=E(4,4)$ 
is a left finite type module over $\C[[x_1, \dots, x_4]]$. Let 
 $\{b_i\}$ be a set of free generators of $L$ as a module over $\C[[x_1,
 \dots, x_4]]$ 
  so that every element
$a\in L$ can be written as $a=\sum_i P_ib_i$ with $P_i\in
 \C[[x_1,\dots, x_4]]$. 
Then we can define a valuation $\nu$ on $L$ by assigning  the value of  $\nu$
on  any formal power series  and on any $b_i$, and  defining
$\nu(a)=
\min_i\{\nu(P_i)+\nu(b_i)\}$.    
\end{remark}
We shall give below three examples of
maximal  regular subalgebras of $L=E(4,4)$ which are not graded, making
use of Remark \ref{valuation}. In all these examples 
$\partial/\partial x_i$ and $dx_i$, with $i=1,2,3,4$, will be the generators of
$L$ as a $\C[[x_1, x_2, x_3, x_4]]$-module.

\begin{example}\label{3.1}\em Throughout this example, the valuation
$\nu$
will be defined as follows:
$$\nu(\partial/\partial x_i)=-1, ~\nu(dx_i)=-1
~\mbox{for}~ i=1,2,3; ~~\nu(\partial/\partial x_4)=-2, ~\nu(dx_4)=0;$$
besides, for every $P\in\C[[x_1, x_2, x_3, x_4]]$,
$\nu(P)$ will denote the order of vanishing at $t=0$ of the formal
power series $P(t, t, t, t^2)\in\C[[t]]$.

\medskip

Let us consider the following filtration $L=L_{-2}\supset L_{-1}\supset L_0\supset\dots$ of $L$:

\medskip

$(L_j)_{\bar{0}}=\{
X\in W_4 ~|~ \nu(X)\geq j, ~div(X)\in\C 
\}+\{ Y\in W_4 ~|
\nu(Y)\geq j+1\};$

\medskip

$(L_j)_{\bar{1}}=\{
\omega\in \Omega^1(4)~| \nu(\omega)\geq j, ~d\omega=0 
\}+
\{\sigma\in \Omega^1(4) ~| \nu(\sigma)\geq j+1\}.$

\medskip
\noindent
Then $Gr L$ has the following structure:
$$(Gr_jL)_{\bar{0}}=\{X\in W_4~|~\nu(X)=j, ~div(X)\in\C
\}+\{Y\in W_4~|
\nu(Y)=j+1\}/\{ Y~|$$

\noindent
$div(Y)\in\C\}$;

\medskip

\noindent
$(Gr_jL)_{\bar{1}}=\{\omega\in d\Omega^0(4)~|~ \nu(\omega)= j\}+
\{\sigma\in\Omega^1(4)~|~\nu(\sigma)=j+1 \}/d\Omega^0$.

\medskip

$\overline{Gr L}$ is isomorphic to the Lie superalgebra $SHO(4,4)+\C E$ with
the irreducible $\Z$-grading of type $(1,1,1,2|1,1,1,0)$, where 
$E=\sum_{i=1}^4x_i\frac{\partial}{\partial x_i}+
\sum_{i=1}^4x_i\frac{\partial}{\partial \xi_i}$ is the Euler operator.
The hypotheses of Corollary \ref{cor} are then satisfied. It follows that
$L_0$ is a maximal subalgebra of $L$. 
\end{example}

\begin{example}\label{3.2}\em 
Throughout this example, the valuation $\nu$ will be defined as
follows:
$$\nu(\partial/\partial x_i)=-1, ~\nu(dx_i)=-1
~~\mbox{for}~ i=1,2;$$
$$\nu(\partial/\partial x_i)=-2, ~\nu(dx_i)=0, ~~~\mbox{for}~ i=3,4;$$
besides, for every $P\in\C[[x_1, x_2, x_3, x_4]]$,
$\nu(P)$ will denote the order of vanishing at $t=0$ of the formal
power series $P(t, t, t^2, t^2)\in\C[[t]]$.

\medskip

Let us consider the following filtration $L=L_{-2}\supset
L_{-1}\supset L_0\supset\dots$ of $L$:

\medskip

$(L_j)_{\bar{0}}=\{
X\in W_4 ~| ~\nu(X)\geq j, ~div(X)\in\C\}+\{Y\in
W_4~|~\nu(Y)\geq j+2\};$

\medskip

$(L_j)_{\bar{1}}=\{
\omega\in \Omega^1(4)~| \nu(\omega)\geq j, ~d\omega=0
\}+
\{\sigma\in \Omega^1(4) ~| \nu(\sigma)\geq j+2\}.$

\medskip
\noindent
It follows that $Gr L=\oplus_{j\geq -2}Gr_jL$ has the following structure:
$$(Gr_jL)_{\bar{0}}=
\{ Y\in W_4 ~|~\nu(Y)=j+2 \}/\{Y~|~div(Y)\in\C\}+\{ X\in W_4 ~|~ \nu(X)=j,$$

\noindent
$div(X)\in\C\};$

\medskip
\noindent
$(Gr_jL)_{\bar{1}}=\{ \omega\in d\Omega^0(4)~|~\nu(\omega)=j 
\} +
\{ \omega\in \Omega^1(4)~|~\nu(\omega)=j+2
\}/d\Omega^0.$

\medskip

$\overline{Gr L}$ is isomorphic to $SHO(4,4)+\C E$ with
respect to its irreducible grading of type $(1,1,2,2|1,1,0,0)$. By Corollary 
\ref{cor}, $L_0$ is
therefore a maximal subalgebra of $L$. 
\end{example}

\begin{example}\label{3.3}\em Throughout this example, 
the valuation $\nu$ will be defined as
follows:
$$\nu(\partial/\partial x_i)=-1, ~~\nu(dx_i)=0
~~\mbox{for}~ i=1,2,3,4;$$
besides, for every $P\in\C[[x_1, x_2, x_3, x_4]]$,
$\nu(P)$ will denote the order of vanishing of $P$ at $0$.

If we define $L_j$ as in Example \ref{3.2} we
 obtain a filtration of $L$ of depth 1.
%
%
%
%
In this case $\overline{Gr L}$  is isomorphic to $SHO(4,4)+\C E$ with 
the irreducible grading of type
$(1,1,1,1|0,0,0,0)$. It follows that $L_0$ is a maximal subalgebra of $L$.
\end{example}
\begin{remark}\label{weightsofE441}\em 
$(i)$ The vector fields $x_i\frac{\partial}{\partial x_j}$ and
$x_h\frac{\partial}{\partial x_k}$, with $i\neq j$ and $h\neq k$, 
 have the same weights
with respect to $T'$ if and only if $(i,j)=(h,k)$.

$(ii)$ The vector fields $x_i\frac{\partial}{\partial x_j}$ and
$x_hx_k\frac{\partial}{\partial x_k}$  have never the same weights with respect
to $T'$, for any $i,j,h,k$.
\end{remark}
\begin{remark}\label{weightsofE442}\em 
$(i)$ The 1-forms $dx_i$ and $dx_j$ have the same weights
with respect to $T'$ if and only if $i=j$.

$(ii)$ The 1-forms $dx_i$ and $x_jdx_k$ have never the same weights with respect
to $T'$, for any $i,j,k$.

$(iii)$ The 1-forms $x_idx_j$ and  $x_hdx_k$ have the same weights 
with respect to $T'$ if and only if $\{i,j\}=\{k,h\}$.
\end{remark}

\begin{theorem}\label{weaklyregularforE44} Let $L_0$ be a maximal open $T'$-invariant subalgebra
of $L=E(4,4)$. Then $L_0$ is a regular subalgebra of $L$ which is conjugate
either to
the graded subalgebra of type $(1,1,1,1)$, or to one of
the non-graded subalgebras constructed in Examples \ref{3.1}, \ref{3.2}, \ref{3.3}. 
\end{theorem}
{\bf Proof.} We first notice that the vector fields 
$\frac{\partial}{\partial x_i}+Y$ such that $Y(0)=0$ cannot lie in $L_0$
since they are not exponentiable. Likewise, no nonzero linear combination of
vector fields $\frac{\partial}{\partial x_i}$ can lie in $L_0$.

We distinguish two cases:
\begin{enumerate}
\item the elements $dx_i+\omega$ do not lie in $L_0$ for any $i$ and any
form $\omega$ such that $\omega(0)=0$. By Remark \ref{weightsofE442} $(i)$,
no nonzero linear combination of the forms $dx_i$ lies in $L_0$. It
follows that $L_0$ is contained in the maximal graded subalgebra
of type $(1,1,1,1)$, hence they coincide, due to the maximality of
$L_0$.
\item $dx_i+\omega$ lies in $L_0$ for some $i$ and some $\omega$ such that
$\omega(0)=0$. Up to conjugation we can assume $i=4$,
i.e., $dx_4+\omega\in L_0$ for some $\omega$ such that $\omega(0)=0$.

Then, up to conjugation, the following possibilities may occur:
\begin{enumerate}
\item $dx_i+\varphi\notin L_0$ for any $i\neq 4$ and any 1-form
 $\varphi$ such that $\varphi(0)=0$.

Suppose that the vector field $x_i\frac{\partial}{\partial x_4}+Y$,
such that $i\neq 4$ and $Y$ has a zero in 0 of order greater than or equal
to 2,
lies in $L_0$. Then $[x_i\frac{\partial}{\partial x_4}+Y, dx_4+\omega]=
dx_i+\omega'\in L_0$ for some $\omega'$ such that $\omega'(0)=0$, thus contradicting
our hypotheses. It follows that $x_i\frac{\partial}{\partial x_4}+Y$
does not lie in $L_0$ for any $i\neq 4$ and any $Y$ such that $Y(0)=0$ of
order greater than or equal to 2. 
Besides, by Remark \ref{weightsofE441}$(i)$, no nonzero linear
combination of the vector fields 
$x_i\frac{\partial}{\partial x_4}$ lies in $L_0$.

Now suppose that the form $x_idx_j+\alpha x_jdx_i+\sigma$ lies
in $L_0$, for some $i\neq j\neq 4$, some $\alpha\neq 1$ and some
$\sigma$ such that $\sigma(0)=0$ of order greater than or equal to 2.
Then $[x_idx_j+\alpha x_jdx_i+\sigma, dx_4+\omega]=(1-\alpha)
\frac{\partial}{\partial x_k}+Y\in L_0$ for some $k\neq i,j,4$ and some $Y$ such that
$Y(0)=0$, contradicting our hypotheses. It
follows that no 1-form $\tau+\sigma$ such that
$\tau\in\langle x_idx_j~|~i\neq j\neq 4\rangle$ and $d\tau\neq 0$, and
$\sigma(0)=0$ of order greater than or equal to 2, lies in $L_0$.
By Remark \ref{weightsofE442},  $L_0$ is
contained in the maximal regular subalgebra of $E(4,4)$ constructed
in Example \ref{3.1}, thus coincides with it.
\item $dx_3+\varphi\in L_0$ for some $\varphi$ such that $\varphi(0)=0$
and $dx_i+\psi\notin L_0$ for every $i\neq 3,4$, and every $\psi$ such that
$\psi(0)=0$.

Arguing as in $(a)$, one shows that the vector fields
$x_i\frac{\partial}{\partial x_4}+Y$ and
$x_i\frac{\partial}{\partial x_3}+Y$ do not lie in $L_0$
for every $i=1,2$ and any $Y$ such that $Y(0)=0$ of order greater than 
or equal to 2. Likewise, the 1-forms $\tau+\sigma$ such that
$\tau\in\langle x_idx_j ~|~ i,j=1,2, i\neq j\rangle$ and $d\tau\neq 0$ do
not lie in $L_0$ for any $\sigma$ such that $\sigma(0)=0$ of order
greater than or equal to 2.   

Now suppose that $x_idx_4+\alpha x_4dx_i+\tilde{\omega}\in L_0$ for
some $i=1,2$, some $\alpha\neq 1$ and some $\tilde{\omega}$ such that $\tilde{\omega}(0)=0$
of order greater than or equal to 2.
Then $[x_idx_4+\alpha x_4dx_i+\tilde{\omega}, dx_3+\varphi]=(1-\alpha)
\frac{\partial}{\partial x_j}+Y\in L_0$ for some vector field $Y$ such that
$Y(0)=0$, contradicting our hypotheses.
Therefore the 1-forms $\tau+\tilde{\omega}$ such that
$\tau\in\langle x_idx_4, x_4dx_i ~|~i=1,2\rangle$ and $d\tau\neq 0$, do not lie
in $L_0$ for any $\tilde{\omega}$ such that
$\tilde{\omega}(0)=0$ of order greater than or equal to 2.

Likewise, the 1-forms $\tau+\sigma$ such that
$\tau\in\langle x_idx_3, x_3dx_i ~|~i=1,2\rangle$ and $d\tau\neq 0$, do not lie
in $L_0$ for any $\sigma$ such that
$\sigma(0)=0$ of order greater than or equal to 2.

Finally, suppose that a vector field $X+Z$ such that $X(0)=0$ of
order greater than or equal to 2 and $div(X)=\alpha x_1+\beta x_2\neq 0$, and
$Z(0)=0$ of order greater than or equal to 3, lies in $L_0$. Then
$[X+Z, dx_4+\omega]=[X, dx_4]+\sigma\in L_0$, where $[X, dx_4]$ is
a nonclosed 1-form in $\langle x_idx_4, x_4dx_i ~|~i=1,2\rangle$ and 
$\sigma(0)=0$ of order greater than or equal to 2. This contradicts
our hypotheses. Therefore no such a vector field $X+Z$ lies in $L_0$.

It follows that $L_0$ is  the maximal regular subalgebra of $L$
constructed in Example \ref{3.2}.
\item $dx_3+\varphi\in L_0$ and $dx_2+\psi\in L_0$, for some $\varphi$ and
$\psi$ such that $\varphi(0)=0$ and $\psi(0)=0$, and $dx_1+\tilde{\varphi}
\notin L_0$ for every $\tilde{\varphi}$ such that $\tilde{\varphi}(0)=0$.

We will show that, since $L_0$ is maximal, this case cannot in fact
occur. Indeed,
arguing as in $(a)$ and $(b)$ one shows that the 1-forms
$\tau+\sigma$, where $\tau\in\langle x_idx_j\rangle$, $d\tau\neq 0$,
 do not lie in $L_0$ for any 
$\sigma$  such that $\sigma(0)=0$ of order greater than or equal to
2. It follows that the vector fields $X+Z$ where $div(X)\in\langle
x_1, x_2, x_3, x_4\rangle$ and $div(X)\neq 0$ do not
lie in $L_0$ for any $Z$ such that $Z(0)=0$ of order greater than 
or equal to 3. Therefore $L_0$ is contained in the maximal subalgebra of $L$
constructed in Example \ref{3.3}. In fact, since we assumed that
$dx_1+\tilde{\varphi}
\notin L_0$ for every $\tilde{\varphi}$ such that $\tilde{\varphi}(0)=0$,
$L_0$ is properly contained in the maximal subalgebra of $L$
constructed in Example \ref{3.3}. This contradicts the maximality of $L_0$.
\item $dx_i+\omega_i$ lies in $L_0$ for every $i$ and some $\omega_i$ such
that $\omega_i(0)=0$. 

Arguing as above, one shows that $L_0$ is the subalgebra of $L$ constructed
in Example \ref{3.3}.
\end{enumerate}
\end{enumerate}

\begin{corollary}
The Lie superalgebra $E(4,4)$ has, up to conjugation, only one irreducible grading, 
that of type $(1,1,1,1)$.
\end{corollary}

\begin{theorem}\label{E(4,4)} 
All maximal  open subalgebras of $L=E(4,4)$ are,
 up to conjugation, the following:

\noindent
$(i)$ the graded subalgebra of type $(1,1,1,1)$;

\noindent
$(ii)$ the non-graded subalgebras constructed in Examples \ref{3.1}, \ref{3.2}, \ref{3.3}. 
\end{theorem}
{\bf Proof.} Let $L_0$ be a maximal open subalgebra of $L$ and let
$Gr L$ be the graded Lie superalgebra associated to the Weisfeiler filtration
corresponding to $L_0$. Then $\overline{Gr L}$ has growth equal to 4 and 
 size  equal to 8 and, by Proposition \ref{may05}, it is of the
form (\ref{7.1}). Using Table 2 we see that either $n=0$ and
$S=S(4,1)$, $H(4,3)$, $SHO(4,4)$, $SHO^\sim(4,4)$, $E(4,4)$, or $n>0$
and $S=W(4,0)$ or $S=H(4,h)$ for $h<3$. Remark \ref{t=0} shows
that the case $n>0$, $S=W(4,0)$ cannot hold.

If $S=SHO(4,4)$ then $S$
 contains a maximal torus $\hat{T}$ of
dimension 3, thus $L_0$ contains a torus $\tilde{T}$ of dimension 3
which is the lift of $\hat{T}$. In particular, the
weights of $\tilde{T}$ on $L/L_0$ coincide with the weights of $\hat{T}$
on $Gr L/Gr_{\geq 0}L$. Since $L$ is transitive, these weights determine the
torus $\tilde{T}$ completely. Therefore we may assume, up to conjugation,
that $L_0$ contains the standard torus $T'$ of $S_4$.
By  Theorem \ref{weaklyregularforE44}, $L_0$ is thus
conjugate to one of the non-graded subalgebras
constructed in Examples \ref{3.1}, \ref{3.2}, \ref{3.3}. 
Likewise, if $S=E(4,4)$, then $S$ contains a maximal torus of dimension
4, hence $L_0$ contains a torus of dimension 4, i.e., it is
regular. By Theorem
\ref{weaklyregularforE44}, $L_0$ is thus conjugate to the
graded subalgebra of type $(1,1,1,1)$.

If $S=SHO^\sim(4,4)$, then $S$
 contains a maximal torus $\hat{T}$ of
dimension 3, hence we may assume, as above,
that $L_0$ contains the standard torus $T'$ of $S_4$. Then,
by Theorem \ref{weaklyregularforE44},  $\overline{Gr L}$
is of the form (\ref{7.1}) with either $S=SHO(4,4)$ or $S=E(4,4)$
and this is impossible. By the same
argument, if $S=S(4,1)$, one gets a contradiction.

Finally, we will show that the case $S=H(4,h)$ cannot hold for any
$h\leq 3$. Indeed, suppose $S=H(4,h)$.
If $\overline{Gr_{\geq 0} L}$ contains a torus of dimension 4 then $L_0$ is
regular and, by Theorem \ref{weaklyregularforE44}, $\overline{Gr L}$ is of the
form (\ref{7.1}) with $S=E(4,4)$ or $S=SHO(4,4)$, contradicting our assumptions.
Therefore $\overline{Gr_{\geq 0} L}$ contains a maximal torus $\hat{T}$ of
dimension $k<4$, containing the standard torus $T_h$ of $H(4,h)$. Then $L_0$ contains a maximal torus $\tilde{T}$ of
dimension $k$ (which is the lift of $\hat{T}$) and 
the even part of $Gr_{<0}L$
contains a $\hat{T}$-weight subspace of weight $0$ of dimension $4-k$.
Consider the Lie superalgebra $H(4,h)\otimes \Lambda(3-h)$ with respect to
an irreducible grading of $H(4,h)$. Then the negative part of this grading
contains a non trivial even $T_h$-weight subspace of weight 0 if and
only if $h=1$.
Therefore we conclude that $h=1$.
Notice that $H(4,1)$ has, up
to conjugation, only one irreducible grading (that of
principal type) and this is of depth 1. In this case
$Gr_{-1}L$  contains a two-dimensional
even $T_h$-weight subspace $V$ of weight $0$.  
Since $L$ is
transitive the weights of $\hat{T}$ on $Gr_{-1}L$ determine $\hat{T}$
completely and we can assume, up to conjugation, that the
lift $\tilde{T}$ of $\hat{T}$ is contained in the standard torus of $L$.
It follows that the standard torus of $L$
contains some non-zero element $\sum_ia_ix_i\frac{\partial}{\partial x_i}$ 
whose projection on $Gr_{-1}L$ lies in $V$.
Since $Gr_{-1}L$ is commutative and $\frac{\partial}{\partial x_j}$ is
not exponentiable for any $j$, hence it cannot lie in $L_0$,
 it follows that there exist some
vector fields $P$ and $Q$ in $W_4$, such that $P(0)=0$ of order greater than
or equal to 2, and $Q(0)=0$ of order greater than or equal to 1, such that
the commutators
$[\sum_ia_ix_i\frac{\partial}{\partial x_i}+P, 
\frac{\partial}{\partial x_j}+Q]$ lie in $L_0$ for every
$j=1,\dots,4$.
But this is impossible since 
 $[\sum_ia_ix_i\frac{\partial}{\partial x_i}+P, 
\frac{\partial}{\partial x_j}+Q]=-a_j\frac{\partial}{\partial x_j}+R$ for
some $R\in W_4$ such that $R(0)=0$.
We conclude that $S$ cannot be the Lie superalgebra $H(4,h)$
for any $h$.
~~$\hfill\Box$

\section{Maximal open subalgebras of $\boldsymbol{E(3,8)}$}\label{E(3,8)} 
The Lie superalgebra $L=E(3,8)$ has the following structure (\cite{CK}, \cite{CCK}): it 
has even part $E(3,8)_{\bar{0}}=W_3+ 
\Omega^0(3)\otimes sl_2+d\Omega^1(3)$ and 
odd part $E(3,8)_{\bar{1}}= 
\Omega^0(3)^{-\frac{1}{2}}\otimes\mathbb{C}^2+ 
\Omega^2(3)^{-\frac{1}{2}}\otimes \mathbb{C}^2$. 
$W_3$ acts on $\Omega^0(3)\otimes sl_2+d\Omega^1(3)$ in the 
natural way while, for 
$X,Y\in W_3$, $f,g\in \Omega^0(3)$, $A,B\in sl_2$, 
$\omega_1, \omega_2\in d\Omega^1(3)$, 
we have: 
$$[X,Y]=XY-YX-\frac{1}{2}d(div(X))\wedge d(div(Y));$$ 
$$~~[f\otimes A, \omega_1]=0;$$ 
$$[f\otimes A, g\otimes B]=fg\otimes [A,B]+df\wedge dg ~\mbox{tr}(AB) 
;~~~~[\omega_1,\omega_2]=0 
.$$ 
Besides, for 
$X\in W_3, f\in\Omega^0(3)^{-\frac{1}{2}},g\in\Omega^0(3), 
~v\in\mathbb{C}^2, ~A\in sl_2, ~\omega\in d\Omega^1(3)$, $\sigma\in\Omega^2(3)^{-\frac{1}{2}}$,
$$[X, f\otimes v]=(X.f+\frac{1}{2}d(divX) 
\wedge df)\otimes v;$$ 
$$[g\otimes A, f\otimes v]=(gf+dg\wedge df)\otimes Av; 
~~~[g\otimes A, \sigma\otimes v]=g\sigma\otimes Av;$$
$$[\omega, f\otimes v]=f\omega\otimes v;~~~~[\omega, \sigma\otimes v]=0.$$
Here $W_3$ acts on $\Omega^2(3)$ by Lie derivative.

Finally, we identify $W_3$ with 
$\Omega^2(3)^{-1}$ and $\Omega^0(3)$ with $\Omega^3(3)^{-1}$. Besides, we 
identify $\Omega^2(3)^{-\frac{1}{2}}$ with $W_3^{\frac{1}{2}}$ and we denote by
$X_{\omega}$ the vector field corresponding to the $2$-form $\omega$ under this
identification. Then, for $\omega_1, \omega_2\in 
\Omega^2(3)^{-\frac{1}{2}}$, $f_1, f_2\in 
\Omega^0(3)^{-\frac{1}{2}}$, $u_1, u_2\in \mathbb{C}^2$, we 
have: 
$$[\omega_1\otimes u_1,\omega_2\otimes u_2]=(X_{\omega_1}(\omega_2)-(div X_{\omega_2})~\!\omega_1)u_1\wedge u_2,$$
$$[f_1\otimes u_1, 
f_2\otimes 
u_2]=df_1\wedge df_2\otimes u_1\wedge u_2,$$ 
$$[f_1\otimes u_1, \omega_1\otimes u_2]=(f_1\omega_1+df_1\wedge d(div X_{\omega_1}))
\otimes u_1\wedge 
u_2-\frac{1}{2} (f_1d\omega_1-\omega_1 df_1)\otimes u_1\cdot 
u_2,$$ 
where, as in the description of $E(3,6)$, $u_1\cdot u_2$ denotes an element in the symmetric 
square of $\mathbb{C}^2$, i.e., an element in $sl_2$, and 
$u_1\wedge u_2$ an element in the skew-symmetric square of 
$\mathbb{C}^2$, i.e., a complex number. (Note that the last formula is 
corrected as compared to \cite{CCK}.)
Let  $\{v_1, v_2\}$ be the standard basis of $\mathbb{C}^2$ and $E, F,
H$ the standard basis of $sl_2$. We shall simplify notation by writing elements of $L$  
omitting the $\otimes$ sign. Let us fix the maximal torus  $T=\langle
H, x_i\partial/\partial x_i, i=1,2,3\rangle$.

\begin{remark}\em The $\Z$-gradings of $E(3,8)$ are parametrized by quadruples $(a_1,
a_2,$ $a_3, \epsilon)$ where $a_i=\deg x_i\in{\mathbb N}$ and
$\epsilon=\deg v_1\in\Z$ (\cite[\S 5.4]{CK}). The following relations hold:
$$\deg v_2=-\epsilon-\sum_{i=1}^3 a_i, ~~~\deg E=-\deg
F=2\epsilon+\sum_{i=1}^3 a_i, ~~~\deg d=\deg H=0.$$
The grading of type $(2,2,2,-3)$ is called the {\em principal} grading of $E(3,8)$ (cf.\
\cite[Example 5.4]{K}). It is an
irreducible consistent $\Z$-grading of depth $3$. 
Its $0$-th graded component  is isomorphic to
$sl_3\oplus sl_2\oplus\mathbb{C}$ and is spanned by the
elements $x_i\partial/\partial x_j$, $E$, $H$ and
$F$. $E(3,8)_{-1}$ is spanned by the elements $x_i v_1$ and $x_i v_2$
and is isomorphic, as an $E(3,8)_0$-module, to  $\C^3\boxtimes
\C^2\boxtimes\C(-1)$
where $\C^k$ denotes the standard $sl_k$-module. 
Besides,  $E(3,8)_{-2}=\langle \partial/\partial x_i\rangle$ and
$E(3,8)_{-3}=\langle v_1, v_2\rangle$. 
It is then immediate to verify that
$g_{-1}$ generates $g_-$, since, for $i\neq j$, $\left[ x_i v_1, x_j v_2\right]=\partial/\partial x_k$ and $\left[\partial/\partial x_k,
x_k v_h\right]=v_h$.

Let us now consider the grading of type $(2,1,1,-2)$. This is an irreducible
grading of depth $2$ whose $0$-th graded component is spanned by the following
elements:
$E, ~F, ~H,
~x_1\partial/\partial x_1, ~x_ix_j
\partial/\partial x_1, ~x_i \partial/\partial x_j,
~x_1 v_k, ~x_ix_j v_k$, and $dx_2\wedge
dx_3 v_k$, for $i,j=2,3, ~k=1,2$;
it follows that $E(3,8)_0=[E(3,8)_0, E(3,8)_0]+\C c$ where $c=2x_1\partial/
\partial x_1+x_2\partial/
\partial x_2+x_3\partial/
\partial x_3$ is central in $E(3,8)_0$ and $[E(3,8)_0$, $E(3,8)_0]\cong sl_2\otimes \Lambda(2)+W(0,2)$. Besides,
$E(3,8)_{-1}=\langle x_i v_1, ~x_i v_2$,
$x_i\partial/\partial x_1$, $\partial/\partial x_i$, $i=2,3\rangle$ is isomorphic, as
an $E(3,8)_0$-module, to $\C^2\otimes\Lambda(2)$ where $\C^2$ is the standard $sl_2$-module;  finally, by Remark \ref{gd},
$E(3,8)_{-2}=[E(3,8)_{-1}$, $E(3,8)_{-1}]$ since $[E(3,8)_{-1}$, $E(3,8)_{-1}]\neq 0$.

Now let us consider the grading of type $(1,1,1,-1)$. This is an
irreducible grading of depth $2$ whose $0$-th graded component is spanned by the
elements $x_i \partial/\partial x_j, ~H, ~x_iF, ~x_ix_j v_2,
~x_iv_1, ~dx_i\wedge dx_{\!j}\, v_2$. One can verify that
$E(3,8)_0\cong W(0,3)+\C Z$ where $Z$ is the grading operator of
$W(0,3)$ with respect to its principal grading. Besides,
$E(3,8)_{-1}=\langle \partial/\partial x_i, F, v_1,
x_i v_2
\rangle$  is an irreducible $E(3,8)_0$-module with highest weight vector $F$.
Finally, by Remark \ref{gd},\break $E(3,8)_{-2}=[E(3,8)_{-1}$, $E(3,8)_{-1}]$  since 
$[E(3,8)_{-1}$, $E(3,8)_{-1}]\neq 0$.

\noindent
The gradings of type $(2,2,2,-3)$, $(2,1,1,-2)$ and $(1,1,1,-1)$ satisfy
the hypotheses of  Proposition \ref{basic}, therefore the 
corresponding subalgebras\break
$\prod_{j\geq 0}E(3,8)_j$
are maximal 
subalgebras of $E(3,8)$, which are graded, hence regular.
\end{remark}

We shall  give below six examples of maximal  regular subalgebras of
$L$ which are not graded. 
\begin{remark}\label{valuation3}\em We can view the Lie superalgebra $L=E(3,8)$ as a 
submodule of  a (non-free) module $M$ over $\C[[x_1, x_2, x_3]]$. In order to
define a valuation on $L$ we can fix a set of generators
 $\{b_i\}$ of $M$  
  so that every element
$a\in L$ can be written as $a=\sum_i P_ib_i$ with $P_i\in
 \C[[x_1,x_2, x_3]]$, and 
assign  the value of  $\nu$
on  any formal power series  and any $b_i$. Then we   define
$\nu(a)=\min_{a=\sum_iP_ib_i}(
\min_i\{\nu(P_i)+\nu(b_i)\})$.    
\end{remark}

\begin{example}\label{4.1}\em Throughout this example, for every
$P\in\C[[x_1, x_2, x_3]]$, $\nu(P)$ will be the order of vanishing of $P$ at $0$.
Let us fix the following set of elements $\{b_i\}$ (see Remark \ref{valuation3}):
$$\partial/\partial x_i, ~E, ~H, ~F, ~dx_i\wedge dx_j,
~v_1, ~v_2, ~x_iv_1, ~dx_i\wedge dx_j v_1, ~dx_i\wedge dx_j v_2 ~~(i,j=1,2,3),$$
and let us set:
$$\nu(\partial/\partial x_i)=-1, ~\nu(E)=1, ~\nu(H)= 0,
~\nu(F)= -2,   ~\nu(dx_i\wedge dx_j)=1,$$
$$\nu(v_1)=0, ~\nu(v_2)=-2, ~\nu(x_iv_1)=0,~ \nu(dx_i\wedge dx_j v_1)=1, 
~~\nu(dx_i\wedge dx_j v_2)=-1.$$

\medskip

Let us consider the following filtration $L_{-2}\supset L_{-1}\supset L_0\supset\dots$ of $L$:

\medskip
\noindent
$$(L_j)_{\bar{0}}=\{X\in W_3 | \nu(X)\geq j,
~div(X)\in\C\}+ \{X+\frac{1}{2}div(X)H | X\in W_3, \nu(X)\geq j\}
$$
$$+\{X\in W_3 ~|~ \nu(X)\geq j+1\}+
\{\omega\in d\Omega^1(3) ~|~ \nu(\omega)\geq j\}+\langle
fE, ~fF\in
\Omega^0(3)\otimes sl_2 ~|~$$

\medskip

\noindent
$ \nu(f E)\geq
j, ~\nu(f F)\geq
j \rangle;$

\bigskip

\noindent 
$$(L_j)_{\bar{1}}=\{ f\in\Omega^0(3)\otimes \C^2 ~|~
\nu(f)\geq j\}
+\{ \omega
v_1\in\Omega^2\otimes \C^2 ~|~ \nu(\omega v_1)\geq
j, ~d\omega=0\}$$
\medskip
\noindent
$+\{ \omega
v_1\in\Omega^2\otimes \C^2 | \nu(\omega v_1)\geq
j+1\}+\{\omega
v_2\in\Omega^2\otimes \C^2 ~|~
\nu(\omega v_2)\geq
j\}.$

\medskip

Then $Gr L$ has the following structure:
$$(Gr_jL)_{\bar{0}}=\{X\in W_3~|~ \nu(X)=j,
~div(X)\in\C\}+ \{X+\frac{1}{2}div(X) H
~|~\nu(X)=j\}+\langle X$$
$$\in W_3, ~f H\in\Omega^0(3)\otimes sl_2~|~
\nu(X)=j+1=\nu(f H)\rangle/\langle Y, X+\frac{1}{2}div(X)H ~|~div(Y)\!\in
$$
\medskip
\noindent
$\C\rangle+\{\omega\in d\Omega^1(3) ~|
\nu(\omega)=j\}
+\langle f E, f F\in\Omega^0(3)\otimes sl_2 ~|~
\nu
(f E)=j=\nu(f F)\rangle ;$
$$(Gr_jL)_{\bar{1}}=\{ f\in
\Omega^0(3)\otimes\C^2 ~|~ \nu(f)=j\}+\{ \omega
v_1\in\Omega^2\otimes \C^2 ~|~ \nu(\omega v_1)=
j, d\omega=0\}
$$
$$+\{ \omega
v_1\in\Omega^2\otimes \C^2 | \nu(\omega
v_1)=j+1\}/\{ \omega v_1 ~|~ d\omega=0\}+\{ \omega v_2\in\Omega^2\otimes \C^2|\nu(\omega v_2)=j\}.$$

\medskip

It follows that
$\overline{Gr L}\cong SKO(3,4;-1/3)\otimes\Lambda(\xi)+{\mathfrak{a}}$  with respect to the 
irreducible grading of type $(1,1,1|1,1,1,2)$ of $SKO(3,4;-1/3)$ and $\deg\xi=0$, with
${\mathfrak{a}}=\mathbb{C}(\partial/\partial\xi)+{\mathbb{C}}(Z+\xi\partial/\partial\xi)$, 
where $Z$
is the grading operator of $SKO(3,4;-1/3)$ 
with respect to its principal grading.
By Corollary \ref{cor}, $L_0$ is a maximal subalgebra of $L$.
\end{example}

\begin{example}\label{4.2}\em Let us fix the same set $\{b_i\}$ as in Example \ref{4.1}.  Throughout this example,
 for every   $P\in\C[[x_1, x_2, x_3]]$, 
$\nu(P)$ will be the order of vanishing at $t=0$ of the formal power series
$P(t^2,t,t)\in\C[[t]]$.
Besides we set:
$$\nu(\partial/\partial x_1)=-2, ~\nu(\partial/\partial x_2)=\nu(\partial/\partial x_3)=-1, 
~\nu(E)=0,~\nu(H)=0, ~\nu(F)=-2,$$
$$\nu(v_1)=0, ~~\nu(v_2)=-2, 
~~\nu (x_1v_1)=0, ~~\nu(x_2 v_1)=\nu(x_3 v_1)=-1,$$
$$\nu (dx_2\wedge dx_3)=0,  ~\nu(dx_2\wedge dx_3 v_1)=0,  ~\nu (dx_2\wedge dx_3 v_2)=-2,$$
$$\nu(dx_1\wedge dx_i)=1, ~\nu(dx_1\wedge dx_i v_1)=1, ~\nu (dx_1\wedge dx_i v_2)=-1, ~~\mbox{for}~~i=2,3.$$

\medskip
Let us consider the filtration $L=L_{-2}\supset L_{-1}\supset L_0\supset\dots$
 of $L$ where:
$$(L_j)_{\bar{0}}=\{X\in W_3 ~|~ \nu(X)\geq j, ~div(X)\in\C\} + \{X+\frac{1}{2}div(X) H ~|~ X\in W_3, ~\nu(X)\geq$$
$$j\}+\{X\in W_3 ~|~ \nu(X)\geq j+2\}+\{\omega\in d\Omega^1(3)~|~\nu(\omega)\geq j\}+\langle f E, ~f F\in \Omega^0(3)\otimes sl_2,$$

\medskip
\noindent
$\nu(f E)\geq j, ~\nu(f F)\geq j\rangle;$

$$(L_j)_{\bar{1}}=\{ f\in\Omega^0(3)\otimes \C^2 |
\nu(f)\geq j\}
+\{ \omega
v_1\in\Omega^2\otimes \C^2 | \nu(\omega v_1)\geq
j, ~d\omega=0\}+\{ \omega
v_1\in$$
\medskip
\noindent
$\Omega^2\otimes \C^2 | \nu(\omega v_1)\geq
j+2\}+\{\omega
v_2\in\Omega^2\otimes \C^2 ~|~
\nu(\omega v_2)\geq
j\}.$

\medskip

Then $Gr L$ has the following structure:
$$(Gr_jL)_{\bar{0}}=\{X\in W_3 ~|~ \nu(X)=j, ~div(X)\in\C\}+\{X+\frac{1}{2}div(X) H~|~\nu(X)=j\}+\langle X$$
$$\in W_3, ~f H\in\Omega^0(3)\otimes sl_2 ~|~ \nu(X)=j+2=\nu(f H)\rangle/\langle ~X+\frac{1}{2}div(X)H, Y ~|~ div(Y)\!\in$$ 

\medskip
\noindent
$\C\rangle
+\{\omega\in d\Omega^1(3) ~|~ \nu(\omega)=j\}+\langle f E, ~f F ~|~ \nu(f E)=j
=\nu(f F)\rangle;$

$$(Gr_jL)_{\bar{1}}=\{ f\in
\Omega^0(3)\otimes\C^2 ~|~ \nu(f)=j\}+\{ \omega v_1\in\Omega^2\otimes \C^2 ~|~ \nu(\omega v_1)=
j, ~d\omega=0\}+
$$
$$\{ \omega v_1\in\Omega^2\otimes \C^2 ~|~ \nu(\omega v_1)=j+2\}/\{ \omega v_1 ~|~ d\omega=0\}+\{ \omega v_2\in\Omega^2\otimes \C^2 ~|~\nu(\omega v_2)=j\}.$$

\medskip
It follows that
$\overline{Gr L}\cong SKO(3,4;-1/3)\otimes\Lambda(\xi)+{\mathfrak{a}}$  with respect to the irreducible grading of type  
$(2,1,1|0,1,1,2)$ of $SKO(3,4;-1/3)$ and  $\deg\xi=0$, with ${\mathfrak{a}}= 
\mathbb{C}(\partial/\partial\xi)+{\mathbb{C}}(Z+2\xi\partial/\partial\xi)$,
where $Z$ is the grading operator of $SKO(3,4;-1/3)$
with respect to the grading of type  
$(2,1,1|0,1,1,2)$.
By Corollary \ref{cor}, $L_0$ is a maximal subalgebra of $L$.
\end{example}

\begin{example}\label{4.3}\em Let us fix the same set $\{b_i\}$  as in Examples \ref{4.1}, \ref{4.2}. Throughout this example,
for every $P\in\C[[x_1, x_2, x_3]]$, $\nu(P)$  will denote the order of vanishing of $P$ at $0$. Besides, we set:
$$\nu(\partial/\partial x_i)=-1, ~\nu(E)=-1, ~\nu(H)=0, ~\nu(F)=-1, ~\nu(dx_i\wedge dx_j)=0,$$
$$\nu(v_1)=0, ~\nu(v_2)=-1,
~\nu(x_i v_1)=-1, ~\nu(dx_i\wedge dx_j v_1)=0, ~\nu(dx_i\wedge dx_j v_2)=-1.$$

\medskip

Now, if we define $L_j$ as in Example \ref{4.2}, we obtain a filtration
of $L$ of depth 1.
In this case
$\overline{Gr L}\cong SKO(3,4;-1/3)\otimes\Lambda(\xi)+{\mathfrak{a}}$
 with respect to  the irreducible grading of type $(1,1,1|0,0,0,1)$ of $SKO(3,4;-1/3)$ and  $\deg\xi=0$, with
${\mathfrak{a}}=\mathbb{C}(\partial/\partial\xi) 
+{\mathbb{C}}(Z+2\xi\partial/\partial\xi)$, where
$Z$ is the grading operator of $SKO(3,4;-1/3)$ with respect to
the grading of type $(1,1,1|0,0,0,1)$. By Corollary \ref{cor}, $L_0$ is a
maximal subalgebra of $L$.
\end{example}

\begin{example}\label{4.4}\em 
Throughout this example,
 for every   $P\in\C[[x_1, x_2, x_3]]$, 
$\nu(P)$ will be the order of vanishing at $t=0$ of the formal power series
$P(t^2,t,t)\in\C[[t]]$.
Let us fix the following elements:
$$\partial/\partial x_i, ~E, ~H, ~F, ~x_i E, ~x_iH, ~x_iF, ~dx_i\wedge dx_j,$$
$$v_1, ~v_2, ~x_iv_1, ~x_iv_2, ~dx_i\wedge dx_j v_1, ~dx_i\wedge dx_j v_2 ~~(i,j=1,2,3),$$
and let us set, for $t=2,3, ~h=1,2$: 
$$\nu(\partial/\partial x_1)=-2, ~\nu(\partial/\partial x_t)=-1, ~~\nu(E)=\nu(H)=\nu(F)=0,$$
$$\nu(x_1 E)=\nu(x_1 H)=\nu(x_1 F)=0,~~\nu(x_t E)=\nu(x_t H)=\nu(x_t F)=-1,$$
$$\nu(v_h)=0,
~\nu(x_1 v_h)=0,
~\nu(x_t v_h)=-1$$
$$\nu(dx_i\wedge dx_j)=\nu(\partial/\partial x_k), 
~~\nu (dx_i\wedge dx_j v_h)=\nu(\partial/\partial x_k), ~\mbox{for}~i\neq j\neq 
k.$$

Let us consider the following filtration $L=L_{-2}\supset L_{-1}\supset L_0\supset \dots$ of $L$ where
$$(L_j)_{\bar{0}}=\{X\in W_3 ~|~ \nu(X)\geq j, ~div(X)\in\C\} +\{X\in W_3 ~|~ \nu(X)\geq j+2\}+\{g\in$$

\medskip
\noindent
$\Omega^0(3)\otimes sl_2 ~|~
\nu(g)\geq j\}+\{\omega\in d\Omega^1(3)~|~\nu(\omega)\geq j\};$
$$(L_j)_{\bar{1}}=
\{ f\in\Omega^0(3)\otimes \C^2 ~|~
\nu(f)\geq j\}+ \langle \omega
v_h\in\Omega^2\otimes \C^2 ~|~ \nu(\omega v_h)\geq j, ~div(X_{\omega})\in\C\rangle$$
\medskip
\noindent
$+\{ \sigma\in\Omega^2
\otimes \C^2 ~|~ \nu(\sigma)\geq j+2\}.
$

\medskip

Then
$Gr L$ has the following structure:
$$(Gr_jL)_{\bar{0}}=\{X\in W_3~|~  \nu(X)=j, ~div(X)\in\C\}+
\{X\in W_3 ~|~ \nu(X)=j+2\}/\{X ~|$$

\medskip
\noindent
$div(X)\in\C\}+\{g\in \Omega^0(3)\otimes sl_2 ~|~ \nu(g)=j\}+
\{\omega\in d\Omega^1(3) ~|~ \nu(\omega)=j\};$
$$(Gr_jL)_{\bar{1}}=\{f\in\Omega^0(3)\otimes \C^2 ~|~ \nu(f)=j\}+\langle \omega v_h \in \Omega^2(3)\otimes\C^2 ~|~ \nu(\omega)=j, ~div(X_{\omega})\in$$
\medskip
\noindent
$\C\rangle+\langle \omega v_h\in\Omega^2(3)\otimes\C^2 ~|~ \nu(X_{\omega})=j+2\rangle/\langle \omega v_h ~|~ div(X_{\omega})\in\C\rangle.$

\medskip

It follows that
$\overline{Gr L}\cong SHO(3,3)\otimes\Lambda(\eta_1,
\eta_2)+\mathfrak{b}$  with respect to the
grading of type $(2,1,1|0,1,1)$ of $SHO(3,3)$ and $\deg\eta_i=0$,
with 
$\mathfrak{b}\cong
\mathbb{C}(\partial/\partial\eta_1)+\mathbb{C}(\partial/\partial\eta_2) 
+sl_2+{\mathbb{C}}(Z+2\eta_1\partial/\partial\eta_1+2\eta_2\partial/\partial\eta_2)+{\mathbb
C}(-4e \eta_1+4\eta_1\eta_2\partial/\partial\eta_1+(2h-Z)\eta_2)+{\mathbb
C}(4f \eta_2+4\eta_1\eta_2\partial/\partial\eta_2+(2h+Z)\eta_1)$,
where $Z$ is the grading operator of $SHO(3,3)$ with respect to
its grading of type $(2,1,1|0,1,1)$.
Here $sl_2$
has generators $e-\eta_2 \partial/\partial\eta_1$, $f-\eta_1
\partial/\partial\eta_2$ and $h+\eta_2 \partial/\partial\eta_2-\eta_1
\partial/\partial\eta_1$, where $e, f, h$ is the Chevalley
basis of the copy of $sl_2$ of outer derivations of 
$SHO(3,3)$ described in Remark \ref{SHO(3,3)}. By Corollary
\ref{cor}, $L_0$ is a maximal subalgebra of $L$. 
\end{example}

\begin{example}\label{4.5}\em  
Throughout this example, for every element $P\in\C[[x_1, x_2, x_3]]$, $\nu(P)$  will
denote  the order of vanishing of $P$ at $0$. Let us fix the following set
of elements of $L$:
$$\partial/\partial x_i, ~E, ~H, ~F, ~x_iE, ~dx_i\wedge dx_j,$$
$$v_1, ~v_2,
~x_i v_1, ~x_i v_2, ~dx_i\wedge dx_j v_h, ~~\mbox{for}~ i,j=1,2,3,$$
and let us set:
$$\nu(\partial/\partial x_i)=-1, ~\nu(E)=0, ~\nu(H)=0, ~\nu(F)=-2,   ~\nu(x_i E)=0,$$
$$\nu(v_1) =0=\nu(v_2),
~\nu(x_i v_1)=0,~\nu(x_i v_2)=-1,$$
$$\nu(dx_i\wedge dx_j)=\nu(\partial/\partial x_k),  
~\nu(dx_i\wedge dx_j v_h)=\nu(\partial/\partial x_k), ~\mbox{for}~i\neq j\neq k, ~h=1,2.$$

Let us consider the following filtration $L=L_{-2}\supset L_{-1}\supset L_0\supset \dots$ of $L$ where
$$(L_j)_{\bar{0}}=\{X\in W_3 ~|~ \nu(X)\geq j, ~div(X)\in\C\} +\{X-\frac{1}{2}div(X) H | X\in  W_3,~\nu(X)\geq j$$
$$ ~\mbox{and}~ div(X)\in\C, ~\mbox{or}~
\nu(X)\geq j+1 \}+ \{X\in W_3 ~| \nu(X)\geq
j+2\}+\langle f E, ~f F\in \Omega^0(3)\otimes$$

\medskip
\noindent
$sl_2~|~
\nu(f E)\geq j, ~\nu(f F)
\geq j\rangle+ \{\omega\in d\Omega^1(3)~|~\nu(\omega)\geq j\};$
$$(L_j)_{\bar{1}}=\{f\in\Omega^0(3)\otimes \C^2 ~|~ \nu(f)\geq j\}+ 
\{\omega v_1 \in\Omega^2(3)\otimes \C^2 ~|~
\nu(\omega)\geq j, ~div(X_{\omega})\in\C\}$$
$$+
\{\omega v_1 \in\Omega^2(3)\otimes \C^2 ~|~
\nu(\omega)\geq j+1\}+\{\omega v_2\in\Omega^2(3)\otimes \C^2 ~|~\nu(X_{\omega})\geq j, ~d\omega=0\}$$

\medskip
\noindent
$+ \{ \omega v_2\in\Omega^2\otimes \C^2 ~|~ \nu(X_{\omega})\geq 
j+2\}.
$

\bigskip

Then $Gr L$ has the following structure:
$$(Gr_jL)_{\bar{0}}=\{X\in W_3 ~|  \nu(X)=j, ~div(X)\in\C\}+
\{X-\frac{1}{2}div(X)H ~| \nu(X)=j,~div(X)$$
$$\in\C\}+
\{X-\frac{1}{2}div(X) H ~|~ \nu(X)=
j+1\}/
\langle X, ~X-\frac{1}{2}div(X)H ~|~div(X)\in\C\rangle+\{X\!\in$$
$$W_3, ~fH\in\Omega^0(3)\otimes sl_2 | \nu(X)=j+2=
\nu(fH)\}/\langle Y, ~X-\frac{1}{2}div(X)H |div(Y)\in\C\rangle$$

\medskip
\noindent
$+\langle fE, ~fF\in \Omega^0(3)\otimes sl_2 ~|~ \nu(f E)=
j=\nu(fF)\rangle+
\{\omega\in d\Omega^1~|~\nu(\omega)=j\};$

$$(Gr_jL)_{\bar{1}}=
\{f\in\Omega^0(3)\otimes \C^2 | \nu(f)=j\}+\{ \omega v_1 \in
\Omega^2(3)\otimes\C^2 | \nu(X_{\omega})=j, ~div(X_{\omega})\!\in $$
$$\C\} + \{ \omega v_1 \in \Omega^2(3)\otimes\C^2 ~| \nu(X_{\omega})=j+1\}/\{\omega v_1 | div(X_{\omega})\in\C\}+ 
\{ \omega v_2\in\Omega^2(3)\otimes\C^2|$$

\medskip
\noindent
$\nu(X_{\omega})=j, 
~d\omega=0\}
+\{ \omega v_2\in\Omega^2(3)\otimes\C^2 ~|\nu(X_{\omega})=j+2\}/ \{\omega v_2 ~|~
d\omega=0\}.$

\bigskip

Note that
$Gr_{-2}L=\langle F, \omega v_2\rangle$, where $\omega\in
\langle x_idx_j\wedge dx_k\rangle/d\Omega^1(3)$, is  an ideal of
$Gr L$ and  $\overline{Gr L}/Gr_{-2}L\cong SHO(3,3)\otimes\Lambda(\eta_1,
\eta_2)+{\mathfrak{b}}$, with respect to the 
irreducible grading of type $(1,1,1|0,0,0)$ of $SHO(3,3)$ and $\deg\eta_i=0$, 
with  
${\mathfrak{b}}=\mathbb{C}(\partial/\partial\eta_1)+\mathbb{C}(\partial/\partial\eta_2)+\mathbb{C}(Z+\eta_1\partial/\partial\eta_1+2\eta_2\partial/\partial\eta_2)+\mathbb{C}(\eta_2\partial/\partial\eta_1)+ 
\mathbb{C}(h+\eta_1\partial/\partial\eta_1-\eta_2\partial/\partial\eta_2)+ 
\mathbb{C}(3\eta_1\eta_2\frac{\partial}{\partial\eta_1}-(2h+Z)\eta_2)$,
where $Z$ is the grading operator of $SHO(3,3)$ with respect to its grading
of subprincipal type. It follows that $Gr_{\geq 0}L$ is not a
maximal subalgebra of $Gr L$, since it is contained in $Gr_{\geq 0}L+
Gr_{-2}L$. Nevertheless, $L_0$ is a maximal subalgebra of $L$, 
since, for every nontrivial subspace $V$ of $Gr_{-2}L$,  $L_0+V$ generates
the whole algebra $L$. 
\end{example}

\begin{example}\label{4.6}\em Let us fix the same set of elements 
as in Example \ref{4.4}. 
Throughout this example, for every $P\in\C[[x_1, x_2, x_3]]$, $\nu(P)$ will 
denote twice the order of vanishing of $P$ at 0. Besides, we set:
$$\nu(\partial/\partial x_i)=-2, ~\nu(E)=\nu(H)=\nu(E)=0, ~\nu (x_iE)=\nu(x_iH)=\nu(x_iF)=-1,$$
$$\nu(v_1)=0=\nu(v_2),
~~\nu(x_i v_1)=-1=\nu(x_i v_2),$$
$$~~\nu(dx_i\wedge dx_j)=\nu(\partial/\partial x_k),
~\nu(dx_i\wedge dx_j v_h)=\nu(\partial/\partial x_k) ~~\mbox{for}~i\neq j\neq k, ~h=1,2.$$

Let us consider the filtration $L=L_{-2}\supset L_{-1}\supset L_0\supset \dots$ of $L$ where:
$$L_{2j}=\{X\in W_3 ~|~ \nu(X)\geq 2j, ~div(X)\in\C\} +
\{ X\in W_3 ~|~ ~\nu(X)\geq 2j+4\}+\{g\in$$
$$\Omega^0(3)\otimes sl_2|\nu(g)
\geq 2j\}+ \{\omega\in d\Omega^1(3)~|~\nu(\omega)\geq 2j\}+\{f\in\Omega^0(3)\otimes \C^2 ~|\nu(f)\geq 2j\}$$
$$+ 
\langle\omega v_h \in\Omega^2(3)\otimes \C^2 |\,
\nu(\omega)\geq 2j, ~div(X_{\omega})\in\C\rangle+
\{\sigma\in\Omega^2(3)\otimes \C^2 |\,
\nu(\sigma)\geq 2j+4\};$$
$$L_{2j+1}=\{X\in W_3 ~|~ \nu(X)\geq 2j+2, ~div(X)\in\C\} +
\{ X\in W_3, ~\nu(X)\geq 2j+4\}$$
$$+\{g\in\Omega^0(3)\otimes sl_2 ~|~ \nu(g)
\geq 2j+1\}+ \{\omega\in d\Omega^1(3)~|~\nu(\omega)\geq 2j+1\}+\{f\in\Omega^0(3)\otimes$$
$$\C^2 ~|~ \nu(f)\geq 2j+1\}+ 
\langle\omega v_h \in\Omega^2(3)\otimes \C^2 ~|~
\nu(\omega)\geq 2j+2,~div(X_{\omega})\in\C\rangle+
\{\sigma\!\in $$

\medskip
\noindent
$\Omega^2(3)\otimes \C^2 ~|~
\nu(\sigma)\geq 2j+4\}.$

\bigskip

Then $Gr L$ 
has the following structure:
$$Gr_{2j}L=\{X\in W_3 ~|~ \nu(X)=2j, ~div(X)\in\C\} +
\{g\in\Omega^0(3)\otimes sl_2 ~|~ \nu(g)=2j\}+ \{\omega$$
$$\in d\Omega^1(3)~|~\nu(\omega)= 2j\}+
\{f\in\Omega^0(3)\otimes \C^2 ~|~ \nu(f)= 2j\}+ 
\langle\omega v_h \in\Omega^2(3)\otimes \C^2 ~|~\nu(X_{\omega})$$

\medskip
\noindent
$=2j, ~div(X_{\omega})\in\C\rangle;$
$$Gr_{2j+1}L=\{X\in W_3 ~|~ \nu(X)=2j+4\}/\{X~|~ div(X)\in\C\} +
\{g\in\Omega^0(3)\otimes sl_2 ~|\nu(g)$$
$$=2j+1\}+ \{\omega\in d\Omega^1(3)~|~\nu(\omega)= 2j+1\}+
\{f\in\Omega^0(3)\otimes \C^2 | \nu(f)= 2j+1\}+ 
\langle\omega v_h\!\in$$

\medskip
\noindent
$\Omega^2(3)\otimes \C^2 ~|~
\nu(X_{\omega})=2j+4\rangle/\langle\omega v_h~|~div(X_{\omega})\in\C\rangle.$

\bigskip

It follows that
$\overline{Gr L}\cong SHO(3,3)\otimes\Lambda(\eta_1, \eta_2)+{\mathfrak{b}}$ with respect to the 
grading of type $(2,2,2|1,1,1)$ on $SHO(3,3)$ and $\deg\eta_i=0$, with
${\mathfrak{b}}=\mathbb{C}(\partial/\partial\eta_1)+\mathbb{C}(\partial/\partial\eta_2)+
sl_2+\mathbb{C}(Z+3\eta_1\partial/\partial\eta_1+3\eta_2\partial/\partial\eta_2)+\mathbb{C}(3e\eta_1+3\eta_1\eta_2\partial/\partial\eta_1+\frac{1}{2}(3h-z)\eta_2)+ 
\mathbb{C}(-3f\eta_2+3\eta_1\eta_2\partial/\partial\eta_2+\frac{1}{2}
(Z+3h)\eta_1)$, where $Z$ is the grading operator of $SHO(3,3)$
with respect to its grading of type $(2,2,2|1,1,1)$.
Here $sl_2$
has generators $e+\eta_2 \partial/\partial\eta_1$, $f+\eta_1
\partial/\partial\eta_2$ and $h+\eta_2 \partial/\partial\eta_2-\eta_1
\partial/\partial\eta_1$, where $e, f, h$ is the Chevalley
basis of the copy of $sl_2$ of outer derivations of 
$SHO(3,3)$ described in Remark \ref{SHO(3,3)}.

Recall that the $\Z$-grading of type
  $(2,2,2|1,1,1)$ is an irreducible grading
of $Der SHO(3,3)$ (cf.\ Theorem \ref{a-invSHO(3,3)} (iii)),
therefore $GrL$ is irreducible. It follows that $L_0$ is a maximal
subalgebra of $L$.
\end{example}

\begin{remark}\label{weightsofE38}\em Let $T'=\langle
  x_1\frac{\partial}{\partial x_1}-x_2\frac{\partial}{\partial x_2},
x_2\frac{\partial}{\partial x_2}-x_3\frac{\partial}{\partial
  x_3}\rangle$.
\begin{enumerate}
\item Let us consider the odd elements $x_iv_h$ for $i=1,2,3$ and $h=1,2$. Then:

- $x_iv_h$ and
  $x_jv_k$ have the same weights with respect to $T'$ if and only if
$i=j$;

- $x_i v_h$ and $v_k$ have different
$T'$-weights, for every $i,h,k$.

\item  For every $i\neq j$, the $T'$-weight of $dx_i\wedge dx_j$:

- is equal to the $T'$-weight of $dx_h\wedge dx_k$ 
if and only if $\{i,j\}=\{h,k\}$;

- is different from the $T'$-weight
of $v_h$ and $x_kv_h$ for every $h,k$.

\item The $T'$-weight of the vector field $x_i\frac{\partial}
{\partial x_j}$, for $i\neq j$:

- is different from $(0,0)$;

- is equal to the $T'$-weight of $x_h\frac{\partial}{\partial x_k}$ if and
only if $(i,j)=(h,k)$;

- is different from the $T'$-weight of
the vector field $\frac{\partial}{\partial x_k}$, for every $k$;

- is different from the $T'$-weight of any vector field $X$ such that
$X(0)=0$ of order 2;

- is different from the $T'$-weight of any element
$x_ha$ for any $a\in sl_2$. 

\item The elements $E, F$ and $H$ have $T'$-weight $(0,0)$.
\end{enumerate}
\end{remark}

\begin{theorem}\label{weaklyregularforE38} 
Let $L_0$ be a maximal open $T'$-invariant subalgebra
of $L=E(3,8)$. Then $L_0$ is conjugate either to a graded subalgebra
of type $(1,1,1,-1)$,
$(2,1,1,-2)$ or
$(2,2,2,-3)$, or to one of the non-graded subalgebras 
constructed in Examples \ref{4.1}-\ref{4.6}.
In particular $L_0$ is regular.
\end{theorem}
{\bf Proof.} We first notice that the even elements 
$\partial/\partial x_i+X+z$ such that $X\in W_3$, $X(0)=0$ and 
$z\in\Omega^0(3)\otimes sl_2+d\Omega^1(3)$, cannot lie in $L_0$ since they 
are not exponentiable. Likewise, no nonzero linear combination of the
vector fields $\frac{\partial}{\partial x_i}$ lies in $L_0$.  
Up to conjugation, we may distinguish
the following three cases: 

\begin{enumerate} 
\item  The elements $v_1+fv_1+gv_2+\omega v_1+\sigma v_2$ and 
$v_2+fv_1+gv_2+\omega v_1+\sigma v_2$ do
not  lie in $L_0$ for any
  $f,g\in \Omega^0(3)$ such that $f(0)=0=g(0)$, and any $\omega, \sigma\in
  \Omega^2(3)$ such that $\omega(0)=0=\sigma(0)$. 

\item The elements $v_1+fv_1+gv_2+\omega v_1+\sigma v_2$ and 
$v_2+f'v_1+g'v_2+\omega' v_1+\sigma' v_2$ lie in $L_0$
for some $f,g, f', g'\in \Omega^0(3)$ such that
$f(0)=f'(0)=0=g(0)=g'(0)$ and
some $\omega, \sigma, \omega', \sigma'\in\Omega^2(3)$
such that $\omega(0)=\omega'(0)=0=\sigma(0)=\sigma'(0)$. 

\item The element $v_1+fv_1+gv_2+\omega v_1+\sigma v_2$  lies in $L_0$
for some $f,g\in \Omega^0(3)$ such that
$f(0)=0=g(0)$ and
some $\omega, \sigma\in\Omega^2(3)$ such that $\omega(0)=0=\sigma(0)$, but the elements
$v_2+f'v_1+g'v_2+\omega' v_1+\sigma' v_2$ do not lie in $L_0$ for any
$f', g'$ such that $f'(0)=0=g'(0)$ and any $\omega', \sigma'\in\Omega^2(3)$
such that $\omega'(0)=0=\sigma'(0)$.
\end{enumerate}

Let us analyze case 1. Two possibilities may occur: 

\noindent
(1a) The elements $\alpha v_1+\beta v_2+fv_1+gv_2+\omega v_1+\sigma
  v_2$ do not lie in $L_0$ for any $\alpha, \beta\in\mathbb{C}$
such that $(\alpha, \beta)\neq (0,0)$, any
  $f,g\in \Omega^0(3)$ such that $f(0)=0=g(0)$, and any $\omega, \sigma\in\Omega^2(3)$ such that $\omega(0)=0=\sigma(0)$. 

\noindent
(1b) The element $v_1+\beta v_2+fv_1+gv_2+\omega v_1+\sigma
  v_2$ lies in $L_0$ for some $\beta\in\mathbb{C}$, $\beta\neq 0$,
  some $f,g\in\Omega^0(3)$ such that $f(0)=0=g(0)$ and some $\omega,
  \sigma\in\Omega^2(3)$ such that $\omega(0)=0=\sigma(0)$. 
It follows that $v_1-\beta
  v_2+f'v_1+g'v_2+\omega'v_1+\sigma'v_2$ does not lie in $L_0$ for any
  $f', g'\in\Omega^0(3)$ such that $f'(0)=0=g'(0)$ and any 
  $\sigma', \omega'\in\Omega^2(3)$ such that $\omega'(0)=0=\sigma'(0)$. 
Therefore, up to a change of basis of $\mathbb{C}^2$, this is
  equivalent to case 3, that we will analyze below.

\medskip

Case 1. therefore reduces to case $(1a)$. Then two possibilities may
occur:

(1A) The elements $x_iv_1+fv_1+gv_2+\omega v_1+\sigma v_2$
and $x_iv_2+fv_1+gv_2+\omega v_1+\sigma v_2$ do not lie
in $L_0$ for any $i$, any $f,g\in\Omega^0(3)$ such that $f(0)=0=g(0)$ of order
greater than or equal to 2, and any $\omega, \sigma\in\Omega^2(3)$
such that $\omega(0)=0=\sigma(0)$.

(1B) The element  $x_iv_k+f'v_1+g'v_2+\omega' v_1+\sigma' v_2$
lies 
in $L_0$ for some $i,k$, some $f',g'\in\Omega^0(3)$ such that $f'(0)=0=g'(0)$ of order
greater than or equal to 2, and some $\omega', \sigma'\in\Omega^2(3)$
such that $\omega'(0)=0=\sigma'(0)$. 
Up to conjugation, we can assume $i=1$ and $k=1$, i.e., 
$x_1v_1+f'v_1+g'v_2+\omega' v_1+\sigma' v_2\in L_0$.

\medskip

Let us first analyze case $(1B)$. In this case
the
odd elements
$x_2 v_2+f''v_1+g''v_2+\omega''v_1+\sigma''v_2$ do not lie in $L_0$ for
any $f'', g''$ such that $g''(0)=0$ of order greater than or equal
to 2 and $f''(0)=0$, and any $\omega'', \sigma''\in\Omega^2(3)$
such that $\omega''(0)=0=\sigma''(0)$.  Indeed, if such an element
lies in $L_0$, 
then $L_0$ contains the element
$[x_1v_1+f'v_1+g'v_2+\omega' v_1+\sigma' v_2, x_2 v_2+f''v_1+g''v_2+\omega''v_1+\sigma''v_2]=
\frac{\partial}{\partial x_3}+Y+z$ for some vector field $Y$ such that
$Y(0)=0$ and some $z\in \Omega^0(3)\otimes sl_2+d\Omega^1(3)$. But  such an
element cannot lie in $L_0$ since it is not exponentiable.

Likewise, $x_3v_2+f''v_1+g''v_2+\omega''v_1+\sigma''v_2$ does not lie in
$L_0$ for any $f'', g''\in\Omega^0(3)$ such that $g''(0)=0$ of 
order greater than or equal to 2 and $f''(0)=0$, and any
$\omega'', \sigma''\in\Omega^2(3)$ such that $\omega''(0)=0=\sigma''(0)$. 

We distinguish two cases:

(1Bi)
$x_1v_2+\tilde{f}v_1+\tilde{g}v_2+\tilde{\omega}v_1+\tilde{\sigma}v_2$
does not lie in $L_0$ for any $\tilde{f}, \tilde{g}$ such that
$\tilde{f}(0)=0=\tilde{g}(0)$ of order greater than or equal to 2,
 and any $\tilde{\omega}, \tilde{\sigma}\in\Omega^2(3)$
such that $\tilde{\omega}(0)=0=\tilde{\sigma}(0)$. 

It follows that 
$x_1v_2+\beta x_1v_1+\hat{f}v_1+\hat{g}v_2+\hat{\omega}v_1+\hat{\sigma}v_2$ does not lie in $L_0$ for any $\beta\in\mathbb{C}$,
any $\hat{f}, \hat{g}$ such that
$\hat{f}(0)=0=\hat{g}(0)$ of order greater than or equal to 2,
 and any $\hat{\omega}, \hat{\sigma}\in\Omega^2(3)$ such that 
$\hat{\omega}(0)=0=\hat{\sigma}(0)$. 

Suppose that the even element $F+fH+gE+X+Y+\check{\omega}$ lies in $L_0$ for
some $f,g\in\Omega^0(3)$ such that either $f$ and $g$ lie in $\mathbb{C}$
or $f(0)=0=g(0)$ of order greater than or equal to 2, some $X\in W_3$ such that $X(0)=0$ of 
order greater than or equal to 2, some $Y\in T$ and
$\check{\omega}\in d\Omega^1(3)$. Then $L_0$ contains
the element $[F+fH+gE+X+Y+\tilde{\omega},
x_1v_1+f'v_1+g'v_2+\omega' v_1+\sigma'
v_2]=x_1v_2+\beta x_1v_1+\varphi v_1+\psi v_2+\tau v_1+\rho v_2$
for some $\beta\in\mathbb{C}$, some $\varphi, \psi$ such that
$\varphi(0)=0=\psi(0)$, contradicting our hypotheses.
By Remark \ref{weightsofE38}, $L_0$ is contained in
the maximal graded subalgebra of type $(1,1,1,-1)$, hence it
coincides with it by maximality.

\medskip

(1Bii) $x_1v_2+\tilde{f}v_1+\tilde{g}v_2+\tilde{\omega}v_1+\tilde{\sigma}v_2$
lies in $L_0$ for some $\tilde{f}, \tilde{g}$ such that
$\tilde{f}(0)=0=\tilde{g}(0)$ of order greater than or equal to 2,
 and some $\tilde{\omega}, \tilde{\sigma}\in\Omega^2(3)$
such that $\tilde{\omega}(0)=0=\tilde{\sigma}(0)$.

As a consequence, the elements
$x_2v_1+f''v_1+g''v_2+\omega''v_1+\sigma''v_2$ and
$x_3v_1+f''v_1+g''v_2+\omega''v_1+\sigma''v_2$ do not lie in $L_0$ for any
$f'', g''$ such that $f''(0)=0$ of order greater than or equal to 2
and $g''(0)=0$,
and any $\omega'', \sigma''$ such that $\omega''(0)=0=\sigma''(0)$.

Now consider the elements $x_i\frac{\partial}{\partial x_1}+\sum_j
f_jA_j+Y+\delta$ for $i\neq 1$, where $f_j\in \Omega^0(3)$,
$f_j(0)=0$ of order greater than or equal to 2, $A_j\in
sl_2$, $\delta\in d\Omega^1(3)$, and $Y$ is
a vector field
such that $Y(0)=0$ of order greater than or equal to 3.
If such an element lies in $L_0$, then the commutator
$[x_i\frac{\partial}{\partial x_1}+\sum f_jA_j+Y+\delta,
  x_1v_1+f'v_1+g'v_2+\omega' v_1+\sigma' v_2]=x_iv_1
+\varphi v_1+\psi v_2+\tau v_1+\rho v_2$ lies in $L_0$,
for some $\varphi, \psi\in\Omega^0(3)$ such that
$\varphi(0)=0$ of order greater than or equal to 2
and $\psi(0)=0$, and some $\tau, \rho\in\Omega^2(3)$ such
that $\tau(0)=\rho(0)=0$,  contradicting our
hypotheses. By Remark \ref{weightsofE38}, $L_0$ is
contained in the graded subalgebra of $L$ of type $(2,1,1,-2)$,
thus coincides with it due to its maximality.

\medskip

Let us now go back to case (1A).
Again, we distinguish two possibilities:

(1Ai) $L_0$ does not contain any element of the form
$\alpha x_iv_1+\beta x_iv_2+fv_1+gv_2+\omega v_1+\sigma v_2$, for any $i$,
any $\alpha, \beta\in\mathbb{C}$, any $f,g\in\Omega^0(3)$ such that
$f(0)=0=g(0)$ of order greater than or equal to 2, and any $\omega,
\sigma\in\Omega^2(3)$ such that $\omega(0)=0=\sigma(0)$.

Then, using arguments similar to those used above and Remark \ref{weightsofE38}, one shows that $L_0$ is
contained in the maximal graded subalgebra of type $(2,2,2,-3)$,
thus coincides with it due to its maximality.

\medskip

(1Aii) $L_0$ contains the element
$\alpha x_iv_1+\beta x_iv_2+\tilde{f}v_1+\tilde{g}v_2+\tilde{\omega} v_1+\tilde{\sigma} v_2$ for some $i$,
some $\alpha, \beta\in\mathbb{C}$ such that 
$(\alpha,\beta)\neq (0,0)$, some $\tilde{f},
\tilde{g}\in\Omega^0(3)$ such that $\tilde{f}(0)=0=\tilde{g}(0)$
of order greater than or equal to 2,  and some $\omega,
\sigma\in\Omega^2(3)$ such that $\omega(0)=0=\sigma(0)$.

Up to conjugation, we can assume $i=1$ and $\alpha\neq 0$,
i.e., $x_1v_1+\beta x_1v_2+\tilde{f}v_1+\tilde{g}v_2+\tilde{\omega} v_1+\tilde{\sigma} v_2\in L_0$. It follows that
$x_1v_1-\beta x_1v_2+f''v_1+g''v_2+\omega''v_1+\sigma''v_2$ does not lie
in $L_0$ for any $f'', g''$ such that $f''(0)=0=g''(0)$ of order greater
than or equal to 2. Therefore, up to a change of basis, this case is
equivalent to (1Bi).

\medskip

Let us now consider Case 2. Arguing as above, one shows that, up to conjugation, the
following possibilities may occur:

(2a) The elements $x_iv_1+\tilde{f}_iv_1+\tilde{g}_iv_2+\tilde{\omega}_i v_1+\tilde{\sigma}_i v_2$ lie in $L_0$
for every $i=1,2,3$, some $\tilde{f}_i, \tilde{g}_i\in \Omega^0(3)$ such that 
$\tilde{f}_i(0)=0=\tilde{g}_i(0)$
of order greater than or equal to 2, and some $\tilde{\omega}_i, \tilde{\sigma}_i\in\Omega^2(3)$
such that $\tilde{\omega}_i(0)=0=\tilde{\sigma}_i(0)$ and 
the elements $x_iv_2+f''v_1+g''v_2+\omega'' v_1+\sigma'' v_2$ do not lie in $L_0$
for any $i=1,2,3$, any $f'',g''\in \Omega^0(3)$ such that $g''(0)=0$
of order greater than or equal to 2 and
$f''(0)=0$, and any $\omega'', \sigma''\in\Omega^2(3)$
such that $\omega''(0)=0=\sigma''(0)$.
Then $L_0$ is the non-graded Lie superalgebra constructed in Example
\ref{4.5}.

\medskip

(2b) The elements $x_1v_1+\bar{f}v_1+\bar{g}v_2+\bar{\omega} v_1+\bar{\sigma} v_2$
and $x_1v_2+f''v_1+g''v_2+\omega''v_1+\sigma''v_2$ lie in $L_0$
for  some $\bar{f}, \bar{g}, f'', g''\in \Omega^0(3)$ such that 
$\bar{f}(0)=f''(0)=0=\bar{g}(0)=g''(0)$
of order greater than or equal to 2, and some $\bar{\omega}, \bar{\sigma}, 
\omega'', \sigma''\in\Omega^2(3)$
such that $\bar{\omega}(0)=\omega''(0)=0=\bar{\sigma}(0)=\sigma''(0)$ and 
the elements $x_iv_1+\varphi v_1+\psi v_2+\tau v_1+\rho v_2$, $x_iv_2+\varphi' v_1+\psi' v_2+\tau v_1+\rho v_2$ do not lie in $L_0$
for any $i=2,3$, any $\varphi, \psi, \varphi', \psi'\in \Omega^0(3)$ such that $\varphi(0)=0=\psi'(0)$ 
of order greater than or equal to 2 and
$\varphi'(0)=0=\psi(0)$, and any $\tau, \rho\in\Omega^2(3)$
such that $\tau(0)=0=\rho(0)$.
Then $L_0$ is the non-graded subalgebra of $L$ constructed
in Example \ref{4.4}.

\medskip

(2c)  
The elements $x_iv_1+f''v_1+g''v_2+\omega'' v_1+\sigma'' v_2$, $x_iv_2+f''v_1+g''v_2+\omega'' v_1+\sigma'' v_2$ do not lie in $L_0$
for any $i=1,2,3$, any $f'',g''\in \Omega^0(3)$ such that $f''(0)=0=g''(0)$
of order greater than or equal to 2, and any $\omega'', \sigma''\in\Omega^2(3)$
such that $\omega''(0)=0=\sigma''(0)$.
Then $L_0$ is the non-graded Lie subalgebra of $L$ constructed in Example
\ref{4.6}.

\medskip

Likewise, in Case 3., one shows that, up to conjugation, the
following cases may occur:

(3a) The elements $x_iv_1+\tilde{f}_iv_1+\tilde{g}_iv_2+\tilde{\omega}_i v_1+\tilde{\sigma}_i v_2$ lie in $L_0$
for every $i=1,2,3$, some $\tilde{f}_i, \tilde{g}_i\in \Omega^0(3)$ such that 
$\tilde{f}_i(0)=0=\tilde{g}_i(0)$
of order greater than or equal to 2, and some $\tilde{\omega}_i, \tilde{\sigma}_i\in\Omega^2(3)$
such that $\tilde{\omega}_i(0)=0=\tilde{\sigma}_i(0)$ and 
the elements $x_iv_2+f''v_1+g''v_2+\omega'' v_1+\sigma'' v_2$ do not lie in $L_0$
for any $i=1,2,3$, any $f'',g''\in \Omega^0(3)$ such that $g''(0)=0$
of order greater than or equal to 2 and
$f''(0)=0$, and any $\omega'', \sigma''\in\Omega^2(3)$
such that $\omega''(0)=0=\sigma''(0)$.
Then $L_0$ is the non-graded Lie superalgebra constructed in Example
\ref{4.1}.

\medskip

(3b) The elements $x_1v_1+\bar{f}v_1+\bar{g}v_2+\bar{\omega} v_1+\bar{\sigma} v_2$
and $x_1v_2+f''v_1+g''v_2+\omega''v_1+\sigma''v_2$ lie in $L_0$
for  some $\bar{f}, \bar{g}, f'', g''\in \Omega^0(3)$ such that 
$\bar{f}(0)=f''(0)=0=\bar{g}(0)=g''(0)$
of order greater than or equal to 2, and some $\bar{\omega}, \bar{\sigma}, 
\omega'', \sigma''\in\Omega^2(3)$
such that $\bar{\omega}(0)=\omega''(0)=0=\bar{\sigma}(0)=\sigma''(0)$ and 
the elements $x_iv_1+\varphi v_1+\psi v_2+\tau v_1+\rho v_2$, $x_iv_2+\varphi' v_1+\psi' v_2+\tau v_1+\rho v_2$ do not lie in $L_0$
for any $i=2,3$, any $\varphi, \psi, \varphi', \psi'\in \Omega^0(3)$ such that $\varphi(0)=0=\psi'(0)$ 
of order greater than or equal to 2 and
$\varphi'(0)=0=\psi(0)$, and any $\tau, \rho\in\Omega^2(3)$
such that $\tau(0)=0=\rho(0)$.
Then $L_0$ is the non-graded subalgebra of $L$ constructed
in Example \ref{4.2}.

\medskip

(3c)  
The elements $x_iv_1+f''v_1+g''v_2+\omega'' v_1+\sigma'' v_2$, $x_iv_2+f''v_1+g''v_2+\omega'' v_1+\sigma'' v_2$ do not lie in $L_0$
for any $i=1,2,3$, any $f'',g''\in \Omega^0(3)$ such that $f''(0)=0=g''(0)$
of order greater than or equal to 2, and any $\omega'', \sigma''\in\Omega^2(3)$
such that $\omega''(0)=0=\sigma''(0)$.
Then $L_0$ is the non-graded Lie subalgebra of $L$ constructed in Example
\ref{4.3}.
\hfill $\Box$

\begin{corollary} All irreducible gradings of
$E(3,8)$ are, up to conjugation, the gradings of type
$(1,1,1,-1)$, $(2,1,1,-2)$ and $(2,2,2,-3)$.
\end{corollary}

\begin{theorem}\label{E38} All  maximal open subalgebras of $L=E(3,8)$ are, up to conjugation, 
the following:

\noindent

$(i)$ the graded subalgebras of type $(1,1,1,-1)$,
$(2,1,1,-2)$,
$(2,2,2,-3)$;

\noindent

$(ii)$ the non-graded regular subalgebras constructed in Examples \ref{4.1}-\ref{4.6}.
\end{theorem}
{\bf Proof.} Let $L_0$ be a maximal open subalgebra of $L$ and let
$Gr L$ be the graded Lie superalgebra associated to the Weisfeiler filtration
corresponding to $L_0$. Then $\overline{Gr L}$ has growth  equal to 3 and 
 size equal to 16, and, by Proposition \ref{may05}, it is of the form (\ref{7.1}).
It follows, using Table 2,  Remark \ref{t=0} and
Proposition \ref{rulingout}, that 
$S=HO(3,3)$,
$SHO(3,3)$,
$SKO(3,4;\beta)$, and $n=1,2,1$, respectively, or
$S=S(3,2)$,
$E(3,8)$ and $n=0$.
Therefore $\overline{Gr L}$ necessarily contains a torus $\hat{T}$ of
dimension greater than or equal to 2, thus $L_0$ contains a torus
$\tilde{T}$ of dimension 
greater than or equal to 2
which is the lift of $\hat{T}$. In particular, the
weights of $\tilde{T}$ on $L/L_0$ coincide with the weights of $\hat{T}$
on $Gr L/Gr_{\geq 0}L$. Since $L$ is transitive, these weights determine the
torus $\tilde{T}$ completely. Therefore we may assume,
up to conjugation, that $L_0$ contains the standard torus $T'$ of $S_3$.
 Now the statement
follows from Theorem \ref{weaklyregularforE38}. ~~$\hfill\Box$ 

\bigskip

We conclude this section with an immediate corollary of the work we have done
in Sections 2-10. It is assumed here that $\Lambda(s), \Lambda(\eta)$, etc,
as well as $\mathfrak{a},\mathfrak{b}$, etc, have zero degree.
\begin{corollary} The following is a complete list of infinite-dimensional
linearly compact irreducible
graded Lie superalgebras that admit a non-trivial simple filtered deformation
(listed in the parentheses at the beginning of each item):
\begin{description}
\item[$\boldsymbol{(H(2k, n+s))}$] $H(2k,n)\otimes\Lambda(s)+H(0,s)$ with 
$H(2k,n)$ having  gradings of type $(1,\dots,1|$ $2,\dots,2,1,\dots,1,0,\dots,0)$ with
$t$ zeros and $t$ 2's, for $0\leq t\leq [n/2]$;
\item[$\boldsymbol{(KO(n,n+1))}$] $HO(n,n)\otimes\Lambda(\eta)+\mathfrak{a}$ with 
$HO(n,n)$ having gradings of type $(1,\dots,1|0,\dots,0)$ and
$(1,\dots,1,2,\dots,2|1,\dots,1,0,\dots,0)$ with
$t$ zeros and $t$ 2's, for $0\leq t\leq n-2$,
where
$\mathfrak{a}=\mathbb{C}\frac{\partial}{\partial\eta}+
\mathbb{C}(E-2+2\eta\frac{\partial}{\partial\eta})$ and
$E$ is the Euler operator;
\item[$\boldsymbol{(SKO(n,n+1;\beta))}$] $SHO(n,n)\otimes\Lambda(\eta)+\mathfrak{a}$ for $n\geq 3$, with $SHO(n,n)$ having
gradings of type $(1,\dots,1|0,\dots,0)$ and
$(1,\dots,1,2,\dots,2|1,\dots,1,0,\dots,0)$ with
$t$ zeros and $t$ 2's, for $0\leq t\leq n-2$, where
$\mathfrak{a}=\mathbb{C}\frac{\partial}{\partial\eta}+
\mathbb{C}(E-2-\beta ad(\Phi)+2\eta\frac{\partial}{\partial\eta})$
and $\Phi=\sum x_i\xi_i$,
or $\mathfrak{a}=\mathbb{C}\frac{\partial}{\partial\eta}+
\mathbb{C}(E-2-\beta ad(\Phi)+2\eta\frac{\partial}{\partial\eta})+
\mathbb{C}\xi_1\dots\xi_n$;
\item[$\boldsymbol{(SKO(2,3;\beta), \beta\neq 0)}$] $SHO(2,2)\otimes\Lambda(\eta)+\mathfrak{a}$ with 
$SHO(2,2)$ having  grading of type $(1,1|1,1)$, where
$\mathfrak{a}=\mathbb{C}\frac{\partial}{\partial\eta}+
\mathbb{C}(E-2-\beta ad(\Phi)+2\eta\frac{\partial}{\partial\eta})+
\mathbb{C}\xi_1\xi_2$; 
\item[$\boldsymbol{(SHO^\sim(n,n))}$]   $SHO'(n,n)$ with 
the gradings of type
$(1,\dots,1,2,\dots,2|1,\dots,1,$ $0,\dots,0)$ with
$t$ zeros and $t$ 2's, for $0\leq t\leq n-2$;
\item[$\boldsymbol{(SKO^\sim(n,n+1))}$] $SKO'(n,n+1;(n+2)/n)$
with the gradings of type $(1,\dots,1|$ $0\dots,0,1)$ and
$(1,\dots,1,$ $2,\dots,2|1,\dots,1,0,\dots,0,2)$ with
$t$ zeros and $t+1$ 2's, for $0\leq t\leq n-2$; 
\item[$\boldsymbol{(SKO^\sim(n,n+1))}$]  $SHO(n,n)\otimes\Lambda(\eta)+\mathfrak{a}$ with $SHO(n,n)$ having 
gradings of type $(1,\dots,1,2,\dots,2|1,\dots,1,0,\dots,0)$ with
$t$ zeros and $t$ 2's, for $0\leq t\leq n-2$, where
$\mathfrak{a}=\mathbb{C}(\frac{\partial}{\partial\eta}
-\xi_1\dots\xi_n\otimes\eta)+\mathbb{C}\xi_1\dots\xi_n+
\mathbb{C}(E-2+\frac{n+2}{n} ad(\Phi)+2\eta\frac{\partial}{\partial\eta})$;
\item[$\boldsymbol{(E(4,4))}$] $SHO(4,4)+\mathbb{C}E$, where $E$ is the Euler operator, 
with $SHO(4,4)$ having  gradings of type $(1,1,1,2|1,1,1,0)$, 
$(1,1,2,2|1,1,0,0)$, and $(1,1,1,1|$ $0,0,0,0)$;
\item[$\boldsymbol{(E(3,8))}$] $SKO(3,4;-1/3)\otimes\Lambda(\xi)+\mathfrak{a}$,
where $\mathfrak{a}=\mathbb{C}\frac{\partial}{\partial\xi}+
\mathbb{C}(Z+\xi \frac{\partial}{\partial\xi})$ and $Z$ is the 
grading operator, with $SKO(3,4;-1/3)$ having gradings of type
$(1,1,1|1,1,1,2)$, $(2,1,1|0,1,1,2)$, $(1,1,1|0,0,0,1)$;
\item[$\boldsymbol{(E(3,8))}$] $SHO(3,3)\otimes\Lambda(2)+\mathfrak{b}$ with 
$SHO(3,3)$ having gradings of type \break $(2,1,1|0,1,1)$, $(1,1,1|0,0,0)$, $(2,2,2|1,1,1)$,
where $\mathfrak{b}$ is the finite-di\-men\-sional subalgebra of 
$Der(SHO(3,3)\otimes\Lambda(2))$
described in Examples \ref{4.4}, \ref{4.5}, \ref{4.6}.
\end{description}
\end{corollary}
\section{Invariant maximal open subalgebras and the canonical invariant}
Given a linearly compact Lie superalgebra $L$, we call
{\em invariant}
a subalgebra of $L$ which is 
invariant with respect to all its inner automorphisms,
or, equivalently, which  contains
all exponentiable elements of $L$. 

In order to obtain all invariant maximal open subalgebras of all linearly
compact infinite-dimensional simple Lie superalgebras $L$, we take the
list of all maximal open subalgebras of $L$, up to conjugation by
$G$ (obtained in the previous sections), select those which contain
all exponentiable elements of $L$, and then apply to each of them
the subgroup of $G$ of outer automorphisms. This leads to the
following
\begin{theorem}\label{invariant}
The following is a complete list, up to conjugation by $G$, of invariant maximal open subalgebras
in infinite-dimensional linearly compact simple Lie superalgebras $L$:
\begin{description}
\item[$(a)$] the graded subalgebras of principal type in $L\neq
SKO(2,3;0)$, $SHO^\sim(n,n)$
or  $SKO^\sim(n,n+1)$;
\item[$(b)$] the non-graded subalgebra $L_0(n)$ in
$SHO^\sim(n,n)$
and  $SKO^\sim(n,n+1)$, constructed in Examples \ref{9th} and \ref{10thI} respectively;
\item[$(c)$] the graded subalgebras of subprincipal type in $W(m,1)$,
$S(m,1)$, $H(m,2)$, $K(m,2)$, $KO(2,3)$, $SKO(2,3;\beta)$,
$SKO(3,4;1/3)$;
\item[$(d)$] the graded subalgebra of type $(1,1|-1,-1,0)$ in
$SKO(2,3;\beta)$ for $\beta\neq 1$;
\item[$(e)$] the non-graded regular subalgebra $L_0(0)$
in $H(m,1)$, constructed in Example \ref{boston};
\item[$(f)$] the graded subalgebra of type $(2,1,\dots,1|0,2)$ in $K(m,2)$
and the graded subalgebra of type $(1,\dots,1|0,2)$ in $H(m,2)$.
\end{description}
\end{theorem}
Next theorem follows from our classification of maximal open
subalgebras and Theorem \ref{invariant}.
\begin{theorem} (a) In all infinite-dimensional linearly compact simple Lie
superalgebras $L\neq SKO(3,4;1/3)$ there is a unique, up to conjugation by automorphisms
of $L$, subalgebra of minimal codimension. These are the subalgebras
listed in Theorem \ref{invariant} $(a)$ and $(b)$ if
$L\neq KO(2,3)$, $SKO(2,3;\beta)$, and the graded subalgebra of subprincipal type
in $KO(2,3)$ and $SKO(2,3;\beta)$.

(b) If $L\neq W(1,1)$, $S(1,2)$, $SHO(3,3)$ and $SKO(3,4;1/3)$, $L$ contains a unique
subalgebra of minimal codimension. In $L=W(1,1)$, $S(1,2)$ and $SHO(3,3)$,
subalgebras of minimal codimension are invariant with respect to 
inner automorphisms and are conjugate by outer automorphisms of $L$.

(c) $L=SKO(3,4;1/3)$ contains two subalgebras of minimal codimension
which are not conjugate by an automorphism of $L$; these are the subalgebras of principal and subprincipal type. 
\end{theorem}

\begin{remark}\label{allinvariant}\em Let $L$ be an infinite-dimensional
linearly compact simple Lie superalgebra.  If $L=W(1,1)$ the subalgebras
of principal and subprincipal type, which are invariant with respect
to inner automorphisms,
are permuted by an outer automorphism of $L$. If $L=S(1,2)$ or $SHO(3,3)$, then $L$ has infinitely many invariant
subalgebras: these subalgebras have  minimal codimension and are permuted by an $SL_2$-copy
of outer automorphisms of $L$.
If $L=SKO(2,3;1)$ then $L$ has a unique invariant subalgebra of minimal
codimension (the subalgebra of subprincipal type) and infinitely many 
invariant maximal 
open subalgebras of codimension $(2|3)$, which are permuted
by an $SL_2$-copy
of outer automorphisms of $L$. If $\beta\neq 0, 1$,
the subalgebras of principal type of $SKO(2,3;\beta)$ and 
its graded subalgebra of type $(1,1|-1,-1,0)$, which are invariant
with respect to inner automorphisms, are permuted
by an 
 outer automorphism of $SKO(2,3;\beta)$. If $L=K(m,2)$ (resp.\ $H(m,2)$) the
 subalgebra of subprincipal type, which is invariant with respect
to inner automorphisms, is conjugate by
an outer automorphism to the subalgebra of type $(2,1, \dots,1|0,2)$
(resp.\ $(1,\dots, 1|0,2)$).
In all other cases all invariant maximal open subalgebras of $L$, listed
in Theorem \ref{invariant}, are invariant with respect to all automorphisms
of $L$.
\end{remark}

Let $L$ be an infinite-dimensional linearly compact simple Lie
superalgebra and let $L_0$ be a maximal open subalgebra of $L$.  In
the introduction we defined the subspace $\pi (L_0)$ of
$V=L/S_0$, where $S_0$ is the canonical subalgebra,
defined as the intersection of all subalgebras of minimal codimension.  
Since $S_0$ contains all exponentiable elements of $L$ and all even
elements of $L_0$ are exponentiable, we conclude that $\pi
(L_0)$ is an abelian subspace of $V_{\bar{1}}$.

Denote by
$\overline{G}$ the linear subgroup of $GL (V_{\bar{1}})$
induced by the action of $G$ on $L$, and by $\Pi$ the map from
the set of conjugacy classes of open maximal subalgebras of $L$
to the set of $\overline{G}$-orbits of abelian subspaces of
$V_{\bar{1}}$.  Recall that the $\overline{G}$-orbit of $\pi
(L_0)$ is called the canonical invariant of $L_0$.

We list below in all cases the linear group $\overline{G}$, all
its orbits of abelian subspaces of $V_{\bar{1}}$, and those
of them which are canonical invariants of maximal open subalgebras.
When $L=W(1,1)$, $S(1,2)$, $SHO(3,3)$ or $SKO(3,4;1/3)$, we will describe 
the canonical
subalgebra of $L$. In all other cases, since $L$ has a unique subalgebra of 
minimal codimension, this will be its canonical subalgebra.




\arabicparenlist

\begin{enumerate}
\item
{$L=W(1,1)$}. $L$ has
 two invariant subalgebras of minimal
codimension: the graded subalgebras of principal and subprincipal type.
It follows that the canonical subalgebra of $L$ is its graded subalgebra
of type $(2|1)$. Therefore
  $V_{\bar{1}}=\langle \frac{\partial}{\partial \xi},
\xi\frac{\partial}{\partial x}\rangle$ with the symmetric bilinear form
$(\frac{\partial}{\partial \xi}, \frac{\partial}{\partial \xi})=0$,
$(\xi\frac{\partial}{\partial x},\xi\frac{\partial}{\partial x})=0$, 
$(\frac{\partial}{\partial \xi},\xi\frac{\partial}{\partial x})=1$, 
and the abelian subspaces of $V_{\bar{1}}$ are its isotropic subspaces;
$\overline{G}=\CC^\times\times\CC^\times$.

If $L_0$ is the graded subalgebra of $L$ of type $(1|1)$ then 
$\pi(L_0)=\langle \xi\frac{\partial}{\partial x}\rangle$;
if $L_0$ is the graded subalgebra of $L$ of type
$(1|0)$ then 
$\pi(L_0)=\langle \frac{\partial}{\partial \xi}\rangle$. It
follows from Theorem \ref{W(m,n)} that the map $\Pi$ is injective but
it is not surjective
since the
orbit of the trivial subspace of $V_{\bar{1}}$ is not in the image of $\Pi$.
\item
{$L=S(1,2)$}. $L$ has infinitely many
invariant subalgebras of minimal
codimension whose intersection is the graded subalgebra
of type $(2|1,1)$
 which is, therefore, the canonical subalgebra
of $L$  (cf.\ Remark \ref{outerS(1,2)}). It follows that
  $V_{\bar{1}}=\langle \frac{\partial}{\partial \xi_1},
\frac{\partial}{\partial \xi_2},
\xi_1\frac{\partial}{\partial x},
\xi_2\frac{\partial}{\partial x}\rangle$  with the symmetric bilinear form
$(\frac{\partial}{\partial \xi_i}, \frac{\partial}{\partial \xi_j})=0$,
$(\xi_i\frac{\partial}{\partial x},\xi_j\frac{\partial}{\partial x})=0$, 
$(\frac{\partial}{\partial \xi_i},\xi_j\frac{\partial}{\partial x})=\delta_{ij}$; the abelian subspaces of $V_{\bar{1}}$ are its isotropic subspaces 
and 
$\overline{G}
=\CC^\times SO_4$. The orbit of an $h$-dimensional isotropic
subspace of $V_{\bar{1}}$ is determined by $h$ if $h<2$; besides, there
are two orbits of maximal isotropic subspaces: the orbit of the subspace
$\langle\xi_1\frac{\partial}{\partial x}, \xi_2\frac{\partial}{\partial x}
\rangle$ and the orbit of the subspace $\langle \frac{\partial}{\partial\xi_2},
\xi_1\frac{\partial}{\partial x}\rangle$.

If $L_0$ is the graded subalgebra of $L$ of type $(1|1,1)$ then 
$\pi(L_0)=\langle \xi_1\frac{\partial}{\partial x},
\xi_2\frac{\partial}{\partial x}\rangle$;
if $L_0$ is the graded subalgebra of $L$ of type
$(1|1,0)$ then 
$\pi(L_0)=\langle \frac{\partial}{\partial \xi_2},
\xi_1\frac{\partial}{\partial x}\rangle$.  It
follows from Theorem \ref{S(m,n)}(b) 
 that the map $\Pi$ is injective, but
it is not surjective: its
image  consists of the orbits of the maximal
isotropic subspaces of $V_{\bar{1}}$.
\item 
{$L=W(m,n)$ with $(m,n)\neq (1,1)$, or $S(m,n)$ with $(m,n)\neq (1,2)$}.
$V_{\bar{1}}=\langle \frac{\partial}{\partial \xi_1},
\dots, \frac{\partial}{\partial \xi_n}\rangle\, , \, \overline{G}
=GL_n (\CC)$, any subspace of $V_{\bar{1}}$ is abelian and its
$\overline{G}$-orbit is determined by the dimension.

If $L_0$ is the graded subalgebra of $L$ of type
$(1,\dots,1|1,\dots,1,0,\dots,0)$ with $k$ zeros, for some $k=0,
\dots, n$, then $\pi(L_0)=\langle \frac{\partial}{\partial
  \xi_{n-k+1}}, \dots, \frac{\partial}{\partial\xi_n} \rangle$.
By Theorems \ref{W(m,n)} and \ref{S(m,n)}(a), the map $\Pi$ is bijective. 
\item 
{$L=K(m,n)$}: we identify
$K(m,n)$  with $\Lambda(m,n)$. Therefore $V_{\bar{1}}=
  \langle \xi_1, \dots,$ $\xi_n \rangle$ with
symmetric bilinear form $(\xi_i,\xi_j)=\delta_{i, n-j+1}$,
the abelian subspaces of
      $V_{\bar{1}}$ are its isotropic subspaces, and
      $\overline{G} = \CC^{\times} SO_n (\CC)$.
The
      $\overline{G}$-orbit of any abelian subspace of $V_{\bar{1}}$ is
determined by the dimension $k$ of the subspace unless $n=2h$ 
and $k=h$. If $n=2h$
 there are two distinct $\overline{G}$-orbits of $h$-dimensional
isotropic subspaces.

Let $L=K(1,2h)$:
%
if $L_0$ is the graded subalgebra of $L$
of type $(1|1, \dots,1, 0,\dots,$ $0)$ with $h$ zeros, then
$\pi(L_0)=\langle \xi_{1}, \dots, \xi_h\rangle$;
if $L_0$ is the graded subalgebra of $L$
of type $(1|1, \dots,1, 0, 1, 0,\dots,0)$ with $h$ zeros, then
$\pi(L_0)=\langle \xi_1, \dots, \xi_{h-1}, \xi_{h+1}\rangle$;
if $L_0$ is the graded subalgebra of $L$
of type $(2|2, \dots, 2,1,\dots, 1, 0,\dots,0)$ with $s+1$
2's and $s$ zeros, for some $s=0, \dots, h-2$,
 then $\pi(L_0)=\langle \xi_{1}, \dots,
\xi_{s}\rangle$.
Therefore, by Theorem \ref{K(2n+1)}$(i)$, all possible images of $\pi$ are
the isotropic subspaces of $V_{\bar{1}}$ except those of dimension $h-1$, and
$\Pi$ is injective.

Let $L=K(2k+1,n)$ where $n$ is odd and $k=0$, or  $n$ is arbitrary
  and $k>0$:
%
if $L_0$ is the graded subalgebra of $L$
of type $(2, 1, \dots,1|2, \dots, 2,$ $1,\dots, 1,$ $0,\dots,0)$ with $s+1$
2's and $s$ zeros, for some $s=0,\dots, [n/2]$,
 then $\pi(L_0)=\langle \xi_{1}, \dots,
\xi_{s}\rangle$.
If $n=2h$ the graded subalgebra of $L$ of type 
$(2, 1, \dots,1|$ $2, \dots, 2,0,2,0,\dots,0)$, with $h$ zeros and $h+1$ 2's,
 is not conjugate to the
graded subalgebra of type $(2, 1, \dots,1|2, \dots, 2,0,\dots,0)$
with $h$ zeros and $h+1$ 2's, and its image through $\pi$ is
the subspace $\langle \xi_{1}, \dots, \xi_{h-1}, 
\xi_{h+1}\rangle$.
By Theorem \ref{K(2n+1)}$(ii)$ and \ref{K(2n+1)}$(iii)$, $\Pi$ is bijective.

\item 
{$L=SHO(3,3)$}: we identify $L$ with the set of elements in
$\{f\in\Lambda(n,n)/\C 1|$ $\Delta(f)=0\}$ not containing the 
monomial $\xi_1\xi_2\xi_3$, with
reversed parity. $L$ has infinitely many invariant subalgebras of minimal
codimension whose intersection is the subalgebra of
type $(2,2,2|1,1,1)$ which is, therefore,
the canonical subalgebra of $L$ (cf.\ Remark \ref{inf.manySHO}).
 It follows that
  $V_{\bar{1}}=\langle x_1, x_2, x_3, \xi_1\xi_2, \xi_1\xi_3,
\xi_2\xi_3\rangle$ and
$\overline{G}=SL_3\times GL_2$. Consider the map $\psi: S^2V_{\bar{1}} 
\longrightarrow \langle \xi_i ~|~ i=1,2,3\rangle$ given by
$\psi(x_j\otimes x_k)=0$, $\psi(\xi_i\xi_j\otimes\xi_h\xi_k)=0$, 
$\psi(x_i\otimes\xi_j\xi_k)=\delta_{ij}\xi_k-
\delta_{ik}\xi_j$. A subspace of $V_{\bar{1}}$ is abelian if and
only if $\psi(a\otimes b)=0$ for any pair of elements $a,b$ of this subspace.
It follows that
  the $\overline{G}$-orbits of the non-trivial abelian subspaces of $V_{\bar{1}}$
are the orbits of the following subspaces: $\langle x_1\rangle$,
$\langle x_1, x_2\rangle$, $\langle x_1, \xi_2\xi_3\rangle$,
$\langle x_1, x_2, x_3\rangle$.

If $L_0$ is the graded subalgebra of $L$ of type $(1,1,1|1,1,1)$,
then $\pi(L_0)= 
\langle 
\xi_1\xi_2, \xi_1\xi_3, \xi_2\xi_3 \rangle$;
if  $L_0$ is the graded subalgebra of $L$
of type 
$(1,1,2|1,1,0)$, then 
$\pi(L_0)= 
\langle 
\xi_1\xi_2, x_3 \rangle$.
By Theorem \ref{HO(n,n)}$(b)$, the map
$\Pi$ is injective but not surjective. Indeed its image does not
contain the orbit of the trivial subspace, that
of the one-dimensional
subspaces and that of the subspace $\langle x_1, x_2\rangle$.
\item 
{$L=HO(n,n)$ (resp.\ $L=SHO(n,n)$ with $n>3$)}:
we identify $HO(n,n)$ with $\Lambda(n,n)/\C 1$ with
reversed parity, and
$SHO(n,n)$ with the set of elements in
$\{f\in\Lambda(n,n)/\C 1~|~ \Delta(f)=0\}$ not containing
the monomial $\xi_1\dots\xi_n$. Then
$V_{\bar{1}}=\langle x_1, \dots, x_n
\rangle \, , \,\overline{G}
=GL_n (\CC)$, any subspace of $V_{\bar{1}}$ is abelian and its
$\overline{G}$-orbit is determined by the dimension.

If $L_0$ is the graded subalgebra of $L$ of type $(1,\dots,1|0,\dots,0)$,
then $\pi(L_0)= 
\langle x_1,
\dots, x_n \rangle$;
if  $L_0$ is the graded subalgebra of $L$
of type 
$(1, \dots, 1,2,\dots,2|1,\dots,1,$ $0,\dots,0)$
with $n-s$ 2's and
$n-s$ zeros, for some
$s=2, \dots, n$, then 
$\pi(L_0)=\langle x_{s+1}, \dots, 
x_n\rangle$.
By Theorem \ref{HO(n,n)}$(a)$,  the image of $\pi$ consists of all
  subspaces of 
$\langle x_1,
\dots, x_n \rangle$
except those of codimension~1, and the map $\Pi$ is injective.
\item 
{$L=H(2k,n)$}:
we identify $L$ with $\Lambda(2k,n)/\mathbb{C}1$. Then
$V_{\bar{1}}=\langle \xi_1, \dots, \xi_{n}
\rangle$ with the bilinear form $(\xi_i,\xi_j)=\delta_{i,n-j+1}$ (cf.\ Example
\ref{boston}), $\overline{G}=\CC^\times SO_n(\mathbb{C})$, and any subspace of
$V_{\bar{1}}$ is abelian.
Let $S$ be a subspace of $V_{\bar{1}}$ and let $S=S^0\oplus S^1$
where $S^0$ is the kernel of the restriction of the bilinear form $(\cdot, \cdot)$ to $S$.
Let $s_i=\dim(S^i)$. Then
the $\overline{G}$-orbit of $S$ is determined by the pair $(s_0, s_1)$
unless $s_1=0$, $n$ is even and $s_0=n/2$. If $n$ is even 
then there are two distinct orbits of maximal isotropic subspaces
of $V_{\bar{1}}$.

%
%

If $L_0=L_0(U)$ is the maximal open subalgebra of $L$ constructed in
Example \ref{boston}, then $\pi(L_0(U))=U^0+(U^1)'$. By
Theorem \ref{H(2k,n)} and Remark \ref{gradedofH}, the map $\Pi$ is
bijective.
\item 
{$L=KO(2,3)$}:
we identify $L$ with $\Lambda(2,3)$ 
 with reversed parity. The canonical subalgebra of $L$ is its
subalgebra of subprincipal type. Therefore
$V_{\bar{1}}=\langle 1, \xi_1\xi_2
\rangle$ and
any subspace of $V_{\bar{1}}$ is abelian. 
$\overline{G}$ is the subgroup of $GL_2(\CC)$ consisting of upper
triangular matrices,
thus there are four $\overline{G}$-orbits of abelian
subspaces in $V_{\bar{1}}$: the orbit of the zero-dimensional
subspace, the orbit of the two-dimensional
subspace, the orbit of the one-dimensional subspace $\langle 1\rangle$
and the orbit of the one-dimensional subspace $\langle \xi_1\xi_2\rangle$.

If $L_0$ is the subalgebra of $L$ of principal type or
the subalgebra of subprincipal type, then
$\pi(L_0)=\langle \xi_1\xi_2\rangle$ or $\pi(L_0)=\langle 0\rangle$,
respectively;
if $L_0$ is the subalgebra constructed in Example \ref{ex1KO}, 
then $\pi(L_0)=\langle 1\rangle$; finally, if 
$L_0(2)$ is the subalgebra constructed in Example \ref{ex2KO}, then
$\pi(L_0(2))=\langle 1, \xi_1\xi_2\rangle$. 
By Theorem \ref{KO}, the map $\Pi$ is bijective.
\item 
{$L=SKO(2,3;\beta)$ with $\beta\neq 0,1$}.
The canonical subalgebra of $L$ is its
subalgebra of subprincipal type. Therefore
$V_{\bar{1}}=\langle 1, \xi_1\xi_2
\rangle$,
any subspace of $V_{\bar{1}}$ is abelian and 
$\overline{G}$ is the subgroup of $GL_2(\CC)$ consisting of diagonal matrices.
It follows that there are five  $\overline{G}$-orbits of abelian
subspaces in $V_{\bar{1}}$: the orbit of the zero-dimensional
subspace, the orbit of the two-dimensional
subspace, the orbit of the one-dimensional subspace $\langle 1\rangle$,
the orbit of the one-dimensional subspace $\langle \xi_1\xi_2\rangle$, and
the orbit of the one-dimensional subspace $\langle 1+\xi_1\xi_2\rangle$.

If $L_0$ is the subalgebra of $L$ of type $(1,1|0,0,1)$, $(1,1|1,1,2)$,
$(1,1|-1,-1,0)$, then $\pi(L_0)=\langle 0\rangle$,
$\pi(L_0)=\langle \xi_1\xi_2\rangle$,
$\pi(L_0)=\langle 1\rangle$, respectively; if 
$S_0(2)$ is the subalgebra of $L$ constructed in Example \ref{ex2sko}, then
$\pi(S_0(2))=\langle 1, \xi_1\xi_2\rangle$. By
 Theorem \ref{SKO}$(a)$, the map $\Pi$ is injective but not surjective,
since its image does not contain the orbit of the subspace
$\langle 1+\xi_1\xi_2\rangle$.
\item 
{$L=SKO(2,3;1)$}.
The canonical subalgebra of $L$ is its
subalgebra of subprincipal type. Therefore
$V_{\bar{1}}=\langle 1, \xi_1\xi_2
\rangle$,
any subspace of $V_{\bar{1}}$ is abelian and 
$\overline{G}=GL_2$.
It follows that the  $\overline{G}$-orbit of an abelian
subspace of $V_{\bar{1}}$ is determined by its dimension.

If $L_0$ is the subalgebra of $L$ of type $(1,1|0,0,1)$ or $(1,1|1,1,2)$, then
$\pi(L_0)=\langle 0\rangle$ or $\pi(L_0)=\langle \xi_1\xi_2\rangle$,
respectively;
if 
$S_0(2)$ is the subalgebra of $L$ constructed in Example \ref{ex2sko}, then
$\pi(S_0(2))=\langle 1, \xi_1\xi_2\rangle$. 
By Theorem \ref{SKO}$(b)$, the map $\Pi$ is bijective.
\item
{$L=SKO(2,3;0)$}.
The canonical subalgebra of $L$ is its
subalgebra of subprincipal type. $V_{\bar{1}}=\langle 1\rangle$ and
any subspace of $V_{\bar{1}}$ is abelian; $\overline{G}=\C^{\times}$.
It follows that there are two $\overline{G}$-orbits of abelian
subspaces in $V_{\bar{1}}$: the orbit of the zero-dimensional subspace
and the orbit of the one-dimensional subspace.

If $L_0$ is the subalgebra of type $(1,1|0,0,1)$, then
$\pi(L_0)=\langle 0\rangle$; if $L_0$ is the subalgebra of type
$(1,1|-1,-1,0)$, then $\pi(L_0)$ is $\langle 1\rangle$.
By Theorem \ref{SKO}$(c)$, $\Pi$ is bijective.
\item 
{$L=SKO(3,4; 1/3)$}. $L$ has 2
subalgebras of minimal codimension: the subalgebras of
principal and subprincipal type. These subalgebras are not
conjugate since the grading of principal type has depth 2 and
the grading of subprincipal type has depth 1.
The canonical subalgebra is, by definition, the graded subalgebra
of type $(2,2,2|1,1,1,3)$, therefore
$V_{\bar{1}}=\langle 1, x_1, x_2, x_3, \xi_2\xi_3, \xi_3\xi_1, \xi_1\xi_2 
\rangle$ with the nontrivial filtration:
$V_{\bar{1}}=V_{-3}\supset V_{-1}$ where $V_{-1}=\langle
x_1, x_2, x_3, \xi_2\xi_3, \xi_3\xi_1,$ $\xi_1\xi_2 
\rangle$.
$\overline{G}=\mathbb{C}^{\times}G'$
where $G'$ consists of matrices
$\left(\begin{array}{cc}
a & c\\
0 & 1
\end{array}
\right)$
where $c$ is an arbitrary
$6\times 1$ matrix
and $a$ belongs to the subgroup of $GL_6(\CC)$ consisting of
matrices $\left\{
\left(\begin{array}{c|c}
A & 0\\
\hline
\sigma A & A
\end{array}
\right)\right\}$ such that $A\in SL_3(\CC)$ and $\sigma\in\CC$. 
Here $\mathbb{C}^{\times}$ acts on $\mathfrak{g}_{-1}
=V_{-1}$ by multiplication by a scalar $\lambda$ 
and on
$\mathfrak{g}_{-3}
=V_{-3}/V_{-1}$ by multiplication by $\lambda^3$.
Consider the map $\psi: S^2V_{\bar{1}}\longrightarrow\langle\xi_i~|~i=1,2,3
\rangle$ given by:
$\psi(1\otimes a)=0$ for $a\in V_{\bar{1}}$, $\psi(x_i\otimes x_j)=0=\psi(\xi_i\xi_j\otimes\xi_k\xi_h)$,
$\psi(x_i\otimes \xi_j\xi_k)=\delta_{ij}\xi_k-\delta_{ik}\xi_j$.
A subspace of $V_{\bar{1}}$ is abelian if and only if $\psi(a\otimes b)=0$
for any pair of elements $a, b$ of this subspace. It follows that the
$\overline{G}$-orbits of the nontrivial abelian subspaces of $V_{\bar{1}}$
are the orbits of the following subspaces: $\langle 1\rangle$,
 $\langle x_1\rangle$,
$\langle \xi_1\xi_2\rangle$, $\langle 1, x_1\rangle$, $\langle 1, 
\xi_1\xi_2\rangle$, $\langle x_3, \xi_1\xi_2\rangle$, $\langle x_1, x_2\rangle$, $\langle \xi_1\xi_2, \xi_1\xi_3\rangle$,
$\langle 1, x_1, x_2\rangle$, $\langle 1, \xi_1\xi_2, \xi_1\xi_3\rangle$,
$\langle x_1, x_2, x_3\rangle$, $\langle \xi_1\xi_2, \xi_1\xi_3, \xi_2\xi_3\rangle$, $\langle 1, \xi_1\xi_2, x_3\rangle$, $\langle 1, x_1, x_2, x_3\rangle$, $\langle 1, \xi_1\xi_2, \xi_1\xi_3, \xi_2\xi_3\rangle$.

If $L_0$ is the subalgebra of  type $(1,1,1|0,0,0,1)$, $(1,1,1|1,1,1,2)$,
\break $(1,1,2|1,1,0,2)$, then
$\pi(L_0)=\langle x_1, x_2, x_3\rangle$,
$\pi(L_0)=\langle \xi_1\xi_2, \xi_1\xi_3, \xi_2\xi_3\rangle$,
$\pi(L_0)=\langle x_3, \xi_1\xi_2\rangle$, respectively;
if $S'_0$ is the subalgebra of $L$ constructed in Example \ref{ex1sko}, then
$\pi(S_0)=\langle 1, x_1, x_2, x_3\rangle$;
if $S_0(2)$ and $S_0(3)$ are  the subalgebras of $L$ constructed in Example \ref{ex2sko}, then
$\pi(S_0(2))=\langle 1, \xi_1\xi_2, x_3\rangle$ and
$\pi(S_0(3))=\langle 1, \xi_1\xi_2, \xi_1\xi_3, \xi_2\xi_3\rangle$.
By Theorem \ref{SKO}$(d)$, the
map $\Pi$ is injective but not surjective. 
\item 
{$L=KO(n,n+1)$ with $n>2$ (resp.\ $L=SKO(n,n+1; \beta)$ with
$n\geq 3$ and
$\beta\neq 1/3$ if $n=3$)}.
$V_{\bar{1}}=\langle 1, x_1, \dots, x_{n}
\rangle$. In this case $V_{\bar{1}}$ has a nontrivial filtration:
$V_{\bar{1}}=V_{-2}\supset V_{-1}$ where $V_{-1}=\langle x_i ~|~ i=1,
\dots,n\rangle$; $\overline{G}=\mathbb{C}^{\times}G'$
where $G'$ consists of matrices
$\left(\begin{array}{cc}
a & c\\
0 & 1
\end{array}
\right)$
with $a\in GL_n(\mathbb{C})$, 
 and where $c$ is an arbitrary
$n\times 1$ matrix. Here $\mathbb{C}^{\times}$ acts on $\mathfrak{g}_{-1}
=V_{-1}$ by multiplication by a scalar $\lambda$ (resp.\
$\sigma^{1-\beta}$) and on
$\mathfrak{g}_{-2}
=V_{-2}/V_{-1}$ by multiplication by $\lambda^2$ (resp.\ $\sigma^2$).
Any subspace of $V_{\bar{1}}$ is abelian. For any $k\in\mathbb{N}$,
$1\leq k\leq n$, there are two $\overline{G}$-orbits of abelian
subspaces of $V_{\bar{1}}$ of dimension $k$: one containing 1 and
the other contained in $\langle x_1, \dots, x_n\rangle$. 

Let $L=KO(n,n+1)$ with $n>2$:
if $L_0$ is the graded subalgebra of type $(1,\dots,1|0,\dots,0,1)$
then $\pi(L_0)=\langle x_1, \dots, x_n\rangle$;
if $L_0$ is the graded subalgebra of type $(1, \dots,1,2,\dots,2|
1,\dots,1,0,\dots,0,2)$ with $n-t+1$ 2's and $n-t$ zeros, for some
$t=2,\dots, n$, then $\pi(L_0)=\langle x_{t+1}, \dots, x_n\rangle$;
if $L_0$ is the subalgebra of $L$ constructed in Example \ref{ex1KO}, then
$\pi(L_0)=\langle 1, x_1, \dots, x_n\rangle$;
if $L_0(t)$ is the subalgebra of $L$ constructed in Example \ref{ex2KO},
for some $t=2, \dots, n$, then
$\pi(L_0)=\langle 1, x_{t+1}, \dots, x_n\rangle$.

By Theorem \ref{KO}
the image of $\pi$ consists of all subspaces
of $\langle x_{1}, \dots, x_n\rangle$ except those of codimension 1,
and of all subspaces
of $\langle 1, x_{1}, \dots, x_n\rangle$ containing 1 except those of codimension 1.
By Theorem \ref{SKO} the same description of the image of $\pi$ holds
for $L=SKO(n,n+1;\beta)$ with $n>2$.
The
map $\Pi$ is therefore injective but not surjective.
\item 
{$L=SHO^\sim(n,n)$}. $V_{\bar{1}}=\langle x_1,
\dots, x_n\rangle$; $\overline{G}=SL_n$; any subspace of
$V_{\bar{1}}$ is abelian and its $\overline{G}$-orbit is determined
by the dimension.

If $L_0$ is the graded subalgebra of type $(1,\dots,1|0,\dots,0)$
then $\pi(L_0)=\langle x_1, \dots,$  $x_n\rangle$;
if $L_0(t)$ is the maximal open subalgebra of $L$ constructed
in Example \ref{9th}, for some $t=2, \dots, n$, then $\pi(L_0(t))=
\langle x_{t+1}, \dots, x_n\rangle$.

By Theorem
\ref{SHOsim} the image of $\pi$ consists of all subspaces of $\langle x_1,
\dots, x_n\rangle$ except those of codimension 1.  
Therefore the map $\Pi$ is injective but not surjective.
\item 
{$L=SKO^\sim(n,n+1)$}. $V_{\bar{1}}=\langle 1, x_1,
\dots, x_n\rangle$. As in the case of $KO(n,n+1)$,
$V_{\bar{1}}$ has a nontrivial filtration:
$V_{\bar{1}}=V_{-2}\supset V_{-1}$ where $V_{-1}=\langle x_i ~|~ i=1,
\dots,n\rangle$; $\overline{G}=\mathbb{C}^{\times}G'$
where $G'$ consists of matrices
$\left(\begin{array}{cc}
a & c\\
0 & 1
\end{array}
\right)$
with $a\in SL_n(\mathbb{C})$, 
 and where $c$ is an arbitrary
$n\times 1$ matrix. Here $\mathbb{C}^{\times}$ acts on $\mathfrak{g}_{-1}
=V_{-1}$ by multiplication by a scalar 
$\sigma^{-2/n}$, and on
$\mathfrak{g}_{-2}
=V_{-2}/V_{-1}$ by multiplication by  $\sigma^2$.
 The description of 
the $\overline{G}$-orbits of the abelian subspaces of $V_{\bar{1}}$ is
the same as for $SKO(n,n+1;(n+2)/n)$ with $n>2$.

If $L_0$ is the subalgebra of $L$ constructed in Example
\ref{10thII}, then $\pi(L_0)=\langle x_1, \dots, x_n\rangle$;
if $L_0(t)$ is the subalgebra of $L$ constructed in Example 
\ref{10thI}, for some $t=2, \dots, n$,  then
$\pi(L_0(t))=\langle x_{t+1}, \dots, x_n\rangle$;
if $S_0(t)$ is the subalgebra of $L$ constructed in Example 
\ref{10thIII}, for some $t=2, \dots, n$,  then
$\pi(S_0(t))=\langle 1, x_{t+1}, \dots, x_n\rangle$.

By Theorem \ref{SKOsim} all possible images of $\pi$ are all subspaces of $\langle x_1, \dots, x_n\rangle$
except those of codimension 1, and all subspaces of
$\langle 1, x_1, \dots, x_n\rangle$ containing 1 except those of codimension 1 and 0.
The map $\Pi$ is therefore injective but
not surjective.
\item 
{$L=E(1,6)$}. $V_{\bar{1}}$,
$\overline{G}$ and  the $\overline{G}$-orbits of the abelian
subspaces of $V_{\bar{1}}$ are the same as for $K(1,6)$.

If $L_0$ is the graded subalgebra of type $(2|1,1,1,1,1,1)$,
 $(1|1,1,1,0,0,0)$,\break $(1|1,1,0,0,0,1)$, $(1|2,1,1,0,1,1)$, then
$\pi(L_0)=\langle 0\rangle$,
$\langle \xi_1, \xi_2, \xi_3\rangle$,
$\langle \xi_1, \xi_2, \eta_3\rangle$, and
$\langle \xi_1\rangle$, respectively.
Therefore, by Theorem \ref{generalE(1,6)},  all
possible images of $\pi$ are (as for $L=K(1,6)$) all isotropic subspaces
of $V_{\bar{1}}$ except those of dimension 2.
 The map $\Pi$ is therefore injective but not 
surjective.
\item 
{$L=E(3,6)$}. $V_{\bar{1}}=\langle a_{ij}:= dx_i v_j ~|~
i=1,2,3 \, , \, j=1,2\rangle$; $\overline{G} = GL_3 (\CC)\times
SL_2 (\CC)$ acting on $V_{\bar{1}} \simeq \CC^3 \otimes \CC^2$.
Consider the map $\psi : S^2 V_{\bar{1}} \to \langle
\frac{\partial}{\partial x_i} | i=1,2,3 \rangle$, given by
$\psi (a_{ij} \otimes a_{rs})= \epsilon (irk) \epsilon (js)
\frac{\partial}{\partial x_k}$, where $\epsilon$ is the sign of
the permutation~$irk$ (resp.~$js$) if all $i,r,k$ (resp.~$j,s$)
are distinct and $\epsilon =0$ otherwise.  A subspace of
$V_{\bar{1}}$ is abelian if and only if $\psi (a \otimes b)=0$ for any pair
of elements $a, b$ of this subspace.  

By Theorem \ref{generalE(3,6)},
all maximal open subalgebras are graded, and they are,
up to conjugation, the subalgebras of type $(2,2,2,0)$,
$(2,1,1,0)$ and $(1,1,1,1/2)$,
 so that the corresponding abelian
subspaces are $0$, $\langle a_{11}, a_{12}\rangle$ and $\langle
a_{11}, a_{21}, a_{31}\rangle$, respectively. 
Therefore all possible non-zero images of $\pi$ are given by all maximal
abelian subspaces of $V_{\bar{1}}$. Thus, the map $\Pi$ is
injective, but not surjective, as the remaining two $\overline{G}$-orbits
 of abelian subspaces, that of $\langle a_{11}\rangle $ and
 $\langle a_{11}, a_{21} \rangle$, are missing.
%
%
%
%
\item 
{$L=E(5,10)$}. $V_{\bar{1}}=\langle q_{ij}:=dx_i\wedge
  dx_j ~|~ i, j=1,2,3,4,5\rangle$, $\overline{G} = GL_5 (\CC)$,
  acting on $V_{\bar{1}} \simeq \Lambda^2 \CC^5$.  Consider the
  map $\varphi : S^2 V_{\bar{1}} \to \langle
  \frac{\partial}{\partial x_i} | i=1,\ldots ,5 \rangle$, given
  by $\varphi (q_{ij} \otimes q_{rs}) =\epsilon (ijrsk)
  \frac{\partial}{\partial x_k}$, where as before, $\epsilon$ is
  the sign of the permutation $ijrsk$ if all $i,j,r,s,k$ are
  distinct and $\epsilon =0$ otherwise.  A subspace of
  $V_{\bar{1}}$ is abelian if and only if $\varphi (a \otimes b) =0$ for any
  pair of elements of this subspace.  

By Theorem \ref{E(5,10)} all maximal
  open subalgebras are graded, of  type $(2,2,2,2,2)$,
  $(3,3,2,2,2)$, $(2,2,2,1,1)$ and $(2,1,1,1,1)$, up to
conjugation, so that the
  corresponding abelian subspaces of $V_{\bar{1}}$ are $0$,
  $\langle q_{12} \rangle$, $\langle q_{12}, q_{13}, q_{23}
  \rangle$ and $\langle q_{1j}| j=2,3,4,5 \rangle$,
  respectively.  Thus the map $\Pi$ is injective, but not
  surjective, as the remaining two $\bar{G}$-orbits of abelian
  subspaces, that of $\langle q_{12}, q_{13} \rangle$ and $\langle q_{12},
  q_{13}, q_{14}  \rangle$,  are missing.
%
%
%
%
\item 
{$L=E(4,4)$}. $V_{\bar{1}}=\langle dx_i~|~i=1,2,3,4\rangle$;
$\overline{G}=GL_4(\mathbb{C})$ acting on $V_{\bar{1}}\cong \mathbb{C}^4$.
Any subspace of $V_{\bar{1}}$ is abelian and its $\overline{G}$-orbit
is determined by the dimension.

If $L_0$ is the graded subalgebra of $L$ of type $(1,1,1,1)$,
then $\pi(L_0)=\langle 0\rangle$;
if $L_0$ is the maximal open subalgebra of $L$ constructed in Examples
\ref{3.1},\ref{3.2}, and \ref{3.3}, then $\pi(L_0)=\langle dx_1\rangle$,
$\pi(L_0)=\langle dx_1, dx_2\rangle$, and
$\pi(L_0)=\langle dx_i ~|~ i=1,2,3,4 \rangle$ respectively.
By Theorem \ref{E(4,4)} all possible images of $\pi$ are all subspaces of $V_{\bar{1}}$
except those of codimension~1. 
Therefore the map $\Pi$ is injective but not surjective.
\item 
{$L=E(3,8)$}. $V_{\bar{1}}=\langle v_1, v_2, x_iv_1, x_iv_2 ~|~ i=1,2,3 \rangle$
has a nontrivial filtration: $V_{\bar{1}}=V_{-3}\supset V_{-1}$
where $V_{-1}=\langle q_{ij}:=x_iv_j ~|~ i=1,2,3, j=1,2\rangle$. 
We can give the following description of abelian subspaces of
$V_{\bar{1}}$: consider the map
$\varphi : S^2 V_{-1} \to \langle
  \frac{\partial}{\partial x_i} | i=1,2,3 \rangle$, given
  by $\varphi (q_{ij} \otimes q_{rs}) =\epsilon (irk)\epsilon(js)
  \frac{\partial}{\partial x_k}$, where, as for $L=E(3,6)$, $\epsilon$ is
  the sign of the permutation $irk$ (resp.\ $js$) if all $i,r,k$ (resp.\
$j, s$) are
  distinct and $\epsilon =0$ otherwise.  A subspace of
  $V_{\bar{1}}$ is abelian if and only if $\varphi (a \otimes b) =0$ for any
  $a,b$ from this subspace.

$\overline{G}=\mathbb{C}^\times(SL_3\times SL_2)$ acts on $V_{\bar{1}}$
as follows: 
$\mathbb{C}^{\times}$ acts
on $\mathfrak{g}_{-1}=V_{-1}$ by multiplication by a scalar $\lambda$ and
on $\mathfrak{g}_{-3}=V_{-3}/V_{-1}$ by multiplication by $\lambda^3$;
$SL_3$ acts trivially on $\mathfrak{g}_{-3}$ and it acts
  on $\mathfrak{g}_{-1}=\mathbb{C}^3\otimes \mathbb{C}^2$
as on the direct sum of two copies of the standard $SL_3$-module;
finally, $SL_2$ acts on $\mathfrak{g}_{-3}$ as on the standard $SL_2$-module and
it acts on $\mathfrak{g}_{-1}$ as on the direct sum of three copies of
the standard $SL_2$-module.

If $L_0$ is the graded subalgebra of type $(2,2,2,-3)$, $(2,1,1,-2)$, or
$(1,1,1,-1)$, then $\pi(L_0)=\langle 0\rangle$,
$\langle x_1v_1, x_1v_2\rangle$, or
$\langle x_iv_1 ~|~ i=1,2,3 \rangle$,
respectively;
if $L_0$ is the maximal subalgebra of $L$ constructed in Example 
\ref{4.1}, \ref{4.2}, \ref{4.3}, \ref{4.4}, \ref{4.5}, or \ref{4.6},  then
$\pi(L_0)=\langle v_1, x_iv_1 ~|~ i=1,2,3 \rangle$,
$\langle v_1, x_1v_1, x_1v_2 \rangle$,
$\langle v_1, x_iv_2 ~|~ i=1,2,3 \rangle$,
$\langle v_1, v_2, x_1v_1, x_1v_2 \rangle$,
$\langle v_1, v_2,  x_iv_1 ~|~ i=1,2,3 \rangle$, or
$\langle v_1, v_2 \rangle$, respectively.

Therefore, by Theorem \ref{E38},
all possible images of $\pi$ are the subspace $\langle v_1, v_2\rangle$
and  every subspace $S$
of $V_{\bar{1}}$ such that 
$S\cap V_{-1}$ is  a maximal
abelian subspace of $V_{-1}$. It follows that  the map
$\Pi$ is injective but not surjective. The
$\overline{G}$-orbits of the following abelian subspaces of $V_{\bar{1}}$ 
are missing: $\langle v_1\rangle$,
$\langle x_1v_1\rangle$, $\langle v_1, x_1v_1\rangle$,
$\langle v_1, x_1v_2\rangle$, $\langle x_1v_1, x_2v_1\rangle$,
$\langle v_1, v_2, x_1v_1\rangle$, $\langle v_1, x_1v_1, x_2v_1\rangle$,
$\langle v_2, x_1v_1, x_2v_1\rangle$,
$\langle v_1, v_2,  x_1v_1, x_2v_1\rangle$.
\end{enumerate}

\medskip

We conclude by listing the maximal  among $\mathfrak{a}_0$-invariant
open subalgebras of $S$ which are not maximal and also those
maximal open subalgebras of $S$, none of whose conjugates is
$\mathfrak{a}_0$-invariant. The lists
follow from Theorems \ref{derS(1,n)}, \ref{a-invS(1,2)}, 
\ref{derSHO}, \ref{a-invSHO(3,3)}, \ref{derSKO}, \ref{SKO(n,n+1;1)}, and
\ref{sl_2-inv}.
\begin{theorem} Let $S$ be an infinite-dimensional linearly
 compact simple Lie superalgebra and
let $\mathfrak{a}_0$ be
a subalgebra of the subalgebra $\mathfrak{a}$ of outer derivations of $S$.

\noindent
$(a)$ A complete list of pairs $(S_0, \mathfrak{a}_0)$ where
$S_0$ is an open, maximal among the $\mathfrak{a}_0$-invariant 
subalgebras of $S$, which is not maximal,
is as
follows: 
\begin{description}
\item[-] $S=S(1,2)$, $S_0$ is the canonical subalgebra, $\mathfrak{a}_0=
\mathfrak{a}\cong sl_2$;
\item[-] $S=SHO(3,3)$, $S_0$ is the canonical subalgebra and $\mathfrak{a}_0
=sl_2$, or $\mathfrak{a}_0=\mathfrak{a}\cong gl_2$;
\item[-] $S=SKO(2,3;0)$, $S_0$ is the subalgebra of principal type or $S_0$ is
the subalgebra $S_0(2)$ constructed in Example \ref{ex2sko}, and
$\mathfrak{a}_0=\C\xi_1\xi_2$ or $\mathfrak{a}_0=\mathfrak{a}$
($\dim\mathfrak{a}=2$).
\end{description}

\noindent
$(b)$ A complete list of pairs $(\mathfrak{a}_0, S_0)$ where $S_0$ is a maximal open 
subalgebra of $S$, none of whose conjugates is $\mathfrak{a}_0$-invariant, is as follows:

\begin{description}
\item[-] $S=S(1,2)$ or $S=SKO(2,3;1)$: $\mathfrak{a}_0=\mathfrak{a}\cong
sl_2$ and $S_0$ is the graded
subalgebra of $S$ of principal type;
\item[-] $S=SHO(3,3)$: $\mathfrak{a}_0=sl_2$ or $\mathfrak{a}_0=
\mathfrak{a}\cong gl_2$ and $S_0$ is the graded subalgebra
of $S$ of principal type;
\item[-] $S=S(1,n)$ with $n\geq 3$: $\mathfrak{a}_0=\CC\xi_1\dots\xi_n\frac{\partial}{\partial x_1}$ or $\mathfrak{a}_0=\mathfrak{a}$ $(\dim\mathfrak{a}=2)$, and $S_0$ is the graded subalgebra of $S$ of type $(1|0,\dots,0)$;
\item[-] $S=SHO(n,n)$ with $n\geq 4$: $\mathfrak{a}_0=\CC\xi_1\dots\xi_n
\rtimes\mathfrak{t}$ where $\mathfrak{t}$ is a torus of
$\mathfrak{a}$ $(\dim\mathfrak{a}=3)$,
 and $S_0$ is the graded subalgebra of $S$ of type $(1, \dots,1|0,\dots,0)$;
\item[-] $S=SKO(2,3;0)$: $\mathfrak{a}_0=\C\xi_1\xi_2$
or $\mathfrak{a}_0=\mathfrak{a}$ $(\dim\mathfrak{a}=2)$, and
$S_0$ is the subalgebra of type $(1,1|0,0,1)$ or the subalgebra
of type $(1,1|-1,-1,0)$;
\item[-] $S=SKO(n,n+1;(n-2)/n)$ with $n>2$: 
$\mathfrak{a}_0=\C\xi_1\dots\xi_n$ or $\mathfrak{a}_0=\mathfrak{a}$ $(\dim\mathfrak{a}=2)$, and $S_0$ is the graded subalgebra of $S$ of type $(1, \dots,1|0,\dots,0,1)$ or the subalgebra $S'_0$ constructed in Example \ref{ex1sko};
\item[-] $S=SKO(n,n+1;1)$ with $n>2$: $\mathfrak{a}_0=\C\xi_1\dots\xi_n\tau$
or $\mathfrak{a}_0=\mathfrak{a}$ $(\dim\mathfrak{a}=2)$, and
$S_0$ is the subalgebra $S'_0$ constructed in Example \ref{ex1sko}.
\end{description}
\end{theorem}

$$$$ 

\end{document}